\documentclass[11pt]{amsart} 
\usepackage[a4paper]{geometry}

\usepackage[english]{babel}
\usepackage{amssymb, amsmath, amsthm}

\usepackage{tikz-cd}
\usepackage{tikz}
\usepackage{mathrsfs}
\usetikzlibrary{arrows,positioning,calc,intersections}
\usepackage{hyperref}

\usepackage[all,cmtip]{xy}

\newtheorem{theorem}{Theorem}[section]

\newtheorem{proposition}[theorem]{Proposition}
\newtheorem{corollary}[theorem]{Corollary}
\newtheorem{claim}[theorem]{Claim}

\theoremstyle{remark}
\newtheorem{remark}[theorem]{Remark}
\newtheorem{example}[theorem]{Example}
\newtheorem{definition}[theorem]{Definition}

\newcommand{\cupprod}{\!\smallsmile\!}

\newcommand{\Lie}{\mathscr{L}}

\DeclareMathOperator{\im}{im}
\DeclareMathOperator{\pr}{pr}
\DeclareMathOperator{\id}{Id}

\DeclareMathOperator{\al}{\mathit{A}}

\DeclareMathOperator{\totdim}{Totdim}
\DeclareMathOperator{\der}{Der}
\DeclareMathOperator{\ann}{Ann}
\DeclareMathOperator{\ider}{End}
\DeclareMathOperator{\Hom}{Hom}
\DeclareMathOperator{\rk}{rank}
\DeclareMathOperator{\Sh}{Sh}

\newcommand{\aaar}{\substack{\longrightarrow\\[-0.85em] \longrightarrow \\[-0.85em] \longrightarrow}}
\newcommand{\cA}{\mathcal{A}}
\newcommand{\cB}{\mathcal{B}}
\newcommand{\cC}{\mathcal{C}}
\newcommand{\cD}{\mathcal{D}iv}
\newcommand{\cE}{\mathcal{E}}
\newcommand{\cF}{\mathcal{F}}
\newcommand{\cG}{\mathcal{G}}

\newcommand{\cI}{\mathcal{I}}

\newcommand{\cK}{\mathcal{K}}
\newcommand{\cL}{\mathcal{L}}
\newcommand{\cM}{\mathcal{M}}
\newcommand{\cN}{\mathcal{N}}
\newcommand{\cO}{\mathcal{O}}
\newcommand{\cP}{\mathcal{P}}

\newcommand{\cS}{\mathcal{S}}
\newcommand{\cT}{\mathcal{T}}
\newcommand{\cV}{\mathcal{V}}

\newcommand{\cZ}{\mathcal{Z}}

\newcommand{\bZ}{\mathbb{Z}}
\newcommand{\bA}{\mathbb{A}}
\newcommand{\bR}{\mathbb{R}}
\newcommand{\bN}{\mathbb{N}}
\newcommand{\bT}{\mathbb{T}}

\newcommand{\g}{\mathfrak{g}}

\newcommand{\fX}{\mathfrak{X}}

\newcommand{\cdo}{\mathbb{D}}
\newcommand{\sdiv}{\mathcal{D}iv}

\newcommand{\bDelta}{\bigtriangleup}

\newcommand{\gs}{\mathtt{s}}
\newcommand{\gt}{\mathtt{t}}

\newcommand{\can}{\mathrm{can}}

\DeclareMathOperator{\sym}{S^\bullet}

\newcommand{\rank} 
{\operatornamewithlimits{rank}}

\newcommand{\man}[1]   {#1\text{-}\mathcal{M}an}
\newcommand{\cob}[1]   {#1\text{-}\mathcal{C}o{A}lg}

\let\oldtocsection=\tocsection

\let\oldtocsubsection=\tocsubsection

\let\oldtocsubsubsection=\tocsubsubsection

\renewcommand{\tocsection}[2]{\bf\hspace{0em}\oldtocsection{#1}{#2}}
\renewcommand{\tocsubsection}[2]{\hspace{1em}\oldtocsubsection{#1}{#2}}
\renewcommand{\tocsubsubsection}[2]{\hspace{2em}\oldtocsubsubsection{#1}{#2}}

\title[The symplectic geometry of graded manifolds and higher Lie groupoids]{Lecture notes on the symplectic geometry of graded manifolds and higher Lie groupoids}

\author{Miquel Cueca}

\address{Departement of Mathematics, KU Leuven. Celestijnenlaan 200B,  Leuven (Heverlee), B-3001, Belgium}

\email{miquel.cuecaten@kuleuven.be}

\author{Antonio Maglio}

\address{ Institute of Mathematics, Polish Academy of Sciences, \'{S}niadeckich 8, 00-656 Warszawa, Poland.}

\email{amaglio@impan.pl}

\author{Fabricio Valencia}

\address{Instituto de Matem\'{a}tica e Estat\'{i}stica, Universidade de S\~{a}o Paulo, Rua do Mat\~{a}o 1010, Cidade Universit\'{a}ria, 05508-090 S\~{a}o Paulo, Brazil.}

\email{fabricio.valencia@ime.usp.br}

\date{} 

\begin{document}

	\begin{abstract} 
		
		In this work, we study symplectic structures on graded manifolds and their global counterparts, higher Lie groupoids. We begin by introducing the concept of graded manifold, starting with the degree 1 case, and translating key geometric structures into classical differential geometry terms. We then extend our discussion to the degree 2 case, presenting several illustrative examples with a particular emphasis on equivariant cohomology and Lie bialgebroids.
		Next, we define symplectic Q-manifolds and their Lagrangian Q-submanifolds, introducing a graded analogue of Weinstein’s tubular neighborhood theorem and applying it to the study of deformations of these submanifolds.
		
		Shifting focus, we turn to higher Lie groupoids and the shifted symplectic structures introduced by Getzler. We examine their Morita invariance and provide several examples drawn from the literature. Finally, we introduce shifted Lagrangian structures and explore their connections to moment maps and symplectic reduction procedures.
		
		Throughout these notes, we illustrate the key constructions and results with concrete examples, highlighting their applications in  mathematics and physics.
		 
	\end{abstract}

	\maketitle

	\tableofcontents 


	\section*{Introduction}
	
	Higher structures such as $L_\infty$-algebras, Courant algebroids, and homotopy Poisson structures, among others, have become ubiquitous in  mathematics and physics over the last fifty years. To deal with them, each community has developed a different framework. In particular, three of those theories are relevant for these lecture notes:
	\begin{enumerate}
		\item[(A)]Using sheaf theory, mathematical physicists have developed the theory of $Q$-manifolds, which is the main topic of Part \ref{partI}. Roughly speaking, a $Q$-manifold $\cM=(M,C_\bullet(\cM))$ is a manifold $M$ whose sheaf of functions is enlarged with new polynomial functions that commute or skew-commute according to a degree, and it is equipped with a differential $Q:C_\bullet(\cM)\to C_{\bullet+1}(\cM)$.
		\item[(B)] Using simplicial methods, homotopy theorists and differential geometers have developed the framework of Lie $n$-groupoids, which is the main topic of Part \ref{partII}.  Roughly speaking, a Lie $n$-groupoid is a simplicial manifold $K_\bullet$ satisfying certain group and smooth-like conditions, known as \emph{Kan conditions}. 
		\item[(C)] Using the theory of descent, algebraic geometers have introduced the theory of $n$-stacks. This will appear while considering the homotopy theory of Lie $n$-groupoids, as in $\S~\ref{sec:morita}$.  
	\end{enumerate}

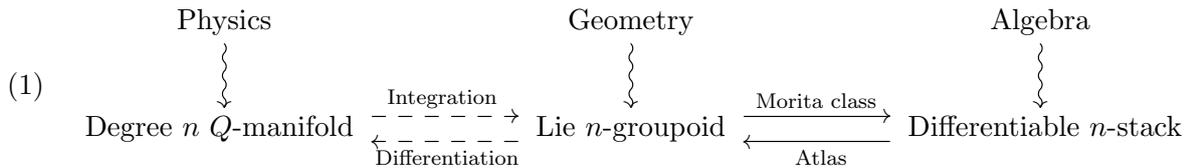
	\begin{equation}\label{diagall}
		\begin{array}{c}
			\xymatrix{\text{Physics} \ar@{~>}[d]&& \text{Geometry} \ar@{~>}[d]&& \text{Algebra}\ar@{~>}[d]\\
				\text{Degree } n  \ Q\text{-manifold }\ar@{-->}@<1ex>[rr]^{\quad\text{Integration}}&& \text{ Lie } n\text{-groupoid }\ar^{\quad \ \text{Differentiation}}@{-->}@<1ex>[ll]\ar@<1ex>[rr]^{\text{Morita class}\quad \ }&&\ar^{\text{Atlas}\quad}@<1ex>[ll]\text{ Differentiable } n\text{-stack}}
		\end{array}
	\end{equation}
	
	These lecture notes have two main goals. On the one hand, we aim to give a gentle introduction to the theory of $Q$-manifolds and higher Lie groupoids, and to explain  the relation\footnote{Or at least introduce the key ideas.} between them as depicted in \eqref{diagall}. On the other hand, we introduce and study the symplectic geometry of $Q$-manifolds and higher Lie groupoids. In particular, we will focus on examples arising from Poisson geometry and mathematical physics.
	
	This work does not contain original results and aims to provide a clear presentation of the different frameworks used to approach higher structures within the Poisson geometry community. In particular, we would like to point out that the previous works \cite{cat:intro, del:not, kap:sup, mnev:book, var:sup} also serve as excellent introductions for Part \ref{partI} (with the main difference being that they work with supermanifolds instead of $\bN$-graded manifolds). The second part is more novel, and we are not aware of any lecture notes dealing with Lie $n$-groupoids and the shifted symplectic geometry on them. Nevertheless, presentations of those topics within the realm of algebraic geometry are available, including \cite{cal:lec3, cal:intro, saf:lec, toe:der}. Finally, we want to advertise the lecture notes by Calaque-Ronchi, included in the same volume \cite{CR:shif}, which offer an alternative perspective on the subject.

	The lecture notes are divided into two parts, which deal with the left and middle columns in \eqref{diagall}, respectively. The content of Part \ref{partI}, the infinitesimal picture, is structured as follows: 
	\begin{itemize}
		\item In \S~\ref{sec:1}, rather than dealing with the general theory, we introduce the category of $\bN$-graded manifolds of degree $1$ and show that it is equivalent to the category of vector bundles. This allows us to establish relations between $\bN$-graded manifolds of degree $1$ equipped with suitably compatible geometric structures and other geometric structures defined over vector bundles, such as Lie algebroids (\S~\ref{sec:Q1}) and bialgebroids (\S~\ref{sec:bia1}).  We also define superdivergences (\S~\ref{sec:div}), which can be used to define algebroid homology. As an application, in \S~\ref{sec:morse} we prove the Morse inequalities by following the elegant approach due to Witten.
		\item  In \S~\ref{sec:2}, we introduce $\bN$-graded manifolds of arbitrary degree and provide the main definitions for Part \ref{partI}. We have gathered the technical aspects of the theory in this section to allow a smoother flow  in the next two sections.
		\item  \S~\ref{sec:3} is devoted to study $\bN$-graded manifolds of degree $2$. First, we present three different ways to characterize them (\S~\ref{sec:geo}). Then, as an application, we show that degree 2-manifolds endowed with homological vector fields are in one-to-one correspondence with Lie 2-algebroids (Theorem \ref{thm:bon}). Finally, we study the tangent and cotangent bundles of a degree $1$ manifold showing their relation to equivariant cohomology (\S~\ref{sec:equico}) and Lie quasi bialgebroids (\S~\ref{sec:liebia}), respectively.
		\item  In \S~\ref{sec:sympQ}, we define symplectic $Q$-manifolds, one of the central concepts of these lecture notes. After deducing some of their basic properties, we show that symplectic $Q$-manifolds of degree $1$ correspond to Poisson manifolds, while those of degree $2$ correspond to Courant algebroids, see also \S~\ref{sec:exaQs} for more examples. We then define Lagrangian $Q$-submanifolds and state a graded version of Weinstein's tubular neighborhood theorem (\S~\ref{sec:lagq}), which is a key tool for studying the deformations of such submanifolds. Finally, we comment on their homotopic versions (\S~\ref{sec:ssa}) and briefly introduce the AKSZ sigma models with target a symplectic $Q$-manifold (\S~\ref{sec:aksz}).
	\end{itemize}

	Let us now describe the content of Part \ref{partII}, the global picture: 
	\begin{itemize}
		\item In \S~\ref{sec:5}, we briefly recall the simplicial language (\S~\ref{sec:bsimp}), introduce simplicial manifolds, their geometric realizations (\S~\ref{sec:georea}) and the generalized de Rham Theorem (Theorem \ref{thm:genderham}). We then define higher Lie groupoids as simplicial manifolds satisfying Kan conditions (\S~\ref{sec:lien}), introduce Morita equivalence via hypercovers (\S~\ref{sec:morita}), and discuss the Lie functor, which relates the category of Lie $n$-groupoids to that of degree $n$ $Q$-manifolds (\S~\ref{sec:liefun}).
		\item In \S~\ref{sec:ss}, we introduce the main notion  of Part \ref{partII}: $m$-shifted symplectic structures on Lie $n$-groupoids. We show that this notion is invariant under Morita equivalence (\S~\ref{sec:defss}) and encodes several well-studied structures in the literature (\S~\ref{sec:fex}). We also provide a series of relevant examples and explain their key properties.  
		\item  Finally, in  \S~\ref{sec:7}, we define shifted Lagrangian structures on an $m$-shifted symplectic Lie $n$-groupoid,  briefly describe their associated ``symplectic category", and establish their connections to various moment map theories and their symplectic reductions (\S~\ref{sec:symcatlag}). We conclude by relating  these reduction procedures with classical field theories (\S~\ref{sec:CTFT}).
	\end{itemize} 
	
	It is worth emphasizing that all the main notions and results discussed throughout these lecture notes are illustrated with numerous examples and applications in mathematics and physics.

	\subsection*{Acknowledgments}
	These lecture notes are based on two mini-courses delivered by the first author at \emph{Geometry in Algebra and Algebra in Geometry VII (2023)} in Belo Horizonte, Brazil, and at the \emph{INdAM Intensive Period: Poisson Geometry and Mathematical Physics (2024)} in Napoli, Italy. We would like to thank the organizers for these wonderful events, and all participants for their questions and contributions, which help improve these notes. Some of the material included here builds upon works by M.~C.~and collaborators, whom we gratefully acknowledge. We also thank the anonymous referee, whose comments and suggestions have been essential in improving the final presentation of this work. M.~C.~was partially supported by FWO grant 1249325N. The research of A.~M.~was funded by the National Science Centre (Poland) within the project WEAVE-UNISONO, No. 2023/05/Y/ST1/00043 and also partially supported by GNSAGA of INdAM. F.~V.~was supported by Grant 2024/14883-6 S\~{a}o Paulo Research Foundation-FAPESP.

	\part{The infinitesimal picture}\label{partI}



	Objects in differential geometry are typically  defined in terms of ``space" or by describing how they act on it. Definitions in algebra are usually formulated in terms of structure or universal properties. Despite these different points of view, there is a way of going from differential geometry to algebra, as illustrated in   Table \ref{tab:dic}.
	
	Indeed, by adding suitable adjectives in each row of Table \ref{tab:dic}, we obtain a one-to-one correspondence. However,  the algebraic characterization of geometric objects is often a non-trivial task, see e.g. \cite{nest:obs} for a more detailed treatment of this topic.
	
	\begin{table}[h!]
		\centering
		\begin{tabular}{|c|c|}
			\hline
			Geometry &  Algebra \\
			\hline
			Manifold M & Algebra $(C^\infty(M), +,\cdot)$\\
			Smooth map $f:M\to N$ & Algebra morphism $f^*:C^\infty(N)\to C^\infty(M)$\\ 
			Submanifold $S\subseteq M$ & Ideal $Z(S)=\{ f_{|S}=0\}\subseteq C^\infty(M)$\\
			Vector field $X\in\fX^1(M)$& Derivation $X\in\der(C^\infty(M))$\\
			Vector bundle $E\to M$ & $C^\infty(M)$-module $\Gamma(E)$\\
			\hline
		\end{tabular}
		\caption{Geometry to algebra dictionary}
		\label{tab:dic}
	\end{table}
	
	In the first half of the past century, the quantum revolution was occurring in physics. The new physical paradigm classified matter particles into two types according to  their spin: bosons and fermions. While bosonic particles commute, fermions behave different: they skew-commute. As a consequence, fermionic particles cannot be represented by coordinates on an ordinary smooth manifold. To provide a rigorous mathematical framework for these new particles, and to attempt a formulation of classical mechanics involving them, in the 1970s, Berezin and Leites \cite{ber:sup}  introduced the concept of \emph{supermanifold}.
	Roughly speaking, they used the above dictionary, adding to the algebra $C^\infty(M)$ new formal Grassmann variables, so that the functions on a supermanifold take the form:
	$$
	C^\infty(M)[\xi^1,\cdots, \xi^k]\quad \text{where}\quad \xi^i\xi^j=-\xi^j\xi^i.
	$$
	
	In the 1990s, the works of Vaintrob \cite{vai:lie}, \v{S}evera  \cite{sev:some} and Roytenberg \cite{roy:on}, among others, refined this idea by introducing an additional $\bN$-grading on the functions of a supermanifold. More precisely, $\cM=(M, C_\bullet(\cM))$ is called an \emph{$\bN$-graded manifold of degree $n$}, or an \emph{$n$-manifold} for short, if $C_\bullet(\cM)$ is a sheaf of graded commutative algebras which locally looks like
	$$
	C_\bullet(\cM)_{|U}=C^\infty(U)\otimes \sym (\mathbf{V}), \;\; \text{ where } \mathbf{V}=\oplus_{i=1}^n V_i.
	$$
	Observe that the bullets on the left and right of the equation do not match. On the left the bullet denotes the internal degree, while on the right the bullet denotes the symmetric power. This abuse of notation will appear in other parts of this document. We hope this will not confuse the reader.
	
	These new objects are the main focus of Part \ref{partI} of these lecture notes. In the following sections, we will explain how they encode many geometric structures appearing in Poisson geometry and mathematical physics. Some important aspects of these notes are the following: we will hide homological algebra by doing differential geometry; usual differential geometric constructions yield  non-trivial results when applied to $n$-manifolds; technicalities are kept to a minimum, for example we have tried to avoid the use of sheaves and just introduced their global sections.

	\section{$\bN$-graded manifolds of degree 1}\label{sec:1}
	
	Here we study the particular case of $\bN$-graded manifolds of degree $1$. We begin by establishing the correspondence between the category of $1$-manifolds and that of vector bundles. Then, we show how vector fields, Poisson structures, and divergences on $1$-manifolds give rise to relevant structures in mathematics and physics. 
	
	\subsection{Definition and the equivalence with vector bundles}
	
	An \emph{$\bN$-graded manifold of degree $1$}  and \emph{dimension $m_0|m_1$} is a pair $\cM=(M, C_\bullet(\cM))$, where $M$ (the \emph{body} of $\cM$) is a 
	manifold of dimension $m_0$ and $ C_\bullet(\cM)$ is a sheaf of graded-commutative algebras on $M$ satisfying the
	following property: any $p\in M$ admits an open neighborhood $U\subseteq M$ such that
	\begin{equation}\label{eq:locmod1man}
		C_\bullet(\cM)_{|U}\cong C^\infty(U)\otimes\wedge^\bullet \bR^{m_1}
	\end{equation}
	where the elements of $\bR^{m_1}$ are of degree $1$. In other words, there exist sections $\{x^i, \xi^\alpha\}$ of $C_\bullet(\cM)$ over $U$, where $\{x^1, \dots, x^i,\dots,  x^{m_0}\}$ are usual coordinates on $U$ of degree $0$, the elements $\{\xi^1,\dots, \xi^\alpha,\dots, \xi^{m_1}\}$ are of degree 1 and every local section of $ C_\bullet(\cM)_{|U}$ can be expressed as a sum of functions that are smooth in $x^i$ and polynomial in $\xi^\alpha$. 
	
	The $\{x^i, \xi^\alpha\}$ are  \emph{local coordinates} of $\cM_{|U}=(U,C_\bullet(\cM)_{|U})$. Global sections of $C_\bullet(\cM)$ are called \emph{functions} on $\cM$. We denote by $C_l(\cM) \subseteq C_\bullet(\cM)$ the subsheaf of homogeneous sections of degree $l$. The degree of a homogeneous section is denoted by $|\cdot|$. Since $\bN$-graded manifold of degree $1$ is rather long, we will usually abbreviate it to \emph{$1$-manifold}.
	
	\begin{remark}[On supermanifolds I]\label{rmk:sup1}
		$\bN$-graded manifolds of degree $1$ should not be confused with the more common notion of a \emph{supermanifold}, see e.g. \cite{ber:sup, witt:int}. The definition looks the same with the only difference that coordinates have parity (i.e. they are $\bZ_2$-graded) instead of an integer degree. The $\bZ_2$-grading makes supermanifolds more flexible but also more technically involved.
	\end{remark}

	A \emph{morphism of $1$-manifolds} $\Psi:\cM\to \cN$ is a  morphism of ringed spaces, given by a pair $\Psi=(\psi, \psi^\sharp)$, where $\psi:M \to N$ is a smooth map and $\psi^\sharp:C_\bullet(\cN)\to\psi_* C_\bullet(\cM)$ is a degree preserving morphism of sheaves of algebras over $N$. Notice that $1$-manifolds with their morphisms form a category, which we denote by $\man{1}$.
	
	\begin{example}\label{ex:1-man}
		Given a vector bundle $E\to M$, we define the $1$-manifold   $$E[1]=\big(M,C_\bullet(E[1])=\Gamma({\wedge^\bullet E^*})\big).$$
		On a chart $U$ of $M$ where $E_{|U}$ is trivializable, if we pick coordinates $\{x^i\}_{i=1}^m$ on $U$ and a frame $\{\xi^\alpha\}_{\alpha=1}^{\rk(E)}$  of $E^*_{|U}$ we have that $\{x^i, \xi^\alpha\}$ 
		are local coordinates of $E[1]_{|U}$, where $|x^i|=0$ and $|\xi^\alpha|=1$.
		Two important examples of 1-manifolds are 
		$$
		T[1]M=\big(M,\Omega^\bullet(M)\big) \qquad \text{and} \qquad T^*[1]M=\big(M,\fX^\bullet(M)\big),
		$$ 
		where $\Omega^\bullet(M)=\Gamma({\wedge^\bullet T^*M})$ and $\fX^\bullet(M)= \Gamma({\wedge^\bullet TM})$ are the differential forms and multivector fields on $M$, respectively.
	\end{example}

	Denote by $\cV ect$ the category whose objects are vector bundles $E\to M$, over different bases, and morphisms are given by vector bundle morphisms covering a smooth map between the bases as the following diagram shows:
		\begin{equation}\label{eq:vbmor}
		\begin{array}{c}
			\xymatrix{E\ar[d]\ar[r]^{\Phi}& F\ar[d]\\ M\ar[r]^{\varphi}& N.}
		\end{array}
	\end{equation}
	 The following result gives a geometric characterization of the category $\man{1}$.
	
	\begin{theorem}\label{thm: 1-manifolds}
		The functor $$[1]:\cV ect \to \man{1}, \qquad (E\to M)\to E[1] $$ is an equivalence of categories.
	\end{theorem}

	\begin{proof}
		In Example \ref{ex:1-man} we defined the functor at the level of objects. Let $\Phi:E\to F$ be a vector bundle morphism covering the smooth map $\varphi:M\to N$. The map $\Phi$ induces the vector bundle morphism $\wedge^l\Phi:\wedge^l E\to \wedge^l F$   for each $l\in\bN$. Therefore, these morphisms define  maps of $C^\infty(N)$-modules at the level of sections given by $$\wedge^l\Phi^\sharp:\Gamma(\wedge^l F^*)\to \varphi_*\Gamma (\wedge^l E^*), \quad s\to \wedge^l\Phi^\sharp(s)=\wedge^l\Phi^*\circ s\circ \varphi.$$
		The collection of all these morphisms satisfies:
		\begin{enumerate}
			\item[$\star$] They preserve degree, because each $\wedge^l\Phi^\sharp$ does it.
			\item[$\star$] $\wedge^l\Phi^\sharp(fs)=(\varphi^*f)\wedge^l\Phi^\sharp(s)$ for each $f\in C^\infty(N)$ and $s\in\Gamma(\wedge^lF^*)$.
			\item[$\star$] $\big(\wedge^k\Phi^\sharp(s_1)\big)\wedge\big(\wedge^l\Phi^\sharp(s_2)\big)=\wedge^{k+l}\Phi^\sharp(s_1\wedge s_2)$.
		\end{enumerate}
		
		The above properties guarantee that $ \wedge^\bullet\Phi^\sharp:\Gamma(\wedge^\bullet F^*)\to\varphi_*\Gamma(\wedge^\bullet E^*)$ is a degree preserving morphism of sheaves of algebras over $N$ and thus $\Phi[1]=(\varphi, \wedge^\bullet\Phi^\sharp):E[1]\to F[1]$  is well-defined. Moreover, two vector bundle morphisms $\Phi,\Psi: E\to F$ are equal if and only if $\Phi^\sharp=\Psi^\sharp$. Therefore, the functor  $[1]$ is faithful because $\Phi[1]=\Psi[1]$ implies $\Phi^\sharp=\wedge^1\Phi^\sharp=\wedge^1\Psi^\sharp=\Psi^\sharp$.
		
		To see that the functor is full, let $\Phi=(\phi,\phi^\sharp):E[1]\to F[1]$ be a morphism between $1$-manifolds. Then 
		$$\phi^\sharp:C_1(F[1])=\Gamma (F^*)\to \phi_\ast C_1(E[1])=\phi_*\Gamma(E^*)$$ is a morphism of $C^\infty(N)$-modules satisfying $\Phi^\sharp(fs)=\phi^*(f)\Phi^\sharp(s)$ for $f\in C^\infty(N)$ and $s\in\Gamma(F^*)$. Hence, $\phi^\sharp$ induces a vector bundle morphism $\Phi:E\to F$ covering $\phi:M\to N$. Moreover, $\Phi[1]=\Phi$ as claimed.
		
		Finally we check that the functor is essentially surjective. Let $\mathcal{M}=(M,C_\bullet(\mathcal{M}))$ be a  $1$-manifold. The multiplication of functions $C_0(\cM)\times C_1(\cM)\to C_1(\cM)$ makes $C_1(\cM)$ into a $C_0(\cM)$-module. Then equation \eqref{eq:locmod1man} has the following implications:
		\begin{enumerate}
			\item[$\star$] $C_0(\cM)=C^\infty(M)$,
			\item[$\star$] $C_l(\cM)_{|U}=(\wedge^l C_1(\cM))_{|U}$,
			\item[$\star$] $C_1(\cM)$ is a locally free   $C^\infty(M)$-module of constant rank, thus it can be identified with the sheaf of sections of some vector bundle $E\to M$, i.e. $\Gamma (E^\ast)\cong C_1(\mathcal{M}) $.
		\end{enumerate}
		So, putting together the above items we get that $E[1]\cong\mathcal{M}$ as desired.
	\end{proof}
	
	\begin{remark}[On supermanifolds II]
		Let  $\mathcal{S}up\mathcal{M}an$ denote the category whose objects are supermanifolds as in Remark \ref{rmk:sup1}. It was shown in \cite{bat:sup} that there is also a functor $$\Pi:\cV ect\to \mathcal{S}up\mathcal{M}an, \quad (E\to M)\mapsto \Pi E=\big(M, \Gamma(\wedge^\bullet E^\ast)\big)$$
		which is faithful and essentially surjective. However, the functor $\Pi$ is far from being full. 
	\end{remark}
	
	The main conclusion of Theorem \ref{thm: 1-manifolds} is that it is enough to study the $1$-manifold $E[1]$ for a given vector bundle $E\to M$. Therefore, in the remainder of this section a generic $1$-manifold will be denoted by  $E[1]$.

	\subsection{Vector fields} Following the dictionary presented in Table \ref{tab:dic}, we say that a \emph{vector field of degree $k$} on the $1$-manifold $E[1]$ is a degree $k$ derivation $X$ of $C_\bullet(E[1])$, i.e.,  an  $\mathbb R$-linear map $X:\Gamma(\wedge^\bullet E^*)\to \Gamma(\wedge^{\bullet+k} E^*)$ with the property that
	\begin{equation}\label{eq:leibvf1}
		X(s_1\wedge s_2)=X(s_1)\wedge s_2+(-1)^{|s_1|k}s_1\wedge X(s_2),\qquad s_i\in\Gamma(\wedge^\bullet E^*).
	\end{equation}
	Vector fields give rise to a $C_\bullet(E[1])$-module. 
	Degree $k$ vector fields are denoted by $\fX^1_k(E[1])$, and all vector fields by $\fX^1_\bullet(E[1])$. The graded commutator of vector fields, defined for homogeneous vector fields $X, Y$ by
	\begin{equation*}
		[X,Y]=XY-(-1)^{|X||Y|}YX,
	\end{equation*}
	makes $\fX^1_{\bullet}(E[1])$ into a graded Lie algebra. Additionally, it is not difficult to show that for $s\in C_\bullet(E[1])$ one has $$[X,sY]=X(s)Y+(-1)^{|X||s|}s[X,Y].$$
	
	Let $U\subseteq M$ be an open chart such that $E_{|U}$ is trivial and consider coordinates $\{ x^i, \xi^{\alpha}\}$  on $E[1]_{|U}$. Define the vector fields  $\frac{\partial}{\partial x^i},\ \frac{\partial}{\partial \xi^{\alpha}}$ on $E[1]_{|U}$ as the derivations acting on coordinates by
	$$\frac{\partial}{\partial x^i}(x^{j})=\delta_{ij}, \quad \frac{\partial}{\partial x^i}(\xi^{\alpha})=0,\quad \frac{\partial}{\partial \xi^{\alpha}}(x^i)=0,\quad \frac{\partial}{\partial \xi^{\alpha}}(\xi^{\beta})=\delta_{\alpha \beta}.$$
	This  definition implies that $\big|\frac{\partial}{\partial x^{i}}\big|=0$ and $\big |\frac{\partial}{\partial \xi^{\alpha}}\big |=-1$.
	
	\begin{proposition}\label{neg-gen-vf}
		Let $U$ be a chart of $E[1]$ with coordinates $\{ x^i, \xi^{\alpha}\}$. Then
		$\fX^1_{\bullet}(E[1])_{|U}$ is freely generated by $\{\frac{\partial}{\partial x^{i}}, \frac{\partial}{\partial \xi^{\alpha}}\}$ as a $C_\bullet({E[1]})_{|U}\text{-module}$.
	\end{proposition}
	\begin{proof} 
		Let  $X$ be a vector field on $E[1]_{|U}$, and consider the vector field
		\begin{equation*}
			X'=X-
			X(x^i)\frac{\partial}{\partial x^{i}}- 
			X(\xi^{\alpha})\frac{\partial}{\partial \xi^{\alpha}}.
		\end{equation*}
	
		It is clear that $X'$ vanishes on the coordinates, hence on arbitrary functions. Thus, $X$ can be expressed as a unique linear combination
		of $\{\frac{\partial}{\partial x^i}, \frac{\partial}{\partial \xi^{\alpha}}\}$ over $C_\bullet(E[1])_{|U}$.
	\end{proof}
	
	A \emph{derivation} of a vector bundle $E\to M$ is a pair $(\cdo,X)$ where $X$ is a vector field on $M$ and  $\cdo: \Gamma(E)\to \Gamma (E)$ is an $\mathbb{R}$-linear operator such that 
	$$\cdo(fe)=X(f)e+f\cdo(e),$$
	for each $f\in C^\infty(M)$ and $e\in \Gamma (E)$. The vector field $X$ is called the \emph{symbol} of $\cdo$ and denoted by $\sigma(\cdo)=X$. One easily shows that derivations form a locally finitely generated and free $C^\infty(M)$-module with the rule $f(\cdo,X)=(f\cdo,fX)$. Therefore, derivations are sections of a vector bundle $\der(E)\to M$
	fitting into the following short exact sequence,
	\begin{equation}\label{es:der}
		0\to \textnormal{End}(E)\xrightarrow{\iota} \textnormal{Der}(E)\xrightarrow[]{\it \sigma} TM
	\end{equation}
	where $\textnormal{End}(E)$ is the vector bundle of endomorphisms and the map $\textnormal{End}(E)\to \textnormal{Der}(E)$ is the inclusion. Also it follows from the definition that $\Gamma\der(E)$ acts on $\Gamma (E^*)$ as
	\begin{equation}\label{eq:deract}
		\cdo(\xi)(e) = \sigma(\cdo)\big(\xi(e)\big)-\xi\big( \cdo(e)\big)\quad \text{with}\quad e\in\Gamma(E),\ \xi\in\Gamma(E^*).
	\end{equation}
	Additionally, there is a canonical Lie algebra structure on $\Gamma \textnormal{Der}(E)$ given by the commutator of derivations $[\cdo,\cdo']_{\der}=\cdo\circ \cdo'-\cdo'\circ \cdo$. 
	
	Recall that $C_\bullet(E[1])$, as an algebra, is generated  by $C_0(E[1])$ and $C_1(E[1])$. Therefore, the Leibniz rule \eqref{eq:leibvf1} implies that vector fields on $E[1]$ are completely determined by their action on functions of degree $0$ and  $1$.  Following \cite[Def. 2.23]{igx:quasi} we say that a pair $(\tau_0, \tau_1)$ is an  \emph{almost $(k+1)$-differential on $E\to M$} if 
	$$\tau_0:C^\infty(M) \to \Gamma( \wedge^k E^\ast)\quad \textnormal{and}\quad \tau_1:\Gamma (E^\ast)\to  \Gamma (\wedge^{k+1}E^\ast)$$
	satisfying the properties 
	$$\tau_0(fg) = \tau_0(f)g+f\tau_0(g) \quad\text{and} \quad\tau_1(f\xi) = \tau_0(f)\wedge \xi+f\tau_1(\xi)$$ for $f,g\in C^\infty(M)$ and $\xi\in\Gamma(E^*)$.
	
	As a direct consequence of the above, we get the following geometric characterization of the vector fields on $E[1]$. 
	\begin{proposition}\label{prop:geovf1}
		Let $E[1]$ be a 1-manifold. The following hold:
		\begin{enumerate}
			\item[1.] There is a bijection between $\fX_k^1(E[1])$ and $\{\text{Almost } (k+1)$-differentials on $ E \}$.
			\item[2.] $\fX^1_{\bullet}(E[1])$ is generated as $C_\bullet(E[1])$-module by $$\Gamma(E)=\fX^1_{-1}(E[1]),\ e\to \iota_e\quad \text{and}\quad \Gamma\der(E^*)=\fX_0^1(E[1]),\ (\cdo,X)\to \widehat{\cdo}.$$ 
			\item[3.] The module structure $C_0(E[1])\cdot\fX^1_{i}(E[1])\subset\fX^1_i(E[1])$ coincides with the one given as sections of the corresponding vector bundle.
			\item[4.] The module structure $C_1(E[1])\cdot\fX^1_{-1}(E[1])\subset\fX^1_0(E[1])$ is given by the map $\iota$ in \eqref{es:der}.
			\item[5.] The bracket $[\cdot,\cdot]:\fX^1_0(E[1])\times\fX_{-1}^1(E[1])\to \fX_{-1}^1(E[1])$ is given by \eqref{eq:deract}.
			\item[6.] The bracket $[\cdot,\cdot]:\fX^1_0(E[1])\times\fX_{0}^1(E[1])\to \fX_{0}^1(E[1])$ is given by the commutator of derivations.
		\end{enumerate}
	\end{proposition}
	\begin{proof}
		Item 1. follows immediately by the Leibniz rule and the fact that $C_\bullet(E[1])$ is generated as an algebra by $C_0(E[1])=C^\infty(M)$ and $C_1(E[1])=\Gamma(E^*)$. For item 2., we now specify the action of vector fields on generators. Given $e\in\Gamma(E)$ and $\cdo\in\Gamma\der(E^*)$ we define 
		$$\iota_e(f)=0,\quad \iota_e(\xi)=\xi(e),\quad \widehat{\cdo}(f)=\sigma(\cdo)(f)\quad \text{and}\quad \widehat{\cdo}(\xi)=\cdo(\xi)$$
		for $f\in C^\infty(M)=C_0(E[1])$ and $\xi\in\Gamma(E^*)=C_1(E[1]).$
		The remaining assertions follow from standard arguments and are left to the reader.
	\end{proof}
	It follows from item $4.$ in the above proposition that the space of vector fields $\fX^1_{\bullet}(E[1])$ is not necessarily freely generated by the spaces listed in item $2.$
	
	\subsection{Q-structures and Lie algebroids}\label{sec:Q1}
	Among all the vector fields on $E[1]$ we now focus on a special class of them. A \emph{$Q$-structure on $E[1]$} is a degree  $1$ vector field $Q\in\fX^1_1(E[1])$ satisfying $$Q^2=\frac{1}{2}[Q,Q]=0.$$
	The pair $(E[1], Q)$ is known as a \emph{$Q$-manifold}.
	
	An immediate consequence of the definition is that the functions on $E[1]$ form a cochain complex  $$C_\bullet(E[1])\xrightarrow{Q}C_{\bullet+1}(E[1]).$$ Its cohomology will be denoted by $H^\bullet_Q(E[1])$.  Therefore, when dealing with $Q$-manifolds we treat homological algebra as differential geometry. 
	
	\begin{remark}
		The equation $\frac{1}{2}[Q,Q]=0$ says that the distribution on $E[1]$ defined by the vector field $Q\in\fX^1_1(E[1])$ is integrable. To rigorously formulate this observation one needs supermanifolds, or to introduce derived manifolds as in \cite{blx}.  In those settings, one can show that $Q$ is the infinitesimal generator of an odd line action, see e.g. \cite[\S 3.2]{sev:gorm} or \cite[Rmk. 3.8]{cat:intro}.   
	\end{remark}
	
	Let us use Proposition \ref{prop:geovf1} to give a more geometric interpretation of $Q$-structures on $E[1].$ In particular, we know that $\fX^1_{-1}(E[1])=\Gamma(E)$. Thus given a $Q\in \fX^1_{1}(E[1])$ we can define the following operations
	\begin{equation}\label{eq:liealg}
		[e,e']_E:=[\iota_e,[\iota_{e'},Q]]\in\Gamma(E)\quad\text{and}\quad \rho(e)(f):=[\iota_e,Q](f) \in C^\infty(M)
	\end{equation}
	for $e,e'\in\Gamma(E)=\fX^1_{-1}(E[1])$ and $f\in C^\infty(M)=C_0(E[1])$, where the bracket on the right is the Lie bracket of vector fields on $\fX^1_\bullet(E[1])$.

	\begin{theorem}[\cite{vai:lie}]\label{thm:Q-lie}
		There is a one-to-one correspondence between  $Q\in\fX_1^1(E[1])$ and pairs $([\cdot,\cdot]_E,\rho)$ where $\rho:E\to TM$ is a vector bundle map and $[\cdot,\cdot]_E:\wedge^2\Gamma (E)\to\Gamma (E)$ satisfying the Leibniz rule
		$$[e,fe']_E=f[e,e']_E+\rho(e)(f)e'\qquad \text{with } \,   e,e'\in\Gamma(E)\text{ and } f\in C^\infty(M).$$
		Moreover $Q^2=0$ if and only if $[e,[e', e'']_E]_E=[[e,e']_E,e'']_E+[e', [e,e'']_E]_E$ for $e, e',e''\in\Gamma(E).$
	\end{theorem}
	\begin{proof}
		Given $Q\in\fX^1_1(E[1])$ define $([\cdot,\cdot]_E,\rho)$ by \eqref{eq:liealg}. Observe that $\rho$ is a vector bundle map because
		$$\rho(ge)(f)=[\iota_{ge},Q](f)=\iota_{ge}(Q(f))+Q(\iota_{ge}(f))=\iota_{ge}(Q(f))=g[\iota_e,Q](f)=g\rho(e)(f),$$
		for $e\in\Gamma(E), f,g\in C^\infty(M)$. The skew-symmetric property of the bracket follows from $\fX_{-2}^1(E[1])=0$, indeed
		$$[e,e']_E=[\iota_e, [\iota_{e'},Q]]=[[\iota_e,\iota_{e'}],Q]-[\iota_{e'},[\iota_e,Q]]=-[\iota_{e'},[\iota_e,Q]]=-[e,e']_E \quad\text{for } e,e'\in\Gamma(E).$$
		Finally, the Leibniz rule is due to the graded Jacobi identity of the Lie bracket of vector fields, more concretely
		\begin{equation}\label{eq:lialgdiff}
			\begin{split}
				[e,fe']_E=&[\iota_e,[\iota_{fe'},Q]]=[[\iota_e,\iota_{fe'}],Q]-[\iota_{fe'},[\iota_e,Q]]=-f[\iota_{e'},[\iota_e,Q]]+[\iota_e,Q](f)e'\\
				=&f[e,e']_E+\rho(e)(f)e'.
			\end{split}
		\end{equation}
		Conversely, given $([\cdot,\cdot],\rho)$ as in the statement we define $Q\in\fX^1_1(E[1])$ by the formula 
		\begin{eqnarray}\label{eq:Qalg}
			Q (s)(e_0,\cdots,e_k) & = & \sum_{i=0}^k(-1)^i\rho(e_i)\Big(s(e_0,\cdots,\hat{e}_i,\cdots,e_k)\Big)\\
			&  &+ \sum_{i<j}(-1)^{i+j}s([e_i,e_j]_E,e_0,\cdots,\hat{e}_i,\cdots,\hat{e}_j,\cdots,e_k)  \notag
		\end{eqnarray}
		for $s\in C_k(E[1])$ and $e_i\in\Gamma(E).$ Observe that $Q(s)$ is skew-symmetric in each entry because so is the bracket and it is linear in the last one, and therefore in all of them, due to the fact that the bracket and the anchor satisfy the Leibniz rule. In summary, $Q$ is well defined. It remains to show that the above formula defines a derivation, i.e.  $$Q(s_1\wedge s_2)=Q(s_1)\wedge s_2+(-1)^{|s_1|}s_1\wedge Q(s_2)$$ for $s_i\in C_\bullet(E[1])$. This follows from the fact that $\rho(e)$ is a vector field, and we leave it to the interested reader. We point out that it is enough to prove it for $s_i\in C_\bullet(E[1])$ with $|s_i|\leq 1.$
		
		It is not difficult to check that formulas \eqref{eq:liealg} and \eqref{eq:Qalg} are inverse of each other, and therefore the above is a one-to-one correspondence.

		Finally we show the moreover part. Given a pair $([\cdot,\cdot]_E,\rho)$ we introduce the bundle morphisms $\lambda:\wedge^2 E\to TM$ and $J:\wedge^3 E\to E$ given by
		\begin{equation}\label{eq:jacobiator}
			\lambda(e, e')=\rho([e,e']_E)-[\rho(e),\rho(e')],\quad J(e, e', e'')=[[e, e']_E, e'']_E-[[e, e'']_E, e']_E+[[e',e'']_E,e]_E
		\end{equation}
		and observe that if $([\cdot,\cdot]_E,\rho)$ satisfies the Leibniz rule then
		$$J(e,e',fe'')-fJ(e,e',e'')=\lambda(e,e')(f)e''.$$
		Therefore, the vanishing of $J$ implies that $\lambda$ also vanishes. To conclude, observe that if $Q^2$ vanishes in the generators it will vanish everywhere. So formula \eqref{eq:Qalg} applied to $f\in C^\infty(M)$ and $\xi\in\Gamma(E^*)$ gives
		\begin{equation*}
			\begin{split}
				Q^2(f)(e,e')=&-\lambda(e,e')(f),\\
				Q^2(\xi)(e,e',e'')=&\xi\big(J(e,e',e'')\big)-\lambda(e,e')(\xi(e''))+\lambda(e, e'')(\xi(e'))-\lambda(e',e'')(\xi(e)),\\
			\end{split}
		\end{equation*}
		for $e,e',e''\in\Gamma(E)$. Hence $Q^2=0$ if and only if $J=0$ as we wanted.
	\end{proof}
	
	Putting together all the conditions in the statement of the previous theorem we get the following definition. A \emph{Lie algebroid $(E\to M, \rho, [\cdot,\cdot]_E)$} is a vector bundle together with a vector bundle map $\rho:E\to TM$, covering the identity, and a Lie algebra structure on its sections satisfying the Leibniz rule
	$$[e,fe']_E=f[e,e']_E+\rho(e)(f)e'\quad \text{for } e,e'\in\Gamma(E),\ f\in C^\infty(M).$$
	\begin{remark}
		One can introduce Lie algebroid and $Q$-manifold morphisms. By defining morphisms appropriately one finds that the functor $[1]$ introduced in Theorem \ref{thm: 1-manifolds} can be promoted to an equivalence of categories between Lie algebroids and degree $1$ $Q$-manifolds \cite{vai:lie}.
	\end{remark}
	
	Let $E[1]$ be a $1$-manifold and choose $U\subseteq M$ a chart and coordinates $\{x^i, \xi^\alpha\}$. By Proposition \ref{neg-gen-vf} we get that $Q\in\fX^1_1(E[1]_{|U})$ is  given by
	$$Q=\rho^i_\alpha(x)\xi^\alpha\frac{\partial}{\partial x^i}+\frac{1}{2}c^\gamma_{\alpha\beta}(x)\xi^\alpha\xi^\beta\frac{\partial}{\partial\xi^\gamma}$$
	for some functions $\rho^i_\alpha(x),c^\gamma_{\alpha\beta}(x)\in C^\infty(U).$ This local expression in coordinates is sometimes useful for computations.
	
	In what follows, we give some fundamental examples of $Q$-structures, i.e. Lie algebroids.
	
	\begin{example}\label{ex:Q-TM-g}
		In a way, the two extreme cases of $Q$-manifolds are given by:
		\begin{itemize}
			\item Tangent bundles $T[1]M=(M, \Omega^\bullet(M))$ with the de Rham differential  $$d:\Omega^\bullet(M)\to \Omega^{\bullet+1}(M).$$
			If we pick a chart $U\subseteq M$ with local coordinates $\{x^i\}$ then  $$\{x^i, \theta^i=dx^i\}\quad \text{with} \quad |x^i|=0\text{ and } |\theta^i|=1$$ are local coordinates on $T[1]M$. Therefore, the vector field given by the de Rham differential on $T[1]M_{|U}$ has the following expression
			$$d=\theta^i\frac{\partial}{\partial x^i}.$$
			\item Lie algebras $\g[1]=(pt, \wedge^\bullet\g^*)$ with the Chevalley-Eilenberg differential $$d_{CE}:\wedge^\bullet\g^*\to \wedge^{\bullet+1}\g^*.$$
			Pick $\{e_\alpha\}$ a basis of $\g$ and denote by $\{\xi^\alpha\}$ its dual basis. Since $\g[1]$ is a $1$-manifold then $\{\xi^\alpha\}$ are coordinates on $\g[1]$ and the Chevalley-Eilenberg differential is 
			$$d_{CE}= \frac{1}{2}c_{\alpha\beta}^\gamma \xi^\alpha\xi^\beta\frac{\partial}{\partial \xi^\gamma}$$
			where $c_{\alpha\beta}^\gamma$ are the structural constants of the Lie algebra $\g$ with respect to the basis $\{e_\alpha\}$, i.e.  they are given by the formula $[e_\alpha,e_\beta]=c_{\alpha\beta}^\gamma e_\gamma$.
		\end{itemize}
	
		A normal form for vector fields on graded manifolds shows that any Lie algebroid is locally the product of these two \cite{vai:vec}, yet the transverse direction is not simply a Lie algebra. This result is also known as the splitting theorem for Lie algebroids, see e.g. \cite{duf:nor}. 
	\end{example}
	
	\begin{example}[Poisson manifolds]\label{ex:poiLA}
		A manifold $M$ together with a bivector $\pi\in\fX^2(M)$ satisfying $[\pi, \pi]_S=0$, where $[\cdot,\cdot]_S:\fX^i(M)\times\fX^j(M)\to\fX^{i+j-1}(M)$ denotes the Schouten bracket, is called a \emph{Poisson manifold}. Then $T^*[1]M=(M,\fX^\bullet(M))$ carries a natural $Q$-structure defined by 
		$$Q_\pi=[\pi,\cdot]_S:\fX^\bullet(M)\to\fX^{\bullet+1}(M).$$
	\end{example}
	
	\begin{example}[Central extensions]
		Let $(E[1], Q)$ be a $Q$-manifold and $r\in C_2(E[1])=\Gamma(\wedge^2E^*)$ with $Q(r)=0$. Then, on the $1$-manifold  $E[1]\times\bR[1]$, see Example \ref{car-pro}, there is a new $Q$-structure given by
		$$\widehat{Q}=Q+r\frac{\partial}{\partial\theta}$$
		where $\theta\in C_1(\bR[1]).$ Moreover, if $r\in C_2(E[1])$ is $Q$-exact then $\widehat{Q}$ is equivalent to the trivial extension. In summary, we claim that $H^2_Q(E[1])$ classifies \emph{central extensions}. The name is due to the fact that, for the induced Lie algebroid structure $(E\oplus\bR_M\to M, \rho, [\cdot,\cdot]_{E\oplus\bR_M})$, the sections of $\bR_M$ are in the kernel of the anchor and central for the bracket.  
		
		Some of the main examples of this construction are:
		\begin{itemize}
			\item Central extensions of Lie algebras. See \cite[\S~3.4]{eti:qua} for the relation with triangular Lie bialgebras.
			\item Given $\omega\in\Omega^2(M)$ a closed $2$-form, then $T[1]M\times\bR[1]$ gives the prequantum Lie algebroid of a closed $2$-form \cite{cra:pre}.
			\item Given a Poisson manifold $(M,\pi)$, then $T^*[1]M\times \bR[1]$ gives the Jacobi Lie algebroid \cite{lich:jac}.
		\end{itemize}
		We will extend this example in Proposition \ref{prop:prequantum}.
	\end{example}
	
	\begin{example}[Derivations]
		Let $E\to M$ be a vector bundle. As shown in the previous section $\textnormal{Der}(E)$ becomes a Lie algebroid with Lie bracket on sections determined by the commutator of derivations and the anchor given by symbol map \eqref{es:der}.  In the literature this is also referred to as the \emph{Atiyah algebroid} of $E\to M$, see e.g. \cite{kosm:der}.
	\end{example}

	\subsubsection{Cartan calculus for Lie algebroids} 
	Note that our approach to Lie algebroids is non-standard. Indeed, we made a long detour to define a bracket and an anchor on a vector bundle. Nevertheless, from our point of view we automatically had Lie algebroid cohomology built in. In fact, this approach provides a full Cartan calculus as follows. 
	
	\begin{proposition}\label{prop:cartcal}
		Given a $Q$-manifold $(E[1], Q)$ and $e\in\Gamma(E)$ the operations
		$$\iota_e\in\fX_{-1}^1(E[1]),\qquad \cL_e=[\iota_e, Q]\in\fX^1_0(E[1])\quad \text{and}\quad  Q\in\fX_1^1(E[1])$$
		act on $\Gamma(\wedge^\bullet E^*)=C_\bullet(E[1])$ by derivations and satisfy the Cartan relations
		\begin{equation*}
			\begin{array}{rclrclrcl}
				[\iota_e, Q]&=&\cL_e,&  [Q,\cL_e]&=&0, &[\cL_e, \cL_{e'}]&=&\cL_{[e,e']_E}, \\
				\left[\iota_e,\iota_{e'}\right]&=&0,& [\cL_e,\iota_{e'}]&=&\iota_{\left[e,e'\right]_E},& [Q,Q]&=&0. 
			\end{array}
		\end{equation*}
	\end{proposition}
	The proof follows directly from the properties of the Lie bracket of vector fields. We will see an application of the Cartan calculus in $\S \ref{sec:equico}$ when dealing with equivariant cohomology.

	\subsubsection{Deformations of Lie algebroids}
	Another advantage of our point of view on Lie algebroids is that $Q$-structures on $E[1]$ have an easy deformation theory. A general principle in deformation theory states that any deformation problem is controlled by a \emph{differential graded Lie algebra}, DGLA for short, $(\g_\bullet, [\cdot,\cdot], d)$ see e.g. \cite{kon:defor, man:def}. More precisely, this means that deformations correspond to \emph{Maurer-Cartan elements}, i.e. elements $$v\in\g_1\quad \text{satisfying}\quad dv+\frac{1}{2}[v,v]=0,$$ modulo \emph{gauge transformations}\footnote{Two Maurer-Cartan elements $v_0, v_1$ in a DGLA are called \emph{gauge equivalent}, if
		there is a curve $v_t$ of Maurer-Cartan elements starting at $v_0$ and ending at $v_1$ together with a curve of degree
		$0$ elements $\lambda_t$, such that $\frac{d}{dt}v_t=d\lambda_t+[v_t,\lambda_t]$, for more detail see \cite[\S 5]{jonas:intro}.}.
	
	\begin{theorem}[\cite{cra:def}]
		Let $E[1]$ be a $1$-manifold with a $Q$-structure given by $Q\in\fX^1_1(E[1])$. Then $(\fX_\bullet^1(E[1]),[\cdot,\cdot], d=[Q,\cdot])$ is the DGLA controlling deformations of $Q$.
	\end{theorem}
	
	\begin{proof}
		Observe that  $(\fX_\bullet^1(E[1]),[\cdot,\cdot], d=[Q,\cdot])$ is a DGLA due to the Jacobi identity for $[\cdot,\cdot].$ Given a smooth path $Q_t\in\fX^1_1(E[1])$ starting at $0$ then the perturbed vector field $Q+Q_t$ defines a new $Q$-structure if and only if 
		$$0=\frac{1}{2}[Q+Q_t, Q+Q_t]=[Q,Q_t]+\frac{1}{2}[Q_t,Q_t]=d(Q_t)+\frac{1}{2}[Q_t,Q_t].$$
		This means $Q_t$ defines a path of Maurer-Cartan elements in our DGLA. Finally we refer to \cite{cra:def} for the fact that gauge transformations correspond to trivial deformations.
	\end{proof}

	\subsection{Poisson algebras and Lie bialgebroids}\label{sec:bia1}
	After the study of vector fields on the $1$-manifold $E[1]$  we focus on Poisson structures. An operation $$\{\cdot,\cdot\}:C_i(E[1])\times C_j(E[1])\to C_{i+j+k}(E[1])$$ is called a \emph{degree $k$ Poisson structure} if  for $s_i\in C_\bullet(E[1])=\Gamma(\wedge^\bullet E^*)$, we have
	\begin{enumerate}
		\item [(P1)] $\{s_1, s_2\}=-(-1)^{(|s_1|+k)(|s_2|+k)}\{s_2, s_1\}$,
		\item [(P2)] $\{s_1, s_2\wedge s_3\}=\{s_1, s_2\}\wedge s_3+(-1)^{(|s_1|+k)|s_2|}s_2\wedge\{s_1, s_3\}$,
		\item [(P3)] $\{s_1,\{s_2,s_3\}\}=\{\{s_1,s_2\},s_3\}+(-1)^{(|s_1|+k)(|s_2|+k)}\{s_2,\{s_1,s_3\}\}.$
	\end{enumerate}
	Moreover, we say that $X\in\fX^1_i(E[1])$ is a \emph{Poisson vector field} if 
	\begin{equation}\label{eq:vfPoi1}
		X(\{s_1,s_2\})=\{X(s_1), s_2\}+(-1)^{(|s_1|+k)i}\{s_1, X(s_2)\}.
	\end{equation}
	A  \emph{degree $k$  PQ-manifold} is $(E[1], Q, \{\cdot,\cdot\})$ where $\{\cdot,\cdot\}$ is a Poisson bracket of degree $-k$ and $Q\in\fX^1_1(E[1])$ is a  $Q$-structure that in addition is a Poisson vector field. Finally, given $s\in C_i(E[1])$ its \emph{hamiltonian vector field} is $X_s=\{s,\cdot\}\in\fX^1_{i+k}(E[1]).$

	Recall that, as an algebra, $C_\bullet(E[1])$ is generated by $C_0(E[1])=C^\infty(M)$ and $C_1(E[1])=\Gamma(E^*)$, see Theorem \ref{thm: 1-manifolds}. Hence, the Leibniz rule (P2) immediately implies that  $\{\cdot,\cdot\}\equiv 0$ for  $k<  -2$. Let us spell out the first two cases.
	
	\begin{proposition}[\cite{gruz:H-twist, ikeda:qp3}]
		A degree $-2$ Poisson structure on $E[1]$ is equivalent to a symmetric pairing $B(\cdot,\cdot):\Gamma(E^*)\times\Gamma(E^*)\to C^\infty(M)$. Moreover, if $([\cdot,\cdot]_E, \rho)$  defines a Lie algebroid then  equation  \eqref{eq:vfPoi1} is equivalent to 
		\begin{equation}\label{-2degreeEq}
		B(\rho^*(df), \xi)=0\quad \text{and}\quad \rho(e)\big(B(\xi, \xi')\big)=B(\cL_e \xi, \xi')+B(\xi,\cL_e\xi')
		\end{equation}
		where $f\in C^\infty(M), \ \xi,\xi'\in\Gamma(E^*)$ and $e\in\Gamma (E)$ with $\cL_e$ given as in Proposition \ref{prop:cartcal}. 
	\end{proposition}
	\begin{proof}
		Since $\{\cdot,\cdot\}$ has degree $-2$ then $(P2)$ implies that the bracket is completely determined by $$\{\cdot,\cdot\}=B(\cdot,\cdot):C_1(E[1])\times C_1(E[1])\to C_0(E[1])$$ that is symmetric by $(P1)$. Observe that, by degree reasons, $(P3)$ is automatically satisfied. The converse is also immediate. 
		
		Let us check the moreover part of the statement. From Theorem \ref{thm:Q-lie} one knows that a $Q$-structure is the same as a Lie algebroid $([\cdot,\cdot]_E, \rho)$. On the one hand, if we plug \[
			s_1=f\in C^\infty(M)=C_0(E[1]),  \quad \text{and} \quad s_2=\xi\in\Gamma(E^*)=C_1(E[1])
		\]
		in Equation \eqref{eq:vfPoi1} then we immediately get that $B(\rho^*(df), \xi)=0$. On the other hand, observe that for any $e\in\Gamma(E)$ we obtain another Poisson vector field  $\iota_e$, so that 
		\begin{equation*}
			\begin{array}{rl}
				\rho(e)\big(B(\xi,\xi')\big)= [\iota_e, Q](\{\xi,\xi'\})=&\{[\iota_e,Q](\xi),\xi'\}+\{\xi,[\iota_e,Q](\xi')\}  \\
				=& B(\cL_e \xi, \xi')+B(\xi,\cL_e\xi'),
			\end{array}
		\end{equation*}
		for $\xi,\xi'\in\Gamma(E^*)=C_1(E[1])$. Similarly, one can use Equations \eqref{-2degreeEq} to show that $([\cdot,\cdot]_E, \rho)$ defines a Lie algebroid.
	\end{proof}

	\begin{proposition}\label{prop:P-liea}
		A degree $-1$ Poisson structure on $E[1]$ is equivalent to a Lie algebroid structure on $E^*\to M.$
	\end{proposition}
	\begin{proof}
		Given a degree $-1$ Poisson structure define 
		$$\rho(\xi)(f)=\{\xi,f\}\quad\text{and}\quad [\xi,
		\xi']_{E^*}=\{\xi,\xi'\}$$
		for $\xi,\xi'\in\Gamma(E^*)=C_1(E[1])$ and $f\in C^\infty(M)=C_0(E[1])$. One can check that indeed equations $(P1),(P2)$ and $(P3)$ imply the Lie algebroid equations. Conversely, given a Lie algebroid $([\cdot,\cdot]_{E^*},\rho),$ one uses the above formulas to define the degree $-1$ Poisson bracket on generators and extend them using $(P2)$  to all the algebra $C_\bullet(E[1])$, that is well-defined due to the Leibniz rule between the bracket and the anchor. The skew-symmetric property of the bracket guarantees $(P1)$, while $(P3)$ follows from the Jacobi identity and the fact that the anchor is bracket preserving.
	\end{proof}
	
	In the literature a  graded algebra with a Poisson bracket of degree $-1$ is sometimes called a \emph{Gerstenhaber algebra}. We also point that one can adopt the Poisson point of view for Lie algebroids and make Lie algebroids into a category with morphisms given by Poisson maps, which leads to \emph{Lie algebroid comorphisms} distinct from standard Lie algebroid morphisms, see e.g. \cite{cat:intpoi} or \cite[\S 4.6]{mac:book} for more details.
	
	\begin{remark}
		If we put together Theorem \ref{thm:Q-lie} with Proposition \ref{prop:P-liea} and include the well known fact that a vector bundle carries a Lie algebroid structure if and only if on its dual  there is a linear Poisson structure, then one get the following result. Given a vector bundle $E\to M$ there is a one-to-one correspondence between the following structures:
		\begin{enumerate}
			\item[1.] Lie algebroid structures on $E\to M$.
			\item[2.] $Q$-structures on $E[1]$.
			\item[3.] Degree $-1$ Poisson structures on $E^*[1].$
			\item[4.] Linear Poisson structures on $E^*\to M.$
		\end{enumerate}
	\end{remark}

	A \emph{Lie bialgebroid} is a Lie algebroid $(E\to M, \rho_E, [\cdot,\cdot]_E)$ whose dual $(E^\ast \to M, \rho_{E^*}, [\cdot,\cdot]_{E^*})$ is also a Lie algebroid and the differential $Q$ of $E$ is a derivation of the Gerstenhaber bracket on $\Gamma( \wedge^\bullet E^\ast)$. That is, the latter satisfies
	\begin{equation*}
		Q \{s,s'\}=\{Q(s),s'\}+(-1)^{\vert s\vert -1}\{s,Q(s')\},\qquad s,s'\in \Gamma (\wedge^\bullet E^\ast).
	\end{equation*}

	Directly from the definition we get the following equivalence.
	
	\begin{corollary}\label{cor:bia}
		There is a one to one correspondence between Lie bialgebroids and degree $1$ PQ-manifolds.
	\end{corollary}
	
	The two main examples of Lie bialgebroids are the following ones:
	
	\begin{example}
		A \emph{Poisson Lie group} $(G,\pi)$ is a Lie group endowed with a bivector $\pi\in\fX^2(G)$ satisfying
		$$[\pi,\pi]_S=0\quad\text{and}\quad \pi_{gh}=(L_g)_*\pi_h+(R_h)_*\pi_g\quad\text{for}\quad g,h\in G.$$
		Note that the latter expression immediately implies that $\pi_e=0$. Therefore, if we see the Poisson structure as a map $\pi:G\to TG\wedge TG$, we get that $T_e\pi:\g\to \g\wedge\g$ makes  $(\g, [\cdot,\cdot])$ into a Lie bialgebra. Poisson Lie groups play a fundamental role in the theory of quantum groups \cite{eti:qua}, cluster algebras \cite{fg:cluster}, meromorphic connections \cite{boa:sto} or field theories \cite{sev:plt}, among others.
	\end{example}
	
	\begin{example}\label{ex:Poibia}
		Let $(M,\pi)$ be a Poisson manifold. Example \ref{ex:poiLA} shows that $(T^*[1]M,Q_\pi)$ is a $Q$-manifold, i.e. $(T^*M,[\cdot,\cdot]_\pi, \pi^\sharp)$ is a Lie algebroid and therefore $(T[1]M, \{\cdot,\cdot\}_\pi)$ is a Poisson manifold by Proposition \ref{prop:P-liea}. 
		\begin{itemize}
			\item On the one hand, the Schouten bracket on multivector fields $[\cdot,\cdot]_S:\fX^i(M)\times \fX^j(M)\to \fX^{i+j-1}(M)$ is a Poisson bracket of degree $-1$ that makes $(T^*[1]M,Q_\pi, [\cdot,\cdot]_S)$ into a degree $1$ $PQ$-manifold.
			\item  On the other hand, $(T[1]M, Q_{dR}=d, \{\cdot,\cdot\}_\pi)$ is also a degree $1$ $PQ$-manifold where the $Q$-structure is given by Example \ref{ex:Q-TM-g}.
		\end{itemize}
	\end{example}
	
	The duality shown in the previous example is a general fact about degree $1$ $PQ$-manifolds,  and it will be explained in \S~\ref{sec:liebia}.
	
	\begin{remark}
		The homotopic version of degree $0$ Poisson structures on $1$-manifolds will appear naturally in   \S~\ref{sec:1sym=poi} while considering deformations of coisotropic submanifolds of Poisson manifolds, see also  \S~\ref{sec:deflag}. 
	\end{remark}
	
	\subsection{Superdivergence and Lie algebroid homology}\label{sec:div}
	Recall that for a smooth manifold, a divergence operator is a rule that assigns to each vector field a divergence, and is compatible with the Lie bracket of vector fields. In the same way, let $E[1]$ be a $1$-manifold. A \emph{superdivergence} is a map $\sdiv:\fX^1_\bullet(E[1])\to C_\bullet(E[1])$  satisfying the following identities:
	\begin{enumerate}
		\item [(D1)] $\sdiv(sX)=s\sdiv(X)+(-1)^{|s||X|}X(s)$,
		\item [(D2)] $\sdiv([X,Y])=X(\sdiv(Y))-(-1)^{|X||Y|}Y(\sdiv(X))$,
	\end{enumerate}
	for $s\in C_\bullet(E[1])$ and $X,Y\in\fX^1_\bullet(E[1])$.

	\begin{proposition}
		A superdivergence $\sdiv$ on the $1$-manifold $E^*[1]$ is equivalent to a morphism $\widetilde{D}:\Gamma\der(E)\to C^\infty(M)$ such that $\widetilde{D}_{|\ider (E)}=-tr$, the usual trace of endomorphisms,
		$$\widetilde{D} (f\cdo)=f\widetilde{D}(\cdo)+\sigma(\cdo)(f)\quad\text{and}\quad \widetilde{D}([\cdo,\cdo']_{\der})=\sigma(\cdo)\big(\widetilde{D}(\cdo')\big)-\sigma(\cdo')\big(\widetilde{D}(\cdo)\big)$$
		for $f\in C^\infty(M)$ and $\cdo,\cdo'\in\Gamma\der(E).$ 
	\end{proposition}
	\begin{proof}
		By Proposition \ref{neg-gen-vf} we get that $\fX^1_\bullet(E^*[1])$ is generated as a $C_\bullet(E^*[1])$-module by $\fX^1_{-1}(E^*[1])$ and by $\fX_0^1(E^*[1])$. Since $\sdiv:\fX^1_\bullet(E^*[1])\to C_\bullet(E^*[1])$ is degree preserving and $C_{-1}(E^*[1])=0$, equation  (D1) implies that $\sdiv$ is completely determined by the map $$\sdiv_0=\widetilde{D}:\fX_0^1(E^*[1])=\Gamma\der(E)\to C^\infty(M)=C_0(E^*[1]).$$
		Clearly equations (D1) and (D2) imply the equations in the statement, while the fact that $\widetilde{D}_{|\ider(E)}=-tr$ is because given $\xi\in\Gamma(E^*)=\fX_{-1}^1(E^*[1])$ and $e\in\Gamma(E)=C_1(E^*[1])$ then
		$$\widetilde{D}_{|\ider(E)}(e\otimes \xi)=\sdiv(e\iota_\xi)=-\iota_\xi(e)=-tr(e\otimes \xi).$$
		Conversely, given $\widetilde{D}$ as in the statement we define the superdivergence on generators $$\sdiv(X)=0\quad\text{and}\quad \sdiv(Y)=\widetilde{D}(Y)$$
		for $X\in\fX_{-1}^1(E^*[1]), \ Y\in\fX^1_0(E^*[1])=\Gamma\der(E)$ and extend it to any degree using (D1).
	\end{proof}
	
	Based on the above result, we define  a \emph{divergence for the Lie algebroid $(E\to M, [\cdot,\cdot]_E, \rho)$} as a map $D:\Gamma E\to C^\infty(M)$ satisfying 
	$$D(fe)=fD(e)+\rho(e)(f)\quad\text{and}\quad D([e,e']_E)=\rho(e)\big(D(e')\big)-\rho(e')\big(D(e)\big)$$
	for $f\in C^\infty(M)$ and $e,e'\in\Gamma(E).$
	
	In order to show why we need superdivergences we first need to give one more definition.  Let $(E^*[1], \{\cdot,\cdot\})$ be a degree $-1$ Poisson structure. A \emph{generating operator} is a map $\partial:C_i(E^*[1])\to C_{i-1}(E^*[1])$ satisfying
	\begin{equation}\label{eq:genop}
		\partial^2=0\quad \text{and}\quad 
		\{s,s'\}=(-1)^{|s|}\Big(\partial(s\wedge s')-\partial(s)\wedge s'-(-1)^{|s|}s\wedge\partial(s')\Big)
	\end{equation}
	for $s,s'\in C_\bullet(E^*[1]).$
	\begin{remark}
		Observe that the first equation in \eqref{eq:genop} implies that $(C_\bullet(E^*[1]), \partial)$ is a chain complex. Its associated homology is denoted by $H^\partial_\bullet(E).$ 
		While the second equation in \eqref{eq:genop} implies that $\partial$ is not a vector field but a second order differential operator.  In the literature, generating operators appear under many different names, including: Batalin-Vilkovisky operators, BV-Laplacians, odd Laplacians, homology operators, right superconnections, among others. 
	\end{remark}
	
	\begin{proposition}[\cite{mon:div}]\label{prop:genop}
		Let $(E^*[1],\{\cdot,\cdot\})$ be a degree $-1$ Poisson structure and $$\sdiv:\fX_\bullet^1(E^*[1])\to C_\bullet(E^*[1])$$ be a superdivergence. Then, the following is a generating operator
		\begin{equation*}
			\begin{array}{cccl}
				\partial: & C_i(E^*[1])&\to& C_{i-1}(E^*[1]) \\
				& s &\rightsquigarrow &  \partial(s):=(-1)^{|s|}\frac{1}{2}\sdiv(\{s,\cdot\}).
			\end{array}
		\end{equation*}
		Moreover, $\partial_1:\Gamma(E)\to C^\infty(M)$ is a divergence for the Lie algebroid given by $\{\cdot,\cdot\}$.
	\end{proposition}
	\begin{proof}
		Let $s,s'\in C_\bullet(E^*[1]).$ Using the identity $X_{s\wedge s'}=sX_{s'}+(-1)^{|s||s'|}s'X_{s}$ one shows that
		\begin{equation*}
			\begin{array}{rl}
				(-1)^{|s|}\partial(s\wedge s') =&(-1)^{|s'|}\dfrac{1}{2}\sdiv\Big(sX_{s'}+(-1)^{|s||s'|}s'X_{s}\Big)\\
				=&  \dfrac{(-1)^{|s'|}}{2}\Big(s\wedge \sdiv(X_{s'})+(-1)^{|s||s'|-|s|}\{s',s\}+ (-1)^{|s||s'|} s'\wedge\sdiv(X_s)\\
				&\qquad\qquad+(-1)^{|s'|}\{s,s'\}\Big)\\
				=&\dfrac{(-1)^{|s'|}}{2}s\wedge \sdiv(X_{s'})+\dfrac{1}{2}\sdiv(X_s)\wedge s'+\{s,s'\}\\
				=&s\wedge\partial(s')+(-1)^{|s|}\partial(s)\wedge s' +\{s,s'\}
			\end{array}
		\end{equation*}
		So we got the second equation in \eqref{eq:genop}. For the first one, observe that by the properties of a superdivergence we get 
		\begin{equation}\label{dereq}
			\begin{array}{rl}\vspace{1mm}
				\partial(\{s,s'\})=& \dfrac{(-1)^{|s|+|s'|-1}}{2}\sdiv(X_{\{s,s'\}})\\ =&\dfrac{(-1)^{|s|+|s'|-1}}{2}\sdiv([X_s,X_{s'}])  \\
				=& \dfrac{(-1)^{|s|+|s'|-1}}{2}\Big(X_s\sdiv(X_{s'})-(-1)^{(|s|-1)(s'-1)}X_{s'}\sdiv(X_s)\Big)\\
				=&\{\partial(s),s'\}+(-1)^{|s|-1}\{s,\partial(s')\}.
			\end{array}
		\end{equation}
		So, using the previous identity, we obtain
		\begin{equation*}
			\begin{array}{rl}
				\partial^2(s\wedge s')=& (-1)^{|s|}\partial\big(s\wedge\partial(s')+(-1)^{|s|}\partial(s)\wedge s'+\{s,s'\}\big)   \\
				=& s\wedge \partial^2(s')+(-1)^{|s|}\partial(s)\wedge\partial(s')+\{s,\partial(s')\}\\
				&+(-1)^{|s|-1}\Big(\partial(s)\wedge\partial(s')+(-1)^{|s|-1}\partial^2(s)\wedge s'+\{\partial(s),s'\}\Big)+(-1)^{|s|}\partial(\{s,s'\})\\ 
				=& \partial^2(s)\wedge s'+s\wedge \partial^2(s').
			\end{array}
		\end{equation*}
		Therefore, since we do not have functions of negative degrees, an easy inductive argument shows that $\partial^2=0.$ Actually, $\partial^2$ is a degree $-2$ vector field on a degree $1$-manifold, which directly guarantees that it must be $0$. For the moreover part notice that since $\partial$ is a generating operator then 
		$$\rho(e)(f)=\{e,f\}=\partial_1(fe)-f\partial_1(e)$$
		and by \eqref{dereq} we also get 
		$$\partial_1([e,e']_E)=\partial_1(\{e,e'\})=\{\partial_1(e),e'\}+\{e,\partial_1(e')\}=\rho(e)\big(\partial_1(e')\big)-\rho(e')\big(\partial_1(e)\big)$$
		Thus, $\partial_1:\Gamma(E)\to C^\infty(M)$ defines a divergence for the Lie algebroid.
	\end{proof}

	\begin{remark}
		In general, degree $-1$ Poisson structures do not admit a canonical generating operator. It was shown in \cite[\S 11]{gor:gcdo} that, given a degree $-1$ Poisson structure $(E^*[1], \{\cdot,\cdot\}),$ the following are equivalent:
		\begin{enumerate}
			\item[1.] Generating operators $\partial:C_i(E^*[1])\to C_{i-1}(E^*[1])$.
			\item[2.] Divergence operators for the Lie algebroid $D:\Gamma(E)\to C^\infty(M)$.
			\item[3.]  Left Lie algebroid  representations on the line bundle $\wedge^{top}E$. 
			\item[4.] Right Lie algebroid  representations on the trivial line bundle $\bR_M$.
		\end{enumerate} 
		In this case, the homology of a generating operator, $H_\bullet^\partial(E)$, is known as \emph{Lie algebroid homology}, and as we are pointing out, it is non-canonical. For a more detailed treatment of this topic see \cite{xu:bv} or \cite[\S 5]{cue:thes}. 
	\end{remark}
	
	An application of our point of view using superdivergences is as follows.
	
	\begin{proposition}\label{prop:divbia}
		Let $(E^*[1],\{\cdot,\cdot\}, Q)$ be a degree $1$ $PQ$-manifold and let $\cD:\fX_\bullet(E^*[1])\to C_\bullet(E^*[1])$ be a superdivergence, then
		$$Q(\cD(Q))=0\quad\text{and}\quad  [Q,\partial]_{Diff}=\frac{1}{2}\{\cD(Q),\cdot\}$$
		where $\partial$ is the generating operator given in Proposition \ref{prop:genop} and $[\cdot,\cdot]_{Diff}$ is the commutator of differential operators.
	\end{proposition}
	\begin{proof}
		For the first identity observe that
		$$0=\cD([Q,Q])=Q\cD(Q)-(-1)^1 Q\cD(Q)=2Q\cD(Q),$$
		while the second follows from
		\begin{equation*}
			[Q,\partial]_{Diff}(s)= Q\partial(s)+\partial Q(s)=(-1)^{|s|}\frac{1}{2}\big(Q(\cD(X_s))-\cD([Q,X_{s}])\big)
			=\frac{1}{2}\{\cD(Q), s\},
		\end{equation*}
		for $s\in C_\bullet(E^*[1]).$
	\end{proof}
	
	The first equation in Proposition \ref{prop:divbia} tells us that $\cD(Q)$ defines a cohomology class on $H^1_Q(E^*[1])$.
	This class is known as the \emph{modular class}, and was introduced in \cite{elw:mod}. In fact, it was suggested in \cite{elw:mod}  that the modular class is the divergence of the $Q$-structure. If $\cD$ and $\cD'$ are two different superdivergences, then by the first property of superdiverngences $$\alpha=\cD-\cD':\fX^1_\bullet(E^*[1])\to C_\bullet(E^*[1])$$ is a $C_\bullet(E^*[1])$-linear map. If one goes to  \S \ref{subs:df-cc}, we immediately see that $\alpha$  defines a degree $0$ differential $1$-form on $E^*[1]$. Moreover, the second property for superdivergences implies that $\alpha$ is a closed form. Hence, by degree reasons we get,  $\alpha\in\Omega^1(M)$ and thus
	$$\cD(Q) - \cD'(Q) = \iota_Q\alpha=\rho_{E^*}^*\alpha\quad \text{for}\quad \alpha\in\Omega^1_{cl}(M).$$
	Consequently, we get that the modular class is independent of the chosen superdivergence and hence defines an invariant of the Lie algebroid $(E^*[1], Q)$.
	
	In general, there are two basic ways to produce a superdivergence on the manifold $E^*[1].$ One is by using superconnections, and the other one is by using \emph{Berezinians}, which are the analogue of volume forms on the manifold $E^*[1]$; details can be found in \cite{mon:div}. In the following, we give three relevant examples. 
	
	\begin{example}
		Let $(\g, [\cdot,\cdot])$ be a Lie algebra. Then $(\g^*[1],\{\cdot,\cdot\})$ is a degree $-1$ Poisson structure. Since $\fX^1_\bullet(\g^*[1])$ is generated by $\fX^1_{-1}(\g^*[1])$ as a $C_\bullet(\g^*[1])$-module, then we have a canonical divergence given by $\cD(X)=0$ for $X\in\fX^1_{-1}(\g^*[1]).$ With this definition, the corresponding generating operator is the usual Lie algebra homology \cite{hue:dual}. In particular, for $e\in\g$ we get $\partial(e)=-\frac{1}{2}tr(\text{ad}_e).$ 
	\end{example}
	
	\begin{example}\label{ex:divtan}
		For an oriented  manifold $M$, consider the $1$-manifold $T[1]M=\big(M, C_\bullet(T[1]M)=\Omega^\bullet(M)\big).$  It was shown in \cite{fn:diff1} that derivations of $\Omega^\bullet(M)$, i.e. vector fields on $T[1]M$, can be written as
		$$\cL_K+\iota_L\qquad\text{for}\quad K,L\in\Omega^\bullet(M;TM).$$
		Using this notation, one shows that $T[1]M$ carries a canonical superdivergence\footnote{We notice here that this canonical superdivergence comes from the canonical Berezinian given by the orientation and it will play an important role in \S~\ref{sec:akszT}.} given by 
		$$\sdiv^c(\iota_L)=(-1)^{|L|}\text{Con}(L),\qquad \sdiv^c(\cL_K)=-d\text{Con}(K),$$
		where $\text{Con}(\cdot)$ denotes contracting the $TM$ part with the form part, see \cite[\S 2.4]{mon:div}.  More concretely, if $$\{x^i, \theta^i=dx^i\}\quad \text{with} \quad |x^i|=0\text{ and } |\theta^i|=1$$ are local coordinates on $T[1]M$ then 
		$$\sdiv^c(\frac{\partial}{\partial \theta^i})=0=\sdiv^c(\frac{\partial}{\partial x^i}).$$
		Since $T[1]M$ carries the $Q$-structure $Q_{dR}=d\in\fX_1^1(T[1]M)$ that is locally given by 
		$$d=\theta^i\frac{\partial}{\partial x^i},\qquad\text{its divergence is}\qquad \sdiv^c(d)=\sdiv^c(\theta^i\frac{\partial}{\partial x^i})=\theta^i\sdiv^c(\frac{\partial}{\partial x^i})+\frac{\partial}{\partial x^i}(\theta^i)=0.$$
		Therefore, we conclude that $Q_{dR}=d$ has a vanishing modular class. 
		
		If, in addition, $(M,\pi)$ is a Poisson manifold, Example \ref{ex:Poibia} implies that $(T[1]M, Q_{dR}=d, \{\cdot,\cdot\}_\pi)$ is a degree $1$ $PQ$-manifold, i.e. a Lie bialgebroid, and the generating operator given by the canonical divergence as in Proposition \ref{prop:genop} is 
		$$\partial^\pi(\omega)=(-1)^{|\omega|}\frac{1}{2}\sdiv^c(\{\omega,\cdot\}_\pi)=[d, \iota_\pi](\omega)\qquad\text{for } \omega\in\Omega^\bullet(M)=C_\bullet(T[1]M).$$
		The map $\partial^\pi$ is known as the \emph{Koszul-Brylinski} operator and gives rise to Poisson homology, see \cite{bry:poih}. By Proposition \ref{prop:divbia}, it satisfies $[d, \partial^\pi]=0.$ Nevertheless, one of the remarkable facts on Poisson homology is that it does not satisfy the Hodge decomposition, see e.g. \cite{mat:har}. For this reason, other homologies and cohomologies have been proposed \cite{yau:hod}.
	\end{example}
	
	\begin{example}\label{ex:divcot}
		For a given oriented manifold $M$, consider the  $1$-manifold $T^*[1]M=\big(M,C_\bullet(T^*[1]M)=\fX^\bullet(M)\big).$ The choice of a volume form $\Lambda\in\Omega^{top}(M)$ on $M$ defines:
		\begin{itemize}
			\item  A (usual) divergence on $M$ given by
			$$div^\Lambda:\fX(M)\to C^\infty(M),\quad X\to div^\Lambda(X)\quad \text{defined by}\quad  \Lie_X\Lambda=div^\Lambda(X)\Lambda,$$ where $\Lie_X$ denotes the usual Lie derivative on $M$.
			\item A star operation on $M$ given by
			$$*:\fX^k(M)\to \Omega^{top-k}(M), \quad \Pi\to *(\Pi)=\iota_\Pi\Lambda.$$
			\item  A superdivergence on $T^*[1]M$ that in local coordinates $\{x^i, \tau_i\}$ is given by 
			$$\sdiv^\Lambda(\frac{\partial}{\partial \tau_i})=0\quad\text{and}\quad  \sdiv^\Lambda(\frac{\partial}{\partial x^i})=div^\Lambda(\frac{\partial}{\partial x^i}).$$
		\end{itemize}
		Moreover, $T^*[1]M$ has a degree $-1$ Poisson bracket given by the Schouten bracket of multivector fields $[\cdot,\cdot]_S$. Then, by Proposition \ref{prop:genop} we get a generating operator of the Schouten bracket given by 
		$$\partial^\Lambda(\Pi)=(-1)^{|\Pi|}\frac{1}{2}\sdiv^\Lambda([\Pi,\cdot]_S)=-*^{-1}d*(\Pi).$$
		If $M$ is, in addition, Poisson, with $\pi\in\fX^2(M)$ satisfying $[\pi,\pi]_S=0$, then Example \ref{ex:Poibia} shows that $(T^*[1]M, Q_\pi, [\cdot,\cdot]_S)$ becomes a degree $1$ $PQ$-manifold. So we can apply Proposition \ref{prop:divbia} and get 
		$$[Q_\pi,\partial^\Lambda]_{Diff}=\frac{1}{2}\Lie_{\partial^\Lambda(\pi)},$$
		where $\partial^\Lambda(\pi)$ is known as the modular vector field and it represents the modular class in Poisson cohomology. The modular class has deep implications in the geometry of Poisson manifolds (including its quantization), see e.g. \cite{bon:quan, elw:mod, fel:deftra, rad:surf}.
	\end{example}
	
	\subsubsection{Morse inequalities}\label{sec:morse}
	
We conclude this section by illustrating the use of superdivergences in the construction of codifferentials associated with Riemannian metrics. As an application, we establish a version of the Morse inequalities, following Witten’s elegant approach in \cite{witt:morse}.
	
Let $(M,g)$ be a connected and compact Riemannian manifold and denote by $\Lambda_g$ the Riemannian volume form on $M$. The metric induces a vector bundle isomorphism $g: TM\to T^*M$ and thus it gives an isomorphism of 1-manifolds $g:T[1]M\to T^\ast[1]M$. Using this isomorphism we can transport the generating operator $\partial^{\Lambda_g}$ constructed in Example \ref{ex:divcot} from $T^*[1]M$ to $T[1]M$ and obtain the usual co-differential and Laplace operators 
	$$\delta:\Omega^\bullet(M)\to \Omega^{\bullet-1}(M),\quad \delta=(-1)^k\star^{-1}d\star\quad\text{and}\quad H:\Omega^\bullet(M)\to \Omega^\bullet(M),\quad H=d\delta+\delta d$$
	where $\star$ is the Hodge star operator.  
	
	Given a function $h\in C^\infty(M)$ one deforms the differential and codifferential as
	$$\textnormal{d}_t=e^{-th} \textnormal{d} e^{th}\quad \textnormal{and}\quad \delta_t=e^{-th} \delta e^{th}\quad \text{for}\quad t\in\bR.$$
	
	It is clear that both operators $\textnormal{d}_t$ and $\delta_t$ square to zero. For each $t$ we also get a deformed Laplacian $H_t=\textnormal{d}_t\delta_t + \delta_t \textnormal{d}_t$, known  as \emph{Witten Laplacian}, which satisfies
	$$\lim_{t\to 0} H_t= H.$$
	The claim is that if $h$ is a \emph{Morse function}, then $\lim_{t\to \infty} H_t$ localizes over critical points of $h$, and therefore it leads to a proof of the Morse inequalities. 
	
	Recall that $h\in C^\infty(M)$ is a \emph{Morse function} if its critical points are all \emph{non-degenerate}, meaning that its Hessian matrix at each of these points, $\textnormal{Hess}_x(f)$, is non-degenerate, see e.g. \cite[\S 1]{aud:morse}. The \emph{index} of a critical point $x$ of $h$ is defined as the number of negative eigenvalues of $\textnormal{Hess}_x(f)$, and is denoted by $\lambda(h,x)$. The \emph{Morse lemma}, see \cite[\S 1.3]{aud:morse}, describes the local behavior of a Morse function $h$ around its critical points by asserting that it looks like the standard quadratic form on $\mathbb{R}^{\dim M}$ of index $\lambda(h,x)$. In particular, this implies that non-degenerate critical points are isolated.
	
	Denote by $$B_k(t)=\dim H^k_{d_{t}}(T[1]M).$$ Since $e^{th}$ is an invertible operator, it follows that $B_k(t)$ is independent of $t$. Let $B_k=B_k(0)$, the \emph{k-th Betti number of $M$}. 
	Moreover, from harmonic decomposition of $H_t$ it follows that 
	$$\sharp\{\text{Zero eigenvalues of } {H_t}_{|\Omega^k}\}=B_k$$
	for all $t\in \bR.$

	Locally around $x\in M$, fix a frame of orthonormal (with respect to $g$) vector fields $\lbrace a^j\rbrace$ and denote by $\lbrace a^{\ast j}\rbrace$ the associated coframe. With these conventions, we obtain that
	\begin{equation}\label{ht}
		H_t=H+t^2g^{ij}\frac{\partial h}{\partial x^i}\frac{\partial h}{\partial x^j}+t\frac{D^2 h}{\partial x^i\partial x^j}[a^{\ast j},a^j].    
	\end{equation}
	Assume $h\in C^\infty(M)$ is a Morse function. This implies that the critical points are isolated and the matrix 
	$$\frac{D^2 h}{\partial x^i\partial x^j}$$
	at a critical point is non-degenerate. Therefore, as $t\to +\infty$ then the second term in \eqref{ht} also goes to infinity except in the vicinity
	of the critical points of $h$. Thus, for
	large $t$ the eigenfunctions of $H_t$ become  concentrated  near the critical points of $h$.
	
	Let $$M_k=\sharp\{ \text{Critical points of } h \text{ with index } k\}.$$  Our goal is to prove the Morse inequalities: $$M_k\geq B_k.$$ 
	
	By the Morse lemma we can rewrite $h$ locally around a critical point $x$ as $$h(y)=h(x)+ \lambda_i(y^i)^2+O(3),$$ where $\lambda_i$ is equal to $-1$ if $i=1,\cdots, \lambda(h,x)$ and equal to $1$ otherwise. Therefore, we approximate $H_t$ around $x$ as
	$$\overline{H}_t=-\frac{\partial^2}{\partial (y^i)^2}+t^2\lambda_i^2(y^i)^2+t\lambda_i [a^{\ast i},a^i]+O(3).$$
	Rewriting  this as $\overline{H}_t=H_i+t\lambda_iK^i$, where $H_i$ is the Hamiltonian of a simple harmonic oscillator and $K^j= [a^{\ast j},a^j]$. In these terms, we get that $H_i$ and $K^j$ mutually commute, and thus can be simultaneously diagonalized. 
	
	It is well known that the eigenvalues of $H_i$ are given by $t\vert \lambda_i\vert (1+2N_i)$ for $N_i=0,1,2,\cdots$ and each of those eigenvalues appears
	with multiplicity one. As the eigenvalues of $K^j$ are $\pm 1$. Therefore, the eigenvalues of $H_t$ are 
	\begin{equation}\label{eq:eigenvalues_morse}
		t\sum_i (\vert \lambda_i\vert (1+2N_i)+\lambda_in_i),\qquad N_i=0,1,2,\cdots,\ n_i=\pm 1
	\end{equation}
	
	If we look at ${H_t}_{|\Omega^k}$ then we know that $\sharp\{\text{positive } n_i\}=k$. Therefore, for \eqref{eq:eigenvalues_morse} to vanish it is needed that $$N_i=0 \quad \forall i\quad \text{and} \quad n_i=+1 \quad \text{if and only if }\quad \lambda_i=-1.$$ Consequently, if we expand $H_t$ around any given critical point of $h$ there is exactly one zero eigenvalue corresponding to a $k$-form if the
	critical point has Morse index $k$. 
	
	That is to say, at most the number of zero energy $k$-forms of $H_t$ equals the number of critical points of Morse index $k$, which means that $M_k\geq B_k$ and we have just established the Morse inequalities.
	
	\section{Basics on $\bN$-graded manifolds}\label{sec:2}
	After our detour on $1$-manifolds, we give the general definitions for $\bN$-graded manifolds of degree $n$ and their fundamental geometric constructions. The main difficulty when working with graded manifolds is that all the constructions must be formulated \emph{without using points on a manifold}. Most of the content of this section consists only of definitions, so we strongly recommend to keep an eye on the other sections for concrete examples that illustrate the objects introduced here.  Our presentation follows the pioneering works \cite{cat:intro, raj:tes, roy:on} and is closely related to \cite{bur:frob, cue:def}.  
	
	\subsection{Graded linear algebra}\label{subsec:GradedVectorSpaces}
	A $\bZ$-graded vector space ${\bf V}$ is a vector space together with a decomposition into a direct sum of subspaces indexed by elements in $\bZ$ that is 
	$${\bf V}=\bigoplus_{i\in\bZ} V_i.$$
	We use $\vert v\vert$ to denote the \emph{degree} of a homogeneous element, i.e. if $v\in V_i$ then $\vert v\vert=i$.
	We say that $f:{\bf V}\to {\bf W}$ is a \emph{morphism} of graded vector spaces if it is a linear map and preserves the degrees, i.e. $f(V_i)\subseteq W_i.$ In this work, $\mathbf{V}[n]$ will denote the graded vector space with degree shifted by $n$, namely $(\mathbf{V}[n])_i=V_{i+n}$.
	
	Most of the operations that we know for vector spaces easily extend to $\bZ$-graded vector spaces including:
	\begin{itemize}
		\item Direct sums, $({\bf V}\oplus{\bf W})_i=V_i\oplus W_i$.
		\item Duals, $(\mathbf{V}^*)_i = (V_{-i})^*$.
		\item Tensor products, $({\bf V}\otimes{\bf W})_i=\oplus_{j\in\bZ} (V_j\otimes W_{i-j})$.
	\end{itemize}
	Nevertheless, one of the main differences is the \emph{sign rule}, see e.g. \cite{del:sign}. This is conventionally defined in such a way that the canonical isomorphism $$\mathbf{V}\otimes \mathbf{W}\cong\mathbf{W}\otimes \mathbf{V} \quad\text{is given by }\quad v\otimes w\to (-1)^{\vert v\vert\vert w\vert}w\otimes v$$  for homogeneous elements $v\in \mathbf{V}$ and $w\in \mathbf{W}$. The rule is that whenever two graded objects are swapped, their degrees show up in the form of a sign.

	In accordance to the previous sign rule, we define the \emph{(graded) symmetric} algebra of a $\bZ$-graded vector space $\mathbf{V}$  as $$\sym \mathbf{V}=T(\mathbf{V})/\langle v\otimes w-(-1)^{\vert v\vert\vert w\vert}w\otimes v\rangle$$  where $T(\mathbf{V})=\bigoplus_{n\geq 0}\mathbf{V}^{\otimes^n}$ is the tensor algebra of $\mathbf{V}$ and  $v,w\in \mathbf{V}$ are homogeneous elements.
	
	\subsection{Definition and first examples}
	Given a non-negative integer $n\in\bN$, let 
	$\mathbf{V}=\oplus_{i=1}^n V_{i}$ be a graded vector space with $\dim V_i=m_i$; let $\sym \mathbf{V}$ denote the graded symmetric algebra of $\mathbf{V}$.
	An \emph{$\mathbb{N}$-graded manifold of degree $n$} (or simply an \emph{$n$-manifold})  is a ringed space $\cM=(M, C_\bullet(\cM))$, where $M$ is a smooth manifold and $C_\bullet(\cM)$ is a sheaf of graded commutative algebras such that any point in $M$ admits a neighborhood $U$ with an isomorphism
	\begin{equation}\label{locality}
		C_\bullet(\cM)|_U\cong C^{\infty}(U)\otimes \sym \mathbf{V},
	\end{equation}
	where the right-hand side is the sheaf of graded commutative algebras on $U$ given by the sheafification of the presheaf $U'\mapsto C^\infty(U')\otimes \sym \mathbf{V}$, for $U'\subseteq U$ open.
	
	We say that $\cM$ has \emph{dimension $m_0|m_1|\cdots|m_n$} and \emph{total dimension} $\totdim(\cM)=\sum_{i=0}^n m_i$, where $m_0=\dim M$.
	The manifold $M$ is known as the \emph{body} of $\cM$. Global sections of $C_\bullet(\cM)$ are called \emph{functions} on $\cM$. We denote by $C_l(\cM)\subseteq C_\bullet(\cM)$ the subsheaf of homogeneous sections of degree $l$. The degree of a homogeneous section $f$ is denoted by $|f|$.
	Given $p\in M$, we denote by $C_\bullet({\cM})_{|p}$ the
	\emph{stalk} of $C_\bullet(\cM)$ at $p$, and given $f\in C_\bullet({\cM})_{|U}$ with $p\in U$, we denote by 
	$\mathbf{f}$ its class in $C_\bullet({\cM})_{|p}$.

	Let $\cM=(M,C_\bullet(\cM))$ be an $n$-manifold of dimension $m_0|\cdots|m_n$. A chart $U\subseteq M$ for which  \eqref{locality} holds is called a \emph{chart} of $\cM$.  We say that 
	$$
	\{  x^{\alpha_i}_i\}  \quad \mathrm{ where } \; i=0,\ldots, n, \, \alpha_i=1, \ldots, m_i, 
	$$   
	are \emph{local coordinates} of $\cM_{|U}=(U,C_\bullet(\cM)_{|_U})$  if $\{x_0^{\alpha_0}\}_{\alpha_0=1}^{m_0}$ are local coordinates on $U$ and $\{x^{\alpha_i}_i\}_{\alpha_i=1}^{m_i}$ form a basis of $V_i$ for all $i =1,\ldots, n$.  Then, any homogeneous section of $C_\bullet(\cM)$ over $U$ can be expressed as a sum of functions that are smooth in $\{x_0^{\alpha_0}\}$ and polynomial in $\{x^{\alpha_i}_i\}$, with $|x^{\alpha_i}_i|=i$.
	We  may occasionally simplify the notation and suppress the sub-indices, denoting local coordinates by 
	$\{x^\alpha\}$.
	
	A \emph{morphism} of $n$-manifolds $\Psi:\cM\to \cN$ is a  morphism of ringed spaces, given by a pair $\Psi=(\psi, \psi^\sharp)$, where $\psi:M \to N$ is a smooth map and $\psi^\sharp:C_\bullet(\cN)\to\psi_* C_\bullet(\cM)$ is a degree preserving morphism of sheaves of algebras over $N$. We will additionally use the pull-back of functions, for $f\in C_\bullet(\cN)$ then $\Psi^\ast f\in C_\bullet(\cM)$ induced by a pair $\Psi=(\psi, \psi^\sharp)$, as we will introduce in Example \ref{ex:pullback}.
	
	For each $n \in \bN$, $n$-manifolds with their morphisms form a category that we denote by $\man{n}$.
	
	An $n$-manifold $\cM=(M,C_\bullet(\cM))$ gives rise to a tower of graded manifolds
	\begin{equation}\label{tower}
		M=\cM_0\leftarrow\cM_1\cdots\leftarrow\cM_{n-1}\leftarrow\cM_{n}=\cM,
	\end{equation}
	where $\cM_r=(M, C_\bullet(\cM_r))$ is an $r$-manifold with $C_\bullet(\cM_r)$ the subsheaf of algebras of $C_\bullet(\cM)$ locally generated by functions of degree $\leq r$, see \cite{roy:on}.

	We illustrate $n$-manifolds with some examples.

	\begin{example}[Linear $n$-manifolds]\label{Ex-linear}
		The most basic example of an ordinary smooth manifold of dimension $m$ is $\bR^m$.  
		For the analogous picture in the graded setting, we consider a collection of non-negative integers $m_0,\cdots,m_n$ and define the graded vector space
		\begin{equation*}
			\mathbf{W}=\oplus_{i=1}^n W_i,
		\end{equation*}
		where $W_i=\bR^{m_i}$. The standard \emph{linear $n$-manifold} of dimension $m_0|\cdots|m_n$ is defined as
		\begin{equation}
			\bR^{m_0|\cdots|m_n}=\big(\bR^{m_0}, C_\bullet({\bR^{m_0|\cdots|m_n}})=C^\infty(\bR^{m_0})\otimes \sym \mathbf{W}\big).
		\end{equation}

		A degree $l$ homogeneous section of $C_\bullet({\bR^{m_0|\cdots|m_n}})$ over an open $U$ is an element in $C^\infty(U)\otimes \mathrm{S}^l \mathbf{W}$. An arbitrary section $f$ over $U$, once restricted to sufficiently small open subsets of $U$ (depending on $f$), is expressed by a finite sum of homogeneous sections. These particular $n$-manifolds are the local models for arbitrary $n$-manifolds, by \eqref{locality}.
	\end{example}

	\begin{example}[Split $n$-manifolds]\label{ex:splitnman}
		The previous example can be generalized as follows. Given a non-positively graded vector bundle  $\mathbf{D}=\oplus_{i=-n+1}^{0} D_i\to M$, we define the graded manifold $$\mathbf{D}[1]=\Big(M,  \Gamma\big({\sym (\mathbf{W})}\big)\Big),$$ where $\mathbf{W} = \oplus_{i=1}^{n} D_{-i+1}^*$ is the graded dual of $\mathbf{D}$ with a degree shifted by $1$ (so it 
		is a graded vector bundle concentrated in degrees $1$ to $n$), $\sym$ denotes, as before, the graded symmetric product, and $\Gamma{\sym (\mathbf{W})}$ is the sheaf of graded algebras given by  $\oplus_{l=0}^\infty \Gamma{\mathrm{S}^l (\mathbf{W})}$. Such $n$-manifolds are known as {\em split}, see e.g. \cite{bon:on}. We point out that, although the split graded manifolds arise from graded vector spaces, the morphisms between them are not necessarily linear maps.
	\end{example}

	\begin{example}[Cartesian product of $n$-manifolds]\label{car-pro}
		Given $n$-manifolds $\cM_i=(M_i,C_\bullet({\cM_i}))$, $i=1,2$, we define an $n$-manifold 
		$\cM_1\times\cM_2=(M_1\times M_2, C_\bullet({\cM_1\times\cM_2}))$, where
		\begin{equation*}
			C_\bullet({\cM_1\times\cM_2})=C_\bullet({\cM_1})\widehat{\otimes}\ C_\bullet({\cM_2}),
		\end{equation*}
		and the hat denotes the usual completion on the product topology,  see e.g. \cite{fio:sup}. Concretely, on open rectangles $U_1\times U_2$, with $U_i\subset M_i$ satisfying $C_\bullet({\cM_i})_{
			|U_i}\cong C^\infty(U_i)\otimes \sym \mathbf{V}_i$ we have
		\begin{equation}\label{prod}
			C_\bullet({\cM_1\times\cM_2})_{|U_1\times U_2}\cong C^\infty(U_1\times U_2)\otimes \sym (\mathbf{V}_1\oplus \mathbf{V}_2),
		\end{equation}
		where $\mathbf{V}_1\oplus \mathbf{V}_2$ denotes the usual graded direct sum. Since open rectangles form a basis of the product topology, \eqref{prod} is sufficient to define $C_\bullet({\cM_1 \times \cM_2})$. 
	\end{example}
	
	\subsection{Vector fields and tangent vectors.}\label{sec:vec}
	Let $\cM=(M,C_\bullet(\cM))$ be an $n$-manifold. A \emph{vector field of degree $k$} on $\cM$ 
	is a degree $k$ derivation $X$ of the graded algebra $C_\bullet(\cM)$ of global sections, i.e.,  
	an  $\mathbb R$-linear map $X\colon C_\bullet(\cM)\to C_\bullet(\cM)$ with the property that, for all 
	$f, g \in C_\bullet(\cM)$ with $f$ homogeneous, $|X(f)| =|f|+k$ and
	\begin{equation}\label{der}
		X(fg)=X(f)g+(-1)^{|f|k}fX(g).
	\end{equation}
	
	Vector fields give rise to a sheaf of $C_\bullet(\cM)$-modules over $M$. The sheaf of degree 
	$k$ vector fields is denoted by $\fX^1_k(\cM)$. The graded commutator of vector 
	fields, defined for homogeneous vector fields $X, Y$ by
	\begin{equation*}
		[X,Y]=XY-(-1)^{|X||Y|}YX,
	\end{equation*}
	turns $\fX^1_\bullet(\cM)$ into a sheaf of graded Lie algebras. If $U\subseteq M$ is a chart of $\cM$ 
	and we have local coordinates $\{x_i^{\alpha_i}\ | \ 0\leq i\leq n, 1\leq \alpha_i\leq m_i\}$ then, as we already showed in Proposition \ref{neg-gen-vf} for $1$-manifolds, the vector 
	fields on $\cM_{|U}$ are generated as $C_\bullet(\cM)_{|U}$-module by
	\begin{equation}\label{eq:basisvf}
		\fX^1_\bullet(\cM_{|U})=\langle \frac{\partial}{\partial x^{\alpha_i}_i}\ | \ 0\leq i\leq n, 1\leq \alpha_i\leq m_i\rangle,
	\end{equation}
	where $\frac{\partial}{\partial x^{\alpha_i}_i}$ are defined as derivations acting on coordinates by 
	$$\frac{\partial}{\partial x^{\alpha_i}_i}(x^{\beta_j}_j)=\delta_{ij}\delta_{\alpha_i\beta_j}.$$
	Notice that this definition implies that $|\frac{\partial}{\partial x^{\alpha_i}_i}|=-i$.

	\begin{example}[The Euler vector field]\label{Ex: Eulervf}
		Given a graded manifold $\cM=(M,C_\bullet(\cM))$ there is always a canonical degree $0$ vector field $\cE_u$, 
		called the \emph{Euler vector field}, defined by $\cE_u(f) = |f| f$. Let $U\subseteq M$ be a chart of $\cM$ 
		with local coordinates $\{x_i^{\alpha_i}\},$ then on $\cM_{|U}$ the Euler vector field has the form
		\begin{equation*}
			{\cE_u}_{|U}=\sum_{i=0}^n\sum_{\alpha_i=1}^{m_i} |x^{\alpha_i}_i|\ x^{\alpha_i}_i\frac{\partial}{\partial x^{\alpha_i}_i}.
		\end{equation*}
	\end{example}

	Vector fields of odd degree  do not need to commute with themselves; if they do, they are of certain interest.

	\begin{definition}[$Q$-manifolds]\label{ex:Qman}
		A degree $1$ vector field on an  $\mathbb{N}$-graded manifold $\cM$, $Q\in\fX^1_1(\cM),$ such that $$[Q,Q]=2Q^2=0,$$ is 
		called a \emph{$Q$-structure or a homological vector field}. The pair $(\cM , Q)$ is called a \emph{$Q$-manifold}.
	\end{definition}

	Note that in the case of Definition \ref{ex:Qman}, the ring of functions 
	$(C_\bullet(\cM), Q)$ is a differential complex, justifying 
	the terminology homological vector field.  Its cohomology will be denoted by $H^\bullet_Q(\cM).$

	A \emph{homogeneous tangent vector of degree $k$ at a point $p\in M$} is 
	a linear map $v\colon C_\bullet({\cM})_{|p}\to \bR$ satisfying
	\begin{equation*}
		v(\mathbf{fg})=v(\mathbf{f})g^0(p)+(-1)^{k|f|}f^0(p)v(\mathbf{g})\quad \forall \mathbf{f, g}\in C_\bullet({\cM})_{|p},
	\end{equation*}
	where $f^0$ and $g^0$ denote the degree zero components of $f$ and $g$, respectively. 
	We denote by $T_p\cM$ the space of tangent vectors at a point $p\in M$. Observe that 
	it is a graded vector space over $\mathbb{R}$ with only non-positive degrees.

	For each open subset $U\subseteq M$, any vector field $X\in\fX^1_\bullet(\cM_{|U})$ defines a 
	tangent vector $X_p$ at each point $p\in U$ by
	\begin{equation}\label{tan-vf-gen}
		X_p(\mathbf{f})=X(f)^0(p)\quad \text{for}\ \mathbf{f}\in C_\bullet(\cM)_{|p}.
	\end{equation}
	For $\mathbb{N}$-manifolds, it is clear that only non-positively graded vector fields 
	can generate non-zero tangent vectors. Therefore, unlike for 
	smooth manifolds, vector fields are not determined by their corresponding tangent 
	vectors. Nevertheless, for all $v_p\in T_p\cM$ there exists a vector field $X\in\fX^1_\bullet(\cM)$, 
	such that $X_p=v_p$. More concretely, let $U$ be a chart of $\cM$ containing $p\in M$ with coordinates 
	$\{ x^{\alpha_i}_i\ |\ 0\leq i\leq n, 1\leq \alpha_i\leq m_i\}$. Then equation \eqref{eq:basisvf} gives a 
	basis of vector fields on $U$ and  we conclude that a basis of  
	the graded vector space $T_p\cM$ is given by 
	$$T_p\cM=\langle \frac{\partial}{\partial x^{\alpha_i}_i}_{|p}\ |\ 0\leq i\leq n, 1\leq \alpha_i\leq m_i\rangle.$$

	\begin{remark}\label{tangentcomplex}
		Let $(\cM, Q)$ be a $Q$-manifold. We immediately deduce from \eqref{tan-vf-gen} that $Q_{|p}=0$ for all $p\in M$, since $Q(f)$ has at least degree $1$ for all functions $f\in C_\bullet(\cM)$.   Thus the linear part of $Q$ induces a 
		differential in $T_p\cM$ and we have that the tangent space at each point is naturally a complex as explained in \cite{aksz}.     
	\end{remark}

	The \emph{differential of a map $\Psi\colon \cM\to\cN$ at $p\in M$} is the linear map
	\begin{equation*}
		\begin{array}{lccl}
			T_p\Psi\colon &T_p\cM&\to& T_{\psi(p)}\cN, \\
			&v&\longrightarrow&T_p\Psi(v)(\mathbf{g})=v(\mathbf{\psi^* g})
			\quad \forall \mathbf{g}\in C_\bullet(\cN)_{|\psi(p)}.
		\end{array}
	\end{equation*}

	We say that a morphism $\Psi\colon \cM\to\cN$ is an 
	\emph{immersion/submersion} if for all  $ p\in M$ the map $T_p\Psi$ is injective/surjective.
	Moreover, $\Psi$ is an \emph{embedding} if it is an injective immersion and $\psi$ 
	is an embedding of smooth manifolds. The following result can be proven in the usual way.

	\begin{proposition}[Local normal form of immersions, cf. \cite{fio:sup}] 
		\label{prop: SubMfldCoo}
		Let $\Psi\colon \cM\to\cN$ be an immersion and $U\subseteq M$ be a chart of $\cM$ around $p$ with local 
		coordinates $\{x_i^{\alpha_i}\}$. Then there exist an open $V$ with $\psi(p)\in V$ and  coordinates 
		$\{\hat{x}^{\alpha_i}_i ,y^{\beta_i}_i\}$ on $\cN_{|V}$ such that  
		$\Psi_{|\psi^{-1}(V)}\colon \cM_{|\psi^{-1}(V)}\to\cN_{|V}$ has the form
		\begin{equation*}
			\psi(x^{\alpha_0}_0)=(\hat{x}^{\alpha_0}_0,0), \quad 
			\psi^* \hat{x}^{\alpha_i}_i= x_i^{\alpha_i}\quad \text{ and }
			\quad \psi^* {y}_k^{\beta_k}= 0.
		\end{equation*}
	\end{proposition}

	\subsection{Submanifolds.}
	\emph{A submanifold of $\cM$} is an $\mathbb{N}$-graded manifold $\cN=(N,C_\bullet(\cN))$ together with an 
	embedding $(j, j^\sharp)\colon \cN\to \cM$ such that the map $j\colon N\to M$ is closed. 
	In addition, we say that it is \emph{wide} if $N=M$.

	\begin{remark}
		Observe that  two submanifolds $(\cN, j)$ and $(\cN', j')$ are \emph{equivalent} 
		if there exists a diffeomorphism $\Psi\colon \cN\to \cN'$ such that $j=j'\circ\Psi$. We 
		will make no distinction between equivalent submanifolds, and we use the term 
		submanifold to refer to an equivalence class.
	\end{remark}

	Let $\cM=(M,C_\bullet(\cM))$ be an $n$-manifold. We say that $\cI\subseteq C_\cM$ is a subsheaf 
	of \emph{homogeneous ideals} if for all $U$ open of $M, \ \cI_{|U}$ is an ideal of 
	$C_\bullet(\cM)_{|U}$, i.e. $$C_\bullet(\cM)_{|U}\cdot \cI_{|U}\subseteq \cI_{|U},$$ and for all 
	$ f\in \cI_{|U}$ its homogeneous components also belong to $\cI_{|U}$.

	Given a subsheaf of homogeneous ideals $\cI\subseteq C_\bullet(\cM)$ we can define a subset of $M$ by
	\begin{equation*}
		Z(\cI)=\bigcap_{f\in \cI\cap C_0(\cM)} f^{-1}(\{0\}).
	\end{equation*}
	We say that $\cI$, a subsheaf of homogeneous ideals, is \emph{regular} if for each $p\in Z(\cI),$ there
	exist $U\subseteq M$ a chart of $\cM$ around $p$
	and local coordinates 
	$\{x_i^{\alpha_i}, y_i^{\beta_i}\}$ of $\cM_{|U}$ for which $\cI_{|U}=\langle y_i^{\beta_i}\rangle$. 
	In this case, $Z(\cI)$ becomes a closed embedded submanifold of $M$. A direct consequence 
	of the local normal form for immersions is the following.

	\begin{proposition}\label{sub-ideal}
		There is a one-to-one correspondence between submanifolds of $\cM$ and subsheaves of 
		regular homogeneous ideals of $C_\bullet(\cM)$.
	\end{proposition}

	\begin{proof}
		We give an idea of the proof, for more details in the context of supermanifolds see \cite[Proposition 5.3.8]{fio:sup}. Given a submanifold $j:\cN\to \cM$ we get the pullback map $j^*:C_\bullet(\cM)\to C_\bullet(\cN)$. Define the corresponding subsheaf of homogeneous regular ideals as $$\cI_\cN=\ker(j^*).$$ Clearly, $\cI_\cN$ is a subsheaf of homogeneous ideals of $C_\bullet(\cM)$. By Proposition \ref{prop: SubMfldCoo} $\cI_\cN$ is regular. Conversely, given $\cI_\cN$ a subsheaf of homogeneous ideals of $C_\bullet(\cM)$ we define the $n$-manifold $$\cN=\big(N=Z(\cI_\cN), j^*(C_\bullet(\cM)/\cI_{\cN})\big)$$ where $j:N\to M$ is the natural embedding. Since $\cI_\cN$ is regular, $N$ is a closed embedded submanifold and $\cN$ is an $n$-manifold. Finally, the natural map $j^*:C_\bullet(\cM)\to j^*(C_\bullet(\cM)/\cI_{\cN})$ makes $j:\cN\to \cM$ into a submanifold of $\cM$.  
	\end{proof}

	Let $(\cN,j)$ be a  submanifold of the $Q$-manifold $(\cM, Q)$ with associated ideal 
	$\cI_\cN$. We say that $(\cN,j)$ is a \emph{$Q$-submanifold} if
	\begin{equation}
		Q(\cI_\cN)\subseteq \cI_\cN.
	\end{equation}

	\begin{example}[Graph of a morphism]
		Let $\Psi=(\psi,\psi^\sharp):\mathcal{M}\to \mathcal{N}$ be a morphism between $n$-manifolds. We have that the graph of $\Psi$ defines a closed embedded submanifold on $\mathcal{M}\times\mathcal{N}$ by the ideal
		\begin{equation*}
			\mathcal{I}_{Gr(\Psi)}=\langle f\otimes g\in C_\bullet(\mathcal{M})\ \widehat{\otimes}\ C_\bullet(\mathcal{N})\ | \ f-\psi^\sharp g=0\rangle.
		\end{equation*}
		
		Moreover, let us verify the following: if $(\mathcal{M}, Q_1)$ and $(\mathcal{N},Q_2)$ are $Q$-manifolds, then $\Psi$ is a $Q$-morphism if and only if $Gr(\Psi)$ is a $Q$-submanifold of $(\mathcal{M}\times\mathcal{N},Q_1\times Q_2)$:
		\begin{equation*}
			\begin{array}{rl}
				\psi^\sharp Q_2 g=Q_1\psi^\sharp g&\Leftrightarrow \Big(Q_1\psi^\sharp g+\psi^\sharp g\Big)\otimes\Big(Q_2 g+g\Big)\in \mathcal{I}_{Gr(\Psi)}\\
				&\Leftrightarrow Q_1\psi^\sharp g\otimes Q_2 g+Q_1 \psi^\sharp g\otimes g+\psi^\sharp g\otimes Q_2 g+\psi^\sharp g\otimes g\in\mathcal{I}_{Gr(\Psi)}\\
				&\Leftrightarrow Q_1\psi^\sharp g\otimes g+\psi^\sharp g\otimes Q_2 g\in\mathcal{I}_{Gr(\Psi)}\\
				&\Leftrightarrow Q_1\times Q_2(\psi^\sharp g\otimes g)\in\mathcal{I}_{Gr(\Psi)}\\
				&\Leftrightarrow Q_1\times Q_2(\mathcal{I}_{Gr(\Psi)})\subseteq \mathcal{I}_{Gr(\Psi)}.
			\end{array}
		\end{equation*}
	\end{example}

	\subsection{Vector bundles.}\label{vector-bundles}
	Following \cite[$\S 2.2.2$]{raj:tes}, we say that a \emph{vector bundle $\pi\colon \cE\to\cM$} is given by two $n$-manifolds $\cE=(E, C_\bullet(\cE))$  
	and $\cM=(M,C_\bullet(\cM))$ and a submersion $\pi=(p,\pi^\sharp)\colon \cE\to\cM$ such that 
	there exists an open cover $\{U_\lambda\}_\lambda$ of $M$ with
	\begin{equation*}
		\cE_{|p^{-1}(U_\lambda)}\cong \cM_{|U_\lambda}\times\bR^{k_0|\cdots|k_n},
	\end{equation*}
	and such that the transition functions between different opens are degree preserving and fiberwise linear.
	The numbers $k_0|\cdots|k_n$ are called the \emph{rank of $\pi\colon \cE\to\cM$}. 
	Therefore, a vector bundle can be defined by specifying trivializations over an open cover
	and gluing them using transition functions, as usual in differential geometry.

	All natural operations such as sums, duals, pull-backs or tensor products are also defined for 
	graded vector bundles\footnote{Many of these operations leave the world of 
		$\mathbb{N}$-graded manifolds but are well defined as $\mathbb{Z}$-graded manifolds, 
		for a more detailed treatment see \cite{raj:tes}.}. In particular, using the tensor 
	product we can define \emph{shifting} by $j\in\bZ$, formally 
	$\cE[j]=\cE\otimes\bR_\cM[j]$ where $\bR_\cM[j]$ denotes the trivial bundle over 
	$\cM$ with just one linear coordinate in degree $j$, thus the new vector bundle has
	fibre coordinates shifted by $j$.

	A \emph{homogeneous section of degree $j$} is a morphism of $\bZ$-graded manifolds 
	$s\colon \cM\to\cE[j]$ such that $\pi\circ s=\id_\cM$. Sections of degree $j$ are denoted by 
	$\Gamma_{j}\cE$, and additionally, $\Gamma\cE:=\oplus_{j\in\mathbb{Z}}\Gamma_j\cE$ is 
	called the space of sections of $\cE$. For us, the sections of a vector bundle always form a left $C_\bullet(\cM)$-module. Therefore, if $\cE_{|U}$ is trivial, for 
	some open $U\subseteq |\cM|$, then
	\begin{align*}
		\Gamma\cE_{|U}=\langle b^i_{\alpha_i} \quad \text{with}\quad |b^i_{\alpha_i}|
		=
		-i \ \text{and} \ 0\leq i\leq n, 1\leq \alpha_i\leq k_i\rangle
	\end{align*}	   
	as a $C_\bullet(\cM)_{|U}$-module. Moreover, as in differential geometry, 
	we establish the following correspondence.

	\begin{proposition}\label{prop:sec-VB}
		Let $\cM=(M, C_\bullet(\cM))$ be an $n$-manifold. The functor 
		$$\Gamma\colon  \{ \text{Vector bundles over } \cM\}\to \{ \text{Locally free finitely generated sheaves of } C_\bullet(\cM)\text{-modules}
		\footnote{If we insist on the total space $\cE$ being an $\bN$-manifold then 
			the local generators of the sheaf must possess non-positive degrees.}\}$$
		is an equivalence of categories.
	\end{proposition}

	\begin{example}[Pull-backs]\label{ex:pullback}
		A standard reference for the technical details on sheaves is $\S$ 2.3 in \cite{sha:she}. Recall that for a sheaf $\cO$ on a manifold 
		$M$ and a map $\psi\colon N\to M$, the \emph{inverse image sheaf} is the sheaf associated to 
		the presheaf $$(\psi^{-1}\cO)(U):=\lim_{\underset{V\supseteq \psi(U)}{\longrightarrow}}\cO(V).$$ 
		Hence, for a vector bundle $\cE\to \cM$ and a map $\Psi\colon \cN\to\cM$, the sheaf of sections of 
		the pull-back bundle  $\Psi^*\cE$ is given by
		$$\Gamma\Psi^*\cE=C_\bullet(\cN)\otimes_{\psi^{-1}C_\bullet(\cM)}\psi^{-1}\Gamma\cE.$$  
		
		This implies that for a given section $s\in\Gamma_j \cE$, we get a section on the pull-back by $\Psi^*s=s\circ \Psi.$
		In particular, for the trivial bundle $\bR_\cM$,  we have that   $\Gamma \bR_{\cM}=\langle 1\rangle$ as $C_\bullet(\cM)$-modules, 
		so we can identify $C_\bullet(\cM)$ with sections of the trivial bundle and $\Psi^*\bR_\cM=\bR_\cN$. Thus, we define the 
		\emph{pull-back of functions} by $$\Psi^*f\in \Gamma \bR_\cN=C_\bullet(\cN),\quad \text{for any}\quad f\in C_\bullet(\cM).$$  
		Note that if the vector bundle $\cE$ is given in terms of (local) transition functions, the transition functions of $\Psi^*\cE$ are given by the pull-backs of the transition functions of $\cE$. 
	\end{example}

	\begin{example}[The tangent bundle]\label{tanget-gm}
		Let $\cM=(M,C_\bullet(\cM))$ be an $n$-manifold of dimension $m_0|\cdots|m_n$. We demostrate that the tangent bundle of $\cM$, 
		$$\pi=(p,\pi^\sharp)\colon T\cM=(TM,C_\bullet({T\cM}))\to\cM,$$  is a graded vector bundle of rank $m_0|\cdots|m_n$. 
		Let $U\subset M$ be a chart of $\cM$ with local coordinates $\{x_i^{\alpha_i}\}$ where
		$0\leq i\leq n, \ 1\leq \alpha_i\leq m_i$. Then $T\cM_{|p^{-1}(U)}$ has linear fiber coordinates given by 
		$$v^{\alpha_i}_i=dx_i^{\alpha_i}\quad \text{with}\quad |v^{\alpha_i}_i|=|x^{\alpha_i}_i|,\quad \text{for}
		\quad 0\leq i\leq n, \ 1\leq \alpha_i\leq m_i.$$
		
		Since we have not introduced  $dx_i^{\alpha_i}$ yet,  we define them by using transition functions of the coordinate vector fields. 
		If $V\subset M$ is another chart of $\cM$ with coordinates $\{\widehat{x}_i^{\alpha_i}\}$ such that on $U\cap V\neq \emptyset$ then we have the change of coordinates 
		$$\widehat{v}_j^{\alpha_j}
		=\sum_{i=0}^j\sum_{\alpha_i=1}^{m_i}v_i^{\alpha_i}\frac{\partial \hat{x}^{\alpha_j}_j}{\partial x_i^{\alpha_i}}.$$
		
		Finally, notice that  sections of the tangent bundle identify with vector fields, i.e. $\Gamma_jT\cM=\fX^1_j(\cM)$, 
		as follows. On the chart $U\subseteq M$, let $b^i_{\alpha_i}\colon \cM_{|U}\to T[-i]\cM_{|p^{-1}(U)}$ be the map defined by 
		$$(b^i_{\alpha_i})^*(f)=f \ \text{for}\ f\in C_\bullet(\cM)_{|U}\ \quad \text{and}
		\quad (b^i_{\alpha_i})^*(v^{\alpha_j}_j)=\delta_{ij}\delta_{\alpha_i\alpha_j}.$$
		Then $b^i_{\alpha_i}$ identifies with the coordinate vector field $\frac{\partial}{\partial x_i^{\alpha_i}}$ on  $C_\bullet(\cM)_{|U}$
		via
		$$(b^i_{\alpha_i})^*(\sum_{j=0}^{n}\sum_{\alpha_j=1}^{m_j}v_j^{\alpha_j}\frac{\partial f}{\partial x_j^{\alpha_j}})=\frac{\partial f}
		{\partial x_i^{\alpha_i}}=\frac{\partial}{\partial x_i^{\alpha_i}}(f)\quad \forall f\in C_\bullet({\cM})_{|U}.$$
	\end{example}

	\subsection{Differential forms and Cartan calculus}\label{subs:df-cc}

	A \emph{differential $k$-form of degree $j$} on the $n$-manifold $\cM=(M,C_\bullet(\cM))$, denoted by $\omega\in\Omega^k_j(\cM)$, is a map
	\begin{equation*}
		\omega\colon \fX^1_{i_1}(\cM)\times\cdots\times\fX^1_{i_k}(\cM)\to C_{i_1+\cdots+i_k+j}(\cM)
	\end{equation*}
	satisfying
	\begin{eqnarray*}
		&&\omega(\cdots, X, Y, \cdots)=-(-1)^{|X||Y|}\omega(\cdots,Y,X,\cdots),\\  
		&&\omega(fX_1,\cdots, X_k)=(-1)^{|f||\omega|}f\omega(X_1,\cdots, X_k),
	\end{eqnarray*}
	for $f\in C_\bullet(\cM).$

	Suppose that 
	$\omega\in\Omega^i_j(\cM)$ and $\eta\in\Omega^k_l(\cM)$. Then, we define a new form 
	$\omega\wedge\eta\in\Omega^{i+k}_{j+l}(\cM)$ by the formula
	\begin{equation*}
		\omega\wedge\eta(X_1,\cdots, X_{i+k})=\sum_{\sigma\in Sh(i,k)} Ksgn(\sigma)
		\omega(X_{\sigma(1)},\cdots, X_{\sigma(i)})\eta(X_{\sigma(i+1)},\cdots,X_{\sigma(i+k)})
	\end{equation*}
	where $X_1,\cdots, X_{i+k}\in\fX^1_\bullet(\cM)$ and $Ksgn(\sigma)$ denotes the 
	signature of the permutation multiplied by the Koszul sign. The differential form
	$\omega \wedge \eta$ is referred to as the \emph{wedge product} of $\omega$ and $\eta$. 
	Observe that, with our 
	sign conventions, we have that $$\omega\wedge\eta=(-1)^{ik+jl}\eta\wedge\omega.$$
	The rule for the wedge product implies that $(\Omega^\bullet_\bullet(\cM),\wedge)$ 
	is generated by $\Omega^1_\bullet(\cM)$ as a sheaf of algebras over $C_\bullet(\cM)$.

	Analogously to the case of vector fields and tangent vectors, we have the following relation between $1$-forms 
	and covectors. The \emph{cotangent space at a point } $p\in M$ is the graded dual of 
	the space $T_p\cM$, i.e. $T^*_p\cM=\underline{\Hom}(T_p\cM, \bR)$, which is a non-negatively graded vector 
	space\footnote{By $\underline{\Hom}$ we denote the internal $\Hom$ in the category of $\bZ$-graded vector spaces, see \S~\ref{sec:akszT}.}.  
	Each differential 1-form of degree $j$, $\alpha\in \Omega^1_j(\cM)$, defines a covector by the formula
	\begin{align*}
		\alpha_p\colon T_p\cM\to \bR[j],\quad  v_p \mapsto \alpha(X)^0(p),
	\end{align*}	 
	where $X\in\mathfrak{X}^1_k(\cM)$ is a vector field, such that $X_p=v_p$. It is easy to show that this 
	map is indeed well-defined. In the same way, for a differential $k$-form $\omega$ one obtains a map 
	$\omega_p\colon \bigwedge^kT_p\cM\to \mathbb{R}$.
	\begin{remark}
		As with vector fields, differential $1$-forms are not determined by the covectors at a point, i.e.  
		there are non-vanishing differential forms $\alpha$ with  $\alpha_p=0$ for all $p\in M$. 
	\end{remark}
	
	Let $U\subseteq M$ be a chart of $\cM$ with coordinates $\{x_i^{\alpha_i}\ | 
	\ 0\leq i\leq n, \ 1\leq \alpha_i\leq m_i\}.$ By \eqref{eq:basisvf},
	$\{ \frac{\partial}{\partial x_i^{\alpha_i}}\}$ is a basis for $\fX^1_\bullet(\cM_{|U})$, thus 
	$$\Omega^1_\bullet(\cM_{|U})=\langle d x_i^{\alpha_i}\ | \ 0\leq i\leq n, \ 1\leq \alpha_i\leq m_i\rangle$$ 
	defines a basis as $C_\bullet(\cM)_{|_U}$-modules as follows  
	\begin{equation}\label{eq:def-1form}
		dx_i^{\alpha_i}(\frac{\partial}{\partial x_j^{\alpha_j}})
		=(-1)^{ij}
		\frac{\partial}{\partial x_j^{\alpha_j}}(x_i^{\alpha_i})=(-1)^{ij}\delta_{ij}\delta_{\alpha_i\alpha_j}.
	\end{equation}
	Moreover, this identity implies that a basis of $T_p^*\cM$ is 
	given by $\{{d x^{\alpha_i}_i}_{|p}\}$. 
	For a function $f\in C_\bullet(\cM)$, we thus have
	\begin{align*}
		df_{|U}= d x_i^{\alpha_i} \frac{\partial f}{\partial x_i^{\alpha_i}}.
	\end{align*}

	\begin{remark}\label{Rem. TransForms}
		For an $n$-manifold $\cM$, in $\S$\ref{vector-bundles}, we considered the sections of vector bundles as left $C_\bullet(\cM)$-modules. However  for differential forms, the more canonical choice  is that of right $C_\bullet(\cM)$-module, as can be seen from the 
		transformation rule of the $d x^{\alpha_i}_i$ above. Nevertheless, we stick to our choice and 
		consider them as a left $C_\bullet(\cM)$-module with local bases $\{dx^{\alpha_i}_i\}$ and the transformation rule 
		\begin{align*}
			d\hat{x}^{\alpha_j}_j
			= (-1)^{i(j-i)}\frac{\partial \hat{x}^{\alpha_j}_j}{\partial x_i^{\alpha_i}} d x_i^{\alpha_i}. 
		\end{align*}
		Note that, up to now, the differential forms are not sections of a vector bundle yet, since $d x^{\alpha_i}_i$ has positive degree, 
		see Proposition \ref{prop:sec-VB}. 
	\end{remark}

	We recall the basics of the Cartan calculus on graded manifolds. As in classical geometry, 
	we define the contraction, Lie derivative, and de Rham differential operators for a graded 
	manifold and prove the usual formulas with some extra signs. In order to know the appropriate 
	sign, just recall that whenever two symbols are transposed, the sign according to their degree appears.

	Let $\cM=(M,C_\bullet(\cM))$ be a graded manifold and consider $X\in\fX^1_r(\cM)$ a vector field. The 
	\emph{contraction} with respect to this vector field is defined as follows:
	\begin{equation*}
		\begin{array}{ll}
			\iota_X\colon \Omega^k_l(\cM)\to\Omega^{k-1}_{l+r}(\cM)\\
			(\iota_X\omega)(X_2,\cdots, X_{k})=(-1)^{|\omega|r}\omega (X, X_2,\cdots, X_k).
		\end{array}
	\end{equation*}
	It acts as a derivation on the wedge product, i.e.
	\begin{equation*}
		\iota_X(\omega\wedge \eta)=\iota_X\omega\wedge \eta+(-1)^{(-1)k+|\omega|r}\omega\wedge \iota_X\eta.
	\end{equation*}
	Therefore, $\iota_X$ is a derivation of the algebra $(\Omega^\bullet_\bullet
	(\cM),\wedge)$ of bidegree $(-1, r)$. We also define the \emph{Lie derivative} with respect to 
	a vector field $X\in\fX^1_r(\cM)$ as
	\begin{equation*}
		\Lie_X f=X(f), \qquad \Lie_X Y=[X,Y],
	\end{equation*}
	and extend it to forms by the following formula
	\begin{equation}\label{lieder}
		\Lie_X \iota_{X_1}\cdots \iota_{X_k}\omega=(-1)^{r\sum\limits_{a=1}^k |X_a|} \iota_{X_1}\cdots \iota_{X_k}
		\Lie_X\omega+\sum_{j=1}^k (-1)^{r\sum\limits_{a=1}^{j-1}|X_a|} \iota_{X_1}\cdots \iota_{[X,X_j]}
		\cdots \iota_{X_k}\omega.
	\end{equation}
	By definition, $\Lie_X$  is a derivation of the wedge product with bidegree $(0,r)$. Finally,
	the de \emph{Rham differential} is defined by the usual Cartan formula
	\begin{equation*}
		\begin{array}{ll}
			\iota_{X_0}\cdots \iota_{X_k}d \omega=&\sum\limits_{i=0}^k (-1)^{i+|X_i|\sum\limits_{a=0}^{i-1}|X_a|}
			\Lie_{X_i}\iota_{X_0}\cdots\widehat{\iota_{X_i}}\cdots \iota_{X_k}\omega\\
			&+ \sum\limits_{i<j} (-1)^{i+1+|X_i|\sum\limits_{a=i+1}^{j-1}|X_a|}\iota_{X_0}\cdots\widehat{\iota_{X_i}}
			\cdots \iota_{[X_i,X_j]}\cdots \iota_{X_k}\omega.
		\end{array}
	\end{equation*}
	This identity implies  that $d$ is a derivation of the wedge product with bidegree $(1,0)$.
	With all these formulas, one gets the usual Cartan  relations.

	\begin{proposition}\label{cartan-cal}
		Let $\cM$ be a graded manifold. The following formulas hold:
		\begin{equation*}
			\begin{array}{l}
				\left[d ,d\right]=2d^2=0,\quad \Lie_X=\left[\iota_X, d \right], 
				\quad \Lie_{\left[X,Y\right]}=\left[\Lie_X, \Lie_Y\right],\\
				\left[\Lie_X, d \right]=0, \quad \iota_{\left[X,Y\right]}=\left[\Lie_X, \iota_Y\right], 
				\quad \left[\iota_X, \iota_Y\right]=0,
			\end{array}
		\end{equation*}
		where the bracket denotes the graded commutator with bidegree, i.e. for example  
		$$[\iota_X,\iota_Y]=\iota_X\circ \iota_Y-(-1)^{(-1)(-1)+|X||Y|}\iota_Y\circ \iota_X.$$ 
		Moreover, for a map $\Psi:\cM\to \cN$, the pull-back of differential forms 
		$\Psi^*:\Omega^\bullet_\bullet(\cN)\to \Omega^\bullet_\bullet(\cM)$ is the unique algebra 
		morphism of bidegree $(0,0)$ that extends the pull-back of functions and commutes with $d$.
	\end{proposition}
	
	\subsection{Multivector fields and Poisson brackets}\label{sec:mult}
	
	Let $\cM$ be an $n$-manifold and $i\in\mathbb{Z}$. We define the \emph{$i$-shifted multivector fields on $\cM$} as the (graded) symmetric algebra 
	\begin{equation}\label{eq:mul}
		^{i}\fX_\bullet^\bullet(\cM)=\bigoplus_{j=0}^\infty \text{S}^j \big(\fX^1_\bullet(\cM)[i]\big)
		\quad\text{with}\quad X\wedge Y=(-1)^{(|X|-i)(|Y|-i)}Y\wedge X\quad X,Y\in \fX^1_\bullet(\cM).
	\end{equation}
	If we extend the Lie bracket using the Leibniz rule, this algebra carries a Poisson bracket of bi-degree 
	$(-1, i)$ that we denote by $[\cdot,\cdot]_S$ and which is known as the \emph{Schouten bracket}. We will make this objects into sections of a vector bundle later on.
	
	Analogously to what happened for $1$-manifolds in the previous section, for $\cM$ an $n$-manifold, an operation $$\{\cdot,\cdot\}:C_i(\cM)\times C_j(\cM)\to C_{i+j+k}(\cM)$$ is called a \emph{degree $k$ Poisson bracket} if  for $f_i\in C_\bullet(\cM)$, we have
	\begin{enumerate}
		\item [(P1)] $\{f_1, f_2\}=-(-1)^{(|f_1|+k)(|f_2|+k)}\{f_2, f_1\}$,
		\item [(P2)] $\{f_1, f_2f_3\}=\{f_1, f_2\}f_3+(-1)^{(|f_1|+k)|f_2|}f_2\{f_1, f_3\}$,
		\item [(P3)] $\{f_1,\{f_2,f_3\}\}=\{\{f_1,f_2\},f_3\}+(-1)^{(|f_1|+k)(|f_2|+k)}\{f_2,\{f_1,f_3\}\}.$
	\end{enumerate}

	An immediate consequence of (P2) is that to each function $f\in C_\bullet(\mathcal{M}),$ one associates a vector field defined by 
	\begin{equation*}
		f\rightsquigarrow X_f=\{f,\cdot\}.
	\end{equation*}
	Such a vector field is known as the \emph{Hamiltonian vector field} of $f$. A vector field $X\in\mathfrak{X}^1_i(\mathcal{M})$ is called a \emph{Poisson vector field} if
	\begin{equation}\label{Q-poisson}
		X(\{f,g\})=\{X(f),g\}+(-1)^{(|f|+k)i}\{f,X(g)\} \qquad \forall f,g\in C_\bullet(\mathcal{M}).
	\end{equation}
	A  \emph{degree $k$ $PQ$-manifold} is a triple $(\mathcal{M},\{\cdot,\cdot\},Q)$ where $\{\cdot,\cdot\}$ is a degree $-k$ Poisson bracket, $Q\in\fX^1_1(\cM)$ is a $Q$-structure and in addition $Q$ is a Poisson vector field.

	\begin{remark}
		It is worth mentioning that	$PQ$-manifolds were introduced by Schwarz in \cite{sch:bv} in the context of supermanifolds in order to provide a more geometric interpretation of the Batalin-Vilkovisky formalism. After that, the works \cite{roy:on,sev:some, vor:dou} among others, made the connections with Poisson and higher geometries that we try to emphasize in these notes.
	\end{remark}
	
	\begin{proposition}
		Given an $n$-manifold $\cM$, there is a one-to-one correspondence between
		\begin{equation*}
			\left\{\begin{array}{c}
				\text{Degree k Poisson}  \\
				\text{ brackets on } \cM
			\end{array}\right\}\\ \leftrightharpoons \left\{\begin{array}{c}
				\pi\in {}^{k-1}\fX_{-k+2}^2(\cM)  \\
				s.t.\  [\pi,\pi]_S=0
			\end{array}\right\} 
		\end{equation*}
	\end{proposition}
	\begin{proof}
		The fact that bivectors correspond to brackets works exactly in the same way as the usual proof. Therefore the only difficulty remains on the degree counting and the desired skew-symmetry. Therefore it is enough to assume that  $\pi=X\wedge Y\in {}^{k-1}\fX^2_{-k+2}(\cM)$ for two given vector fields. With respect to the original degree on $\cM$, we obtain that
		$$-k+2=|\pi|=|X|-k+1+|Y|-k+1\quad\text{so}\quad |X|+|Y|=k.$$
		Therefore the bracket defined by $\{f,g\}:=\pi(df,dg)$ has degree $k$. Finally the definition of the wedge product in \eqref{eq:mul} implies $(P1).$
	\end{proof}

	The correspondence established in the above proposition allows us to introduce the homotopic versions of a Poisson bracket as follows. Let $(\cM, Q)$ be a $Q$-manifold. Following \cite{raj:poi}, we say that a \emph{degree $m$ homotopic Poisson structure\footnote{$m$-Shifted Poisson in perhaps a better terminology coming from \cite{cptvv}.} on $(\cM, Q)$} is a collection $$\pi^\bullet=\sum_{i=2}^\infty \pi^i\in {}^{-m-1}\fX_{m+2}^i(\cM)\quad\text{satisfying}\quad [Q,\pi^\bullet]_S+\frac{1}{2}[\pi^\bullet,\pi^\bullet]_S=0.$$
	Therefore, if $\pi^i=0$ for $i>2$, the  definition above yields  a degree $m$ $PQ$-manifold, i.e. a degree $-m$ Poisson bracket together with a $Q$-structure that is a Poisson vector field. Further examples and applications will be provided in the following sections, see e.g. \S \ref{sec:deflag}.

	\section{$\bN$-graded manifolds of degree 2}\label{sec:3}
	In this section, we focus on $\bN$-graded manifolds of degree $2$. Just as vector bundles encode $1$-manifolds, here we provide three distinct geometric descriptions for $2$-manifolds and explore some of their applications. We then turn to several relevant examples of degree $2$ manifolds (primarily tangent and cotangent bundles) and examine their relation to higher structures.    
	
	\subsection{Definition and geometric descriptions}\label{sec:geo}
	
	It follows from the general definition that a $2$-manifold $\cM$ of dimension $m_0|m_1|m_2$ is a pair $(M, C_\bullet(\cM))$, where $M$ is a dimension $m_0$ manifold and $C_\bullet(\cM)$ is a sheaf of graded-commutative algebras on $M$ satisfying the following property: any $x \in M$ admits an open neighborhood $U \subseteq M$ such that
	$$C_\bullet(\cM)_{|U}\cong C^\infty(U)\otimes \wedge^\bullet\bR^{m_1}\otimes \sym\bR^{m_2},$$
	where the elements of $\bR^{m_q}$ are of degree $q$. In other words, there exist sections $\{x^i, \xi^\alpha, \nu^I \}$ of $C_\bullet(\cM)$ over $U$, where $\{x^i\}$ are usual coordinates on $U$ , the elements $\{\xi^\alpha\}$ and $\{\nu^I \}$ are respectively of degree $1$ and $2$, and every local section can be
	expressed as a sum of functions that are smooth in $x^i$ and polynomial in $\xi^\alpha$ and $\nu^I$.
	
	\begin{example}
		The $2$-manifolds with dimension of special type $m_0|0|0$ are naturally identified with usual smooth manifolds of dimension $m_0$, while  those of dimension type $m_0|m_1|0$ are $1$-manifolds. 
	\end{example}
	
	In what follows we give three different geometric characterizations of $2$-manifolds that extend the one given for $1$-manifolds in Theorem \ref{thm: 1-manifolds}.

	\subsubsection{Graded bundles}\label{sec:graded_bundles}
	The geometric description that we are ready to give was fully explained for all $n$-manifolds in \cite{bon:on}. Therefore, we briefly present the key ideas here. Recall that $\man{2}$ denotes the category of $2$-manifolds and let us denote by $2$-$\cV ect$ the category whose objects are graded vector bundles ${\bf E}\to M$ of the form ${\bf E}=(E_{-1}\oplus E_0)$ and whose morphisms from ${\bf E}\to M$ to ${\bf F}\to N$ are vector bundle maps $\Phi:\wedge^\bullet{\bf E}\to {\bf F}$ covering a smooth map $\varphi:M\to N$. More concretely, $\Phi$ decomposes into three components, which are vector bundle maps covering $\varphi:M\to N$, given by
	\begin{equation}\label{eq:2vectmor}
		\Phi_{0,1}\colon E_0\to F_0,\quad  \Phi_{-1,1}\colon E_{-1}\to F_{-1}\quad \text{ and }\quad \Phi_{0,2}\colon \wedge^2E_0\to F_{-1}.    
	\end{equation}
	\begin{theorem}[\cite{bon:on}]\label{thm:gem2-1}
		The functor 
		\[
		[1]\colon 2\text{-}\cV ect \to \man{2},\quad ({\bf E}\to M)\rightsquigarrow {\bf E}[1]=\Big(M, \Gamma\sym\big(({\bf E}[1])^*\big)\Big)
		\]
		is an equivalence of categories.
	\end{theorem}
	
	\begin{remark}
		We hope that there is no confusion between the $2$-manifold $E[1]$ on the left and the graded vector bundle $E[1]$ on the right. Observe that both things are essentially the same, the graded vector bundle is the ``space" while the $2$-manifold is ``the functions on that space".
	\end{remark}
	
	Since the image of this functor is given by the split $2$-manifolds, as introduced in Example \ref{ex:splitnman}, this result extends Theorem \ref{thm: 1-manifolds}. It is clear that the above functor is well defined and faithful, to show that it is full and essentially surjective we  refer the reader to the original work \cite{bon:on} and \S~\ref{sec:geo2vec}. Once again, we emphasize  that the above result holds for any $n$ as shown in \cite{bon:on}.

	\subsubsection{2-Coalgebra bundle}\label{sec:geo2vec}
	The main disadvantage of Theorem \ref{thm:gem2-1} is that the image of the functor $[1]\colon 2\text{-}\cV ect \to \man{2}$ is too small, which leads to the need for non-canonical choices, as we will explain later. Therefore it is convenient to introduce a new category with an equivalence to $\man{2}$ whose  image is larger. The idea behind the new category is that $2$-manifolds are completely determined by their degree $0,1$ and $2$ functions. Instead of merely retaining the degree $2$ coordinates, we keep all degree $2$ functions. This also requires keeping track of which of these functions arise as products of degree $1$ functions.

	Following \cite{bur:red}, we introduce the category $\cob{2}$ whose objects are triples $(E\to M,\widetilde{E}\to M,\mu)$ where $E\to M$ and $\tilde{E}\to M$ are vector bundles over $M$ and $\mu\colon \tilde{E}\twoheadrightarrow \wedge^2 E$ is a surjective vector bundle morphism covering the identity map on $M$. A \emph{morphism} in $\cob{2}$ is a pair  $$(\Phi,\widetilde{\Phi}):(E\to M, \tilde{E}\to M, \mu)\to (F\to N, \tilde{F}\to N, \mu')$$ where $\Phi\colon E \to F$ and $\tilde{\Phi}\colon \tilde{E}\to \tilde{F}$ are vector bundle morphisms covering the same smooth map $\varphi:M\to N$ and satisfying $\wedge^2\Phi \circ \mu = \mu' \circ \tilde{\Phi}$.
	
	\begin{theorem}[\cite{bur:red}]\label{thm:gem2-2}
		The functor 
		\[
		\cF\colon \cob{2} \to \man{2},\quad (E\to M, \widetilde{E}\to M, \mu)\rightsquigarrow \cM=\Big(M, \dfrac{\Gamma\big(\wedge^\bullet E^*\otimes\sym\widetilde{E}^*\big)}{\langle(\xi\wedge \xi')\otimes 1-1\otimes\mu^\sharp(\xi\wedge \xi')\rangle}\Big)
		\]
		with $\xi,\xi'\in\Gamma(E^*),$
		is an equivalence of categories.
	\end{theorem}

	A full proof of Theorem \ref{thm:gem2-2} can be found in \cite[\S 2.2]{bur:red}, and its extension to $n$-manifolds is one of the main results in \cite{bur:frob}. In the following sections, we will see that the image of the functor $\cF$ includes more than just split manifolds. 
	
	Now, let us compare the categories $2$-$\cV ect$ and $\cob{2}$. If one introduces the functor $\cF':2\text{-}\cV ect\to \cob{2}$ given by
	\begin{eqnarray*}
		\cF'({\bf E}\to M)&=&(E_{0}\to M, \wedge^2E_{0}\oplus E_{-1}\to M, \pr_{\wedge^2E_{0}}),\\
		\cF'(\Phi_{0,1},\Phi_{-1,1},\Phi_{0,2})&=&\big(\Phi_{0,1},
		(\wedge^2\Phi_{0,1}\oplus\Phi_{0,2})\circ\text{Diag}_{\wedge^2 E_0}\oplus\Phi_{-1,1}\big),
	\end{eqnarray*}
	one gets that $\cF'$ is an equivalence of categories and satisfies $[1]=\cF\circ\cF'$. To produce an inverse for $\cF'$ at the level of objects presents no difficulty: given 
	$$(E\to M, \widetilde{E}\to M, \mu)\quad \text{define} \quad E_0=E\quad \text{and}\quad E_{-1}=\ker(\mu).$$
	Nevertheless, given a morphism $(\Phi,\widetilde{\Phi}):(E\to M, \widetilde{E}\to M, \mu)\to (F\to M, \widetilde{F}\to M, \mu')$ one needs to define three maps as in \eqref{eq:2vectmor}. The first two are given by 
	$$\Phi_{0,1}=\Phi\quad\text{and}\quad  \Phi_{-1,1}=\widetilde{\Phi}_{|\ker\mu},$$
	which are well defined because $\wedge^2\Phi \circ \mu = \mu' \circ \tilde{\Phi}$. However, defining $\Phi_{0,2}$ requires choosing splittings of the exact sequences
	$$0\to\ker(\mu)\to \widetilde{E}\xrightarrow{\mu} \wedge^2E\to 0\quad\text{and}\quad 0\to\ker(\mu')\to \widetilde{F}\xrightarrow{\mu'} \wedge^2F\to 0.$$
	These splittings are the non-canonical choices that will appear later in concrete examples.
	

	\subsubsection{Involutive double vector bundles}
	In Theorems \ref{thm:gem2-1} and \ref{thm:gem2-2} we used the shape of a chain complex to give a geometric characterization of $2$-manifolds. However, mathematicians working with higher categories know that it is sometimes useful to use square shapes to model $2$-categories. Building on this idea, an alternative geometric characterization of $2$-manifolds was originally presented in \cite{libla:thes}. The key aspects are as follows.

	Recall that a \emph{double vector bundle} is a vector bundle object in the category of vector bundles. More concretely, it  is a commutative square \eqref{eq:dbv}, where all four sides are vector bundles, all the structural maps are vector bundle morphisms and the two additions on $D$ satisfy the exchange law, see e.g. \cite{raj:lie} for more details. 
	\begin{equation}\label{eq:dbv}
		\begin{array}{c}
			\xymatrix{D\ar[d]\ar[r]& F\ar[d]\\ E\ar[r]& M}
		\end{array}
	\end{equation}
	
	The vector bundles $E$ and $F$ are known as the \emph{side bundles} and we will denote double vector bundles by $(D; E,F;M).$ There is a third relevant vector bundle $C\to M$, known as the \emph{core}, defined as the intersection of the kernels of the bundle map projections of $D$.
	
	Following \cite{fer:geo, libla:thes} we introduce the category $\cI DVB$ whose objects are \emph{involutive double vector bundles} $\mathfrak{A}=((D;E, E;M),J)$, where $(D;E, E;M)$ is a double vector bundle with equal sides and 
	\[
	J=(J; -\operatorname{id}, \operatorname{id}; \operatorname{id})\colon (D;E,E;M)\to (D^f;E,E; M)
	\]
	is a double vector bundle morphism from $D$ to its flip\footnote{This just means that we flip the vertical and horizontal directions.} $D^f$ such that the morphism induced on the cores is the identity. 
	
	A \emph{morphism} in $\cI DVB$ is a double vector bundle morphism $\Phi:(D;E,E;M)\to (D';E',E';M')$ such that 
	\[
	J'\circ \Phi = \Phi^f\circ J.
	\]
	
	Recall that, for a double vector bundle $(D;E,F;M)$, we can define the double polynomial functions $C_{i,j}(D)$ that tells how much linear they are with respect to each vector bundle structure, see \cite[\S 4.2]{mei:wei}. Moreover, given an involutive double vector bundle $\mathfrak{A}=((D;E, E;M),J),$ we say that $f\in C^\infty(D)$ is \emph{involutive} if $J^*f=f$. We denote the space of involutive functions by $C^\infty(D)^J.$
	
	\begin{theorem}[\cite{libla:thes}]\label{thm:gem2-3}
		The functor 
		\[
		\cF_I \colon \cI DVB \to \man{2},\quad \mathfrak{A}=((D;E, E;M),J)\rightsquigarrow \cM=\big(M,C_\bullet(\cM)\big)
		\]
		where
		\begin{equation*}
			\left\{\begin{array}{l}
				C_0(\cM)=C^\infty(M )= C_{0,0}(D) \cap C^\infty(D)^J  \\
				C_1(\cM) =\Gamma(E^*) = C_{0,1}(D) \cap C^\infty(D)^J\\
				C_2(\cM) = C_{1,1}(D) \cap C^\infty(D)^J 
			\end{array}\right.
		\end{equation*}
		is an equivalence of categories.
	\end{theorem}
	
	There is a direct way of going from $\cob{2}$ to $\cI DVB$ by using the characterization of double vector bundles as sequences given in \cite{che:ond}, for this approach see \cite{cue:thes}. We also point out that the above result was extended to arbitrary $n$-manifolds in \cite{man:new, eli:geo}.
	
	\subsection{Examples}
	We already saw that vector bundles and ${\bf E}=E_{-1}\oplus E_0\to M$, i.e. $1$-manifolds and split $2$-manifolds,  give examples of $2$-manifolds. In this section we provide examples of $2$-manifolds that are non-canonically split. 
	
	\subsubsection{The Keller-Waldmann algebra}\label{ex:kw}
	Let $(E \to M, \langle\cdot,\cdot\rangle)$ be a vector bundle equipped with a non-degenerate symmetric bilinear form. To such a structure we can associate a graded algebra as defined by Keller and Waldmann \cite{kw:def}.
	
	A \emph{$k$-cochain} on $E$, $k \geq 1$, is a map 
	$\omega: \Gamma(E) \times \overset{k)}{\cdots} \times \Gamma(E)\to C^\infty(M)$
	that is $C^\infty(M)$-linear in the last entry and, for $k \geq 2$, there exists a map $\sigma_\omega: \Gamma(E) \times\overset{k-2)}{\cdots} \times \Gamma(E)\to \fX(M),$
	called the \emph{symbol}, such that
	\begin{align*}
		&\omega(e_1, \dots, e_i, e_{i+1}, \dots, e_k) + \omega(e_1, \dots, e_{i+1}, e_i, \dots, e_k) 
		= \sigma_\omega (e_1, \dots, \widehat{e}_i, \widehat{e}_{i+1}, \dots, e_k)(\langle e_i,e_{i+1}\rangle)
	\end{align*}
	for all $e_j \in \Gamma(E)$ and $1 \leq i \leq k-1$. We write $C_k(E,\langle\cdot,\cdot\rangle)$ to denote the space of $k$-cochains.
	
	In low degrees, we can give simple descriptions of $C_k(E,\langle\cdot,\cdot\rangle)$: $$C_0(E,\langle\cdot,\cdot\rangle)=C^\infty(M)\quad \text{and}\quad C_1(E,\langle\cdot,\cdot\rangle)=\Gamma(E^*).$$  For each $\omega \in C_2(E,\langle\cdot,\cdot\rangle)$, let $\widehat{\omega}: \Gamma(E) \to \Gamma(E)$ be given by $\langle\widehat{\omega}(e_1),e_2\rangle=\omega(e_1,e_2).$
	Then $\widehat{\omega}\in\Gamma\der(E)$ with symbol $\sigma_\omega$, and it is skew-symmetric, i.e.\ it satisfies
	\[ \sigma_\omega \langle e_1,e_2\rangle = \langle\widehat{\omega}(e_1),e_2\rangle + \langle e_1,\widehat{\omega}(e_2)\rangle\]
	for all $e_1, e_2 \in \Gamma(E)$. The skew-symmetric derivations are the sections of a Lie subalgebroid, which we denote by $\bA_E^{\langle\cdot,\cdot\rangle}$, and the map $\omega \mapsto \hat{\omega}$ gives an isomorphism $C_2(E,\langle\cdot,\cdot\rangle) \cong \Gamma(\bA_E^{\langle\cdot,\cdot\rangle})$.

	The space of cochains $C_\bullet(E,\langle\cdot,\cdot\rangle)$ is a graded-commutative algebra, where the product is given by
	\[
	(\omega\cupprod\tau)(e_1,\cdots, e_{k+m})=\sum_{\pi\in \Sh(k,m)}\text{sgn}(\pi)\ \omega(e_{\pi(1)},\cdots,e_{\pi(k)})\tau(e_{\pi(k+1)},\cdots, e_{\pi(k+m)})
	\]
	for $\omega\in C_k(E, \langle\cdot,\cdot\rangle)$ and $\tau\in C_m(E, \langle\cdot,\cdot\rangle)$; see \cite[Corollary 3.17]{kw:def}.

	It was shown in \cite[Corollary 5.11]{kw:def} that
	$C_\bullet(E,\langle\cdot,\cdot\rangle)$ is generated as an algebra by the elements in degrees $0, 1$ and $2$ and it satisfies the local condition \eqref{locality}. Therefore, given $(E\to M, \langle\cdot,\cdot\rangle)$, we get that  $$\cM=\big(M, C_\bullet(E,\langle\cdot,\cdot\rangle)\big) $$
	defines a $2$-manifold. Moreover, the isomorphism $C_2(E,\langle\cdot,\cdot\rangle) \cong \Gamma(\bA_E^{\langle\cdot,\cdot\rangle})$ implies that $\cM$ is not a split manifold. Nevertheless, a linear connection $\nabla$ on $E\to M$ preserving $\langle\cdot,\cdot\rangle$ induces an isomorphism between $\cM$ and the split $2$-manifold given by  the graded bundle ${\bf E}=(T^*M)_{-1}\oplus (E)_0\to M,$ see \cite[Proposition 5.2]{kw:def}.
	
	\subsubsection{Double vector bundles}\label{ex:vor}
	Take as input an exact sequence of the form
	\begin{equation}\label{eq:dvbs}
		0\to C\to \Omega \stackrel{\varphi}{\to} E\otimes F \to 0,
	\end{equation}
	where $A$, $B$, $C$ and $\Omega$ are vector bundles over $M$, and the maps are over the identity map. We will refer to such exact sequences as {\em dvb-sequences}, since, as shown in \cite{che:ond}, the natural category they form is equivalent to the category of double vector bundles; more concretely, a double vector bundle
	\begin{equation}\label{dbv}
		\xymatrix{D\ar[d]\ar[r]& F\ar[d]\\ E\ar[r]& M}
	\end{equation}
	with core $C$  gives rise to a dvb-sequence \eqref{eq:dvbs} where $\Omega\to M$ is the dual of the vector bundle whose space of sections are the double linear functions on $D$ (see  \cite{raj:lie} and \cite{mei:wei}), and this sequence fully recovers $D$ up to isomorphism (see also \cite{fer:geo} for a detailed discussion).
	
	For the double vector bundles given by \emph{tangent and cotangent prolongations} of a given vector bundle $E\to M$, 
	\begin{equation}\label{eq:tgcotg}
		\begin{array}{c}
			\xymatrix{TE\ar[d]\ar[r]& TM\ar[d]\\ E\ar[r]& M}
		\end{array}\qquad \text{and}\qquad\begin{array}{c}
			\xymatrix{T^*E\ar[d]\ar[r]& E^*\ar[d]\\ E\ar[r]& M,}
		\end{array} 
	\end{equation}
	the corresponding dvb-sequences are the dual sequences of
	\begin{equation}\label{eq:j1der}
		0\to \Hom(TM,E^*) \to J^1E^*\to E^*\to 0  \qquad\text{and}\qquad 0\to \text{End}(E) \to \der(E)\to TM\to 0,
	\end{equation}
	where $J^1E^*$ is the \emph{first jet prolongation} of $E^*$.
	
	For $D$ a double vector bundle (as in \eqref{dbv}) with corresponding dvb-sequence \eqref{eq:dvbs}, we define
	$$
	E^D=E\oplus F,\quad \widetilde{E}^D=\wedge^2 E\oplus \Omega \oplus \wedge^2 F,\quad \mu^D(\alpha+ \omega+\tau)=\alpha+\varphi(\omega)+\tau,
	$$
	where $\alpha\in\wedge^2E,\ \omega\in \Omega,\ \tau\in\wedge^2 F$. 
	
	For a geometric description of the $2$-manifold
	$$
	\cM^D:= \cF(E^D\to M, \widetilde{E}^D\to M, \mu^D),
	$$
	we consider the $1$-manifolds $E[1]$ and $D[1]$ corresponding to the vector bundles $E\to M$ and $D\to F$. The additional vector bundle structure of $D$ over $E$ makes $D[1]\to E[1]$ into a vector bundle in the category of $1$-manifolds. Denoting by $\mathcal{R}\to E[1]$ the trivial vector bundle with  $\mathcal{R}=\bR^{0|1}\times E[1]$, the $2$-manifold $\cM^D$ can be regarded as the total space of the graded vector bundle $D[1]\otimes \mathcal{R}\to E[1]$. 
	
	When $D$ is the tangent or cotangent prolongation of $E\to M$ (see \eqref{eq:tgcotg} and \eqref{eq:j1der}), one can check that  
	$$
	\cM^{TE}=T[1]E[1] \quad \text{and}\quad \cM^{T^*E}=T^*[2]E[1].
	$$ 
	
	Let us make a couple of final remarks on the above construction:
	\begin{itemize}
		\item The manifold $\cM^D$ is denoted in \cite{vor:dou} by  $D[1]_F[1]_E$ that indicates the shift of the fibre coordinates. In \S~\ref{sec:doua} we endow $\cM^D$ with a $Q$-structure.
		\item The choice of a splitting for the exact sequence \eqref{eq:dvbs} is equivalent to give an isomorphism between $\cM^D$ and the split $2$-manifold defined by ${\bf E}=C_{-1}\oplus (A\oplus B)_0$.
		\item Notwithstanding the procedure described here involves double vector bundles, it is different from the functor described in Theorem \ref{thm:gem2-3} and must not be confused with it.
		\item Applying the same ideas one can  produce examples of interesting $n$-manifolds, for more details see \cite[Example 3.15]{bur:frob}.
	\end{itemize}

	\subsection{Lie 2-algebroids and Q-manifolds}
	
	A \emph{Lie $n$-algebroid} is a graded vector bundle ${\bf E}=\oplus_{i=-n+1}^0 E_i\to M$  together with a bundle map $\rho:E_0\to TM$ and graded skew-symmetric  brackets $$[\cdot,\dots, \cdot]_i:\Gamma{\bf E}\times\overset{i)}{\cdots}\times \Gamma{\bf E}\to \Gamma{\bf E},\  i \in  \{1, \cdots, n + 1\},$$ of degree $2-i$, such that 
	$$[e_0, fe]_2=f[e_0, e]_2+\rho(e_0)(f)e\quad \text{for}\quad e_0\in\Gamma(E_0),\ e\in\Gamma{\bf E},\ f\in C^\infty(M),$$
	all the brackets are $C^\infty(M)$-linear   and the following Jacobi identity is satisfied: For any $k\ge 1$
	\begin{equation}\label{eq:higjac}
		\sum_{i+j=k+1}\sum_{\sigma\in\Sh(i,j-1)}(-1)^{i(j-1)}Ksign(\sigma) [[e_{\sigma(1)},\cdots, e_{\sigma(i)}]_i, e_{\sigma(i+1)},\cdots, e_{\sigma(r)}]_j=0
	\end{equation}
	where $Ksign(\sigma)$ is the signature multiplied by the Koszul sign with respect to the grading of ${\bf E}$.
	
	Motivated by Theorem \ref{thm:Q-lie} of \cite{vai:lie} and \cite[\S~4.3]{kon:def}, it was suggested in \cite{sev:some, vor:qman} that $Q$-manifolds are related to ``higher Lie algebroids". Following \cite{sta:shla}, the article \cite{cc:extensions} introduced the above definition of Lie $n$-algebroids. Finally its equivalence with $Q$-manifolds was proved in \cite{bon:on}. Here we give the key ideas for $n=2.$

	If we specify the above definition for $n=2$ we obtain that  a Lie $2$-algebroid is given by:
	\begin{itemize}
		\item Two vector bundles ${\bf E}=E_{-1}\oplus E_0\to M$,
		\item Two vector bundle maps over the identity $\rho:E_0\to TM$ and $\partial=[\cdot]_1:E_{-1}\to E_0$,
		\item Two maps $[\cdot,\cdot]:\Gamma(E_0\wedge E_0)\to \Gamma (E_0)$ and $\nabla:\Gamma (E_0\otimes E_{-1})\to \Gamma (E_{-1})$ coming from the 2-bracket that satisfy the following Leibniz rules
		$$[e,fe']=f[e,e']+\rho(e)(f)e',\quad \nabla_e(fv)=f\nabla_ev+\rho(e)(f)v,\quad \nabla_{fe}v=f\nabla_ev,$$
		for $e,e'\in\Gamma (E_0),\ v\in\Gamma (E_{-1})$ and $f\in C^\infty(M),$
		\item A vector bundle map $[\cdot,\cdot,\cdot]_3:\wedge^3 E_0\to E_{-1}$ given by the 3-bracket.
	\end{itemize}
	The higher Jacobi identity \eqref{eq:higjac}  is equivalent to the following equations
	\begin{eqnarray}
		&\rho\circ\partial=0,\quad  \partial(\nabla_e v)=[e,\partial v], \quad \nabla_{\partial v}w=-\nabla_{\partial w}v,\label{eq:J1}\\
		&J(\cdot,\cdot,\cdot)=\partial\circ [\cdot,\cdot,\cdot]_3, \quad F_\nabla(e,e')(v)=[\partial v, e,e']_3,\label{eq:J2}\\
		&\big(\nabla_e[e',e'',e''']+c.p.\big)+[[e,e''],e',e''']_3+[[e',e'''],e,e'']_3=\big([[e,e'],e'',e''']_3+c.p.\big),\label{eq:J3}
	\end{eqnarray}
	where $J:\wedge^3 E_0\to E_0$ is given as in \eqref{eq:jacobiator}, $F_\nabla(e,e')(v)=\nabla_{[e,e']}v-[\nabla_e,\nabla_{e'}](v)$ is the curvature, $e,e',e'',e'''\in\Gamma (E_0), \ v,w\in\Gamma (E_{-1})$ and $c.p.$ denotes cyclic permutations.
	
	Roughly speaking, equation \eqref{eq:J1} tells us that $\rho$ and $\partial$ form a chain complex and that the bracket and the connection descend to the homology of the complex, \eqref{eq:J2} says that on the homology we get a graded Lie algebra and \eqref{eq:J3} is a higher coherence condition, see e.g. \cite{bc:lie2} for a more conceptual explanation.
	
	\begin{theorem}[\cite{bon:on}]\label{thm:bon}
		The functor $[1]$, given in Theorem \ref{thm:gem2-1}, provides an equivalence of categories between 
		$$[1]:\{\text{Lie 2-algebroids} \}\to \{\text{Degree 2 Q-manifolds}\}$$
	\end{theorem}
	\begin{proof}
		Given a Lie $2$-algebroid $({\bf E}\to M, \rho, \partial,  [\cdot,\cdot],\nabla, [\cdot,\cdot,\cdot]_3),$ define the graded manifold ${\bf E}[1]$ as in Theorem \ref{thm:gem2-1}. To define the $Q$-structure we modify formula \eqref{eq:Qalg} to include the new brackets as follows: for $\eta\in C_{k+2l}({\bf E}[1])=\wedge^k E_0^*\otimes S^l E^*_{-1}$ then $Q(\eta)\in C_{k+2l+1}({\bf E}[1])$ and it decomposes into  three components  given by
		\begin{align*}
			(Q\eta)(e_0,\cdots, e_{k-2}; v_0,\cdots, v_{l})& = (-1)^{l+1}\sum_{i=0}^{l} \eta(e_0,\cdots, e_{k-2}, \partial(v_i); v_0,\cdots, \widehat{v}_i,\cdots, v_{l})\\
			(Q\eta)(e_0,\cdots, e_{k}; v_0,\cdots, v_{l-1})& = \sum_{i=0}^k(-1)^i\rho(e_i)\Big(\eta(e_0,\cdots,\hat{e}_i,\cdots,e_k; v_0,\cdots, v_{l-1})\Big)\\
			&\quad +\sum_{i,j}(-1)^{i+1} \eta(e_0,\cdots, \widehat{e}_i,\cdots, e_{k}; \nabla_{e_i}v_j, v_0,\cdots,  \widehat{v}_j,\cdots, v_{l-1})\\
			&\quad+\sum_{i<j}(-1)^{i+j} \eta([e_i,e_j], e_0,\cdots, \widehat{e}_i,\cdots, \widehat{e}_j,\cdots, e_{k}; v_0,\cdots, v_{l-1}) \\
			(Q\eta)(e_0,\cdots, e_{k+2}; v_0,\cdots, v_{l-2})& =\sum_{i<j<k}(-1)^{i+j+r} \eta(e_0,\cdots, \widehat{e}_i, \cdots,\widehat{e}_j,\cdots,\widehat{e}_r,\cdots, e_{k+2};\\
			&\hspace{6.5cm}[e_i,e_j, e_r]_3, v_0,\cdots, v_{l-2}) 
		\end{align*}
		for $e_i\in\Gamma (E_0)$ and $v_j\in\Gamma (E_{-1}).$
		
		To check that $Q^2=0$ if and only if \eqref{eq:J1},\eqref{eq:J2} and \eqref{eq:J3} are satisfied is a bit more cumbersome but it can be done in the same way as in Theorem \ref{thm:Q-lie}, see also \cite{bon:on}. 
		
		For the converse we just point out the following. Let $(\cM, Q)$ be a degree $2$ $Q$-manifold, then Theorem \ref{thm:gem2-1} implies that there exists a diffeomorphism of graded manifolds $\Phi:\cM\overset{\simeq}{\to}{\bf E[1]}$ for some graded bundle ${\bf E}=E_{-1}\oplus E_0$. Under this isomorphism a $Q\in\fX^1_1({\bf E}[1])$ in coordinates looks like
		$$Q=\rho^i_\alpha(x)\xi^\alpha\frac{\partial}{\partial x^i}+\Big(\frac{1}{2}c_{\alpha\beta}^\gamma(x)\xi^\alpha\xi^\beta+\partial^\gamma_I(x)\nu^I\Big)\frac{\partial}{\partial \xi^\gamma}+\Big(\nabla_{\alpha I}^J(x)\xi^\alpha\nu^I+\frac{1}{3}R_{\alpha\beta\gamma}^J(x)\xi^\alpha\xi^\beta\xi^\gamma\Big)\frac{\partial}{\partial \nu^J}$$
		and thus, inside $Q$ one identifies all structural constants of the objects defining a Lie $2$-algebroid. Again, we direct the reader to \cite{bon:on} for  the fact that $Q^2=0$ is equivalent to the aforementioned equations.
	\end{proof}
	
	\begin{remark}\label{rmk:bon}
		We want to point out again that the above result follows for arbitrary $n$, see \cite{bon:on}. Indeed it establishes an equivalence of categories between the category of Lie $n$-algebroids  and the category of degree $n$ $Q$-manifolds.
	\end{remark}
	
	\subsubsection{Double Lie algebroids}\label{sec:doua}
	In \S~\ref{ex:vor}, we saw that a double vector bundle defines a $2$-manifold.  Here, we add more structures to the double vector bundle in order to obtain a $Q$-manifold. 
	
	A \emph{double Lie algebroid} \cite{mac:dou} is a double vector bundle as in \eqref{eq:dbv}, with core denoted by $C$, and with Lie algebroid structures on each of the four vector bundles such that each pair of parallel Lie algebroids endows $D$ with the structure of a $VB$-algebroid,
	and such that a suitable dual\footnote{To be precise $(D^{\uparrow\to}, D^{\to\uparrow})$.} with the induced Lie algebroid structures
	as defined above, is a Lie bialgebroid. 
	
	This complicated definition, involving several duals of $D$, can be reformulated in many different ways, see e.g. \cite{bch:vec, raj:dou, vor:macken}. Examples of double Lie algebroids are given by the tangent and cotangent prolongations of any Lie bialgebroid \cite{mac:dou}, linking to the $Q$-manifold we will present in $\S~\ref{sec:liebia}$. We will see in $\S~\ref{sec:liefun}$ that double Lie algebroids appear as the infinitesimal object of a double Lie groupoid. 
	
	The following result is an immediate consequence of \cite[Theorem 1]{vor:macken}. Moreover, this proposition provides the infinitesimal analogue of the Artin–Mazur construction, which we will present in  \S~\ref{sec:doug}.
	
	\begin{proposition}
		Let $D$ be a double Lie algebroid. Then the $2$-manifold $\cM^D$ constructed in \S~\ref{ex:vor} is a $Q$-manifold.
	\end{proposition}
	\begin{proof}
		The idea is as follows. Since $D\to E$ is a Lie algebroid then $D[1]_E$ has a $Q$-structure $Q_E$. Moreover, $Q_E$ is linear with respect to the graded vector bundle $D[1]_E\to F[1]$. Analogously we also get on $D[1]_F$ a $Q$-structure $Q_F$ that is linear with respect to $D[1]_F\to E[1]$. Therefore, $\cM^D =D[1]_E[1]_F$ carries a degree $1$-vector field given by $Q=Q_E+Q_F$. The compatibility condition of a double Lie algebroid implies that $[Q_E,Q_F]=0$ so $Q^2=0$ as we wanted.
	\end{proof}
	
	\subsection{$T[1]\g[1]$ and equivariant cohomology}\label{sec:equico}
	Let $G$ be a compact and connected Lie group with Lie algebra $(\g,[\cdot,\cdot]),$ and let $M$ be a smooth manifold with an action $\psi:G\times M\to M$. The goal is to give a homological model for the leaf space $M/G$ using graded geometry. In this section we will mostly follow \cite{gui:sup} and references therein.
	
	Since $(\g,[\cdot,\cdot])$ is a Lie algebra, Theorem \ref{thm:Q-lie} implies that $(\g[1], d_{CE})$ is a $Q$-manifold, see Example \ref{ex:Q-TM-g} and consider $\Omega^\bullet_\bullet(\g[1])$ the differential forms on $\g[1]$ as constructed in \S\ref{subs:df-cc}. Then, as it happens in Example \ref{ex:1-man}, the $2$-manifold $T[1]\g[1]$ constructed in \S~\ref{ex:vor} satisfies that
	$$C_i(T[1]\g[1])=\bigoplus_{j+k=i}\Omega^k_j(\g[1])=\bigoplus_{a+2b=i}\wedge^a\g^*\otimes S^b\g^*.$$
	The algebra $\text{W}(\g)=\wedge^\bullet\g^*\otimes S^\bullet\g^*$ is called the \emph{Weil algebra of $\g$}. The graded geometric interpretation goes back to \cite{kalk:brst}, see \cite{raj:tes} for the Lie algebroid case.

	Pick $\{e_\alpha\}$ a basis of $\g$ and denote by $\{\xi^\alpha\}$ its dual basis. Since $\g[1]$ is a $1$-manifold then $\{\xi^\alpha\}$ are coordinates on $\g[1]$ and the generators of its differential forms are 
	$$\Omega^\bullet_\bullet(\g[1])=\langle \xi^\alpha,\  \nu^\alpha=d\xi^\alpha\rangle_{C_\bullet(\g[1])}.$$
	Therefore, the $2$-manifold $T[1]\g[1]$ has local coordinates given by $$\{\xi^\alpha,\nu^\alpha\}\quad \text{with}\quad |\xi^\alpha|=1\text{ and } |\nu^\alpha|=2.$$
	On the one hand, the $2$-manifold $T[1]\g[1]$ has the following relevant vector fields: 
	\begin{equation*}
		\begin{array}{|c|c|c|}
			\hline
			\text{Degree}&\text{Vector fields} & \text{Coordinate expression}\\
			\hline
			-1& 
			\begin{array}{c}
				\text{Contractions}  \\
				\iota_\alpha=\Lie_{\frac{\partial}{\partial\xi^\alpha}}
			\end{array} & \dfrac{\partial}{\partial\xi^\alpha}\\
			\hline
			0&
			\begin{array}{c}
				\text{Lie derivatives}  \\
				L_\alpha=[\iota_\alpha, d_W]=\Lie_{[d_{CE},\frac{\partial}{\partial\xi^\alpha}]}
			\end{array}& c_{\alpha\beta}^\gamma\xi^\beta\dfrac{\partial}{\partial \xi^\gamma}+c_{\alpha\beta}^\gamma\nu^\beta\dfrac{\partial}{\partial \nu^\gamma}\\
			\hline
			0&
			\begin{array}{c}
				\text{An automorphism}  \\
				\iota_{d_{CE}}
			\end{array} &\frac{1}{2}c_{\alpha\beta}^\gamma\xi^\alpha\xi^\beta\dfrac{\partial}{\partial \nu^\gamma}\\
			\hline
			1&
			\begin{array}{c}
				\text{The differential}  \\
				d_W=d+\Lie_{d_{CE}}
			\end{array}   & \nu^\alpha\dfrac{\partial}{\partial\xi^\alpha}+\frac{1}{2}c_{\alpha\beta}^\gamma\xi^\alpha\xi^\beta\dfrac{\partial}{\partial \xi^\gamma}+c_{\alpha\beta}^\gamma\xi^\alpha\nu^\beta\dfrac{\partial}{\partial \nu^\gamma} \\
			\hline
		\end{array}
	\end{equation*}
	where $[e_\alpha,e_\beta]=c_{\alpha\beta}^\gamma e_\gamma$ are the structural constants of the Lie algebra. 
	
	One can easily compute that $H^\bullet_{d_W}(T[1]\g[1])=H^0_{d_W}(T[1]\g[1])=\bR.$ Combined with the existence of the vector fields $\iota_\alpha, L_\alpha$ and $d_W$ the latter implies that   $C_\bullet(T[1]\g[1])$  gives a model for $EG$, the universal bundle over the classifying space $BG$, see \cite{gui:sup} for more details.

	On the other hand, denote by $\rho:\g\to\fX^1(M)$ the infinitesimal action of $\psi:G\times M\to M$ and consider the $1$-manifold $T[1]M$ with local coordinates given by $\{x^i, \theta^i\}$ as in Example \ref{ex:Q-TM-g}. Due to the Lie algebra action we also get the following vector fields on $T[1]M$:
	\begin{equation*}
		\begin{array}{|c|c|c|}
			\hline
			\text{Degree} & \text{Vector field}& \text{Coordinate expression}\\
			\hline
			-1&
			\widehat{\iota}_\alpha=\iota_{\rho(e_\alpha)}
			& \rho^i_\alpha(x)\dfrac{\partial}{\partial\theta^i} \\
			\hline
			0&
			\widehat{L}_\alpha=\Lie_{\rho(e_\alpha)}
			& \rho^i_\alpha(x)\dfrac{\partial}{\partial x^i}+\dfrac{\partial\rho^i(x)}{\partial x^j}\theta^j\dfrac{\partial}{\partial \theta^i}\\
			\hline
			1 &
			d
			& \theta^i\dfrac{\partial}{\partial x^i}\\
			\hline
		\end{array}
	\end{equation*}
	Finally, consider the $2$-manifold $\cM=T[1]\g[1]\times T[1]M$ with the $Q$-structure given by $Q=d_W\times d$.
	A function $\omega\in C_\bullet(\cM)$ is called \emph{basic} if
	$$(\iota_\alpha\times\widehat{\iota}_\alpha)(\omega)=0\quad\text{and}\quad (L_\alpha\times \widehat{L}_\alpha)(\omega)=0$$
	One can check that $Q$ restricts to basic functions. Indeed $(C_\bullet(\cM)_{|\text{basic}}, Q)$ is called the \emph{Weil model for equivariant cohomology} and is the right cohomology for the leaf space $M/G$.
	
	Nevertheless there is another famous model for the equivariant cohomology known as the \emph{Cartan model}. The relation between the Weil and the Cartan models can be summarized in the following result.
	
	\begin{theorem}[\cite{kalk:brst, mq}]
		The $2$-manifold $\cM$ has an automorphism given by
		$$\Phi=e^{(\xi^\alpha\widehat{\iota}_\alpha)}:\cM\to \cM$$
		and it satisfies
		$$\Phi Q\Phi^{-1}=Q-\nu^\alpha\widehat{\iota}_\alpha+\xi^\alpha\widehat{L}_\alpha=:Q_{BRST}\quad \text{and}\quad \Phi(\iota_\alpha\times\widehat{\iota}_\alpha)\Phi^{-1}=\iota_\alpha.$$
		Moreover, $$C_\bullet(\cM)_{|\text{basic}}\cong\big(\sym\g^*\otimes\Omega^\bullet(M)\big)^G\quad\text{and}\quad Q_{|\text{basic}}=1\otimes d-\nu^\alpha\widehat{\iota}_\alpha$$ that is known as the \emph{Cartan} model for equivariant cohomology.
	\end{theorem}
	The proof of the above result can be found in \cite[\S~4]{gui:sup}. As far as we know, the generalization of this result for Lie groupoids or Lie $n$-groups is still an open problem.   
	
	\subsection{$T^*[2]E[1]$ and Lie bialgebroids}\label{sec:liebia}
	In the previous section we studied the tangent of a $1$-manifold; here we do the opposite, we focus on the cotangent bundle $T^*[2]E[1]$ of the $1$-manifold $E[1]$. 
	
	One possible option is to describe $T^*[2]E[1]$ as we did in \S~\ref{ex:vor}. However, we want to relate the cotangent bundles with the shifted multivector fields of \S~\ref{sec:mult}. Therefore, for a given $k\geq 1$ we define 
	\begin{equation}\label{eq:defcot2}
		C_\bullet(T^*[k]E[1]):={}^{-k}\fX_\bullet^\bullet(E[1])=\bigoplus_{j=0}^\infty S^j \big(\fX^1_\bullet(E[1])[-k]\big).
	\end{equation}
	Since $E[1]$ is a $1$-manifold, Proposition \ref{neg-gen-vf} implies that $T^*[2]E[1]$ is indeed a $2$-manifold whose functions in low degrees are\footnote{In this section we eliminate the left index $-2$ on ${}^{-k}\fX_\bullet^\bullet(E[1])$ and follow the usual convention that sub-indices denote the original degree  on the manifold $E[1]$ instead of the one on $T^*[2]E[1]$. We hope this will not create confusion since one should take care of the correct symmetric properties of each object.}
	\begin{enumerate}
		\item[1.] $C_0(T^*[2]E[1])=C^\infty(M),$
		\item[2.] $C_1(T^*[2]E[1])=C_1(E[1])\oplus \fX^1_{-1}(E[1])=\Gamma(E^*\oplus E),$
		\item[3.] $C_2(T^*[2]E[1])=C_2(E[1])\oplus \fX_0^1(E[1])\oplus \fX^2_{-2}(E[1])=\Gamma(\wedge^2 E^*\oplus \der(E^*)\oplus\wedge^2 E).$
	\end{enumerate}
	Therefore, it matches with the description given in Example \ref{ex:vor}. Moreover, observe that the vector bundle $E\oplus E^*$ carries a canonical non-degenerate and symmetric pairing. Under this pairing we can identify 
	$$\bA_{E\oplus E^*}^{\langle\cdot,\cdot\rangle}=\wedge^2 E^*\oplus \der(E^*)\oplus\wedge^2 E.$$
	Hence, this manifold also admits a description in terms of the Keller-Waldmann algebra given in \S~\ref{ex:kw}. 
	
	If $\{x^i,\xi^\alpha\}$ denotes a set of local coordinates for the $1$-manifold $E[1]$ then we get that 
	$$\{x^i, \xi^\alpha, \nu_i=\frac{\partial}{\partial x^i}, p_\alpha=\frac{\partial}{\partial \xi^\alpha}\}\quad \text{with}\quad |x^i|=0, |\xi^\alpha|=1, |\nu_i|=2, |p_\alpha|=1$$
	are a set of local coordinates for $T^*[2]E[1].$ In particular the $\nu_i$ are sections of $TM$ corresponding to the symbols for $\der(E^*)$. As we already mentioned in \S~\ref{ex:vor}, a linear connection on $E$ will identify the $2$-manifold $T^*[2]E[1]$ with the split $2$-manifold corresponding to the graded bundle ${\bf E}=(T^*M)_{-1}\oplus(E\oplus E^*)_0$. To complete the picture, see e.g. \cite{fer:geo, libla:thes} for details, we point out that the involutive double vector bundle corresponding to $T^*[2]E[1]$ under the functor of Theorem \ref{thm:gem2-3} is
	\begin{equation}\label{quadradito}
		\begin{array}{c}
			\xymatrix{T^*(E\oplus E^*)\ar[d]\ar[r]& E\oplus E^*\ar[d]\\ E\oplus E^*\ar[r]& M.}
		\end{array}
	\end{equation}
	In this example we see the difference between the construction given in \S~\ref{ex:vor}, that use the cotangent prolongation of $E\to M$, and the functor described in Theorem  \ref{thm:gem2-3} that needs \eqref{quadradito}.

	It follows immediately from \eqref{eq:defcot2} that the Schouten bracket on multivector fields  defines a degree $-2$ Poisson bracket on the functions $C_\bullet(T^*[2]E[1]).$ A geometric characterization of this bracket in terms of the low degree functions is as follows
	
	\begin{proposition}
		The degree $-2$ Poisson bracket on $C_\bullet(T^*[2]E[1])$  given by the Schouten bracket of multivector fields is completely encoded in the following operations:
		\begin{enumerate}
			\item[1.] The canonical pairing $\langle\cdot,\cdot\rangle$ on $E\oplus E^*.$
			\item[2.] The Lie algebroid structure on $\bA_{E\oplus E^*}^{\langle\cdot,\cdot\rangle}=\wedge^2 E\oplus \der(E^*)\oplus \wedge^2 E^*.$
			\item[3.] The Lie algebroid action on $E\oplus E^*$ as skew-symmetric derivations.
		\end{enumerate}
	\end{proposition}
	\begin{proof}
		Since $T^*[2]E[1]$ is a $2$-manifold, the Leibniz rule implies that the Poisson bracket is totally determined by their action on functions of degree $0,1$ and $2$. Since the bracket has degree $-2$ the only non-trivial brackets are 
		$$[\cdot,\cdot]_S:C_1(T^*[2]E[1])\times C_1(T^*[2]E[1])\to C_0(T^*[2]E[1])$$
		that gives the canonical symmetric pairing; the map
		$$[\cdot,\cdot]_S:C_2(T^*[2]E[1])\times C_2(T^*[2]E[1])\to C_2(T^*[2]E[1])$$
		that gives the Lie bracket on $\Gamma(\bA_{E\oplus E^*}^{\langle\cdot,\cdot\rangle})$, 
		$$[\cdot,\cdot]_S:C_2(T^*[2]E[1])\times C_0(T^*[2]E[1])\to C_0(T^*[2]E[1])$$
		corresponding to the anchor of the Lie bracket, that is the symbol map, and finally, the last item on the statement corresponds to the bracket between degree $2$ and degree $1$ functions. 
	\end{proof}
	
	We finish by studying the $Q$-structures on $T^*[2]E[1]$ that are compatible with the Schouten bracket. A way to guarantee that a vector field is Poisson is by assuming it is Hamiltonian. In the next section we will see that all degree $1$ Poisson vector fields are Hamiltonian for this bracket. Therefore, since we want to find a degree $1$ vector field and the bracket has degree $-2$ we need to study the degree $3$ functions on $T^*[2]E[1].$ A direct consequence from \eqref{eq:defcot2} is that
	$$C_3(T^*[2]E[1])=C_3(E[1])\oplus \fX_1^1(E[1])\oplus\fX^2_{-1}(E[1])\oplus\fX^3_{-3}(E[1]).$$
	Hence, we get that $\Theta\in C_3(T^*[2]E[1])$ decomposes into $$\Theta=H+ Q+\pi^2+\pi^3$$ with $H\in C_3(E[1])=\Gamma(\wedge^3 E^*)$, $Q\in\fX_1^1(E[1]),$ $\pi^2\in\fX^2_{-1}(E[1])$ and $\pi^3\in\fX^3_{-3}(E[1])=\Gamma(\wedge^3 E).$
	Moreover, the vector field $Q_{\Theta}=[\Theta,\cdot]_S$ will define a $Q$-structure if and only if 
	$$0=\frac{1}{2}[Q_\Theta,Q_\Theta]=Q_{[\Theta,\Theta]_S}\leftrightharpoons [\Theta,\Theta]_S=0.$$
	The above equation is called \emph{the classical master equation (CME)}. In this case we say that the vector bundle $E\to M$ is a \emph{proto-Lie bialgebroid}, see \cite{ive:quasi, roy:on}. We will see in \S~\ref{sec:cou2} that $(T^*[2]E[1], [\cdot,\cdot]_S, Q_\Theta)$ corresponds to the Courant algebroid known as the \emph{Drinfeld double} of the  proto-Lie bialgebroid $E\to M$.
	
	The more common proto-Lie bialgebroids are the \emph{quasi Lie bialgebroids} that correspond to $H=0$, in this case $\Theta$ satisfying the CME encodes a $Q$-structure on $E[1]$ (given by $Q$) together with a degree $1$ homotopic Poisson structure  on $E[1]$ (given by $\pi^2+\pi^3$), see \S~\ref{sec:mult}. If, in addition, we also have $\pi^3=0$ then we get that $E[1]$ is a degree $1$ $PQ$-manifold and by Corollary \ref{cor:bia} the vector bundle $E$ has a Lie bialgebroid structure on it.

	\section{Symplectic Q-manifolds}\label{sec:sympQ}
	Here, we accomplish two goals simultaneously: we study symplectic $Q$-manifolds and their Lagrangian $Q$-submanifolds, objects that lie in the core of these lecture notes, thereby introducing relevant  examples of $n$-manifolds for arbitrary $n$. Several higher structures relevant in mathematics and physics will appear naturally.
	
	\subsection{Definition and first properties}
	
	A \emph{symplectic $m$-manifold} $(\cM, \omega)$ is 
	an  $m$-manifold $\cM=(M, C_\bullet(\cM))$ together with $\omega\in\Omega^2_m(\cM)$ 
	such that $d\omega=0$ 
	and $\omega$ is non-degenerate, i.e. the map
	\begin{align*}
		\omega^\flat\colon \mathfrak{X}^1_\bullet(\cM)\to \Omega^1_{\bullet+m}(\cM), \quad  X\to \iota_X\omega,
	\end{align*}
	is an isomorphism of $C_\bullet(\cM)$-modules.    
	
	\begin{remark}
		Since we just consider $\bN$-graded manifolds, the fact that the degree of the manifold is equal to the 
		degree of the symplectic form is not a restriction at all, see e.g. \cite{roy:on}.
	\end{remark}

	For every function $f\in C_{k}(\cM)$ we denote by $X_f$ the unique degree $k-m$ vector field, such that
	\begin{equation}\label{Eq: HamVF}
		\omega^\flat(X_f)=df.
	\end{equation} 
	$X_f$ is called the \emph{Hamiltonian} vector field of $f$. 
	Similar to the case $m=0$, Hamiltonian vector fields allow us to define the 
	associated Poisson bracket $\{\cdot,\cdot\}\colon C_i(\cM)\times C_j(\cM)\to C_{i+j-m}(\cM)$ 
	given by the formula
	\begin{equation}
		\{f,g\}=\iota_{X_f}\iota_{X_g}\omega.
	\end{equation}
	
	A vector field $X\in\fX^1_k(\cM)$ is called \emph{symplectic}  if $\Lie_X\omega=0$.
	
	\begin{proposition}\label{bra}
		Let $(\cM,\omega)$ be a  symplectic $m$-manifold. The following equalities hold:
		$$\Lie_{X_f}g=\{f,g\},\qquad \Lie_{X_f}\omega=0\quad \text{ and }\quad X_{\{f,g\}}=[X_f,X_g].$$
		Moreover, $\{\cdot,\cdot\}$ is a Poisson bracket of degree $-m$.
	\end{proposition}
	
	\begin{proof}
		We use Proposition \ref{cartan-cal} and formula \eqref{lieder}:
		\begin{eqnarray*}
			\Lie_{X_f}g&=&(\iota_{X_f}d +d \iota_{X_f})g=\iota_{X_f}d g=\iota_{X_f}\iota_{X_g}\omega=\{f,g\}.\\
			\Lie_{X_f}\omega&=& d \iota_{X_f}\omega+\iota_{X_f}d\omega=0.\\
			\iota_{[X_f,X_g]}\omega&=&[\Lie_{X_f},\iota_{X_g}]\omega=\Lie_{X_f}\iota_{X_g}\omega
			=\Lie_{X_f}d g=d \Lie_{X_f}g=d \{f,g\}=\iota_{X_{\{f,g\}}}\omega.
		\end{eqnarray*}
		Next, we show that the bracket is a Poisson bracket. Since $0=[\iota_X,\iota_Y]=\iota_X\iota_Y-(-1)^{1+|X||Y|}\iota_Y \iota_X$,
		we obtain that $\{\cdot,\cdot\}$ is skew-symmetric.
		The graded Leibniz identity:
		\begin{equation*}
			\begin{array}{rl}
				\{f,gh\}=&\iota_{X_f}\iota_{X_{gh}}\omega=\iota_{X_f}d(gh)=\iota_{X_f}d g \ h+\iota_{X_f}g\ d h\\
				=&\{f,g\}h+(-1)^{|g||X_f|}g\ \iota_{X_f}d h=\{f,g\}h+(-1)^{|g|(|f|-m)}g\{f,h\}.
			\end{array}
		\end{equation*}
		The graded Jacobi identity:
		\begin{align*}
			\{\{f,g\},h\}+(-1)^{(|f|-m)(|g|-m)}\{g,\{f,h\}\}=&\iota_{X_{\{f,g\}}}\iota_{X_h}
			\omega+(-1)^{|X_f||X_g|}\iota_{X_g}\iota_{X_{\{f,h\}}}\omega\\
			=&\iota_{[X_f,X_g]}\iota_{X_h}\omega+(-1)^{|X_f||X_g|}\iota_{X_g}\iota_{[X_f,X_h]}\omega\\
			=&\Lie_{X_f}\iota_{X_g}\iota_{X_h}\omega=\Lie_{X_f}\{g,h\}\\
			=&\{f,\{g,h\}\}.\qedhere
		\end{align*}
	\end{proof}
	
	Two special properties of symplectic $m$-manifolds are stated in the next proposition, which is due to Roytenberg.

	\begin{proposition}[see \cite{roy:on}]\label{sym=ham}
		Let $(\cM,\omega)$ be a symplectic $m$-manifold. The following 
		statements hold:
		\begin{enumerate}
			\item[1.] If $m\geq 1$, then $\omega$ is exact. Moreover $\omega=d(\frac{1}{m} 
			\iota_{\cE_u}\omega)$, where $\cE_u$ is the Euler vector field of $\cM$.
			\item[2.] Let $X\in\fX^1_l(\cM)$ be a symplectic vector field. If $l+m\neq 0$,	then $X$ is the Hamiltonian of $$f=\frac{1}{l+m} \iota_{\cE_u}\iota_X\omega.$$
		\end{enumerate}
	\end{proposition}
	
	\begin{proof}
		Both formulas are straightforward:
		\begin{equation*}
			m\ \omega=\Lie_{\cE_u}\omega=d \iota_{\cE_u}\omega\Rightarrow \omega
			=d (\frac{1}{m}\iota_{\cE_u}\omega).
		\end{equation*}
		\begin{equation*}
			d \iota_{\cE_u}\iota_X\omega=\Lie_{\cE_u}i_X\omega-\iota_{\cE_u}d \iota_X\omega
			=(l+m)\iota_X\omega\Rightarrow \iota_X\omega=d (\frac{1}{m+l}\iota_{\cE_u}\iota_X\omega).\qedhere
		\end{equation*}
	\end{proof}

	There is also a graded version of the classical Darboux Theorem for symplectic $m$-manifolds (see \cite[\S 2.6]{cue:def} for a general definition of the graded cotangent bundles and \cite[Cor 3.8]{cue:def} for a proof) and compare with Theorem 5.3 in \cite{kos:gra} for supermanifolds.
	
	\begin{corollary}[Darboux coordinates]\label{Cor: Darboux}
		Let $(\cM,\omega)$ be a  symplectic $m$-manifold of dimension 
		$k_0|\dots|k_m$ and $p\in M$. Then there exist $U\subseteq M$ a chart of $\cM$ around $p$ 
		and a symplectomorphism $\Phi\colon \cM_{|U}\to \cN$ where
		\begin{equation*}
			\cN=\left\{\begin{array}{ll}
				\text{for}\quad m=2l+1,& (T^*[m]\bR^{k_0|\cdots|k_l}, \omega_{\can});   \\
				\text{for}\quad m=4l,& 
				(T^*[m]\bR^{k_0|\cdots|k_{2l-1}|\frac{k_{2l}}{2}}, \omega_{\can});\\
				\text{for}\quad m=4l+2, &
				(T^*[m]\bR^{k_0|\cdots|k_{2l}}, \omega_{\can})\times (\bR^{0|\cdots|0 | k_{2l+1}}, 
				\sum_{j=1}^{k_{2\ell+1}} \epsilon_jd y^j\wedge d y^j) 
			\end{array}\right.
		\end{equation*}
		with $\epsilon_j=\pm 1.$
	\end{corollary}

	\subsubsection{Symplectic Q-manifolds and their prequantization}
	As it happens for $1$ and $2$-manifolds we are interested in having and addition $Q$-structure on it that must be compatible with the symplectic form as follows.
	
	A \emph{degree $m$ symplectic $Q$-manifold} is a triple $(\cM,\omega,Q)$ where $(\cM,\omega)$ 
	is a symplectic $m$-manifold and $(\cM, Q)$ is a $Q$-manifold such that 
	\begin{align*}
		\Lie_Q\omega=0.
	\end{align*}

	\begin{proposition}
		Let $(\cM,\omega,Q)$ be a degree $m$ symplectic $Q$-manifold with $m\geq 1$. The following hold:
		\begin{enumerate}
			\item[1.] The vector field $Q$ is Hamiltonian, so is\begin{align*}
				Q=X_\Theta \ \quad  \text{ with } \quad \Theta \in C_{m+1}(\cM) 
				\quad  \text{ and } \quad \{\Theta,\Theta\}=0. 
			\end{align*}	
			\item[2.] The differential forms on $\cM$ are a double complex $(\Omega^\bullet_\bullet(\cM),\Lie_Q, d)$ with total differential $D=\Lie_Q+(-1)^id$, where $i$ denotes the internal degree. Moreover 
			$$\omega=D((-1)^m\lambda-\Theta)\quad\text{with }\quad \lambda=\frac{1}{m}\iota_{\cE_{u}}\omega.$$
		\end{enumerate}
	\end{proposition}
	
	\begin{proof}
		From Proposition \ref{sym=ham}, we know that $Q$ is Hamiltonian with $\Theta\in C_{m+1}(\cM)$. Moreover, 
		by Proposition \ref{bra},  $0=[Q, Q]=X_{\{\Theta,\Theta\}},$ thus $d\{\Theta,\Theta\}=0$ and 
		due to the degree $\{\Theta,\Theta\}=0$.
		
		By the Cartan calculus of Proposition 
		\ref{cartan-cal} we obtain that $(\Omega^\bullet_\bullet(\cM), \Lie_Q, d)$ becomes a 
		double complex and since 
		\begin{align*}
			\Lie_Q\lambda=&\Lie_Q(\frac{1}{m}\iota_{\cE_u}\omega)=\frac{1}{m} d\iota_Q\iota_{\cE_u}\omega+\frac{1}{m}\iota_Qd\iota_{\cE_u}\omega\\
			=&\frac{-1}{m}d\iota_{\cE_u}\iota_Q\omega+\frac{1}{m}\iota_Q\omega=-\frac{m+1}{m}d\Theta+\frac{1}{m}d\Theta=-d\Theta.
		\end{align*}
		We get that $\omega=D((-1)^m\lambda-\Theta).$
	\end{proof}
	
	An immediate interesting consequence of the previous results is that for $m\geq 1$, degree $m$ symplectic $Q$-manifolds always admit a \emph{prequantization $\bR[m]$-bundle with connection}, see e.g. \cite{frs:hig, fss:hig}. 
	
	\begin{proposition}\label{prop:prequantum}
		Let $(\cM, \omega=d(-1)^m\lambda, Q=X_\Theta)$ be a degree $m$ symplectic $Q$-manifold with $m\geq 1$. Then a prequantum $Q$-bundle with connection is given by $\cP=\cM\times \bR[m]$ together with $$\tau=(-1)^m(\lambda+d\hbar)\in\Omega^1_m(\cP)\quad\text{and}\quad \widetilde{Q}=Q\times\Theta\frac{\partial}{\partial \hbar}$$ where $\hbar$ denotes the coordinate in $\bR[m]$. 
	\end{proposition}
	\begin{proof}
		First observe that $(\cP,\widetilde{Q})\to(\cM,Q)$ is a $Q$-bundle \cite{ks:qbun} because $\widetilde{Q}^2=Q^2+Q(\Theta)=0.$ Moreover $\tau$ is a connection $1$-form with curvature $\omega$ because
		\begin{equation*}
			d\tau=d\big((-1)^m(\lambda+d\hbar)\big)=\omega\quad\text{and}\quad \Lie_{\widetilde{Q}}\tau=(-1)^m\big(\Lie_Q\lambda+d\Lie_{\widetilde{Q}}\hbar\big)=(-1)^m(-d\Theta+d\Theta)=0.\qedhere
		\end{equation*}
	\end{proof}
	
	We finish by mentioning that the above prequantum bundle is connected to the \emph{graded contact $Q$-manifolds} introduced in \cite{raj:con}.

	\subsection{Lagrangian Q-submanifolds}\label{sec:lagq}
	
	Just as happens in standard symplectic geometry, Lagrangian submanifolds play an important role in describing fine geometrical and algebraic features of graded symplectic manifolds. This is due to the fact that graded symplectic manifolds also enjoy a very simple normal form around Lagrangian submanifolds, as evidenced by a graded version of the Weinstein’s Lagrangian tubular neighborhood.

	Let $(\mathcal{M}, \omega)$ be a symplectic $m$-manifold  and let $j\colon \mathcal{N}\hookrightarrow \mathcal{M}$ be a submanifold with 
	vanishing ideal $\mathcal{I}_\mathcal{N}$. We say that $\mathcal{N}$ is:
	\begin{enumerate}
		\item[1.] \emph{isotropic} if $j^*\omega=0$,
		\item[2.] \emph{coisotropic}  if 	$\{\mathcal{I}_\mathcal{N},\mathcal{I}_\mathcal{N}\}\subseteq \mathcal{I}_\mathcal{N}$, and
		\item[3.] \emph{Lagrangian} if it is coisotropic with $\totdim(\mathcal{N})=\frac{1}{2}\totdim(\mathcal{M})$. 
	\end{enumerate}
	Furthermore, if $(\mathcal{M}, \omega, Q)$ is a degree $m$ symplectic $Q$-manifold then we say that $\mathcal{N}$ is a 
	\emph{Lagrangian $Q$-submanifold} if $Q$ preserves its vanishing ideal.
	
	\begin{remark}
		
		There are two special issues in the context of graded symplectic manifolds which do not pop up in standard symplectic geometry. Firstly, the conditions on tangent spaces for isotropic and coisotropic submanifolds are not sufficient. For instance, let $\mathcal{M}=\mathbb{R}^{1|4|1}$ be the $2$-manifold with coordinates
		$$\{x, q^1, q^2, p^1, p^2, y\} \quad \textnormal{of degree}\quad |x|=0, |q^1|=|q^2|=|p^1|=|p^2|=1, |y|=2$$
		together with  the symplectic form given by $\omega=\textnormal{d} y\wedge\textnormal{d} x+ \textnormal{d} p^1\wedge\textnormal{d} q^1+\textnormal{d} p^2\wedge\textnormal{d} q^2.$ Consider 
		the submanifold $j\colon \mathcal{N}=\mathbb{R}^{1|2|0}\hookrightarrow \mathcal{M}$ given by the vanishing ideal $\mathcal{I}_\mathcal{N}=\langle y-p^1q^2,\ p^2,\ q^1\rangle$. The submanifold $\mathcal{N}$ is not isotropic, because $j^*\omega=(q^2\textnormal{d} p^1-p^1 \textnormal{d} q^2)\wedge \textnormal{d} x\neq 0$ 
		but we get $T_p\mathcal{N}\subseteq T_p\mathcal{N}^\omega$ for all $p\in N=\mathbb{R}$. Additionally, $\mathcal{N}$ is also not coisotropic, since $\{ y-p^1q^2, q^1\}=\pm q^2$. However, $T_p\mathcal{N}^\omega\subseteq T_p\mathcal{N}$ for all $p\in M=\mathbb{R}$.
		
		Secondly, not all graded symplectic manifolds admit Lagrangian submanifolds. Indeed, let us consider the graded $2$-manifold $\mathcal{M}=\mathbb{R}^{0|1|0}$ with coordinate $\{ e\}$ of degree 
		$1$ and symplectic form $\omega= \textnormal{d} e\wedge \textnormal{d} e$. Since $\totdim{\mathcal{M}}=1$ we deduce that $(\mathcal{M}, \omega)$ does not admit Lagrangian submanifolds. As a trivial consequence of the local description provided in Corollary \ref{Cor: Darboux} we  also deduce that if $(\mathcal{M},\omega)$ is a symplectic $m$-manifold that does not admit any (local) Lagrangian submanifolds then $m=4l+2$.   
	\end{remark}
	
	We are now ready to state the graded version of Weinstein's tubular neighborhood theorem for Lagrangian submanifolds. In the context of supermanifolds the statement already appeared  in \cite{aksz, sch:bv}.

	\begin{theorem}[\cite{cue:def}]\label{Thm: dgWeinstein}
		Let $(\mathcal{M}=(M,C_\bullet (\mathcal{M})), \omega)$ be a symplectic $m$-manifold and let $i\colon \mathcal{L}=(L,C_\bullet(\mathcal{L}))\hookrightarrow \mathcal{M}$ be a Lagrangian submanifold. Then, there exist an 
		open neighborhood $U$ of $L$ in $M$ and a symplectomorphism $\Psi\colon T^\ast[m] \mathcal{L}\to \mathcal{M}|_{U}$ 
		such that the following diagram commutes:
		\begin{center}
			\begin{tikzcd}
				T^\ast[m]\mathcal{L} \arrow[rr, "\Psi"]&&\mathcal{M}|_{U} \\
				&\mathcal{L}\ar[ur,hookrightarrow, "i"]\ar[lu,hookrightarrow, "0_{\mathcal{L}}"']&
			\end{tikzcd}
		\end{center}
	\end{theorem}
	
	The proof of the previous theorem relies on applying two key results which are interesting in their own right:
	\begin{enumerate}
		\item[1.] A graded version of the tubular neighborhood theorem for general graded submanifolds, saying that if $j\colon \mathcal{N}\to \mathcal{M}$ is a submanifold  then there is a tubular neighborhood $N\subseteq U\subseteq M$ and a diffeomorphism $\Psi\colon \mathcal{M}|_{U}\to j^* T \cM/T \cN$ such that $j$ is identified with the zero-section.
		\item[2.] A Darboux theorem along submanifolds.  Let $j\colon N\hookrightarrow M$ be a submanifold of the body, $\omega_0,\omega_1\in \Omega^2_m(\mathcal{M})$ two symplectic structures 
		and   $\omega_t\in \Omega^2_m(\mathcal{M})$ a smooth path of closed $2$-forms connecting $\omega_0$ and $\omega_1$. If $(\omega_t)_p$ is non-degenerate for all $p\in N$ then there exist an open neighborhood $N\subseteq U\subseteq M$ and a symplectomorphism  $\Phi\colon (\mathcal{M}|_U,\omega_0)\to (\mathcal{M}|_U,\omega_1)$. In the context of symplectic supermanifolds this is partially discussed 
		in \cite{sch:bv, sch:semi}.
	\end{enumerate}

	\subsubsection{Deformations and homotopic Poisson structures}\label{sec:deflag}
	An application of the above result is the study of deformations of Lagrangian $Q$-submanifolds done in \cite{cue:def}. Let $(\cM,\omega,X_\Theta)$ be a degree $m$ symplectic $Q$-manifold and $i:\cL\to(\cM,\omega,X_\Theta)$ be a Lagrangian $Q$-submanifold. In this case one constructs an $L_\infty$-algebra on $C_\bullet(\cL)[m-1]$ whose Maurer-Cartan elements modulo gauge  control the formal deformations of $\cL$ as a Lagrangian  $Q$-submanifold modulo flows of $Q$-exact Hamiltonian vector fields.
	
	Let us outline how to construct the $L_\infty$-algebra on $C_\bullet(\cL)[m-1]$. 
	For simplicity, assume $\cL$ is a wide submanifold (so the only degree $0$ coordinates on $T^*[m]\cL$ are the ones coming from $\cL$). By Theorem \ref{Thm: dgWeinstein} there exists a symplectomorphism $\Psi:(T^*[m]\cL,\omega_{\can})\to(\cM, \omega)$ and since 
	$$C_\bullet(T^*[m]\cL)={}^{-m}\fX_\bullet^\bullet(\cL),$$
	we can look into
	$$\Psi^*\Theta=Q_{\cL}+\pi^2+\cdots\pi^{m+1}\in {}^{-m}\fX_{m+1}^\bullet(\cL)=C_{m+1}(T^*[m]\cL).$$
	The classical master equation for $\Theta$ implies that $$Q_\cL^2=0\quad\text{and} \quad [Q_\cL,\pi^\bullet]_S+\frac{1}{2}[\pi^\bullet,\pi^\bullet]_S=0$$ so $\pi^\bullet$ is an $(m-1)$ homotopic Poisson structure on $\cL$, as defined in \S~\ref{sec:mult}, and its functions, $C_\bullet(\cL)[m-1]$, inherits an  $L_\infty$-algebra structure. Observe that this result also gives a geometric way of recovering the $(m-1)$-shifted Poisson structure underlying a Lagrangian $Q$-submanifold\footnote{That is a particular kind of shifted Lagrangian structure as in \S~\ref{sec:ssa}.}, see \cite{saf:lec}.
	
	\begin{remark}
		Notice that the $L_\infty$-algebra on $C_\bullet(\cL)[m-1]$ is a Lie $m$-algebra. This feature is due to the fact that we started with a wide submanifold. In the general case, one usually gets all the brackets, see \cite{cue:def} for details. 
	\end{remark}
	
	\subsection{Examples}\label{sec:exaQs}
	In the following we give several relevant examples of degree $m$ symplectic $Q$-manifolds for several $m$. Observe that for $m=0$ we recover the usual symplectic manifolds $(M,\omega)$ and this case is kind of special due to the fact that the form does not need to be exact. 
	
	\subsubsection{Poisson manifolds}\label{sec:1sym=poi}
	
	Recall that a $1$-manifold must be of the form $E[1]$ for some vector bundle $E\to M$, see Theorem \ref{thm: 1-manifolds}. Now assume that $\omega\in\Omega^2_1(E[1])$ is a symplectic form, then 
	$$\omega^\flat:\fX^1_{-1}(E[1])\overset{\sim}{\to}\Omega^1_0(E[1])$$
	is a vector bundle isomorphism. Proposition \ref{neg-gen-vf} implies that $\fX^1_{-1}(E[1])=\Gamma E$ and a degree counting argument shows that $\Omega^1_0(E[1])=\Omega^1(M)$,  therefore $E=T^*M$ and $\omega=\omega_\can\in\Omega^2_1(T^*[1]M).$ Moreover, by Proposition \ref{sym=ham} we get that any symplectic $Q$-structure is given by a degree $2$ function satisfying the classical master equation. Hence, since $\omega=\omega_\can$ the Poisson bracket of $\omega$ is given by the Schouten bracket on $\fX^\bullet(M)$ and $\Theta=\pi\in C_2(T^*[1]M)=\fX^2(M)$ satisfying $[\pi,\pi]_S=0$ is exactly a Poisson structure on $M$.
	
	Submanifolds of $T^*[1]M$ are the same as  subbundles $E\to N\subseteq T^*M\to M$. The submanifold $E[1]\subseteq T^*[1]M$ is Lagrangian for the canonical symplectic form if and only if $E=\ann(TN)$. Finally, $Q_\pi$ will preserve $\ann(TN)[1]$ if and only if $\ann(TN)\to N\subset T^*M\to M$ is a Lie subalgebroid and this, by definition, means that $N\subseteq (M, \pi)$ is a \emph{coisotropic submanifold for the Poisson structure}, see \cite[\S 5]{cat:intpoi} for details. In summary we have proven the following result
	
	\begin{proposition}[\cite{cat:intpoi}]\label{prop:poi1sym}
		There is a one-to-one correspondence between
		\begin{equation*}
			\left\{\begin{array}{c}
				\text{Poisson manifolds}  \\
				(M,\pi)
			\end{array}\right\} \leftrightharpoons \left\{\begin{array}{c}
				\text{Deg 1 symplectic Q-manifolds}  \\
				(T^*[1]M, \omega_{\can}, Q_\pi)
			\end{array}\right\} 
		\end{equation*}
		Moreover, this correspondence also establishes an equivalence between coisotropic submanifolds $C\subseteq M$ and Lagrangian $Q$-submanifolds $\ann(TC)[1]\subseteq T^*[1]M.$
	\end{proposition}
	
	The prequantum $Q$-bundle corresponding to $(T^*[1]M, \omega_{\can}, Q_\pi)$ as in Proposition \ref{prop:prequantum} gives exactly the Lie algebroid structure on $T^*M\oplus\bR_M\to M$ corresponding to seeing the Poisson structure $\pi$ as a Jacobi manifold. The study of those extensions and their relation to prequantization was carried out in \cite{cc:jacobi}.
	
	Finally, we note that the $L_\infty$-algebra governing deformations of coisotropic submanifolds in Poisson manifolds was constructed in \cite{cf:def} and it give a degree $0$ homotopic Poisson structure on the Lie algebroid of the coisotropic submanifold.
	
	\subsubsection{Courant algebroids}\label{sec:cou2}
	In  \S~\ref{ex:kw} we have shown how given a vector bundle with a symmetric and non-degenerate pairing, $(E\to M, \langle\cdot,\cdot\rangle)$, one constructs a $2$-manifold whose functions are given by the Keller-Waldmann algebra $\cM=(M, C_\bullet(E, \langle\cdot,\cdot\rangle)).$ It was shown is \cite[Theorem 3.18]{kw:def} that the Keller-Waldmann algebra carries a Poisson bracket of degree $-2$ that indeed is non-degenerate. 
	
	Conversely, given $(\cM, \omega)$ a symplectic $2$-manifold we know that $C_0(\cM)=M, \ C_1(\cM)=\Gamma (E)$ for some vector bundle and 
	$$\{\cdot,\cdot\}:C_1(\cM)\times C_1(\cM)\to C_0(\cM) $$
	gives the symmetric and non-degenerate pairing. Indeed, the above procedures are inverses of each other, and we recover the following result.
	
	\begin{proposition}[Theorem 3.3 in \cite{roy:on}]
		There is a one-to-one correspondence between 
		
		\begin{equation*}
			\left\{\begin{array}{c}
				\text{Vector bundles with a symmetric}\\
				\text{and non-degenerate pairing}  \\
				(E\to M, \langle\cdot,\cdot\rangle)
			\end{array}\right\} \leftrightharpoons \left\{\begin{array}{c}
				\text{Symplectic 2-manifolds}  \\
				\cM=\big(M, C_\bullet(E, \langle\cdot,\cdot\rangle)\big)
			\end{array}\right\} 
		\end{equation*}   
	\end{proposition}
	A full proof of this result using the functor $\cF$ of Theorem \ref{thm:gem2-2} is in \cite[\S 3]{bur:red}. The next step is to add the $Q$-structure. For that we need the following definition introduced in \cite{lwx:ca}.
	
	A \emph{Courant algebroid} $(E\to M,\langle \cdot, \cdot \rangle, [\![\cdot,\cdot]\!], \rho)$ is a vector bundle $E\to M$ equipped with a symmetric and non-degenerate pairing $\langle \cdot , \cdot\rangle\colon \Gamma(E)\times \Gamma(E)\to C^\infty(M)$, a bracket $[\![-,-]\!]\colon \Gamma(E)\times \Gamma(E)\to \Gamma(E)$ and a vector bundle morphism $\rho\colon E\to TM$ called the anchor satisfying the following properties:
	\begin{enumerate}
		\item[(C1)] $  [\![e_1, fe_2]\!] = f [\![ e_1, e_2]\!] + \rho(e_1)(f)e_2$;
		\item[(C2)] $\rho(e_1)\langle e_2, e_3\rangle = \langle  [\![ e_1, e_2]\!], e_3\rangle + \langle e_2,  [\![ e_1, e_3]\!] \rangle $;
		\item[(C3)] $ [\![ e_1,  [\![ e_2, e_3 ]\!] ]\!] =  [\![  [\![ e_1, e_2 ]\!], e_3 ]\!] +  [\![ e_2,  [\![ e_1, e_3 ]\!] ]\!]$;
		\item[(C4)] $ [\![e_1, e_2 ]\!] +  [\![ e_2, e_1 ]\!] = D d\langle e_1, e_2 \rangle$,
	\end{enumerate}
	where $e_1, e_2, e_3\in \Gamma(E)$, $f\in C^\infty(M)$ and $D= \langle -, - \rangle ^\flat \circ \rho^\ast \colon T^\ast M \to E$. Given a Courant algebroid $(E\to M,\langle \cdot, \cdot \rangle, [\![\cdot,\cdot]\!], \rho)$, a \emph{Dirac structure with support} is a subbundle $L\to N\subset E\to M$ satisfying
	\begin{enumerate}
		\item [(D1)] $L$ is isotropic for the pairing and $\rank L=\frac{1}{2}\rank E$;
		\item [(D2)] $\rho (L)\subseteq TN$ and given $e_1, e_2\in\Gamma (E)$ such that ${e_1}_{|N}, {e_2}_{|N}\in\Gamma(L)$ then $[\![ e_1, e_2]\!]_{|N}\in\Gamma(L).$
	\end{enumerate}
	
	\begin{theorem}[\cite{roy:on, sev:some}]\label{thm:sevroycorr}
		There is a one-to-one correspondence between
		\begin{equation*}
			\left\{\begin{array}{c}
				\text{Courant algebroids}  \\
				(E\to M,\langle \cdot, \cdot \rangle, [\![\cdot,\cdot]\!], \rho)
			\end{array}\right\} \leftrightharpoons \left\{\begin{array}{c}
				\text{Deg 2 symplectic Q-manifolds}  \\
				(\cM,\omega,Q_\Theta)
			\end{array}\right\} 
		\end{equation*}
		Moreover, this correspondence also establishes an equivalence between Dirac structures with support and Lagrangian $Q$-submanifolds.
	\end{theorem}
	\begin{proof}
		Given a Courant algebroid $(E\to M,\langle \cdot, \cdot \rangle, [\![\cdot,\cdot]\!], \rho),$ two different ways for constructing the $Q$-structure are as follows:
		\begin{enumerate}
			\item[1.] Define $\Theta\in C_3(E,\langle\cdot,\cdot\rangle)$ by 
			\begin{equation}\label{eq:thcou}
				\Theta(e_1,e_2,e_3)=\langle[\![e_1,e_2]\!],e_3\rangle,\quad e_i\in\Gamma (E)
			\end{equation}
			and show that $\{\Theta,\Theta\}=0$.
			\item[2.] Define directly $Q:C_k(E,\langle\cdot,\cdot\rangle)\to C_{k+1}(E,\langle\cdot,\cdot\rangle)$ by the Cartan-like formula
			\begin{align*}
				(Q \eta)(e_0, \dots, e_k) =& \sum_{i=0}^k (-1)^i \rho(e_i)\eta(e_0, \dots, \widehat{e}_i, \dots, e_k)  \\
				&- \sum_{i < j} (-1)^i \eta(e_0, \dots, \widehat{e}_i, \dots, e_{j-1}, [\![e_i,e_j]\!], e_{j+1}, \dots, e_k)
			\end{align*}
			and then check that it is a derivation of the cup product and $Q^2=0.$
		\end{enumerate}
		For the converse direction one uses the derived bracket formalism, as in the original proof, and obtains
		\begin{equation}\label{eq:derCA}
			\rho(e)(f)=-\{f,\{e,\Theta\}\},\quad\text{and}\quad [\![e,e']\!]=\{e',\{e,\Theta\}\}
		\end{equation}
		for $f\in C^\infty(M)=C_0(E,\langle\cdot,\cdot\rangle)$ and $e,e'\in\Gamma(E)=C_1(E,\langle\cdot,\cdot\rangle).$
		Taking some care, one shows that the classical master equation implies the axioms (C1)-(C4) of a Courant algebroid.
		
		Let us give some ideas for the moreover part, we point to \cite{roy:on} and to \cite{bur:red} for details. Recall from Proposition \ref{sub-ideal} that a submanifold $\cN\subseteq\cM$ is the same as a subsheaf of regular ideals $\cI_\cN\subseteq C_\bullet(\cM)=C_\bullet(E,\langle\cdot,\cdot\rangle).$ In particular, the regular condition implies
		\begin{itemize}
			\item  $(\cI_\cN)_0=Z(N)=\{ f\in C^\infty(M)\ |\ f_{|N}=0\}$ for some submanifold $N\subseteq M$,
			\item $(\cI_\cN)_1=\{e\in\Gamma (E)\ | e_{|N}\in\Gamma (L)\}$ for some subbundle $L\to N\subseteq E\to M,$ 
			\item $(\cI_\cN)_2=\{\cdo\in\Gamma(\bA_{E}^{\langle\cdot,\cdot\rangle})\ |\ \cdo_{|N}\in\Gamma (A)\}$ for some subbundle $A\to N\subseteq \bA_{E}^{\langle\cdot,\cdot\rangle}\to M$ fitting in the exact sequence
			$$0\to L\wedge E_{|N}\to A\to F\to 0$$
			with $F\subseteq TM_{|N}$ a subbundle.
		\end{itemize}
		The submanifold $\cN$ is coisotropic if and only if $\cI_\cN$ is a Poisson ideal, i.e. $\{\cI_\cN,\cI_\cN\}\subseteq\cI_\cN$. The coisotropic condition is equivalent to (see \cite[Theorem 4.5]{bur:red} for details):
		\begin{itemize}
			\item $L\to N$ is an isotropic subbundle,
			\item $F\subseteq TN$ is a regular and involutive distribution,
			\item $A$ is determined by a  flat, metric $F$-connection $\nabla$ on the vector bundle $L^\perp/L\to N.$
		\end{itemize}
		And $\cN$ will be Lagrangian if, in addition, $$\totdim\cN=\frac{1}{2}\totdim\cM,$$ which is equivalent to $F=TN$ and $L^\perp=L$ (so the connection disappear) and $\rank L=\frac{1}{2}\rank E$, so we get (D1).
		
		Finally, $\cN$ is a $Q$-submanifold if and only if $$\Big(Q(\cI_\cN)=\{\Theta,\cI_\cN\}\Big)\subseteq\cI_\cN.$$
		Using the formulas in \eqref{eq:derCA} and the fact that $L$ is maximally isotropic, one immediately recognizes that the above condition is equivalent to (D2) as we want. Observe that the implications are if and only if, therefore given a Dirac structure with support one defines a sheave of homogeneous regular ideals $\cI_\cN\subseteq C_\bullet(E,\langle\cdot,\cdot\rangle)$ that will define a Lagrangian $Q$-submanifold.
	\end{proof}
	
	\begin{example}\label{ex:CA}
		In these lecture notes we encounter two main examples of Courant algebroids:
		\begin{enumerate}
			\item[1.] {\bf Quadratic Lie algebras:} Recall that a Lie algebra $(\g,[\cdot,\cdot])$ is called \emph{quadratic} if it has a symmetric and non-degenerate pairing $\langle\cdot,\cdot\rangle:\g\times\g\to \bR$ satisfying
			\begin{equation}\label{eq:adinv}
				\langle[e_1,e_2],e_3\rangle+\langle e_2,[e_1,e_3]\rangle=0
			\end{equation}
			for $e_i\in\g.$ One can directly check that a Courant algebroid for which the base manifold $M=pt$ is the same as a quadratic Lie algebra. The corresponding degree $2$ symplectic $Q$-manifold is
			$$(\g[1], \varpi, Q=d_{CE})$$
			where $Q$ is the Chevalley-Eilenberg differential (see \ref{ex:Q-TM-g}) and the degree $2$ symplectic form $\varpi\in\Omega^2_2(\g[1])=\text{S}^2\g^*$ is given by the pairing. Notice that $d\varpi=0$ is equivalent to the pairing having constant coefficients, the non-degeneracy of $\varpi$ is  equivalent to the  non-degeneracy for the pairing and, finally, the equation $\Lie_Q\varpi=0$ is equivalent to \eqref{eq:adinv}.
			\item[2.] {\bf Drinfeld doubles of proto-Lie bialgebroids:} Recall from \S~\ref{sec:liebia} that a proto-Lie bialgebroid structure on a vector bundle $E\to M$ is a function $\Theta\in C_3(T^*[2]E[1])$ satisfying the classical master equation 
			$$[\Theta,\Theta]_S=0.$$
			Therefore, $(T^*[2]E[1],\omega_\can, Q=X_\Theta)$ defines a degree $2$ symplectic $Q$-manifold. Using the results from \S~\ref{sec:liebia} together with Theorem \ref{thm:sevroycorr}, it is easy to see that one obtains a Courant algebroid  
			$$(E\oplus E^*\to M, \langle\cdot,\cdot\rangle, \rho, [\![\cdot,\cdot]\!])$$ where the pairing is the canonical one, and anchor and bracket are given by the formulas in \eqref{eq:derCA}. This Courant algebroid is known as the \emph{Drinfeld double of the proto-Lie bialgebroid $E\to M$}, see \cite{lwx:ca, roy:on}. 
			
			Moreover, the degree $2$ symplectic manifold $(T^*[2]E[1],\omega_\can, X_\Theta)$ has two natural Lagrangian submanifolds: the zero section, $0_{E[1]}:E[1]\to T^*[2]E[1]$, and the fiber at zero $j:E^*[1]\to T^*[2]E[1]$. However, the vector field $X_\Theta$ is not, in general, tangent to either them. If $X_\Theta$ is tangent to both submanifolds, then $E\to M$ is a Lie bialgebroid; if it is tangent only to the zero section, then  $E\to M$ is a quasi-Lie bialgebroid. 
		\end{enumerate}  
		Other relevant examples of Courant algebroids appearing in the literature are transitive Courant algebroids \cite{bh:tranca, gfrt:hol} or action Courant algebroids \cite{lbm:ca}, among others.
	\end{example}
	
	The above point of view for Courant algebroids and Dirac structures has some important consequences including:
	\begin{enumerate}
		\item[1.] Courant algebroids have a cohomology given by $H_Q(\cM).$ Indeed, one can also construct a Cartan calculus for Courant algebroids on the Keller-Waldmann algebra, see \cite{cue:cou}.
		\item[2.] As one performs reduction of symplectic manifolds using Hamiltonian actions, one can do the same for Courant algebroids, see \cite{bur:red}.
		\item[3.] The prequantum bundle of Proposition \ref{prop:prequantum} was studied for Courant algebroids in \cite{roy:shl}. In particular, when the Courant algebroid is a quadratic Lie algebra $(\g,[\cdot,\cdot],\langle\cdot,\cdot\rangle)$ the prequantum bundle is given by the string Lie $2$-algebra 
		\begin{equation}\label{eq:stringalg}
			\text{String}(\g)=(\bR\to\g,\partial=0, [\cdot,\cdot], [\cdot,\cdot,\cdot]_3=\langle[\cdot,\cdot],\cdot\rangle).
		\end{equation}
		Observe that in this case $C_\bullet(\cP=\g[1]\times\bR[2])=C_\bullet(\g[1])[[\hbar]]=\wedge^\bullet\g[[\hbar]]$ and a deformation quantization of the Poisson algebra $(\wedge^\bullet\g,\wedge,\langle\cdot,\cdot\rangle)$  in the direction of the pairing is given by the Clifford algebra, \cite{lei:cli, mei:clif}. Indeed one can also quantize the $Q$-structure in this picture as in \cite[\S 7]{mei:clif}. There are several other ways of performing a quantization of Courant algebroids using the symplectic $Q$-manifold point of view, including \cite{xu:weyl, kw:def, roy:aksz}.
		\item[4.] The Lie $3$-algebra controlling deformations of Dirac structures (without support) appeared in \cite{lwx:ca} and was extensively studied in \cite{gua:def}. For Dirac structures with support see \cite{cue:def}.
	\end{enumerate}
	
	\begin{remark}
		The case of degree $3$ symplectic $Q$-manifolds is less well studied than the two previous ones. Nevertheless some works dealing with it are \cite{gruz:H-twist, ikeda:qp3, she:clwx}. 
	\end{remark}

	\subsubsection{Higher Courant algebroids as graded cotangent bundles}\label{sec:higcou}
	
	The classification of degree $m$ symplectic $Q$-manifolds stops at $m=3$, nevertheless other examples of higher symplectic $Q$-manifolds had appeared in the literature. In what follows we explain the case of  higher Courant algebroids and their interpretation as cotangent bundles. 
	
	Consider the triple $\mathfrak{A}=(A[1], Q, H)$ where $A[1]$ is a $1$-manifold, $Q\in\fX^1_1(A[1])$ is a $Q$-structure and $H\in C_{m+1}(A[1])$ with $Q(H)=0.$ Then we can form the symplectic $Q$-manifold 
	$$\mathfrak{B}=(T^*[m]A[1],\omega_\can, X_{\Theta}) \quad \text{with} \quad \Theta=H+Q\in C_{m+1}(T^*[m]A[1])={}^{-m}\fX^\bullet_{m+1}(A[1]).$$
	Observe that $\Theta$ satisfies the classical master equation because $Q^2=0$ and $Q(H)=0.$ Moreover, if $m>3$ a degree counting argument shows that any $\Theta\in C_{m+1}(T^*[m]A[1])$ satisfying the classical master equation is as above.
	
	By Theorem \ref{thm:Q-lie} the data encoded in $\mathfrak{A}$ is a Lie algebroid $(A\to M,\rho,[\cdot,\cdot]_A)$ with an $m+1$ Lie algebroid cohomology class represented by $H$. Using the general version of the functor $\cF$ given in Theorem \ref{thm:gem2-2}, it was shown in \cite{cue:cot} that the data encoded in $\mathfrak{B}$ corresponds to the \emph{higher Courant algebroid}
	$$(A\oplus\wedge^{m-1}A^*,\langle\cdot,\cdot\rangle,\rho\circ \pr_A, [\![\cdot,\cdot]\!])$$
	where $\langle a\oplus\omega,b\oplus\eta\rangle=\iota_a\eta+\iota_b\omega\in\wedge^{m-2}A^*$ and the bracket is given by
	\begin{equation}\label{eq:brahigc}
		[\![a\oplus\omega,b\oplus\eta]\!]_H=[a,b]_A\oplus\cL_a\eta- \iota_bQ(\omega) - \iota_b\iota_aH\in\Gamma(A\oplus\wedge^{m-1}A^*),
	\end{equation}
	where $\cL_a$ and $Q$ are defined as in Proposition \ref{prop:cartcal}.
	
	The correspondence between $\mathfrak{B}$ and higher Courant algebroids works in the same way as in Theorem \ref{thm:sevroycorr}. Roughly speaking, it proceeds as follows: the direction from the higher Courant algebroid to $\mathfrak{B}$ is straightforward, since all the necessary data are readily available. Conversely, given $\mathfrak{B}$, we obtain the higher Courant algebroid by looking into
	\begin{eqnarray*}
		C_{m-1}(T^*[m]A[1])&=&\fX_{-1}^1(A[1])\oplus C_{m-1}(A[1])=\Gamma(A\oplus \wedge^{m-1}A^*), \\ 
		C_{m-2}(T^*[m]A[1])&=&C_{m-2}(A[1])=\Gamma(\wedge^{m-2}A^*).
	\end{eqnarray*}
	Hence, using the same formulas as in the Courant algebroid, one gets the pairing anchor and bracket by
	\begin{eqnarray*}
		\langle \cdot,\cdot\rangle&=&[\cdot,\cdot]_S:C_{m-1}(T^*[m]A[1])\times C_{m-1}(T^*[m]A[1])\to C_{m-2}(T^*[m]A[1]);\\
		\rho(a+\omega)(f)&=&[f,[a+\omega,\Theta]_S]_S;\\
		\left[\![a\oplus\omega,b\oplus\eta \right]\!]_H &=&[a\oplus\omega,[b\oplus\eta,\Theta]_S]_S;
	\end{eqnarray*}
	where $f\in C^\infty(M)=C_0(T^*[m]A[1])$ and $a\oplus\omega,b\oplus\eta\in C_{m-1}(T^*[m]A[1]).$

	If one studies Lagrangian $Q$-submanifolds of $\mathfrak{B}$ one gets, see \cite[Theorem 4.8]{cue:cot}, that they are the same as 
	$(L\to N)\subseteq A\oplus\wedge^{m-1}A^*$ satisfying
	\begin{enumerate}
		\item [(HD0)] $\pr_A(L) \subseteq A$ is a subbundle and $L \cap \wedge^{m-1} A^*_{|N}= \ann(\pr_A(L)) \wedge(\wedge^{m-2}A^*)_{|N};$ 
		\item [(HD1)] $\langle L,L\rangle\subseteq \ann(\pr_A(L))\wedge(\wedge^{m-3}A^*)_{|N}$;
		\item [(HD2)] $\rho(L) \subseteq T N$ and if $e, e'\in\Gamma(A\oplus\wedge^{k-1}A^*)$ such that $e_{|N}, e'_{|N}\in\Gamma (L)$ then $[\![e,e']\!]_{H|N}\in\Gamma(L).$
	\end{enumerate}
	
	Compared with the definition of a Dirac structure with support we get that $(HD1)$ is similar to $(D1)$ and $(HD2)$ is analogous to $(D2)$ but the condition $(HD0)$ makes this new structure more rigid. The above conditions appeared previously in the work \cite{hag:nam} and when $A=TM$ and $H=0$ they encode graphs of closed $m$-forms and graphs of Nambu tensors (that are higher analogues of Poisson structures).
	
	Finally let us mention that the connection between higher Courant algebroids and graded cotangent bundles was hinted in the works \cite{bou:aksz, zam:hdir} but without giving a systematic study.
	
	\subsubsection{Extended Courant algebroids} In the recent work \cite{ghu:std} a new kind of Courant algebroids appeared in connection with spherical $T$-duality. Here, we sketch their relation with symplectic $Q$-manifolds.

	Let $\pi: E\to M$ be an oriented $\mathbb{S}^{2n-1}$-bundle and $H\in \Omega^{2(n+k)-1}(E)$ with $d_EH=0$. As in  \cite{ghu:std}, we use the following notation:
	\begin{enumerate}
		\item[1.] A \emph{representative of the Euler class of $E$} is $\varepsilon\in \Omega^{2n}(M)$;
		\item[2.] A \emph{global angular form for $\varepsilon$} is $\psi\in\Omega^{2n-1}(E)$ such that $d_E\psi=\varepsilon$.
		\item[3.] A \emph{decomposition of $H$} is $H_0\in\Omega^{2(n+k)-1}(M)$ and $H_1\in\Omega^{2k}(M)$ such that $H=\pi^*H_0+\psi\wedge\pi^*H_1$ as in \cite[Thm 4.2]{ghu:std}. Observe that $d_E H=0$ becomes
		\begin{equation*}
			d_M H_1=0,\qquad d_MH_0+\varepsilon\wedge H_1=0.
		\end{equation*}
	\end{enumerate}
	With the above data we introduce the $Q$-manifold  
	$$\cP=T[1]M\times \bR[2n-1]\qquad\text{with}\qquad Q_\varepsilon=d_M+\varepsilon\frac{\partial}{\partial u},$$
	where $u$ is the coordinate in $\bR[2n-1]$,  and the $Q$-morphism  
	$$\Phi_\psi:(T[1]E, d_E)\to (\cP, Q_\varepsilon), \quad \Phi^\sharp(\omega_0+u\omega_1)=\pi^*\omega_0+\psi\wedge\pi^*\omega_1\quad \text{for}\quad \omega_0+u\omega_1\in C_k(\cP),$$
	that induces an isomorphism in the $Q$-cohomology and satisfies $\Phi^\sharp(H_0+uH_1)=H$.
	
	This leads directly to the following results:
	\begin{enumerate}
		\item[1.] By \S~\ref{sec:higcou}, the higher Courant algebroid $TE\oplus \wedge^{2(k+n)-3}T^*E$ with $H$-twisted bracket is given by the symplectic $Q$-manifold $(\cM=T^*[2(k+n-1)]T[1]E, \omega_\can, X_{\Theta_\cM})$  with 
		$$\Theta_\cM= H+d_E\in C_{2(k+n)-1}(\cM).$$
		\item[2.] The higher cotangent bundle $(\cN=T^*[2(k+n-1)]\cP, \omega_{\can}, X_{\Theta_\cN})$ with  
		$$\Theta_\cN=H_0+uH_1+Q_\varepsilon\in C_{2(k+n)-1}(\cN)$$
		is also a symplectic $Q$-manifold. 
		\item[3.] The conormal bundle to the graph of $\Phi_\psi$ $$j:N^*[2(k+n-1)]\Gamma_{\Phi_\psi}\to \cM\times \cN= T^*[2(k+n-1)](T[1]E\times\cP)$$
		is a Lagrangian $Q$-submanifold.
	\end{enumerate}
	Using the general version of Theorem \ref{thm:gem2-2} given in \cite{bur:frob}, the symplectic $Q$-manifold $(\cN, \omega_\can, X_{\Theta_\cN})$ is equivalent to the \emph{extended Courant algebroid} $(C_\psi, \langle\cdot,\cdot\rangle, [\![\cdot,\cdot]\!]_H)$ with $H$-twisted bracket as introduced in \cite[\S 8]{ghu:std}. More concretely, one gets $C_{2(k+n)-3}(\cN)=\Gamma(C_\psi)$ where
	$$C_\psi=T M \oplus (\wedge^{2n-2}T^*M \otimes\langle\partial\psi\rangle)\oplus \wedge^{2(n+k)-3}T^*M\oplus (\langle\psi\rangle \otimes \wedge^{2k-2}T^*M).
	$$
	The pairing on $C_\psi$ corresponds to the Poisson bracket
	$$[\cdot,\cdot]_S:C_{2(k+n)-3}(\cN)\times C_{2(k+n)-3}(\cN)\to C_{2(k+n)-4}(\cN)$$
	and the Courant-like bracket in Theorem \cite[Theorem 8.2]{ghu:std} is given by the derived bracket formula
	$$[\![e,e']\!]_H=[e',[e,\Theta_\cN]_S]_S\quad\text{for}\quad e,e'\in C_{2(k+n)-3}(\cN).$$
	
	\subsection{On the Homotopic version}\label{sec:ssa}
	Let $(\cM, Q)$ be a $Q$-manifold of degree $n$. As explained in Remark 
	\ref{tangentcomplex}, we have that $T_p\cM$ and $T_p^*\cM$ are chain complexes. In 
	\cite{pym:shift} an \emph{$m$-shifted symplectic Lie  $n$-algebroid} is defined as a degree $n$ $Q$-manifold $(\cM, Q)$ 
	endowed with $\omega\in\Omega^2_m(\cM)$ such that 
	$$\Lie_Q\omega=0,\quad d\omega=D(\sum_{i=1}^m \beta^i)\quad \text{for some}\ \beta^i\in\Omega^{2+i}_{m-i}(\cM)
	\quad \text{and the map}\quad \omega_p\colon T_p\cM\to (T^*_p\cM)[m]$$
	is a quasi-isomorphism for all $p\in M$. Therefore, the symplectic $Q$-manifolds that we studied in this section are a strict version of the shifted symplectic Lie $n$-algebroids introduced in \cite{pym:shift}. This generalization will play an important role in Part \ref{partII}. Observe that one can unify the two form with the higher forms $\beta^i$ into a single element $\omega^\bullet$ and then the two first equations simplify to $D(\omega^\bullet)=0.$
	
	Analogously to the above, one introduces \emph{shifted Lagrangian structures}. Let $(\cM, Q, \omega^\bullet)$ be an $m$-shifted symplectic Lie $n$-algebroid and let $\Phi:(\cL,Q_\cL)\to (\cM,Q)$ be a $Q$-morphism. A \emph{Lagrangian structure on $\Phi$} consists of $\eta^i\in\Omega^{2+i}_{m-i-1}(\cL)$ for $i=\{0,\cdots, m-1\}$ satisfying
	$$\Phi^*\omega^\bullet=D\eta^\bullet\quad\text{and}\quad \eta^0_p\oplus (T^*_p\Phi\circ\omega^0_p):T_p\cL\oplus T_p\cM[-1]\to T^*_p\cL[m-1]$$
	is a quasi-isomorphism for all $p\in L,$ the body of $\cL$. The study of shifted Lagrangian structures for $m=1,2$ appeared in \cite{cro:red, pym:shift}.
	
	As an illustration, let us show the result found in \cite{pym:shift} that $1$-shifted symplectic structures on the $Q$-manifold $(E[1],Q)$ are in one-to-one correspondence with Dirac structures on exact Courant algebroid of the base. Therefore, this explicitly exhibits the difference between symplectic $Q$-manifolds (which, in degree $1$, correspond to Poisson manifolds) and their homotopic versions.  
	
	\begin{example}\label{ex:dir}
	  By Theorem \ref{thm:Q-lie} we know that $(E[1],Q)$ is the same as a Lie algebroid structure on $(E\to M,[\cdot,\cdot],\rho)$. Also observe that for each $p\in M,$ the tangent complex is given by $$T_pE[1]= E_p\xrightarrow{\rho_p}T_pM.$$
	 Then, a closed $1$-shifted symplectic structure on $(E[1],Q)$ is the same as an
	 	$$\omega\in\Omega^2_1(E[1])\quad\text{together with}\quad H\in\Omega^3_0(E[1])=\Omega^3(M)$$ satisfying
	 	\begin{equation}\label{imformsdir}
	 		\Lie_Q\omega=0,\quad \Lie_QH=d\omega,\quad dH=0\quad\text{and}\quad \omega_p:T_pE[1]\xrightarrow{\sim} T^*_pE[1].
	 \end{equation} 
	 It was shown in \cite{bcwz:dir}, see also \cite{bc:2form}, that the first two equations implies $\omega$ is an IM $2$-form completely determined by the evaluations $\omega_p$. In other words, if we denote the induced vector bundle map by $\mu:E\to T^*M$, then it satisfy 
	 \begin{equation}	\label{imformsdir2}
	 	\iota_{\rho(e_1)}\mu(e_2)+\iota_{\rho(e_2)}\mu(e_1)=0\quad\text{and}\quad \mu([e_1,e_2])=\Lie_{\rho(e_1)}\mu(e_2)-\iota_{\rho(e_2)}d\mu(e_1)-\iota_{\rho(e_1)}\iota_{\rho(e_2)}H.
	 \end{equation}
	 	 The last condition in \eqref{imformsdir} implies that
	 	 $$0\to E\xrightarrow{(\rho,\mu)}TM\oplus T^*M\xrightarrow{-\mu^*+\rho^*}E^*\to 0$$
	 	 is an exact sequence. While condition $dH=0$ together with the equations \eqref{imformsdir2} implies that $(\rho,\mu):E\to TM\oplus T^*M$ becomes a Dirac structure on the standard Courant algebroid twisted by $H$, see Example \ref{ex:CA}.
	\end{example}

	\subsection{AKSZ sigma models}\label{sec:aksz}
	
	Perhaps the most important aspect of symplectic $Q$-manifolds is that given a degree $m$ symplectic $Q$-manifold $(\cM,\omega, Q),$ there is an associated $(m+1)$-dimensional Topological Quantum Field Theory (TQFT) constructed in  \cite{aksz} and today known as the \emph{AKSZ sigma model with target $(\cM,\omega, Q)$}. This field theory is defined using the \emph{Batalin-Vilkovisky formalism}, see \cite{mnev:book} for a comprehensive introduction. In the following we shortly recall the AKSZ transgression procedure as well as the examples of the Poisson sigma model \cite{stro:poi} and $3d$-Chern-Simons theory \cite{wit:jon}. For a more detailed exposition see e.g. \cite{cmr:cla, mnev:book}. In this section we consider $\bZ$-graded manifolds instead of $\bN$-graded.

	The \emph{BV formalism} is a device oriented towards quantization
	of field theories with generalized gauge symmetries through path integrals with additional superfields. In the BV setting, we have a space of superfields $\mathfrak{F}$ (typically a formal infinite dimensional $\bZ$-graded manifold) endowed with $\omega_{\mathfrak{F}}\in\Omega^2_{-1}(\mathfrak{F})$ a symplectic structure and a (formal) Berezinian volume $\mu$ satisfying the compatibility condition encoded in the notion of \emph{SP-structure}, see \cite{sch:bv}. The
	key quantum computations are of the form
	\begin{equation}\label{eq:correlation}
		\langle O\rangle=\int_{\mathfrak{L}\subset\mathfrak{F}}\sqrt{\mu}\ O\  e^{\frac{i}{\hbar}S}
	\end{equation}
	where $\mathfrak{L}\subset\mathfrak{F}$ is a Lagrangian submanifold implementing a \emph{gauge fixing} condition and  $S,O\in C_0(\mathfrak{F})[\![\hbar]\!]$ are functions representing the \emph{BV-action} and an observable. 
	
	The key fact, which can be rigorously proven when $\mathfrak{F}$ is a finite-dimensional $SP$-manifold (see \cite{sch:bv}), is that the integral $\langle O\rangle$ is stable under deformations of the gauge fixing $\mathfrak{L}\subset\mathfrak{F}$ when
	the following two conditions are met:
	\begin{enumerate}
		\item[1.]  $S$ satisfies the Quantum Master Equation (QME): $\{S, S\} - 2i\hbar\Delta S = 0$
		\item[2.] $O$ is a quantum observable: $\{S, O\} - i\hbar\Delta O = 0.$
	\end{enumerate}
	Here $\Delta$ is the so-called \emph{BV-laplacian} defined by the $SP$-structure as $\Delta f = (-1)^{|f|} 2 \cD_\mu(X_f )$, c.f. with the generating operator in Proposition \ref{prop:genop}. It is worth noting that the above two conditions combine
	to give $\Delta(O \ e^{\frac{i}{\hbar}S})=0$ which is the general requirement for the stability of the integral under
	deformations of $\mathfrak{L}$.
	
	The leading order in $\hbar$ in the QME yields the Classical Master Equation, $\{S, S\} = 0$. Forgetting about the (formal) measure $\mu$, the triple $(\mathfrak{F}, \omega_\mathfrak{F}, S)$ thus defines a $PQ$-structure that is referred to as a \emph{classical BV-manifold}.

	\subsubsection{The AKSZ transgression}\label{sec:akszT}
	
	In \cite{aksz} it was shown how using a transgression procedure one can obtain a classical $BV$-manifold from the following data:
	\begin{itemize}
		\item A degree $m$ symplectic $Q$-manifold $(\cM,\omega, X_\Theta)$.
		\item A closed and oriented $(m+1)$-dimensional manifold $N$.
	\end{itemize}
	The manifold $N$ is known as the \emph{source manifold or spacetime} while $(\cM,\omega, X_\Theta)$ is called the \emph{target manifold}. In the following we give the construction of the AKSZ procedure, see also \cite{aksz, mnev:book}.

	Recall that the category of $\bZ$-graded vector spaces is monoidal with respect to the tensor product given by
	$$({\bf V}\otimes{\bf W})_i=\bigoplus_{j+k=i}V_j\otimes W_k.$$
	Therefore, one defines a \emph{mapping space or internal hom} as a right adjoint to the tensor product functor, i.e. for each pair of $\bZ$-graded vector spaces ${\bf V}, {\bf W}$ there is another $\bZ$-graded vector space $\underline{\Hom}({\bf V},{\bf W})$ defined by the property
	$$\Hom({\bf U}\otimes{\bf V}, {\bf W})\cong \Hom({\bf U},\underline{\Hom}({\bf V},{\bf W})) \qquad \forall\  {\bf U}\ \bZ\text{-graded vector space},$$
	that is natural in all the three arguments. In this case, we can explicitly compute it by the formula
	$$\underline{\Hom}_i({\bf V},{\bf W})=\Hom({\bf V},{\bf W}[i]).$$

	We see in Example \ref{car-pro} that the category of $\bZ$-graded manifolds has a cartesian product, i.e. given $\cM$ and $\cN$ one gets $\cM\times\cN$ which is again a $\bZ$-graded manifold. Therefore we construct a mapping space in the category of $\bZ$-graded manifolds using the universal property
	$$\Hom(\cZ\times \cN, \cM)\cong \Hom(\cZ,\text{Maps}(\cN,\cM)) \qquad \forall\  \cZ\ \bZ\text{-graded manifold}.$$
	In our examples $\cN=T[1]N$ for some manifold $N$. Thus, in the special case when $\cM={\bf V}[1]$ is a graded vector space, one has 
	$$\text{Maps}(T[1]N, {\bf V}[1])=C_\bullet(T[1]N)\otimes {\bf V}[1]=\Omega^\bullet(N)\otimes {\bf V}[1].$$ For a general $\bZ$-graded manifold, mappings whose range fits in a single coordinate chart can be described by this construction. More care is required to describe mappings whose
	body is a non-trivial morphism with range outside a single coordinate chart.

	The other important feature of the mapping space that we need is the \emph{evaluation map}
	$$ev:\cN\times\text{Maps}(\cN,\cM)\to\cM.$$
	Choose $\{u^i, \theta^i=du^i\}$ local coordinates on $T[1]N$ and $\{x_j^{\alpha_j}\}$ coordinates on $\cM$. Then, the evaluation map is completely determined by the formula
	$$ev^*x_j^{\alpha_j}=\mathbb{X}^{\alpha_j}=\sum_{r=0}^{\dim N}\sum_{i_1<\cdots<i_r}X^{\alpha_j}_{i_1\cdots i_r}(u)\theta^{i_1}\cdots\theta^{i_r},$$
	where $X^{\alpha_j}_{i_1\cdots i_r}(u)\in C^\infty(N).$ The  $\mathbb{X}^{\alpha_j}$ is the \emph{superfield} associated to the
	target coordinate $x_j^{\alpha_j}$ and  the functions $X^{\alpha_j}_{i_1\cdots i_r}(u)$ give coordinates on the
	mapping space of degree $j-r$. 
	
	Recall that when $N$ is a closed and oriented $(m+1)$-dimensional manifold, then $T[1]N$ has a canonical Berezinian defined by the orientation, and let $\mathfrak{F}=\text{Maps}(T[1]N,\cM).$ Then \cite{aksz} shows how to use the evaluation map and the Berezinian on $T[1]N$ to produce the transgression map
	\begin{equation*}
		\begin{array}{c}
			\xymatrix{ T[1]N\times \text{Maps}(T[1]N,\cM)\ar[d]^\pr\ar[r]^{\qquad\qquad ev}& \cM\\\text{Maps}(T[1]N,\cM) }
		\end{array}
		\qquad \bT_N:\Omega^\bullet_\bullet(\cM)\to \Omega^\bullet_{\bullet-m-1}(\mathfrak{F}), \quad \bT_N(\omega)=\int_{T[1]N}ev^*\omega
	\end{equation*}
	satisfying the following.
	\begin{proposition}[\cite{aksz}]
		Let $(\cM,\omega=d\lambda, Q_\Theta)$ be a degree $m$ symplectic $Q$-manifold and $N$ be a closed and oriented $(m+1)$-dimensional manifold. Then $\mathfrak{F}=\text{Maps}(T[1]N,\cM)$ is a classical BV-manifold with
		$$\omega_{\mathfrak{F}}=\bT_N(\omega)\in\Omega^2_{-1}(\mathfrak{F})\quad\text{and}\quad  S=\bT_N(\Theta)+i_{\widehat{d_N}}\bT_N(\lambda)\in C_{0}(\mathfrak{F}),$$
		where $\widehat{d_N}$ is the lift of the de Rham vector field on $T[1]N$ to $\mathfrak{F}$, seen as infinitesimal
		transformations of the source.
	\end{proposition}
	
	\begin{remark}
		If $\partial N\neq\emptyset$ the above proposition fails to hold. A way of restoring the result is by imposing boundary conditions as in \cite{cat:kon}, another option is to implement the BV-BFV formalism \cite{cmr:cla}. It was shown in  \cite{cat:kon} that appropriate boundary conditions are given by  $$i:\cL\to \cM \quad \text{Lagrangian }Q\text{-submanifold with}\quad  i^*\lambda=0.$$ 
		For a more general setting see \cite{cal:lag}.
	\end{remark}
	
	\subsubsection{The Poisson sigma model and star products}\label{sec:PSM}
	Let $(M, \pi)$ be a Poisson manifold, by Proposition \ref{prop:poi1sym} we get a degree $1$ symplectic $Q$-manifold $(T^*[1]M,\ \omega_{\can}=d\lambda_{\can},\ \Theta=\pi)$. Using coordinates $\{y^i\}$ on $M$, a coordinate description of the degree $1$ symplectic $Q$-manifold is
	$$C^\infty(T^*[1]M)\supset \{y^i, v_i\} \quad \text{ with degree }\quad  |y^i|=0, \ |v_i|=1 \text{ then }$$ 
	$$ \omega_{\can}=dv_i dy^i,\quad \lambda_{\can}=v_i d y^i,\quad Q_\Theta=\pi^{ij}v_i\frac{\partial}{\partial y^j}+\frac{1}{2}\partial_i\pi^{jk}v_jv_k\frac{\partial}{\partial v_i}, \quad \Theta=\frac{1}{2}\pi^{ij}v_iv_j.$$
	
	Let $\Sigma$ be a closed and oriented 2-dimensional manifold.   The induced space of fields for the Poisson sigma model $\mathfrak{F}=Map(T[1]\Sigma, T^*[1]M)$ can be described by the coordinate superfields 
	\begin{equation*}
		\begin{array}{llc}
			ev^*y^i=\mathbb{Y}^i=Y^i_0+Y^i_1+Y^i_2, &|Y^i_j|=-j,& Y_0:\Sigma\to M, \quad Y_j\in\Omega^j(\Sigma; Y^*_0TM),\\
			ev^*v_i=\mathbb{V}_i=
			V^0_i+ V^1_i+V^2_i, &|V^j_i|=1-j, &V^j\in\Omega^j(\Sigma; Y^*_0T^*M).
		\end{array}
	\end{equation*}
	Let us describe the components. The fields with $|\cdot|=0$ are called \emph{classical fields}, here $Y_0$ and $V_1$. \emph{Gauge symmetries} are encoded by fields with degree $|\cdot|=1$, which in this case correspond to $V^0$. The other three are called \emph{antifields}, they are conjugate to the former with respect to the symplectic form 
	$$ \omega_\mathfrak{F}=\int_{T[1]\Sigma}d \mathbb{V}_id \mathbb{Y}^i=\int_\Sigma dV_i^0 dY^i_2+ dV_i^1 dY^i_1+dV_i^2 dY^i_0.$$  
	The action $S_{PSM}$ is locally given by
	\begin{equation*}
		\begin{split}
			S_{PSM}(\mathbb{Y, V})=&\int_{T[1]\Sigma} \mathbb{V}_i  d_{\Sigma}\mathbb{Y}^i+\frac{1}{2}\pi^{ij}(\mathbb{Y})\mathbb{V}_i\mathbb{V}_j\\
			=&\int_\Sigma\Big( V^1_i  d_{\Sigma}Y^i_0+\frac{1}{2}\pi^{ij}(Y_0)V_i^1V_j^1\Big)
			+ V_j^2\pi^{ij}(Y_0)V_j^0+Y^i_1\Big(d_{\Sigma}V^0_i+\partial_i\pi^{jk}(Y_0)V_j^0V_k^1\Big)\\
			&+\frac{1}{2}Y^i_2\partial_i\pi^{jk}(Y_0)V_j^0V_k^0+\frac{1}{4}Y^i_1Y^j_1\partial_i\partial_j\pi^{kl}(Y_0)V_k^0V_l^0.\end{split}
	\end{equation*}
	The first two terms give the classical action of the Poisson sigma model while the other terms give the gauge transformations, the bracket of gauge transformations and the last term is  the homotopy for the defect of integrability of the gauge
	symmetry as a distribution on the space of classical fields.

	Perhaps the main application of the Poisson sigma model is that the Kontsevich star product\footnote{See \cite{bffls} for the deformation quantization program.}  deforming the point-wise product on $C^\infty(M)$ in the direction of the Poisson bracket \cite{kon:def} is given by an observable in the Poisson sigma model described in \cite{cat:kon}.
	
	\begin{theorem}[\cite{cat:kon}]
		Let $(M,\pi)$ be a Poisson manifold and $p\in M$. Consider a disk $D$ and denote by $\mathfrak{F}=\{ \text{Maps}(T[1]D, T^*[1]M) \ | \ \mathbb{V}_{|\partial D}=0\}$. Then the formula for the Kontsevich star product is 
		$$(f\star g)(p)=\int_{\mathfrak{L}\subset\mathfrak{F}}f(\mathbb{Y}(0))g(\mathbb{Y}(1))\delta_p(\mathbb{Y}(\infty)) \ e^{\frac{i}{\hbar}S_{PSM}(\mathbb{Y, V})}$$
		where $0, 1$ and $\infty$ are three distinct points on the boundary of $D$ and $\mathfrak{L}$ is an appropriate Lagrangian gauge fixing.
	\end{theorem}


	\subsubsection{Chern-Simons and knot invariants}\label{sec:CS}
	In Theorem \ref{thm:sevroycorr} we saw that degree $2$ symplectic $Q$-manifolds are the same as Courant algebroids. We saw in Example \ref{ex:CA} that a Courant algebroid with $M=pt$  corresponds to a \emph{quadratic Lie algebra} $(\g,[\cdot,\cdot],\langle\cdot,\cdot\rangle)$. The corresponding degree $2$ symplectic $Q$-manifold is $(\g[1],\varpi=d\lambda, \Theta)$ where $\varpi\in\Omega^2_2(\g[1])$ is given by the pairing and $\Theta\in C_3(\g[1])$ is as in \eqref{eq:thcou}.
	
	Using a base $\{a_i\}$ of $\g$, a coordinate description of the degree $2$ symplectic $Q$-manifold is
	$$C^\infty(\g[1])\supset \{a^i\} \quad \text{ with degree }\quad  |a^i|=1 \text{ then }$$ 
	$$ \varpi=\frac{1}{2}g_{ij}da^i da^j,\quad \lambda=\frac{1}{2}g_{ij}a^i d a^j,\quad Q_\Theta=\frac{1}{2}c_{jk}^ia^ja^k\frac{\partial}{\partial a^i}, \quad \Theta=\frac{1}{6}T_{ijk}a^ia^ja^k.$$
	where $\langle a_i, a_j\rangle=g_{ij},\ [a_j, a_k]=c^i_{jk}a_i$ and $T_{ijk}=g_{il}c^l_{jk}.$
	
	Let $N$ be a closed and oriented $3$-dimensional manifold.   The induced space of fields for Chern-Simons is $\mathfrak{F}=Map(T[1]N, \g[1])=\Omega^\bullet(N)\otimes \g[1]$ and can be described by the coordinate superfield 
	\begin{equation*}
		ev^*a^i=\mathbb{A}^i=A^i_0+A^i_1+A^i_2+A^i_3,\quad  |A^i_j|=1-j, \quad A_j\in\Omega^j(N; \g),\\
	\end{equation*}
	Let us describe the components. The classical field is $A_1=A$, gauge symmetries are encoded by fields with degree $|\cdot|=1$, which in this case correspond to $A^0=c$. The other two are the antifields and they are conjugate to the former with respect to the symplectic form, denoting $A_3=c^+$ and $A_2=A^+$ we get
	$$ \omega_\mathfrak{F}=\int_{T[1]N} \frac{1}{2}g_{ij}\ d \mathbb{A}^id \mathbb{A}^j=\int_N \frac{1}{2}g_{ij}\big(dA_0^i dA^i_3+ dA^i_1 dA^i_2\big)=\int_N \frac{1}{2}g_{ij}\big(dc^i dc^{+ i}+ dA^i dA^{+i}\big).$$  
	The action $S_{CS}$ can be written globally as
	\begin{equation*}
		\begin{split}
			S_{CS}(\mathbb{A})=&\int_{T[1]N} \frac{1}{2}\langle\mathbb{A},  d_{N}\mathbb{A}\rangle+\frac{1}{6}\langle\mathbb{A},[\mathbb{A},\mathbb{A}]\rangle\\
			=&\int_N\frac{1}{2}\langle A,  d_{N}A\rangle+\frac{1}{6}\langle A,[A,A]\rangle
			+\langle A^+, d_Nc+[A,c]\rangle+\frac{1}{2}\langle c^+,[c,c]\rangle.\end{split}
	\end{equation*}
	The first two terms give the classical action of Chern-Simons theory  while the other terms give the gauge transformations and the bracket of gauge transformations.
	
	In $\S$ \ref{sec:bg} we will show how Chern-Simons theory helps to integrate the symplectic $Q$-manifold $(\g[1],\varpi,Q_\Theta)$ to a $2$-shifted symplectic Lie $2$-group. The above is a shadow from the relation between Chern-Simons and the symplectic structure on the moduli space of flat connections on a closed and oriented surface explained in \cite{ab:ym}. Finally, we want to mention that \cite{wit:jon} found a way of computing knot invariants by using  observables in Chern-Simons known as \emph{Wilson loops}. Given $\gamma:S^1\to N$, and $R$, an irreducible representation of $G$, one defines an observable on Chern-Simons by 
	$$W_R(\gamma)= Hol^R_\gamma(A)=Tr_R Pexp(\int_\gamma A)
	.$$
	It was shown in \cite{wit:jon} how \eqref{eq:correlation}, with $O$ given by $W_R(\gamma)$ and $S=S_{CS}$, can be used to compute knot invariants such as the Jones polynomial.

	\part{The global picture}\label{partII}
	Topologists of the early 20th century envisioned a generalization  of the non-abelian fundamental group to higher dimensions. In order to study relative homotopy groups, the work \cite{whi:on} introduced the term \emph{crossed
		module} to describe the properties of the boundary map
	$$\partial: \pi_2(X,A,x)\to  \pi_1(A,x).$$
	
	In the 1970s, several works, see e.g. \cite{bro:dou, ehr:cat, higg:cat, pra:jet}, introduced groupoids and double groupoids in geometry and  topology and established connections with  the crossed modules of \cite{whi:on}. The elements of a double groupoid can be visualized as squares that can be composed in two distinct ways, vertically and horizontally, and both compositions must satisfy an exchange law, as suggested by the diagram \eqref{diagdou}.
	\begin{equation}\label{diagdou}
		\begin{array}{c}
			\begin{tikzpicture}
				\draw (0,0) rectangle (1,1);
				\draw (0,1.2) rectangle (1,2.2);
				\draw (1.2,0) rectangle (2.2,1);
				\draw (1.2,1.2) rectangle (2.2,2.2);
				\draw[->] (2.4,0.5) -- (2.8,0.5);
				\draw[->] (2.4,1.7) -- (2.8,1.7);
				\draw[->] (0.5,-0.2) -- (0.5,-0.6);
				\draw[->] (1.7,-0.2) -- (1.7,-0.6);
				\draw (3,1.2) rectangle (5,2.2);
				\draw (3,0) rectangle (5,1);
				\draw (1.2,-0.8) rectangle (2.2,-2.8);
				\draw (0,-0.8) rectangle (1,-2.8);
				\draw[->] (2.4,-1.8) -- (2.8,-1.8);
				\draw[->] (4,-0.2) -- (4,-0.6);
				\draw (3,-0.8) rectangle (5,-2.8);
			\end{tikzpicture}
		\end{array}
	\end{equation}
	In the 1990s, a categorical approach introduced non-strict and higher dimensional versions of double groupoids \cite{ bc:lie2, bc:2grou, sta:shla}. However, it soon became clear that this point of view presented two main problems. On the one hand, it was difficult to determine the appropriate definition among the many weaker notions available in the literature. On the other hand, the  definitions for higher versions involved too many axioms to manage effectively. 
	
	Parallel to those developments, homotopy theorists soon realized, see e.g. \cite{lane:book, may:simp}, that simplicial methods provide a good framework to deal with higher structures. Roughly speaking, a simplicial set $X_\bullet$ is a collection of sets $X_0, X_1, X_2,\cdots,$ together with maps $$d_i:X_l\to X_{l-1}\quad \text{and}\quad s_i:X_l\to X_{l+1},$$ 
	called face and degeneracy maps, respectively. One thinks of $X_0$ as the set of points on the space, $X_1$ as arrows, $X_2$ as oriented triangles, and, in general, $X_l$ as the $l$-simplices of the space. The maps $d_i$ specify the faces of these simplices, while the degeneracy maps $s_i$ describe how lower-dimensional simplices are included as degenerate simplices within higher-dimensional ones.

	In these notes, we follow  \cite{get:int, hen:int, cc:n-grou} and introduce Lie $n$-groupoids as simplicial manifolds satisfying Kan conditions (\S~\ref{sec:lien}). A priori, this definition is not very intuitive, and it may be not immediately clear why the Kan conditions provide a meaningful generalization of Lie groupoids. To motivate the above definition, let us highlight a few key properties that complete Diagram \eqref{diagall}:
	\begin{itemize}
		\item The theory of double Lie groupoids is subsumed by that of  Lie $2$-groupoids (\S~\ref{sec:doug}).
		\item Lie $n$-groupoids admit a Lie functor whose target is the category of  degree $n$ $Q$-manifolds introduced in Part \ref{partI} (\S~\ref{sec:liefun}). This functor serves as a bridge to the formalism used in mathematical physics.
		\item A homotopy theory of Lie $n$-groupoids was developed in \cite{bg:hig, rcc:hom}, establishing a connection   with the theory of differentiable $n$-stacks (\S~\ref{sec:morita}) which is commonly employed in algebraic geometry \cite{cptvv, ptvv, toe:der}.
	\end{itemize}
	
	After briefly introducing Lie $n$-groupoids, we go for the main topic of Part \ref{partII}: $m$-shifted symplectic structures on Lie $n$-groupoids, as introduced in \cite{get:slides}, along with their shifted Lagrangian structures. We will show that $m$-shifted symplectic structures can be transported along Morita equivalences (see Theorem \ref{thm:sme}). Therefore, the $m$-shifted symplectic structures on Lie $n$-groupoids studied in these lecture notes provide concrete presentations of the smooth analog of the shifted symplectic structures introduced in \cite{ptvv}.
	
	Rather than developing the general theory, our focus is on concrete examples and applications in Poisson geometry, higher structures, and their connection with topological field theories. We also emphasize the relation between symplectic double groupoids and $2$-shifted symplectic structures, as well as the  importance of having multiple presentations for the same shifted symplectic stack.
	
	In the final section, following \cite{cal:lag, ptvv}, we introduce shifted Lagrangian structures for the above  $m$-shifted symplectic Lie $n$-groupoids. Our presentation is based on \cite{bur:lag}. One of the main applications of shifted Lagrangian structures is that they encode many of the moment map theories previously developed in the Poisson geometry community (\S~\ref{sec:mome}). Furthermore, the reduction procedures associated with these moment map theories can be understood as Lagrangian compositions in the ``symplectic category" proposed in \cite{wei:symcat}.  
	
	In this second part, we omit many proofs and refer the reader to the original sources. Our goal is to offer a comprehensive introduction to the emerging field of shifted symplectic structures on higher Lie groupoids and their connections with diverse areas of mathematics and physics.

	\section{Higher Lie groupoids}\label{sec:5}
	
	In this section, we briefly review the theory of simplicial manifolds and their geometric realizations. In particular, they  allow us to present a generalized version of the de Rham theorem, as stated in Theorem \ref{thm:genderham}. We then adopt the approach outlined in \cite{get:int, hen:int, cc:n-grou} to introduce higher Lie groupoids as simplicial manifolds that satisfy the Kan conditions. Finally, we discuss the arrows shown in Diagram \eqref{diagall}. More specifically, we introduce the Morita theory of Lie $n$-groupoids, which provides a connection to differentiable $n$-stacks, and present the Lie functor that relates Lie $n$-groupoids to the degree $n$ $Q$-manifolds introduced in Part \ref{partI}.

	\subsection{Simplicial objects}\label{sec:bsimp} For a more detailed introduction see \cite{gh:simp, hov:mod, may:simp}. Denote by $\Delta$ the category of finite ordinals
	\[
	[0]=\{0\}, \quad [1]=\{0,1\}, \quad \dots ,\ [l]=\{0,1, \dots, l\}, \quad \dots, 
	\]
	with order-preserving maps. Among all of the maps appearing in $\Delta$ there are special ones, namely for $0\leq i\leq l$ we define the \emph{cofaces} and \emph{codegeneracies}
	$$\delta^i_l: [l - 1] \to [l]  \quad\text{and}\quad \sigma^i_{l+1} : [l + 1] \to  [l]$$
	where $\delta^i_l$ skips the $i$-th position and $\sigma^i_{l+1}$ repeats the $i$-th position. It
	is an exercise to show that these maps satisfy a list of identities, called the \emph{cosimplicial identities}\footnote{Dual to the ones in \eqref{eq:simpid}.}. The maps $\delta^i_l$, $\sigma^i_{l+1}$ and these relations can be viewed as a set of generators and
	relations for $\Delta$, see \cite[\S 2]{eml:on}.

	Let $\mathsf{C}$ be a category. A \emph{simplicial object  in $\mathsf{C}$} is a  functor $X_\bullet:\Delta^{\text{op}}\to \mathsf{C}$. By the above description of the category $\Delta$, such a functor $X_\bullet$ consists of a tower of objects $X_l$ in $\mathsf{C}$ given by the image of $[l]$ under the functor, and for $0\leq i\leq l$, morphisms in $\mathsf{C}$ $$d_i^l\colon X_l\to X_{l-1}\quad \text{and}\quad s_i^l\colon X_l\to X_{l+1}$$  satisfying the following \emph{simplicial identities}
	\begin{equation}\label{eq:simpid}
		\left\{\def\arraystretch{1.4}
		\begin{array}{ll}
			d_i^{l-1}d_j^l= d_{j-1}^{l-1}d_i^l &\text{if}\quad i< j; \\
			s_i^ls_j^{l-1} = s_{j+1}^l s_i^{l-1} & \text{if}\quad  i \leq j;\\
			d_i^ls_j^{l-1} = 
			s_{j-2}^{l-1}d_i^{l-1} & \text{if}\quad i<j; \\
			d_j^ls_j^{l-1}=\operatorname{id}=d_{j+1}^ls_j^{l-1} &  \\
			d_i^ls_j^{l-1}=s_j^{l-2} d_{i-1}^{l-1} & \text{if}\quad  i>j+1.
		\end{array}\right.
	\end{equation}
	
	In what follows we will omit the upper indices and just write $d_i$ and $s_i$ if there is no risk of confusion. The morphisms $d_i$ and $s_i$ are called the \emph{face} and \emph{degeneracy maps} respectively. 
	
	A morphism $\Phi_\bullet:X_\bullet\to Y_\bullet$ between two simplicial objects in $\mathsf{C}$ is a collection of maps $\Phi_l:X_l\to Y_l$ that commute with all the faces and degeneracy maps. We denote the set of simplicial morphisms by $\Hom(X_\bullet, Y_\bullet).$

	Some relevant examples of simplicial sets, i.e. simplicial objects in the category of sets, are the following.
	\begin{example}[Simplicial simplex]
		Let $l$ be a non-negative integer. The \emph{$l$-simplex} $\Delta[l]$ is given by
		\[
		(\Delta[l])_k=\left\{u_0 \dots u_k \, | \, u_j\leq u_{j+1}, \, u_j\in \{0, \dots, l\}\right\},
		\]
		with faces and degeneracies given by
		\begin{align*}
			d_j(u_0 \dots u_k)&= u_0\dots \hat{u}_j\dots u_k\quad \text{and}\quad
			s_j(u_0\dots u_k)= u_0\dots u_ju_j\dots u_k
		\end{align*}
		where the hat means omission. For instance, when $l=2$, $$(\Delta[2])_0=\{0, 1, 2\},\quad  (\Delta[2])_1=\{00, 01, 02, 11, 12, 22\},\ \cdots.$$
		Observe that for another simplicial set $X_\bullet$ we get 
		$$\Hom(\Delta[l]_\bullet,X_\bullet)=X_l\quad\text{given by}\quad  \Phi_\bullet\to \Phi_l(01\cdots l).$$
	\end{example}

	\begin{example}[The boundary]
	The \emph{boundary $\partial\Delta[l]$ of the simplex $\Delta[l]$} is the simplicial subset of $\Delta[l]$ obtained by removing the interior of $\Delta[l]$, i.e. we remove the element	$v=02\cdots l$ and all their degeneracies.
	\end{example}
	
	\begin{example}[The horns]
		Let $0\leq j\leq l$ be integers. The \emph{horn} $\Lambda[l,j]_\bullet\subseteq \Delta[l]_\bullet$ is the simplicial subset of $\Delta[l]$ obtained by removing the interior of $\Delta[l]$ and its $j$-th face, i.e. we remove the elements $$v=01\cdots l\qquad\text{and}\qquad d_j(v)=01\cdots\widehat{j}\cdots l$$ and all their degeneracies.
	\end{example}

	\subsection{Simplicial manifolds}
	From now on $X_\bullet$ will denote a \emph{simplicial manifold}, i.e. a simplicial object in the category of smooth manifolds. In general, we will work with finite dimensional manifolds. Nevertheless, in some examples we will encounter infinite dimensional manifolds. In those cases, we assume they are Banach manifolds.
	
	\subsubsection{Geometric realizations}\label{sec:georea}
	One can  view  a simplicial manifold as a ``resolution" for a more complicated topological space. Here we make  this idea precise. Let $\bDelta^n$ be the \emph{geometric $n$-simplex},  this is $$\bDelta^n=\{(t_0,\cdots,t_n)\in[0,1]^{n+1} \ |\ t_0+\cdots+t_n=1\}.$$ Following \cite{seg:cat}, for a simplicial manifold $X_\bullet$ one defines the topological spaces: 
	\begin{itemize}
		\item The \emph{fat geometric realization}
		$$||X_\bullet||=\coprod_{l=0}^\infty X_l\times\bDelta^l/\sim\quad\text{with}\quad (d_i(x),t_0,\cdots t_{l-1})\sim (x,t_0,\cdots t_{i-1},0,t_i,\cdots,t_{l-1}),$$
		for $x\in X_l$ and $(t_0,\cdots,t_{l-1})\in\bDelta^{l-1}.$
		\item The \emph{geometric realization}
		$$|X_\bullet|=\coprod_{l=0}^\infty X_l\times\bDelta^l/\approx\quad\text{with}\quad \left\{\begin{array}{l}
			(d_i(x),t_0,\cdots t_{l-1})\approx (x,t_0,\cdots t_{i-1},0,t_i,\cdots,t_{l-1}),  \\
			(s_i(x),t'_0,\cdots t'_{l+1})\approx (x,t'_0,\cdots t'_{i-1},t'_i+t'_{i+1},\cdots,t'_{l+1}),
		\end{array}\right.$$
		for $x\in X_l, \ (t_0,\cdots,t_{l-1})\in\bDelta^{l-1}$ and $(t'_0,\cdots,t'_{l+1})\in\bDelta^{l+1}.$ 
	\end{itemize}
	
	Segal proves \cite[Prop A.1]{seg:cat} that if $X_\bullet$ is a simplicial manifold for which the degeneracies $s_i:X_l\to X_{l+1}$ are closed embeddings\footnote{In his terminology, such simplicial manifolds are called \emph{good}.}, then there is a weak homotopy equivalence between $|X_\bullet|$ and $||X_\bullet||.$

	\subsubsection{Differential forms and the generalized de Rham theorem}
	Perhaps the most relevant aspect of simplicial manifolds for these notes is the analogue of differential forms. Let $X_\bullet$ be a simplicial manifold, the differential forms on $X_\bullet$ form a double complex  $(\Omega^\bullet(X_\bullet), \delta, d)$ known as  the  \emph{Bott-Shulman-Stasheff double complex}, see \cite{bss:on}, where $\delta$ is the \emph{simplicial differential} given by the alternating sum of the pullback along the face maps $d_i:X_l\to X_{l-1}$ and $d$ is the de Rham differential on each manifold, i.e. 
	$$\delta=\sum_i (-1)^i d_i^\ast:\Omega^k(X_{l-1})\to\Omega^k(X_{l})\quad\text{and}\quad d:\Omega^j(X_l)\to \Omega^{j+1}(X_l).$$
	We can depict this double complex as  
	\begin{equation*}
			\begin{tikzcd}
				&\vdots &\vdots &\vdots\\
				0\arrow[r] &\Omega^2(X_0) \arrow[u, "d"] \arrow[r, "\delta"] &\Omega^2(X_1)\arrow[u,"d"] \arrow[r,"\delta"] &\Omega^2(X_2)\arrow[u,"d"] \arrow[r,"\delta"] &\cdots
				\\
				0\arrow[r] &\Omega^1(X_0) \arrow[u, "d"]\arrow[r,"\delta"] &\Omega^1(X_1)\arrow[u,"d"] \arrow[r,"\delta"] & \Omega^1(X_2)\arrow[u,"d"]\arrow[r,"\delta"] &\cdots
				\\
				0\arrow[r] &C^\infty(X_0) \arrow[u, "d"]\arrow[r,"\delta"] &C^\infty(X_1)\arrow[u,"d"] \arrow[r,"\delta"]& C^\infty(X_2)\arrow[u,"d"]\arrow[r,"\delta"] &\cdots
				\\
				&0 \arrow[u] &0\arrow[u] &0\arrow[u]
			\end{tikzcd}
	\end{equation*}
	We denote the associated differential on the total complex of the double complex $(\Omega^\bullet(X_\bullet), \delta, d)$   by $D=\delta + (-1)^l d$.
	Observe that the simplicial differential just uses the face maps. Following \cite[\S~4]{eml:on}, we introduce the  \emph{normalized double subcomplex} as 
	\begin{equation*}
		\hat{\Omega}^\bullet(X_\bullet)=\{\alpha \in \Omega^\bullet(X_\bullet) \, | \, s_i^\ast \alpha = 0  \text{ for all possible } i\},
	\end{equation*}
	in other words, a differential form $\alpha$ is called \emph{normalized} if it vanishes on degeneracies. Nevertheless, the normalized double subcomplex is quasi-isomorphic to the Bott-Shulman-Stasheff double complex as shown in \cite[Thm. 4.1]{eml:on}.
	
	To fix notation, we also say that $\alpha_\bullet$ is an \emph{$m$-shifted $k$-form} if
	\begin{equation*}
		\alpha_\bullet = \sum_{i=0}^{m}\alpha_i \quad \text{with} \quad \alpha_i\in \hat{\Omega}^{k+m-i}(X_i).
	\end{equation*}
	\begin{remark}
		Following \cite{ptvv}, one can interpret $m$-shifted $k$-forms as elements in a filtration for the normalized subcomplex.    
	\end{remark}
	
	With all these definitions in place, we are ready to state the generalized de Rham theorem for simplicial manifolds, see \cite[\S~1]{bss:on}.
	\begin{theorem}\label{thm:genderham}
		Let $X_\bullet$ be a simplicial manifold. Then 
		$$H^\bullet(\text{Tot}(\widehat{\Omega}^\bullet(X_\bullet)),D)\cong H^\bullet_{sing}(||X_\bullet||).$$
	\end{theorem}
	
	This result shows precisely how we can think that simplicial manifolds give a resolution for complicated topological spaces. 
	
	\begin{remark}
		A point that we will not explore in these notes is that the Bott-Shulman-Stasheff double complex carries additional structure. In particular one would like to have an analogue of the wedge product on the differential forms on a manifold. Indeed, one can define a product by the formula
		\begin{equation}\label{eq:cuppro}
			(\alpha\cupprod\beta)_p=\sum_{q=0}^p (d^*_{q+1})^q\alpha_{p-q}\wedge(d^*_0)^{p-q}\beta_q.
		\end{equation}
		This product is not graded commutative, for more details see \cite{get:tra}.
	\end{remark}

	\subsection{Lie n-groupoids and first examples}\label{sec:lien}
	\emph{Kan simplicial} sets \cite[\S~1]{may:simp}, an important special class of simplicial sets, are the simplicial sets $S_\bullet$ such that any map $\Lambda[l, j]_\bullet \to S_\bullet$ extends to a map $\Delta[l]_\bullet \to S_\bullet$. Kan simplicial sets are \emph{fibrant objects}, see e.g. \cite{gh:simp,  hov:mod}, which makes calculations easier in some cases.    Similarly to the case of simplicial sets, we restrict our attention to a subclass of simplicial manifolds by imposing a version of the Kan conditions.

	A \emph{Lie $n$-groupoid} \cite{get:int, hen:int, cc:n-grou} is a simplicial manifold $K_\bullet$ where the natural projections 
	\begin{equation}\label{eq:plj}
		p^l_j\colon K_l=\Hom(\Delta[l],K_\bullet)\to \operatorname{Hom}(\Lambda[l,j],K_\bullet)
	\end{equation}
	are surjective submersions for all $ l\geq 1$ and $0\leq j\leq l$ and diffeomorphisms for all $l\geq n$ and $0\leq j\leq l$. When $K_0$ is a point we say that it is a \emph{Lie $n$-group}.
	
	The fact that $\operatorname{Hom}(\Lambda[l,j],K_\bullet)$ is again
	a manifold is not immediate but can be proven inductively, for example 
		\begin{equation*}
			\begin{array}{rcl}
			\operatorname{Hom}(\Lambda[1,0],K_\bullet)&=&K_0=\operatorname{Hom}(\Lambda[1,1],K_\bullet)\\
			\operatorname{Hom}(\Lambda[2,0],K_\bullet)&=&K_1\times_{d_0,K_0,d_0} K_1\\ \operatorname{Hom}(\Lambda[2,1],K_\bullet)&=&K_1\times_{d_1,K_0,d_0} K_1,\\\ \operatorname{Hom}(\Lambda[2,2],K_\bullet)&=&K_1\times_{d_1,K_0,d_1} K_1,\\
			\operatorname{Hom}(\Lambda[3,0],K_\bullet)&=&\big(K_2\times_{d_1,K_1,d_1} K_2\big)\times_{d_2\times d_2,\operatorname{Hom}(\Lambda[2,2],K_\bullet),p^2_2}K_2.
			\end{array}
		\end{equation*}
		We refer to \cite[Corollary 2.5]{hen:int} or \cite[\S~2]{cc:n-grou} for more details on this fact.
	

	\begin{remark}
		In \cite{dus:hig, dus:sim, glen:rea} higher groupoids were defined as Kan complexes because Kan conditions guarantee some invertibility and composition of simplices, see \S~\ref{sec:lie1grou} for a concrete example. In analogy with the above, the works \cite{get:int, hen:int, cc:n-grou} introduced  Lie $n$-groupoids in the smooth setting.
	\end{remark}
	
	Notice also that for a Lie $n$-groupoid $K_\bullet$ all the face maps are surjective submersions, see e.g. \cite[Lemma 2.44]{du:thes} who built on  \cite{joy:qua}, therefore the degeneracy maps are closed embeddings (because they are sections of a submersion). Hence, the work of Segal implies that the geometric realization $|K_\bullet|$ and the fat geometric realization $||K_\bullet||$ are weak homotopy equivalent. 
	
	\subsubsection{Manifolds as Lie 0-groupoids}\label{sec:0grou}
	Given a manifold $M$, we define the simplicial manifold $\underline{M}_\bullet$ with $\underline{M}_l=M$ for all $l\in\bN$ and $d_i=s_i=\id.$ This simplicial manifold is a Lie $0$-groupoid because all the horn projections $p^l_j$ are diffeomorphisms. 
	
	Conversely if $K_\bullet$ is a Lie $0$-groupoid, then $K_l\cong K_0$ and all the face maps must be diffeomorphisms, so $K_\bullet\cong (\underline{K_0})_\bullet$.

	\subsubsection{Lie groupoids as Lie 1-groupoids}\label{sec:lie1grou}
	
	Recall that a \emph{Lie groupoid} $K\rightrightarrows M$ consists of two manifolds $M$ \emph{(objects)} and
	$K$ \emph{(arrows)} together with surjective submersions $\gs,\gt:K\to M$ \emph{(source and target)} and smooth maps $\cdot:K{}_{\gs}\times_{\gt}K\to M$ \emph{(multiplication)}, ${\cdot}^{-1}:K\to K$ \emph{(inversion)} and $1:M\to K$ \emph{(unit)}
	and which satisfy
	\begin{eqnarray}
		&\gs(g\cdot h)=\gs(h),\qquad \gt(g\cdot h)=\gt(g),\qquad \gs(1_p)=p=\gt(1_p),\label{eq:lg1}\\
		&((g\cdot h)\cdot k)=(g\cdot (h\cdot k)),\qquad g\cdot 1_{\gs(g)}=g=1_{\gt(g)}\cdot g,\label{eq:lg2}\\
		&\gs(g^{-1})=\gt(g),\qquad \gt(g^{-1})=\gs(g),\qquad g^{-1}\cdot g=1_{\gs(g)},\qquad g\cdot g^{-1}=1_{\gt(g)},\label{eq:lg3}
	\end{eqnarray}
	for $g,h,k\in K$ and $p\in M$.
	
	The nerve $N_\bullet K$ of a Lie groupoid $K\rightrightarrows M$ is given by the simplicial manifold 
	\[
	N_0K=M\quad \text{and}\quad N_lK = \{(k_1, \dots, k_l)\in K^{\times l} \, | \, \gs(k_i)= \gt(k_{i+1})\}
	\]
	with face and degeneracy maps given by
	\begin{align*}
		d^1_0(k)&= \gs(k), \quad d^1_1(k)=\gt(k),\quad 
		d^l_0(k_1, \dots, k_l)= (k_2, \dots, k_l), \\
		d_i^l(k_1, \dots, k_l)&= (k_1, \dots, k_i\cdot k_{i+1}, \dots, k_l), \quad 
		d_l^l(k_1, \dots, k_l)= (k_1, \dots, k_{l-1}),\\
		s_0^0(p)&=1_p,\quad s_i^l(k_1, \dots, k_l)= (k_1, \dots, 1_{\gs(k_i)}, \dots, k_l).
	\end{align*}
	Since $K\rightrightarrows M$ is a Lie groupoid, then $\gs$ and $\gt$ are surjective submersions and the next level face maps give a diffeomorphism between $K_2$ and $K{}_{\gs}\times_
	{\gt}K$. Therefore $N_\bullet K$ is a Lie $1$-groupoid. 
	
	Conversely, given a Lie $1$-groupoid $K_\bullet$ we define a Lie groupoid as follows
	$$M=K_0,\quad K=K_1,\quad \gs=d_0^1,\quad \gt=d_1^1\quad \text{and}\quad 1=s_0^0,$$
	that satisfy \eqref{eq:lg1} by the simplicial indentities. 
	Since $p^2_1:K_2\to K{}_{\gs}\times_\gt K$ is a diffeomorphism we can use this map to define the multiplication  
	$$\cdot=d_1^2\circ (p_1^2)^{-1}$$
	and the inverse is given by the fact that $p^2_0$ and $p^2_2$ are diffeomorphisms, so given $g\in K$ there must exist $(g,g^{-1})\in K{}_{\gs}\times_\gt K\cong K_2$ filling the horn. Observe that \eqref{eq:lg2} are satisfied due to $p^3_2$ and $p^2_1$ being diffeomorphisms while the equations in \eqref{eq:lg3} hold  because we are identifying $K{}_{\gs}\times_\gt K\cong K_2$ using the map $p^2_1.$ For more detail see \cite[$\S~1$]{cc:n-grou}.
	
	\subsubsection{Double Lie groupoids and the Artin-Mazur construction}\label{sec:doug}
	\begin{equation}\label{eq:dlg}
		\begin{array}{c}
			\xymatrix{D\ar@<-.5ex>[r]\ar@<.5ex>[r]\ar@<-.5ex>[d] \ar@<.5ex>[d]& H\ar@<-.5ex>[d] \ar@<.5ex>[d]\\ G\ar@<-.5ex>[r]\ar@<.5ex>[r]&M}
		\end{array}
	\end{equation}
	The square \eqref{eq:dlg} is a \emph{double Lie groupoid} if each of the sides is a Lie groupoid and the following conditions hold:
	\begin{enumerate}
		\item[1.] $\gs^v,\gs^h,\gt^v, \gt^h$ are Lie groupoid morphisms.
				\item[2.] The interchange law
				$$(\alpha_{11}\cdot^h\alpha_{12})\cdot^v(\alpha_{21}\cdot^h\alpha_{22})=(\alpha_{11}\cdot^v\alpha_{21})\cdot^h(\alpha_{12}\cdot^v\alpha_{22})$$
				holds for all $\alpha_{11},\alpha_{12},\alpha_{21},\alpha_{22}\in D$ such that $\gs^h (\alpha_{i1}) = \gt^h(\alpha_{i2})$ and $\gs^v (\alpha_{1i}) =\gt^v (\alpha_{2i})$ for $i = 1, 2$.
				\item[3.] The double-source map $(\gs^v , \gs^h ) : D \to G{}_{\gs}\times_{\gs}H$ is a submersion. If in addition it is surjective we say that the double Lie groupoid is \emph{full}.
			\end{enumerate}
			
			In order to avoid confusion, we denote a double Lie groupoid by $D_{\bullet,\bullet}$. This notation is motivated by the fact that for any double Lie groupoid \eqref{eq:dlg}, there is an associated \emph{bisimplicial manifold}, see \cite[Prop 3.10]{mt:from}, given as follows: For $p, q > 0$, let $D_{p,q}$ consist of all $pq$-tuples $\{\alpha_{ij} \}$, where $\alpha_{ij}\in D $ for
			$1 \leq i \leq p$ and $1 \leq j \leq q$, and the elements satisfy the following compatibility conditions:
			$$\gs^h (\alpha_{ij}) = \gt^h (\alpha_{i(j+1)}),\qquad \gs^v (\alpha_{ij}) = \gt^v (\alpha_{(i+1)j}).$$
			We take $D_{0,q} = N_qH$ and $D_{p,0} = N_pG$ for $p, q > 0$, and we set $D_{0,0} = M$.
			
			In the 1960s, it was shown in \cite{am:on} how taking codiagonals on a bisimplicial
			set $S_{\bullet,\bullet}$ one obtains a simplicial set $\overline{W}_\bullet S$ described
			as follows.
			For each $l$, let $\overline{W}_lS$ be the fiber product
			$$\overline{W}_lS=S_{0,l}\ {}_{d^v_0}\times_{d^h_1}S_{1,l-1}\ {}_{d^v_0}\times_{d^h_2}\cdots {}_{d^v_0}\times_{d^h_l}S_{l,0}$$
			with  face and degeneracy maps defined by
			$$d_i(x_0, \cdots, x_l) = (d^v_i(x_0), d^v_{i-1}(x_1),\cdots , d^v_1(x_{i-1}), d^h_i(x_{i+1}), \cdots , d^h_i(x_l )),$$
			$$s_i(x_0, . . . , x_l) = (s^v_i(x_0), s^v_{i-1}(x_1),\cdots , s^v_0(x_i), s^h_i(x_i), \cdots , s^h_i(x_l )) .$$
			
			This method is known as the \emph{Artin-Mazur construction} and it was used in \cite[Theorem 4.5]{mt:from} to prove the following result:
			
			\begin{theorem}\label{thm:raj}
				Let $D_{\bullet,\bullet}$ be a double Lie groupoid. Then $\overline{W}_\bullet D$ is a local Lie $2$-groupoid. Furthermore, if $D_{\bullet,\bullet}$ is a full double Lie groupoid, then $\overline{W}_\bullet D$ is a Lie
				$2$-groupoid.
			\end{theorem}
			
			
			\subsubsection{Cochain complexes and Dold-Kan}

			Another way of constructing Lie $n$-groupoids is by using the Dold-Kan correspondence. Here we give the basic ideas on the simpler case of cochain complexes,  see e.g. \cite[\S 8.4]{weib:intro} or \cite[\S~1.2.3]{lurie:ha}. Given a cochain complex concentrated in non-positive degrees $({\bf V},d)$ define a simplicial vector space $DK_\bullet({\bf V})$ given by 
			$$DK_l({\bf V})=\bigoplus_{[l]\twoheadrightarrow[k]}V_{-k}=\bigoplus_{k} (V_{-k})^{\binom{l}{k}}\qquad \text{for}\ l,k\in\bN$$
			and for a morphism in $\Delta$, $\beta : [l'] \to [l]$, the induced map $\beta^* : DK_{l}({\bf V})\to DK_{l'}({\bf V})$
			is given by the matrix of morphisms $\{f_{\alpha,\alpha'} : V_{-k} \to V_{-k'} \}$, with $\alpha:[l]\to[k]$  and  $\alpha':[l']\to[k']$ surjective maps on  $\Delta$, where the map $f_{\alpha,\alpha'}$ is zero unless
			\begin{itemize}
				\item  $k = k'$ and the following diagram commutes
				\begin{equation*}
					\begin{array}{c}
						\xymatrix{\left[l'\right]\ar[d]_{\alpha'} \ar[r]^\beta &\left[l\right]\ar[d]^{\alpha}\\ \left[k'\right]  \ar[r]^\id & \left[k\right],
						}
					\end{array}
				\end{equation*}
				then the map $f_{\alpha,\alpha'}=\id$ or 
				\item  $-k' = -k + 1$ and the following diagram commutes
				\begin{equation*}
					\begin{array}{c}
						\xymatrix{\left[l'\right]\ar[d]_\alpha \ar[rr]^\beta &&\left[l\right]\ar[d]^{\alpha'}\\ \left[k'\right]  \ar[r]^{\sim\qquad} & \{0, \cdots, k-1\}\ar[r]& \left[k\right],
						}
					\end{array}
				\end{equation*}
				then the map $f_{\alpha,\alpha'}=d$, the differential.
			\end{itemize}
			The fact that if ${\bf V}$ is concentrated in degree $0$ to $n$ then $DK_\bullet({\bf V})$ is a Lie $n$-groupoid goes back to \cite[Thm. 3]{moo:hom}.
			
			Since we will need it later, let us analyze the converse statement. Given  a simplicial vector space $V_\bullet$, its \emph{normalized complex} $({\bf N}_\bullet V, \partial)$ is given by
			\[
			{\bf N}_{-l}V:=\bigcap_{i=0}^{l-1} \ker d_i\subset V_l, \quad \partial := (-1)^l d_l\colon {\bf N}_{-l}V \to {\bf N}_{-l+1}V.
			\]
			The Dold-Kan correspondence for simplicial vector spaces says that the functor sending $V_\bullet$ to $({\bf N}_\bullet V, \partial)$ is an equivalence of categories between simplicial vector spaces and cochain complexes of vector spaces concentrated in non-positive degrees. 
			
			\begin{remark}
				The above Dold-Kan correspondence can be extended to treat the case of vector bundle over Lie $n$-groupoids as in \cite{ mat:vec, raj:vbgrou}.
			\end{remark}
			
			\subsection{Morita equivalence and differentiable stacks}\label{sec:morita}

			We started by saying that simplicial manifolds give ``resolutions" for topological spaces via the geometric realizations explained in \S~\ref{sec:georea}. Moreover, the work \cite{seg:cat} suggests that one should care just about the homotopy type of the geometric realization. From that perspective, it is natural to ask when two different Lie $n$-groupoids give weakly homotopic geometric realizations. A refined version of that question is addressed by the following definitions.
			
			The \emph{boundary space} of a simplicial manifold $X_\bullet$ is
			\begin{equation}\label{eq:bun}
				\partial_i(X_\bullet):=\Hom(\partial\Delta[i],X_\bullet)=\{(x_0,\cdots,x_i)\in X^{\times i+1}_{i-1}|\ d_j(x_k)=d_{k-1}(x_j) \ \forall j<k\}.
			\end{equation}
			Let $K_\bullet, J_\bullet$ be two Lie $n$-groupoids. We say that $\Phi_\bullet: K_\bullet\to  J_\bullet$ is a \emph{hypercover} if $\Phi_0$ is a surjective submersion and the maps
			\begin{equation}\label{eq:hyp}
				q_i := ((d_0, · · · , d_i), \Phi_i) : K_i \to \partial_i(K_\bullet)\times_{\partial_i(J_\bullet)}J_i
			\end{equation}
			are surjective submersions for $1 \leq i < n$ and an isomorphism for $i = n$.  Moreover, $K_\bullet$ and $J_\bullet$ are \emph{Morita equivalent} if there is a third Lie $n$-groupoid $Z_\bullet$ with hypercovers $\Phi_\bullet:Z_\bullet \to K_\bullet$ and $\Psi_\bullet : Z_\bullet \to J_\bullet$.
			
			\begin{remark}
				For details on the definition of hypercover see \cite[\S~3]{bg:hig}. In general, $\partial_i(X_\bullet)$ might not be a manifold any more, even if $X_\bullet$ is a Lie
				$n$-groupoid, and the right-hand side of \eqref{eq:bun} is a set-theoretical description of it. Nevertheless, the right-hand-side of \eqref{eq:hyp}  is always a manifold if $\Phi_\bullet$ is a hypercover, even though
				it seems to be logically dependent, this is something one can prove level-wise inductively, see \cite[Lemma 9.4]{dhi:hyp} or \cite[\S~2.1]{cc:n-grou}.  Observe also that  $q_i$ is automatically
				an isomorphism for $i \geq n + 1$, see \cite[Lemma 4.11]{dhi:hyp} or \cite[Lemma 2.5]{cc:n-grou}. 
			\end{remark}
			
			\begin{remark}[Weak equivalences]
				Instead of defining Morita equivalences using hypercovers one can say that two groupoids are Morita equivalent if there is a zig-zag of \emph{weak equivalences}, see \cite[\S~5]{rcc:hom} for a definition. It was shown in  \cite[\S~6]{rcc:hom} that the two definitions give the same equivalence relation. 
			\end{remark}
			
			For $n=1$, it was shown in \cite[\S~2]{bx:sta} that the Morita class of a Lie groupoid encodes a \emph{differentiable stack}, that is a category fibred in groupoids satisfying descent conditions and having a presentation, for a precise definition see \cite[Def 2.10 and 2.15]{bx:sta}. For higher groupoids, the full homotopy theory was developed in \cite{bg:hig, rcc:hom}. Therefore, we think of a Lie $n$-groupoid $K_\bullet$ as a particular atlas presenting a differentiable $n$-stack. 
			
			One of the important facts of hypercovers that is relevant to these notes is the following
			
			\begin{proposition}\label{prop:cohomdes}
				Let $\Phi_\bullet:K_\bullet\to J_\bullet$ be an hypercover between Lie $n$-groupoids. Then 
				$$\Phi^*:H^*(Tot(\widehat{\Omega}^\bullet(J_\bullet)),D)\to H^*(Tot(\widehat{\Omega}^\bullet(K_\bullet)),D) $$
				is an isomorphism. Moreover, the map $\Phi^*$ is an isomorphism on the cohomology of the filtration given by shifted forms.
			\end{proposition}
			\begin{proof}
				The case for Lie groupoids follows from \cite[Prop.~2]{beh:coh}. For general Lie $n$-groupoids a proof can be found in \cite{get:slides}.
			\end{proof}

			\subsection{On the Lie functor}\label{sec:liefun}
			
			On the one hand, in Part \ref{partI} we introduced $Q$-manifolds and showed in Theorem \ref{thm:bon} (see also Remark \ref{rmk:bon}) that they are the ``same thing" as Lie $n$-algebroids. On the other hand, in Part \ref{partII} we introduced Lie $n$-groupoids as simplicial manifolds satisfying Kan conditions. Nevertheless, we did not give any relation between these two objects apart from the obvious semantic. 
			
			If we specify the above discussion when $n=1$, one obtains Lie algebroids (by Theorem \ref{thm:Q-lie}) and  Lie groupoids (by $\S~\ref{sec:lie1grou}$). In \cite{pra:lie}, it was shown that there is a differentiation functor between these two objects extending the classical one between Lie groups and Lie algebras.
			
			Inspired by the above Lie functor, together with ideas coming from homotopy theory \cite{sul:inf} and mathematical physics \cite{aksz} the work \cite{sev:some} proposed that $Q$-manifolds should be integrated into a kind of simplicial manifolds\footnote{The article \cite{sev:some} uses balls instead of simplices.}. Later on \cite{get:int, hen:int} introduced Lie $n$-groupoids as simplicial manifolds satisfying Kan conditions and used them to integrate $L_\infty$-algebras. Using the above definition of Lie $n$-groupoids, the works \cite{sev:gorm, sev:jet, sev:diff} gave a conceptual Lie functor, see \cite[Proposition 9.2]{sev:jet}. 
			
			Given a Lie $n$-groupoid $K_\bullet$, define the graded vector bundle 
			\begin{equation}\label{eq:nalg}
				{\bf A}_K=\bigoplus_{i=-n+1}^0 A_i\to K_0\quad\text{with}\quad A_i=(\ker Tp_{-i+1}^{-i+1})_{|K_0},
			\end{equation}
			where $p_{-i+1}^{-i+1}$ denotes the horn projection as introduced in \eqref{eq:plj}.
			Using ${\bf A}_K$, a concrete model  for the above conceptual Lie functor was given in  \cite[Theorem 3.3]{ccl:diff}. 
			
			We state  \cite[Proposition 9.2]{sev:jet} and  \cite[Theorem 3.3]{ccl:diff} together, in the following result whose proof is beyond the scope of these notes.

			\begin{theorem}\label{thm:liefun}
				There is a differentiation functor 
				$$\text{Lie}:\{ \text{Lie $n$-groupoids}\}\to \{\text{Degree $n$ $Q$-manifolds}\}, \quad K_\bullet\to \text{Lie}(K_\bullet)=(\cK, Q).$$
				Moreover, $\cK$ is non-canonically isomorphic to ${\bf A}_K[1]$ as $n$-manifolds.
			\end{theorem}
			
			Even though Theorem \ref{thm:liefun} is conceptually clear, in concrete examples it is difficult to compute the $Q$-structure. For recent progress in clarifying these points see \cite{ch:diff, dor:diff, ccl:diff}.  After this discussion several important remarks are in order:
			\begin{enumerate}
				\item[1.] {\bf Integrability:} The problem of integrability, i.e. understanding if the functor Lie is essentially surjective, remains open. The obstructions for $n=1$ were computed in \cite{rui:int}, see also \cite{cc:ints} for a different solution. The study of higher $n$ started with the string group, see e.g. \cite{baez:from}, and the study of double Lie algebroids and groupoids, see e.g. \cite{luwe, ste:laint, bch:vec}. For Lie $n$-algebras the integration problem is treated in \cite{get:int, hen:int, rw:int} among others. The integration of Courant algebroids also received much attention including \cite{lis:int,  mt:from, cc:extensions, zhu:bg, alv:tran}. Finally, a first attempt for the general (local) theory can be found in \cite{sev:int}.
				\item[2.] {\bf The van Est map:} The above theorem explains how to differentiate the space. However, we would like to apply it to differential forms as well. When $n=1$, this can be done using the van Est map constructed in \cite{ac:ve}.  Given a Lie groupoid $K_\bullet$ denote by  $(A_K[1], Q)$ its degree $1$ $Q$-manifold, i.e. its Lie algebroid. Then one has a map of double complexes
				\begin{equation}\label{eq:VEmap}
					\text{VE}:(\widehat{\Omega}^\bullet(K_\bullet),\delta,d)\to (\Omega^\bullet_\bullet(A_K[1]), \Lie_Q, d)
				\end{equation}
				that allows us to differentiate shifted forms on $K_\bullet$ to forms on the $Q$-manifold $A_K[1].$ For strict Lie $2$-groups a version of the van Est map can be found in \cite{cue:ve}. For higher $n$ see the work in progress \cite{ch:diff}, that also gives an alternative proof for Theorem \ref{thm:liefun}. 
				\item[3.] {\bf Double Lie groupoids:} The differentiation of double Lie groupoids (\S~\ref{sec:doug}) was constructed in \cite{mac:dou} as a two step process: first from double Lie groupoids to $LA$-groupoids and then to double Lie algebroids (\S~\ref{sec:doua}). In the respective sections we already showed how double Lie groupoids give examples of Lie $2$-groupoids (via the $\overline{W}$ construction) and double Lie algebroids of $Q$-manifolds (via the [1][1] shifting construction). We believe that the following diagram commutes, however we have not found a proof in the literature
				\begin{equation}\label{diag:diff}
					\begin{tikzcd}
						&& {} \\
						\begin{array}{c} \begin{array}{c}\text{Double Lie}\\\text{groupoids} \end{array} \end{array} && {LA\text{-groupoids} } && \begin{array}{c} \begin{array}{c}\text{Double Lie}\\\text{algebroids} \end{array} \end{array} \\
						\\
						{\text{Lie 2-groupoids}} &&&& \begin{array}{c} \begin{array}{c}\text{Degree 2}\\ Q\text{-manifolds.} \end{array} \end{array} \\
						&&&& {}
						\arrow["{\text{Lie}}", from=2-1, to=2-3]
						\arrow["{\overline{W}}"', from=2-1, to=4-1]
						\arrow["{\text{Lie}}", from=2-3, to=2-5]
						\arrow["{[1][1]}", from=2-5, to=4-5]
						\arrow["{\text{Lie}}"', from=4-1, to=4-5]
					\end{tikzcd}
				\end{equation}
			\end{enumerate}
			
			\section{Shifted symplectic structures}\label{sec:ss}
			In this section, following \cite{get:slides}, we introduce $m$-shifted symplectic structure on Lie $n$-groupoids. Our focus is on examples and on illustrating some of their main properties. The concept of shifted symplectic structures has two primary motivations. On the one hand, it aims to provide a concrete differential geometric treatment of the shifted symplectic structures on higher stacks appearing in \cite{ptvv}. On the other hand, it serves as a global counterpart to the symplectic $Q$-manifolds introduced in $\S~\ref{sec:sympQ}$.

			\subsection{Tangent complex}
			
			Given  a Lie $n$-groupoid $K_\bullet$, we get a simplicial vector bundle over $K_0$, given by
			\begin{equation}\label{eq:tan}
				\begin{tikzcd}
					\cdots TK_2|_{K_0} \arrow[r] \arrow[r, shift left=1ex] \arrow[r, shift right=1 ex] & TK_1|_{K_0} \arrow[r, shift left=0.5ex] \arrow[r, shift right=0.5 ex] & TK_0,
				\end{tikzcd} 
			\end{equation}
			where we used that $K_0$ is a submanifold of $K_l$ via $s_0\circ \dots \circ s_0\colon K_0\to K_l$ and $TK_l|_{K_0}$ denotes the restriction of the tangent bundle $TK_l$ to this submanifold. The face and degeneracy maps are the restrictions of $Td_i$ and $Ts_i$, respectively. 
			
			The \emph{tangent complex} $(\cT^K_\bullet, \partial)$ of a Lie $n$-groupoid $K_\bullet$ is the normalized complex of \eqref{eq:tan}, more concretely, is  the cochain complex of vector bundles over $K_0$ given by 
			\begin{equation*}
				\cT_l^K := 
				\begin{cases}
					(\ker Tp^{-l}_{-l})_{|K_0} & \text{if } l<0, \\
					TK_0 & \text{if } l=0, \\
					0 & \text{if } l>0,
				\end{cases}\qquad \text{with }\qquad \partial:= (-1)^l Td_{}^{-l}.
			\end{equation*}
			We denote by $H^\bullet(\cT^K)$ the cohomology groups of the tangent complex $(\cT^K_\bullet , \partial)$. We also introduce the \emph{cotangent complex} $(\cT_\bullet^{\ast K}, \partial^\ast)$ as the dual to the tangent complex, this means
			\[
			\cT_l^{\ast K}= (\cT_{-l}^K)^\ast \qquad \text{and} \qquad \partial^*= (-1)^{l+1}\partial^t.
			\]

			\begin{remark}
				The notion of tangent complex works well for local Lie $n$-groupoids, and the tangent complex of a local Lie $n$-groupoid $K_\bullet$ is still a complex of vector bundles over $K_0$, since only the submersion condition is  used in the definition, see \cite{dor:diff} for more details. 
			\end{remark}
			From the Kan conditions for a Lie $n$-groupoid, we get that $\cT_{-l}^K=0$ if $l>n$ and one can show that the rank of $\cT_{-l}^K$ is 
			\begin{equation*}
				\dim K_l - \sum_{i=0}^{l-1} (-1)^i \binom{l}{i+1}\dim K_{l-1-i}.
			\end{equation*}

			\begin{example}\label{ex:tcxgr}
				Let $K\rightrightarrows M$ be a Lie groupoid. Its nerve is a Lie $1$-groupoid and the tangent complex is concentrated in degrees $-1$ and $0$. In particular we get
				\begin{equation*}
					\cT_l^K := 
					\begin{cases}
						A_K & \text{if } l=-1, \\
						TM & \text{if } l=0, \\
						0 & \text{otherwise},
					\end{cases}
				\end{equation*}
				where $A_K$ is the Lie algebroid of $K$, and the differential is just the anchor map $\rho\colon A_K\to TM$. This complex is, sometimes, called the \emph{core complex} associated to the Lie groupoid $K$.
			\end{example}
			
			\begin{remark}
				For a Lie $n$-groupoid, the non-positive part of $(\cT^K_\bullet,\partial)$ is (up to a shift by $1$) what we called ${\bf A}_K$ in \eqref{eq:nalg}. Therefore, by Theorems \ref{thm:liefun} and $\ref{thm:bon}$ (see also Remark \ref{rmk:bon}) it carries the structure of a Lie $n$-algebroid whose unary bracket is $\partial$. Nevertheless, to give an explicit construction for the brackets (as one does when $n=1$) is not an easy task.
			\end{remark}
			
			We finish this section by stating the relation between hypercovers and the tangent complex appearing in \cite[Lemma 2.27]{zhu:bg}. 
			
			\begin{proposition}\label{prop:hyp}
				A hypercover $\Phi_\bullet : K_\bullet \to J_\bullet$ between Lie $n$-groupoids induces a map on their
				tangent complexes $\cT_\bullet \Phi : (\cT_\bullet^K,\partial) \to (\cT_\bullet^J,\partial)$ that is a quasi-isomorphism. In fact, each $\cT_i\Phi:\cT_i^K\to \cT_i^J$ is surjective.
			\end{proposition}
			
			\subsection{Definition and Morita invariance}\label{sec:defss}
			
			Let $K_\bullet$ be a Lie $n$-groupoid and let $\omega_\bullet$ be an $m$-shifted $2$-form, recall that this means 
			\begin{equation*}
				\omega_\bullet = \sum_{i=0}^{m}\omega_i \quad \text{with} \quad \omega_i\in \hat{\Omega}^{2+m-i}(X_i).
			\end{equation*}
			We say that $\omega_\bullet$ is \emph{closed} if it is closed for the total differential, i.e. $D(\omega_\bullet)=0$. One of the main achievements of \cite{get:slides} was to find an appropriate way to express the non-degeneracy condition for $m$-shifted $2$-forms. 
			
			The \emph{IM-pairing} $\lambda^{\omega_m}$ associated to $\omega_\bullet$ is a graded symmetric pairing in the tangent complex defined as follows: for $l\in \mathbb{N}$ and $v\in (\cT_{-l}^K)_x\subset T_xK_l$, $u\in (\cT_{-m+l}^K)_x\subset T_xK_{m-l}$ in the tangent complex at $x\in K_0$, then
			\begin{equation}\label{eq:im-form}
				\lambda^{\omega_m}(v,u):= \sum_{\sigma\in \Sh(l,m-l)} \text{sgn}(\sigma)\ \omega_m \left(T(s_{\sigma(m-1)}\dots s_{\sigma(m+l)})v, T(s_{\sigma(m+l+1)}\dots s_{\sigma(0)})u\right).
			\end{equation}
			For an interpretation of the $IM$-pairing in terms of the Eilenberg-Zilber map see \cite[\S 2.2]{ccs:vs}. As stated in \cite{get:slides}, it was checked in \cite[App. E]{zhu:bg} that the $IM$-pairing has the following properties.
			
			\begin{proposition}\label{prop:IMpairing}
				Let $K_\bullet$ be a Lie $n$-groupoid and $\omega_\bullet$ be a closed $m$-shifted $2$-form with associated $IM$-pairing $\lambda^{\omega_m}$. Then the following hold:
				\begin{enumerate}
					\item[1.] For $u\in \cT_{-l-1}^K$ and $v\in \cT_{-m+l}^K$ we get $$\lambda^{\omega_m}(\partial u, v) + (-1)^{-l-1}\lambda^{\omega_m}(u, \partial v)=0.$$
					\item[2.] The $IM$-pairing is invariant under a gauge transformation. That is
					$$\lambda^{\omega_m+\delta\eta_{m-1}}(u,v)-\lambda^{\omega_m}(u,v)=\lambda^{\eta_{m-1}}(\partial u,v)+(-1)^l\lambda^{\eta_{m-1}}(u,\partial v)$$
					for any $(m-1)$-shifted $2$-form $\eta_\bullet$ and $u\in\cT^K_{-l}$ and $v\in\cT^K_{-m+l}$.
				\end{enumerate}
			\end{proposition}
			
			The above proposition implies that for closed $m$-shifted $2$-forms  $$(\lambda^{\omega_\bullet})^\flat:(\cT^K,\partial)\to(\cT^{*K}[m],\partial^*)$$ 
			is a cochain map. Therefore, one can make the following definition.

			Following \cite{get:slides}, we say that an \emph{$m$-shifted symplectic structure} on a Lie $n$-groupoid $K_\bullet$ is a closed and non-degenerate $m$-shifted $2$-form $\omega_\bullet$, meaning that  
			$$D(\omega_\bullet)=0\quad\text{and}\quad  (\lambda^{\omega_m})^\flat:(\cT^K,\partial)\to (\cT^{*K}[m],\partial^*)$$
			is a quasi-isomorphism. The pair $(K_\bullet, \omega_\bullet)$ is called an \emph{$m$-shifted symplectic Lie $n$-groupoid}. 
			
			Before going on, let us comment on the definition for a fixed Lie $n$-groupoid $K_\bullet$.
			\begin{enumerate}
				\item[1.] The cohomology groups $H^\bullet(\cT^K)$ form a bundle of vector spaces over $K_0$, but this bundle may not have constant rank and therefore is not, in general, a vector bundle. Nevertheless, the non-degeneracy condition can be equivalently expressed by  stating that the induced pairing on the cohomology groups $H^\bullet(\cT^K)$ is pointwise non-degenerate. 
				\item[2.] Only $m$-shifted $2$-forms with $m\in\bN$ exist on $K_\bullet$, and therefore we do not recover the full picture described in \cite{ptvv}. This limitation arises from the absence of derived directions on Lie $n$-groupoids. In other words, Lie $n$-groupoids present differentiable stacks but not derived differentiable stacks.
				\item[3.] The only possible non-zero cohomology groups are concentrated between $-n$ and $0$. Suppose that $$H^0(\cT^K) \neq 0 \neq H^{-n}(\cT^K)$$  then $K_\bullet$  admits only $n$-shifted symplectic structures. This will be the generic situation but other cases may occur. The cases when $m > n$ are included in $m$-shifted symplectic $m$-groupoids, because any Lie $n$-groupoid is also a Lie $m$-groupoid for $m > n$. The cases when $m < n$ can be understood as “adding singularities” to the case of $m$-shifted symplectic $m$-groupoids, as we will see in $\S~\ref{subsec:sss1}$.
				\item[4.] The $IM$-pairing $\lambda^\omega_\bullet$ of a closed $m$-shifted $2$-form  should be understood as part of the van Est map \eqref{eq:VEmap} for Lie $n$-groupoids. Indeed, this interpretation is true when $n=1$, notwithstanding it just gives the leading term of the full $IM$-form and does not contain all the infinitesimal information \cite{bcwz:dir, bur:mul}.
			\end{enumerate}
			
			\begin{remark}[On differentiation I]\label{rmk:diff1}
				The Lie functor introduced in Theorem \ref{thm:liefun} allows us to produce a $Q$-manifold out of a Lie $n$-groupoid, i.e. 
				$$\text{Lie}(K_\bullet)=(\cK, Q).$$
				However, to differentiate an $m$-shifted symplectic Lie $n$-groupoid, the Lie functor alone is not sufficient; one also needs a  van Est map for Lie $n$-groupoids, namely
				$$\text{VE}:(\widehat{\Omega}^\bullet(K_\bullet),\delta, d)\to (\Omega^\bullet_\bullet(\cK), \Lie_Q, d).$$
				Such a map is currently work in progress, see \cite{ch:diff}, and it is expected that if $(K_\bullet, \omega_\bullet)$ is an $m$-shifted symplectic Lie $n$-groupoid, then  $(\cK, Q, \text{VE}(\omega_\bullet))$ is an $m$-shifted symplectic Lie $n$-algebroid as in $\S~\ref{sec:ssa}$. When $m=n=1$, this was proved in \cite{bcwz:dir}, see Proposition \ref{prop:diffn=1} below. 
			\end{remark}
			
			Following the ideas exposed in \cite{get:slides} one gets the Morita invariance of $m$-shifted symplectic structures.
			
			\begin{theorem}\label{thm:sme}
				Let $\Phi_\bullet:K_\bullet\to J_\bullet$ be a hypercover between Lie $n$-groupoids. Then $K_\bullet$ admits an $m$-shifted symplectic structure $\omega_\bullet^K$ if and only if $J_\bullet$ admits an $m$-shifted symplectic structure $\omega_\bullet^J$ satisfying 
				$$\omega^K_\bullet-\Phi_\bullet^*\omega^J_\bullet=D\eta_\bullet$$
				for some $(m-1)$-shifted $2$-form $\eta_\bullet$ in $K_\bullet$.
			\end{theorem}
			\begin{proof}
				Assume $\omega_\bullet^J$ is an $m$-shifted symplectic structure on $J_\bullet$. Since $\Phi_\bullet$ is a simplicial morphism we have that $D\Phi^*_\bullet=\Phi^*_\bullet D$ and $\Phi^*s_i^*=s_i^*\Phi^*$ so 
				$$\omega_\bullet^K=\Phi_\bullet^*\omega_\bullet^J$$
				is closed and normalized,  thus we only need to show that it is non-degenerate. As $\Phi_0:K_0\to J_0$ is a surjective submersion, we take $x_0 \in J_0$ and one of its preimages $y_0 \in K_0$. Since $\Phi_\bullet$ is a simplicial morphism, it commutes with all face and degeneracy maps and we have
				$$\lambda^{\omega_\bullet^K}_{y_0}(v,u)=\lambda^{\omega_\bullet^J}_{x_0}(\cT_{y_0}\Phi( v),\cT_{y_0}\Phi(u)).$$
				By Proposition \ref{prop:hyp}, the induced map  $\cT\Phi^*: H^\bullet(\cT^K) \to H^\bullet(\cT^J)$ is an isomorphism, thus the induced pairings on the cohomology groups are the same under this isomorphism. Therefore $\omega_\bullet^K$ is
				non-degenerate if and only if $\omega_\bullet^J$ is non-degenerate.
				
				Conversely, assume $\omega_\bullet^K$ is an $m$-shifted symplectic structure. In particular $D\omega_\bullet^K=0$, hence Proposition \ref{prop:cohomdes} implies that there exist a closed (and normalized) $m$-shifted $2$-form $\omega^J_\bullet$ on $J_\bullet$ and a (normalized) $(m-1)$-shifted $2$-form $\eta_\bullet$ in $K_\bullet$ satisfying 
				$$\omega^K_\bullet-\Phi_\bullet^*\omega^J_\bullet=D\eta_\bullet.$$
				By the second item in Proposition \ref{prop:IMpairing} we get that  $\omega^K_\bullet-D\eta_\bullet$ is also an $m$-shifted symplectic form because the $IM$-pairings $\lambda^{\omega^K_\bullet+D\eta_\bullet}$ and $\lambda^{\omega^K_\bullet}$ are homotopic. Therefore by the previous discussion we get that   $\omega^J_\bullet$ is non-degenerate if and only if $\omega^K_\bullet$ so is.
			\end{proof}
			
			The above result allows us to make the following definition.
			Let $(K_\bullet, \omega^K_\bullet)$ and $(J_\bullet, \omega^J_\bullet)$ be two $m$-shifted symplectic Lie $n$-groupoids. We say that $(K_\bullet, \omega^K_\bullet)$ and $(J_\bullet, \omega^J_\bullet)$ are \emph{symplectic Morita equivalent} if there exist another Lie
			$n$-groupoid $(Z_\bullet, \eta_\bullet)$, with $\eta_\bullet$ an $(m-1)$-shifted $2$-form and hypercovers $\Phi_\bullet, \Psi_\bullet$ satisfying
			\begin{equation*}
				K_\bullet\xleftarrow{\Phi_\bullet}Z_\bullet\xrightarrow{\Psi_\bullet}J_\bullet \quad\text{and}\quad  \Psi^*_\bullet\omega^J_\bullet-\Phi_\bullet^*\omega^K_\bullet=D\eta_\bullet.
			\end{equation*}
			
			Almost immediately from the definition, it follows that symplectic Morita equivalence is an equivalence relation, see \cite{zhu:bg}. We also indicate here that symplectic Morita equivalence for $m=n=1$ was introduced in \cite{xu:mome} using principal bibundles.  The fact that these two notions define the same equivalence relation is proved in \cite[Thm. 2.37]{zhu:bg}.
			
			\subsection{First examples}\label{sec:fex}
			With the above machinery we are ready to give the first examples of $m$-shifted symplectic Lie $n$-groupoids. Here we treat the cases when $n=0,1,2.$
			
			\subsubsection{On Lie 0-groupoids}\label{subsub:sss0}
			From $\S~\ref{sec:0grou}$ we know that a Lie $0$-groupoid is a manifold $M$, seen as the constant simplicial manifold $\underline{M}_\bullet$. The tangent complex is concentrated in degree $0$ and given by $\cT_0^{\underline{M}} = TM$. Hence we get two cases:
			\begin{enumerate}
				\item [a)] When $\dim M=0$ (so a bunch of points): Its tangent complex is $0$ everywhere, therefore the $0$ is an $m$-shifted symplectic $2$-form for all $m$! Even though this example looks trivial, we will show it is fundamental.
				\item [b)] When $\dim M\neq 0$: then $\underline{M}_\bullet$  admits only $0$-shifted symplectic structures. It is straightforward to see that $\omega_\bullet=\omega\in\Omega^2(M)$ is a $0$-shifted symplectic structure if and only if it is a usual  symplectic $2$-form. Thus we obtain the following correspondence
				$$\{\text{Symplectic manifolds}\}\leftrightharpoons\{\text{0-Shifted symplectic Lie 0-groupoids}\}.$$
			\end{enumerate}
			
			\subsubsection{On Lie 1-groupoids}\label{subsec:sss1}
			From $\S~\ref{sec:lie1grou}$, we know that a Lie $1$-groupoid is a usual Lie groupoid $K\rightrightarrows M$. Its tangent complex was computed in Example \ref{ex:tcxgr} and it is given by 
			$$A_K\xrightarrow{\rho}TM,$$
			where $A_K$ is the Lie algebroid of the Lie groupoid and $\rho$ is the anchor map. Thus, the only possible non-vanishing cohomology groups are
			$$H^{-1}(\cT^K)=\ker\rho\quad\text{and}\quad H^0(\cT^K)=\text{coker}\rho.$$
			We get four different situations depending on which cohomology groups vanish:
			\begin{enumerate}
				\item [a)] When both cohomology groups vanish. This means that our Lie groupoid is Morita equivalent to a discrete Lie groupoid. As in \ref{subsub:sss0} a), we get that the zero form is an $m$-shifted symplectic $2$-form for all $m$. Two important examples are: the pair groupoid of a manifold and discrete groups, the quantization of this last when $m=2$ is the main subject of study in \cite{dw:top}.
				\item [b)] When $H^{-1}(\cT^K)=0$ and $H^{0}(\cT^K)\neq 0$. This means we get a Lie groupoid with injective anchor map, such groupoids are called {\em foliation groupoids}. In this case $m=0$ and a closed $0$-shifted $2$-form is given by $\omega\in\Omega^2(M)$ satisfying
				\begin{equation}\label{eq:c02f}
					d\omega=0\quad \text{and}\quad \gs^*\omega=\gt^*\omega.
				\end{equation}
				The non-degeneracy condition simply says that the restriction $\omega
				_{|\text{coker}\rho}$ is non-degenerate, so $\omega$ is transversally non-degenerate.  In particular,
				in case the orbit space of $K$ is smooth, it inherits a usual symplectic form. Thus, these Lie groupoids are a singular version of \ref{subsub:sss0} b). The study of these objects was done in \cite{rey:0sym, lm:ham} and model $0$-shifted symplectic stacks.
				\item [c)] When $H^{-1}(\cT^K)\neq 0$ and $H^{0}(\cT^K)=0$. This means we get a Lie groupoid with surjective anchor map, such groupoids are called {\em transitive groupoids} and, when $M$ is connected, they are Morita equivalent to Lie groups. In this case $m=2$, we explain more features of this example in \S~\ref{sec:bg}.
				\item [d)] When $H^{-1}(\cT^K)\neq 0$ and $H^{0}(\cT^K)\neq 0$. This is the generic case and therefore $m=1$. A $1$-shifted symplectic structure on $K\rightrightarrows M$ consists of $\omega\in\Omega^2(K)$ and $H \in\Omega^3(M )$ satisfying
				\begin{equation}\label{eq:1qs}
					\delta\omega=0,\quad d\omega = \delta H,\quad dH = 0 \quad \text{and}\quad \ker \omega_x\cap\ker T_x\gs \cap \ker T_x\gt = 0 \quad \forall x \in M.
				\end{equation} 
				The three equations from the left are equivalent to $D(\omega_\bullet) = 0$, while the non-degeneracy condition for $\omega_\bullet$ can be reformulated to be the last condition. A Lie groupoid with $\omega+H$ satisfying \eqref{eq:1qs} was called a  \emph{twisted presymplectic Lie groupoid} in \cite{bcwz:dir} and \emph{quasi symplectic Lie groupoid} in \cite{xu:mome}.
				So we obtain the following correspondence
				$$\{ \text{Twisted presymplectic Lie groupoids} \}\leftrightharpoons \{ 1\text{-shifted symplectic Lie 1-groupoids} \}.
				$$
				In particular, the classical notion of \emph{symplectic groupoid} \cite{cdw:sym} is a particular case of a $1$-shifted
				symplectic Lie $1$-groupoid.
			\end{enumerate}
			
			Since for Lie groupoids we have the van Est map given in \eqref{eq:VEmap} we can produce a full differentiation of $m$-shifted symplectic structures on Lie groupoids. Let $K_\bullet=(K\rightrightarrows M)$ be a Lie groupoid and denote by $\text{Lie}(K)=(A_K[1], Q)$ its corresponding degree $1$ $Q$-manifold, i.e. its Lie algebroid.
			
			\begin{proposition}\label{prop:diffn=1}
				If $(K_\bullet, \omega_\bullet)$ is an $m$-shifted symplectic Lie $1$-groupoid, then $(A_K[1], Q, \text{VE}(\omega_\bullet))$  is an $m$-shifted symplectic Lie $1$-algebroid as introduced in $\S~\ref{sec:ssa}$.
			\end{proposition}
			\begin{proof}
				The van Est map is a map of double complexes, hence $\text{VE}(\omega_\bullet)$ is a closed $m$-shifted $2$-form on $(A_K[1], Q)$. Moreover, $$(\cT^K_\bullet)_{|p}=(A_p\xrightarrow{\rho}T_pM)=T_pA[1]\quad \forall p\in M.$$
				Therefore, the global and infinitesimal non-degeneracy conditions coincide. Hence, we proved the result.
			\end{proof}
			
			For $m=1$ this theorem appeared in \cite{bcwz:dir}. See also Example \ref{ex:dir}

			\subsubsection{On Lie 2-groupoids}
			
			Observe that Lie $2$-groupoids also have a finite-data description in terms of truncated simplicial manifolds, see \cite{cc:n-grou}. Therefore, we denote a Lie $2$-groupoid and its tangent complex by $$K_\bullet=H\aaar K\rightrightarrows M\quad \text{and}\quad \cT^K_\bullet=B\xrightarrow{\partial}A\xrightarrow{\rho}TM,$$ 
			where $A$ sits in degree $-1$ and $B$ in degree $-2$. A similar analysis as in the previous case shows that Lie $2$-groupoids admit $m$-shifted symplectic structures with $m\in \{0,\cdots,4\}$, or any number $m$ if the tangent complex has no cohomology. The generic case is $m=2$; we interpret $m=0$ and $m=1$ as adding singularities, corresponding to cases (b) and (c) in \ref{subsec:sss1}, respectively, while $m=3$ and $4$ correspond to simpler  shifted symplectic structures on Lie $3$ and $4$-groupoids.  Instead of giving a full classification as before, we just compute the case $m=0$ and compare it with \ref{subsec:sss1} b). 
			
			Let $K_\bullet$ be a Lie $2$-groupoid.  A closed $0$-shifted $2$-form on $K_\bullet$ is given by  $\omega\in \Omega^2(K_0)$ satisfying \eqref{eq:c02f}. The non-degeneracy condition is equivalent to
			\begin{equation*}
				\begin{tikzcd}
					0 \arrow[r] & B \arrow[r, "\partial"] \arrow[d] & A \arrow[r, "\rho"] \arrow[d] & TM \arrow[r] \arrow[d, "\omega^{\flat}"] & 0 \arrow[r] \arrow[d] & 0 \arrow[r] \arrow[d] & 0 \\
					0 \arrow[r] & 0\arrow[r] & 0 \arrow[r] & T^\ast M \arrow[r, "\rho^\ast"'] & A^* \arrow[r, "\partial^\ast"'] & B^* \arrow[r] &0
				\end{tikzcd}
			\end{equation*}
			defining  a quasi-isomorphism. This can be rephrased as
			\begin{equation*}
				\ker (\partial) =0, \quad \ker (\rho )=\im (\partial)
			\end{equation*}
			and $\omega_{|\text{coker}\rho}$ is non-degenerate. Therefore we observe that the only difference with \ref{subsec:sss1} b) is that we do not need an injective anchor and the failure of injectivity is controlled by the next level.
			
			\subsection{Integration of Poisson manifolds}   
			A Poisson manifold $(P,\pi)$ may not always integrates to a symplectic groupoid, see \cite{rui:intp} for the obstructions, but always integrate to a Lie $2$-groupoid $K_\bullet$ with a symplectic form $\omega\in\Omega^2(K_1)$ as demonstrated in \cite{tz:int}. It was shown in \cite{zhu:bg} that $(K_\bullet, \omega)$ gives a $1$-shifted symplectic Lie $2$-groupoid. Here we reproduce their argument.
			
			For a Poisson manifold $(P, \pi)$ the \emph{$A$-path space} $P_0 T^*P$ is
			\[
			P_0 T^*P=\{ a\in C^1(I, T^*P) | \pi ( a)= \gamma', \text{with}\; \gamma=p(a) \; \text{the based path of}\; a, a|_{\partial I}=a'|_{\partial I}=0 \}, 
			\]where $p: T^*P\to P$ is the natural projection. This space has a foliation $F\subseteq T(P_0 T^*P)$ of codimension $2\dim P$ and  take complete local transversals\footnote{This means $\sqcup U_i$ intersects with all the leaves.} $U_i$ with respect to $F$. Denote by $Mon$  the monodromy groupoid of $F$, by $Mon|_{P_0 T^*P \times K_1}$ its source fibre over the submanifold $K_1=\sqcup U_i$ and by $\odot$ the concatenation of two paths. It was shown in \cite{tz:int} that 
			$$K_\bullet=\cdots K_1\times_{K_0} K_1 \times_{\odot, P_0 T^*P, \gt} Mon|_{P_0 T^*P \times K_1}\aaar\sqcup U_i\rightrightarrows P$$
			defines a Lie 2-groupoid with tangent complex given by $\cT^K_\bullet=0\to T^*P\xrightarrow{\pi^\sharp}TP.$
			 
			On the path space $P_0 T^*P$, we have the 2-form 
			\[
			\omega(X, Y ):= \int^1_0 \omega_\can(X(t), Y(t)) \ dt, \quad \forall X, Y \in T(P_0T^*P),
			\]where $\omega_\can$ is the canonical 2-form on $T^*P$. This 2-form, restricted to $K_1$, becomes a symplectic 2-form satisfying $\delta \omega =0$ (\cite[Theorem 1.1]{tz:int}). It is clear that $\omega$ is normalized. The fact that $\omega$ is non-degenerate follows from  the behavior of $\omega$ in a small neighborhood $U$ of $K_0$ in $K_1$, where $(U, \omega)$ forms a symplectic local Lie groupoid.  
			
			\begin{remark}\label{rmk:intpoi}
				Observe that the integration presented here gives a de-singularization of the \emph{Weinstein groupoid} given by $A$-paths modulo $A$-homotopies, see e.g. \cite{rui:int}. It was shown in \cite{cat:sym} that the Weinstein groupoid can be obtained via symplectic reduction from the classical phase space of the Poisson sigma model, see \ref{sec:PSM}, on a square $I\times I$. 
			\end{remark}

			\subsection{The AMM groupoid}\label{sec:amm}
			The \emph{AMM groupoid} was introduced in \cite{xu:mome} to codify the information of group-valued moment maps studied in \cite{amm:lie}. Moreover, the article  \cite{xu:mome} realized that the equivalence between quasi Hamiltonian $G$-spaces and Hamiltonian $LG$-spaces, found in \cite{amm:lie}, is due to a symplectic Morita equivalence between $1$-shifted symplectic Lie groupoids. Here we recall such a construction. 
			
			In this section we assume that $G$ is a Lie group whose Lie algebra is quadratic  $(\g, [\cdot,\cdot], \langle\cdot,\cdot\rangle)$ and denote by $(\g[1],\varpi, Q=d_{CE})$ the corresponding  degree $2$ symplectic $Q$-manifold as in Example \ref{ex:CA}.

			Denote the Maurer-Cartan $1$-forms on $G$ by $\theta^l,\theta^r\in\Omega^1(G;\g).$ The Lie group $G$ acts on itself by conjugation and we denote its action groupoid by $$G^{AMM}_\bullet=G\ltimes G\rightrightarrows G.$$ 
			The degree $2$ symplectic structure $\varpi\in\Omega^2_2(\g[1]),$ i.e. the pairing $\langle\cdot,\cdot\rangle:\g\times\g\to\bR$, induces a $1$-shifted $2$-form on $G^{AMM}_\bullet$ as follows: 
			\begin{eqnarray}\label{eq:AMM}
				\omega_{(g,x)}&=& -\frac{1}{2} \left(\langle Ad_x\operatorname{pr}_1^\ast \theta^l, \operatorname{pr}_1^\ast \theta^l\rangle + \langle \operatorname{pr}_1^\ast\theta^l , \operatorname{pr}_2^\ast(\theta^l +\theta^r)\rangle \right)\in \Omega^2(G^{AMM}_1=G\times G),\\
				H&=& \frac{1}{12} \langle\theta^l, [\theta^l, \theta^l]\rangle\in \Omega^3(G)=\Omega^3(G^{AMM}_0=G),\label{eq:c3f}
			\end{eqnarray}
			where $\operatorname{pr}_{1,2}\colon G\times G\to G$ are the natural projections, $(g,x)\in G\times G$, and $Ad$ denotes the adjoint action of $G$ on $\mathfrak{g}$.
			
			\begin{proposition}[{\cite[Proposition 2.8]{xu:mome})}]\label{prop:amm1}
				The form $\omega_\bullet=\omega+H$ defines a $1$-shifted symplectic structure on $G^{AMM}_\bullet$.
			\end{proposition}
			
			We will give a different proof for this result in \S~\ref{sec:csf}. 	It is worth noticing that:
			\begin{enumerate}
				\item[1.] The infinitesimal counterpart of the AMM groupoid is given by the \emph{Cartan-Dirac structure} $L\subseteq(TG\oplus T^*G, \langle\cdot,\cdot\rangle, \pr_{T}, [\![\cdot,\cdot]\!]_H)$, see e.g. \cite{bcwz:dir} or \cite{abm:pure}.
				\item[2.] If we denote by $\cB G$ the differentiable stack presented by $G\rightrightarrows pt$. Then,  the groupoid $G^{AMM}_\bullet$ presents the \emph{inertia stack} of $\cB G$.
			\end{enumerate}
			
			By using actual loops, see \cite[\S 4.1]{mw:loop} or \cite[Example 3.11]{xu:mome}, one constructs another presentation of the inertia stack of $\cB G$ that is infinite dimensional.  Recall from \cite[\S 4]{seg:loop} that the Lie algebra $L\g$ has a central extension given by $\widehat{L\g} = L\g \times\bR$ with bracket
				$$\left[(\xi_1,t_1),(\xi_2,t_2)\right] =\big ([\xi_1,\xi_2], \sigma(\xi_1,\xi_2) \big) \quad\text{where}\quad \sigma(\xi_1,\xi_2)=\int_{S^1}\langle \xi_1,d\xi_2\rangle.$$
			When $\sigma/2\pi$ represents an integral chomology class,  the Lie algebra $\widehat{L\g}$ integrates to a Lie group $\widehat{LG}$ that is an $S^1$-central extension of $LG$, see \cite[Thm 4.4.1]{seg:loop}.
			
			The coadjoint action of $LG$ on $\widehat{L\g^*} := L\g^* \times \bR$ is given by
			$$g\cdot(\xi,\lambda) = (Ad_g\xi +\lambda dg g^{-1},\lambda ).$$
			When $L\mathfrak{g}^*\cong \Omega^1(S^1,\g)$ is identified with the space of connections on the trivial $G$-bundle over $S^1$, this action on the affine hyper-plane $(\widehat{L\mathfrak{g}^*})_1:=L\g^* \times \{1\} \subseteq \widehat{L\g^*}$  becomes the action by gauge transformations. 
			
			The action groupoid $J_\bullet=LG\ltimes (\widehat{L\mathfrak{g}^*})_1\rightrightarrows (\widehat{L\mathfrak{g}^*})_1$ is a symplectic groupoid with symplectic form $\omega^J\in\Omega^2(J_1)$ given by 
			$$(\pi\times\id)^*\omega^J=(\id\times i)^*\omega_{\can}\quad\text{with}\quad \omega_{\can}\in\Omega^2(T^*\widehat{LG}),$$
			here $\pi:\widehat{LG}\to LG$ denotes the projection of the central extension and $i:(\widehat{L\mathfrak{g}^*})_1\to \widehat{L\mathfrak{g}^*}$ the natural inclusion.
			\begin{remark}
				Perhaps an important observation is that the Poisson manifold $P=(\widehat{L\mathfrak{g}^*})_1$ appears naturally from the AKSZ sigma model perspective as explained in \S~\ref{sec:aksz}. Indeed, there is a natural identification $$(T^*[1]P, \omega_{\can}, Q_\pi)=\big(Maps(T[1]S^1, \g[1]), \bT_{S^1}(\varpi), X_S)$$
				as degree $1$ symplectic $Q$-manifolds.
			\end{remark}
			
			It remains to show that the two models are symplectic Morita equivalent, this is the content of \cite[Prop. 4.26]{xu:mome} and was also implicit in \cite[Prop. 8.1]{amm:lie}. The main ideas of the proof were the following.

			Let $\operatorname{Hol}\colon L\mathfrak{g}\to G$ be the \emph{holonomy map}, i.e. the solution at time $1$  of the differential equation
			\[
			\operatorname{Hol}_s(r)^{-1} \frac{\partial}{\partial s}\operatorname{Hol}_s (r)= r, \quad \operatorname{Hol}_0(r)=e.
			\]
			The pullback groupoid of $G^{AMM}_\bullet$ along the holonomy map is isomorphic to $J_\bullet$ because
			\[
			\operatorname{Hol}^\ast (G\times G)=\{(r_1(s), r_2(s), g) \in L\mathfrak{g}\times L\mathfrak{g}\times G\, | \, g^{-1} \operatorname{Hol}(r_1)g = \operatorname{Hol}(r_2)\}
			\]
			and for any $(r_1(s), r_2(s), g)\in \operatorname{Hol}^\ast(G\times G)$ pick $g(s)\in LG$ given by
			\[
			Ad_{g(s)^{-1}}r_1(s) - g(s)^{-1}\frac{dg(s)}{ds} = r_2(s), \quad g(0)=g,
			\]
			thus the map
			\[
			\operatorname{Hol}^\ast(G\times G)\to LG\times L\mathfrak{g}, \quad (r_1(s), r_2(s), g)\mapsto (g(s), r_1(s))
			\]
			defines  a  diffeomorphism. Moreover, it exchanges the groupoid structures. In consequence, the map $f\colon LG\times L\mathfrak{g}\to G\times G$ defined by
			\[
			f(g(s), r(s))= (g(0), \operatorname{Hol}(0)), \quad f(r(s))= \operatorname{Hol}(r)
			\]
			for all $(r(s))\in L\mathfrak{g}$ and $g(s)\in LG$, provides a hypercover. Finally, define $\eta\in\Omega^2(J_1)$  by
			\[
			\eta= \frac{1}{2} \int_{0}^{1} \langle \operatorname{Hol}_s^\ast \theta^r, \frac{\partial}{\partial s} \operatorname{Hol}_s^\ast \theta^r\rangle ds.
			\]
			Then, it is shown in \cite[Prop 8.1]{amm:lie} that
			\[
			d\eta=-\text{Hol}^*H\quad\text{and}\quad  \omega^J=f^\ast\omega+ \delta(\mu),
			\]
			and therefore this concludes the proof showing that $(G^{AMM}_\bullet, \omega+H)$ and $(J_\bullet, \omega^J)$ are  symplectic Morita equivalent.

			\subsection{Double symplectic groupoids as 2-shifted symplectic}\label{sec:dousym}
			
			In $\S~\ref{sec:doug}$, we explained how from a full double groupoid $D_{\bullet,\bullet}$ one obtains a Lie $2$-groupoid $\overline{W}_\bullet(D)$. In \cite{alv:tran} (after the works \cite{mt:from} and \cite{zhu:bg}) it was shown that Theorem \ref{thm:raj} can be extended to the realm of symplectic geometry. Indeed, this that was the original motivation of \cite{mt:from} but the 
			appropriate framework was missing at that time.
			\begin{equation*}
				\begin{array}{c}
					\xymatrix{D\ar@<-.5ex>[r]\ar@<.5ex>[r]\ar@<-.5ex>[d] \ar@<.5ex>[d]& H\ar@<-.5ex>[d] \ar@<.5ex>[d]\\ G\ar@<-.5ex>[r]\ar@<.5ex>[r]&M}
				\end{array}
			\end{equation*}
			
			For a double Lie groupoid $D_{\bullet,\bullet}$ the analogue of the Bott-Shulman-Stasheff double complex is a triple complex $(\Omega^\bullet(D_{\bullet,\bullet}),\delta^v, \delta^h, d)$ and denote the differential on the total complex by $$\widetilde{D}=\delta^v+(-1)^i\delta^h+(-1)^{i+j}d.$$ It was shown in \cite[Lemma 4.20]{alv:tran} that the $\overline{W}$ functor intertwines the total differentials between the above triple complex of $D_{\bullet,\bullet}$ and the Bott-Shulmann-Stasheff double complex of $\overline{W}_\bullet D$.

			Following \cite{alv:tran}, we say that a double Lie groupoid $D_{\bullet,\bullet}$ is a \emph{double quasi-symplectic groupoid} if there are normalized differential forms $$\omega_{1,1} \in \Omega^2(D=D_{1,1}),\quad \omega_{2,0}\in \Omega^2(H_2=D_{2,0})\quad \text{and}\quad \omega_{1,0} \in \Omega^3(H=D_{1,0})$$ such that $\widetilde{D}(\omega_{1,1}+\omega_{1,0}+\omega_{2,0})=0,$ and the following non-degeneracy condition is satisfied
			\begin{equation}\label{eq:nonddou}
				\ker{\omega_{1,1}}_{|x} \cap \ker T_{ |x} \gt^h \cap \ker T_{ |x} \gs^h= 0 \quad \forall x \in H\quad \text{and}\quad  \dim D = 2 \dim H.
			\end{equation}
			With these definitions, we can state \cite[Theorem 4.19]{alv:tran}.
			
			\begin{theorem}\label{thm:daniel}
				If $(D_{\bullet,\bullet}, \omega_{1,1}+\omega_{1,0}+\omega_{2,0})$ is a full double quasi-symplectic groupoid then $$\big(\overline{W}_\bullet D, \overline{W}^*(\omega_{1,1}+\omega_{1,0}+\omega_{2,0})\big)$$ is a $2$-shifted symplectic Lie $2$-groupoid.
			\end{theorem}
			\begin{proof}
				By Theorem \ref{thm:raj}, we know that $\overline{W}_\bullet D$ is a Lie $2$-groupoid. Since the forms are normalized and closed for the total differential, then  \cite[Lemma 4.20]{alv:tran} implies that $\overline{W}^*(\omega_{1,1}+\omega_{1,0}+\omega_{2,0})$ is a closed $2$-shifted $2$-form. It remains to show that \eqref{eq:nonddou} implies the  non-degeneracy  condition for the $2$-shifted $2$-forms, and that is exactly the content of  \cite[Lemma 4.21]{alv:tran}.
			\end{proof}
			
			When $\omega_{2,0}=0=\omega_{1,0}$, one recovers the usual notion of a \emph{double symplectic groupoid} as in \cite{luwe}. If one applies the horizontal differentiation in \eqref{diag:diff} to a  double symplectic groupoid $(D_{\bullet,\bullet}, \omega)$, the result is a double Lie algebroid on a cotangent prolongation, see \cite{luwe, mac:sym}, that codifies a Lie bialgebroid given by the $PQ$-manifold $(A_G[1], Q, \pi)$ as in Corollary \ref{cor:bia}.  It was explained in \cite{mt:from} that $(\overline{W}_\bullet D, \overline{W}^*(\omega))$ integrates the Courant algebroid $A\bowtie A^*$, i.e. the Drinfeld double of the Lie bialgebroid $(A_G[1], Q, \pi)$ \cite{lwx:ca}. Therefore, using $\S~\ref{sec:liebia}$ we get
			$$\text{Lie}(\overline{W}_\bullet D, \overline{W}^*(\omega))=(T^*[2]A_G[1],\omega_{\can},X_{Q+\pi}).$$
			Hence, we have shown the commutativity of \eqref{diag:diff} in this particular example.  
			
			For more examples of double quasi-symplectic groupoids see \cite{alv:tran}. We just point out that the definition of a double quasi-symplectic groupoid can be generalized to include additional differential forms. Its differentiation  then produces examples of  $2$-shifted symplectic Lie $2$-algebroids, see $\S~\ref{sec:ssa}$, which  encode Courant algebroids twisted by $4$-forms \cite{stro:ca4form, pym:shift}.

			\subsection{The classifying space of a Lie group}\label{sec:bg}
			
			As in $\S~\ref{sec:amm}$, let $G$ be a Lie group whose Lie algebra is quadratic and denote by $(\g[1], \varpi, Q=d_{CE})$ the corresponding degree $2$ symplectic $Q$-manifold, see Example \ref{ex:CA}.  We still denote the Maurer-Cartan $1$-forms on $G$ by $\theta^l,\theta^r\in\Omega^1(G;\g).$
			
			Denote by $\cB G$ the differentiable stack presented by the Lie group $G\rightrightarrows pt$. In this section we give three different models for the $2$-shifted symplectic structure on $\cB G$ discussed in \cite{ptvv}. Our presentation follows \cite[\S~3 and \S~4]{zhu:bg}.
			
			A clear way of presenting $\cB G$ is using $N_\bullet G$, the nerve of the Lie group $G$, as explained in $\S~\ref{sec:lie1grou}$. 
			
			\begin{proposition}\label{prop:BG}
				The Lie $1$-group $N_\bullet G$ carries a $2$-shifted symplectic structure $\omega_\bullet$ given by 
				\[
				\omega_2= \frac{1}{2}\langle d_2^\ast \theta^l , d_0^\ast \theta^r\rangle \in \Omega^2(N_2G=G\times G), \quad \omega_1= -\frac{1}{12} \langle\theta^l , [\theta^l , \theta^l ]\rangle\in \Omega^3(N_1G=G),
				\]
				Moreover, $\text{VE}(\omega_\bullet)=\varpi\in\Omega^2_2(\g[1]).$
			\end{proposition}
			The fact that $\delta\omega_2=0$ already appeared in \cite{ polwie}. A full proof for $D(\omega_\bullet)=0$ can be found in  \cite[\S~6.4]{bry:loop} or in \cite[Cor. 3.5]{wei:mod}. The non-degeneracy condition was explicitly computed in \cite[Theorem 3.1]{zhu:bg}.
			
			The relevance of this structure comes from the fact that it gives the \emph{Atiyah-Bott symplectic structure} \cite{ab:ym} on the moduli space of flat $G$-connections on a closed surface $\Sigma$, see \cite{lis:sym, wei:mod}. The geometric realization of $N_\bullet G$ gives the (topological) classifying space $BG$ of $G$, see e.g. \cite{bss:on}. The $4$-class $[\omega_\bullet]\in H^4_{sing}(BG)$ is known as the \emph{first Pontryagin class} and is a fundamental subject of study in mathematics and physics, see e.g. \cite{baez:from, bre:pon, bry:loop, dw:top}.
			
			\begin{remark}[On differentiation II]
				In this example, we see that the differentiation of a $2$-shifted symplectic Lie $1$-groupoid yields a degree $2$ symplectic $Q$-manifold. Or from the opposite perspective, in general, the integration of a symplectic $Q$-manifold  generally results in a shifted symplectic higher Lie groupoid. Thus, even if one is primarily interested in  symplectic $Q$-manifolds, integrating them necessarily involves introducing shifted symplectic Lie $n$-groupoids.
			\end{remark}
			
			\subsubsection{The 2-truncation}
			One can think of $N_\bullet G$ as the $1$-truncation of the universal integration of the Lie algebra $\g$, see \cite{hen:int} for details. What we explain now is the $2$-truncation $\mathbb{G}_\bullet$, that is an infinite dimensional Lie $2$-group.   
			
			The spaces $\mathbb{G}_k$ are given by:
			\begin{itemize}
				\item $\mathbb{G}_0= pt$;
				\item $\mathbb{G}_1= P_eG$, consists of paths\footnote{To make $\mathbb{G}$ into a Banach Lie $2$-group one needs to fix a Sobolev class. We will ignore such analytic problems here, see \cite{zhu:bg} for a proper definition.} on $G$ starting at the identity;
				\item $\mathbb{G}_2=\Omega G$, consists of loops on $G$ that start and end at the identity;
				\item $\mathbb{G}_3= \{(\tau_0, \tau_1, \tau_2)\in \mathbb{G}_2^{\times 3} \, | \, \tau_0(t)= \tau_1(t) , \, \tau_1(\tfrac{2}{3}+ t)= \tau_2(\tfrac{2}{3}+t) , \, \tau_0(1-t)= \tau_2(t), \, t\in [0, \tfrac{1}{3}]\}$.
			\end{itemize}
			In other words, elements in $\mathbb{G}_3$ are empty tetrahedrons on $G$ based at the identity.

			The face maps $d_i\colon \Omega G\to P_eG$ for $i=0,1,2$ are defined by
			\[
			d_0(\tau)(t)= \tau\left(\frac{t}{3}\right), \quad d_1(\tau)(t)=\tau\left(1-\frac{t}{3}\right), \quad d_2(\tau)(t)= \tau\left(\frac{1}{3}+ \frac{t}{3}\right)\cdot \tau\left(\frac{1}{3}\right)^{-1},
			\]
			with $t\in [0,1]$, and $d_i\colon \mathbb{G}_3\to \Omega G$ for $i=0,1,2,3$ are given by
			\[
			d_i(\tau_0, \tau_1, \tau_2)= \tau_i \quad\text{where}\quad \tau_3(t)= 
			\begin{cases}
				\tau_0\left(t+\frac{1}{3}\right) \cdot \tau_0\left(\frac{1}{3}\right)^{-1} & t\in [0, \frac{1}{3}],\\
				\tau_2(t)\cdot \tau_0\left(\frac{1}{3}\right)^{-1} & t\in [\frac{1}{3}, \frac{2}{3}], \\
				\tau_1\left(\frac{4}{3}- t\right)\cdot \tau_0\left(\frac{1}{3}\right)^{-1} & t\in [\frac{2}{3}, 1].
			\end{cases}
			\]
			The degeneracy maps $s_i\colon P_eG \to \Omega G$ for $i=0,1$ are given by
			\[
			s_0(\gamma)(t)= 
			\begin{cases}
				\gamma(3t) & t\in [0, \frac{1}{3}],\\
				\gamma(1) & t\in [\frac{1}{3}, \frac{2}{3}], \\
				\gamma(3-3t) & t\in [\frac{2}{3}, 1],
			\end{cases}\quad\text{and}\quad 
			s_1(\gamma)(t)= 
			\begin{cases}
				e & t\in [0, \frac{1}{3}],\\
				\gamma(3t-1) & t\in [\frac{1}{3}, \frac{2}{3}], \\
				\gamma(3-3t) & t\in [\frac{2}{3}, 1].
			\end{cases}
			\]
			A $2$-form on the based loop group $\omega\in\Omega^2(\mathbb{G}_2=\Omega G)$ was given in \cite{seg:loop} as follows
			\begin{equation}\label{eq:segalsform}
				\omega_\tau(a,b)= \int_{S^1}\langle \hat{a}'(t), \hat{b}(t)\rangle dt \quad\text{for}\quad \tau\in \Omega G, \ a, b\in T_\tau \Omega G
			\end{equation}
			where  $\hat{a}(t)=L_{\tau(t)^{-1}}a(t)$, $\hat{b}(t)=L_{\tau(t)^{-1}}b(t)$, and the prime denotes time derivative. As a consequence of Theorem \ref{thm:SMEbg}, that we  will see later, we obtain that  $(\mathbb{G},\omega)$ is a $2$-shifted symplectic Lie $2$-group. 
			
			To our knowledge, the multiplicativity of  Segal's $2$-form $\omega$ was not identified until \cite{zhu:bg}. In \cite{seg:loop}, it was shown that $(\Omega G, \omega)$ is a symplectic manifold. If we compare with the definition of a shifted symplectic structure, the fact that \eqref{eq:segalsform} is symplectic is much stronger. Moreover, if we further have an integral condition, that is $[\omega]\in H^2(\Omega G,\mathbb{Z})$, then the canonical prequantum $S^1$ bundle $\widetilde{\Omega G}$ also gives rise to a Lie $2$-group. Explicitly, the $S^1$ principal bundle over $\Omega G$ with curvature $\omega$ is the $S^1$-central extension $\widetilde{\Omega G}$ of $\Omega G$ classified by $[\omega]$. As proven in \cite{baez:from, hen:int}, there is a Lie 2-group structure on the simplicial manifold
			\begin{equation} \label{eq:stringG}
				\text{String(G)}_\bullet=... \widetilde{\Omega G} \aaar P_eG\rightrightarrows  pt,
			\end{equation}
			and it gives a model of the string group $\text{String}(G)$, associated to $G$. Thus, in this way, the string Lie $2$-group provides a prequantization of $(\mathbb{G}_\bullet, \omega)$. Moreover, it was shown in \cite{baez:from, hen:int} that the Lie $2$-algebra $\text{Lie(String(G))}$ is quasi isomorphic to $\text{String}(\g)$ as computed in \eqref{eq:stringalg}.

			As shown in \cite[Section 2.1]{rem:cs} and \cite{mei:mod}, the manifold $(\Omega G, \omega)$ (the second level of a 2-shifted Lie 2-groupoid) is obtained via symplectic reduction as the reduced phase space of the $3d$ Chern-Simons  sigma model with source $\bDelta^2\times [0,1]$, cf. with Remark \ref{rmk:intpoi}. As explained in \S~\ref{sec:CS}, the Chern-Simons action on a $3$-dimensional manifold $N$ is given by $$S(A)=\int_N \frac{1}{2}\langle A, dA\rangle+\frac{1}{6}\langle A, [A,A]\rangle\quad \text{with}\quad A\in\Omega^1(N,\g).$$
			For a surface with boundary $\Sigma$, the Chern-Simons theory on $N=\Sigma\times [0,1]$ gives a symplectic manifold $\cA_{\Sigma}=\{$Connections on $G$-principal bundles with base $\Sigma\}$ with a Hamiltonian action of the gauge group $\cG_\Sigma=\{g:\Sigma\to G\ |\ g(\partial\Sigma)=e\}$ and moment map $\mu:\cA_\Sigma\to(\g_\Sigma)^*$
			\begin{equation*}
				\mu(A)(\alpha)=\int_\Sigma \langle\alpha, dA+\frac{1}{2}[A,A]\rangle, \qquad A\in\Omega^1(\Sigma;\g), \ \alpha\in\Omega^0(\Sigma,\g).
			\end{equation*}
			The symplectic quotient $\cA_\Sigma//\mu^{-1}(0)$ is the moduli space of flat connections on $\Sigma$. If  the surface is a disk, i.e. $\Sigma=D$, it was shown in \cite[Example 3.1]{mei:mod} that the moduli space of flat connections is the based loop group with the Segal's $2$-form, more precisely $$\cA_D//\mu^{-1}(0)=(\Omega G,\omega).$$

			\subsubsection{Manin triples}
			The last model for $\cB G$ that we present is using Manin triples. A \emph{Manin triple} is a triple $(\mathfrak{g}, \mathfrak{h}_+, \mathfrak{h}_-)$ where $\mathfrak{g}$ is a quadratic Lie algebra and $\mathfrak{h}_+, \mathfrak{h}_-\subseteq \mathfrak{g}$ are two Lie subalgebras such that
			\[
			\langle \mathfrak{h}_+,\mathfrak{h}_+\rangle = \langle \mathfrak{h}_-, \mathfrak{h}_-\rangle=0, \quad \mathfrak{h}_+\oplus \mathfrak{h}_-= \mathfrak{g}.
			\]
			There is a one-to-one correspondence between Manin triples and Lie bialgebras, that are the infinitesimal counterparts of Poisson-Lie groups \cite{dri:qua}.

			Let $(\mathfrak{g}, \mathfrak{h}_+, \mathfrak{h}_-)$ be a Manin triple. Denote with $H_+$, $H_-$ and $G$ the connected and simply-connected Lie groups integrating $\mathfrak{h}_+$, $\mathfrak{h}_-$ and $\mathfrak{g}$ respectively. Moreover, we denote by $\varphi_\pm\colon H_\pm\to G$ the integration of the inclusions $\mathfrak{h}_\pm \to \mathfrak{g}$. It was shown in \cite{luwe} how to produce a  double symplectic Lie groupoid (see \S~\ref{sec:dousym}) that integrates the given Manin triple. Its construction is as follows: consider
			\[
			\Gamma=\{(h_2, a_2,a_1, h_1)\in H_+\times H_-\times H_-\times H_+ \, | \, \varphi_+(h_2) \varphi_-(a_1)= \varphi_-(a_2) \varphi_+(h_1)\}
			\]
			and define the double symplectic groupoid by
			\begin{equation}
				\label{eq:dlg_gamma}
				\begin{tikzcd}
					\Gamma \arrow[r, shift left=0.5ex, ] \arrow[r, shift right=0.5ex, ] \arrow[d, shift right=0.5ex, '] \arrow[d, shift left=0.5ex, ] & H_+ \arrow[d, shift right=0.5ex] \arrow[d, shift left=0.5ex] \\
					H_- \arrow[r, shift right=0.5ex] \arrow[r, shift left=0.5ex] & \ast
				\end{tikzcd}
			\end{equation}
			where the structural maps are 
			\begin{eqnarray*}
				s^h(\xi)= h_1, \quad t^h(\xi) = h_2, \quad m^h(\xi, \xi')= (h_2, a_2a_2', a_1a_1', h_1'),\\
				s^v(\xi)= a_1, \quad t^v(\xi)= a_2, \quad m^v(\xi, \xi')= (h_2h_2', a_2, a_1', h_1h_1'),
			\end{eqnarray*}
			for $\xi=(h_2, a_2, a_1, h_1)$ and $\xi'=(h_2', a_2', a_1', h_1')\in \Gamma$ and the symplectic form $\omega_\Gamma\in\Omega^2(\Gamma)$, as shown in \cite[Thm. 3]{am:slp}, is given by
			\[
			\omega_\Gamma :=\frac{1}{2}\Big( \langle (t^h)^\ast \theta^l_{H_+} , (s^v)^\ast \theta^r_{H_-}\rangle - \langle (t^v)^\ast \theta^l_{H_-} , (s^h)^\ast \theta^r_{H_+}\rangle\Big) \in \Omega^2(\Gamma). 
			\]
			An immediate consequence of Theorem \ref{thm:daniel} is that $(\overline{W}_\bullet(\Gamma),\overline{W}^*\omega_\Gamma)$ is a $2$-shifted symplectic (local) Lie $2$-group.

			\subsubsection{The equivalences}
			The evaluation map $ev_1\colon P_eG\to G$ extends to a simplicial morphism $ev_\bullet \colon \mathbb{G}_\bullet\to N_\bullet G$ with  $ev_2\colon \Omega G\to G\times G$  defined by
			\[
			ev_2(\tau)=\left(\tau\left(\frac{2}{3}\right), \tau\left(\frac13\right)^{-1}, \tau\left(\frac13\right)\right),
			\]
			
			One of the main results in \cite{zhu:bg} is to show that the above three models are symplectic Morita equivalent, and thus they provide three different presentations of the $2$-shifted symplectic structure on $\cB G$ studied in \cite{ptvv}. The following result summarizes Theorems 3.15 and 4.8 in \cite{zhu:bg}.
			
			\begin{theorem}\label{thm:SMEbg}
				The following are symplectic Morita equivalences:
				\begin{enumerate}
					\item[1.] For $G$ connected and simply connected, 
					$$(\mathbb{G}_\bullet, \omega)\xleftarrow{\id}(\mathbb{G}_\bullet, \eta^P)\xrightarrow{ev_\bullet}(N_\bullet G,\omega_\bullet)$$
					where $\eta^P\in\Omega^2(\mathbb{G}_1=P_eG)$ is given by
					\[
					\eta^P_\gamma(u,v)= \frac12 \int_{0}^{1} (\hat{u}'(t), \hat{v}(t))- (\hat{u}(t), \hat{v}'(t))dt,
					\]
					for all $\gamma\in P_eG$, $u,v\in T_\gamma P_eG$ and $\hat{u}'(t)= L_{\gamma(t)^{-1}}u(t)$, $\hat{v}'(t)= L_{\gamma(t)^{-1}}v(t)$.
					\item[2.] For $H_+$ a complete Poisson Lie group,
					$$(\overline{W}_\bullet(\Gamma),\overline{W}^*\omega_\Gamma)\xleftarrow{\id}(\overline{W}_\bullet(\Gamma),\beta)\xrightarrow{\Phi_\bullet}(N_\bullet G,\omega_\bullet)$$
					where $\beta\in\Omega^2(\overline{W}_1(\Gamma)=H_- \times H_+)$ is $\beta=\langle\theta_{H_-}^l , \theta^r_{H_+}\rangle$ and $\Phi_\bullet$ is given by 
					\[
					\Phi_1(a,h)\mapsto \varphi_-(a)\varphi_+(h)\quad\text{and}\quad  \Phi_2 (a_3, (h_3, a_2, a_1, h_2), h_1)\mapsto (\varphi_-(a_3)\varphi_+(h_3), \varphi_-(a_1)\varphi_+(h_1)).
					\]
				\end{enumerate}
				
			\end{theorem}
			
			The proof of the above result is not difficult but long; therefore, we suggest the interested reader to look into the original work \cite{zhu:bg}. We also mention that in \cite{cgm:dir}, using Dirac structures,  a relation between the forms $\omega\in\Omega^2(\Omega G)$ and $\omega_2\in\Omega^2(G\times G)$ was found.

			\subsection{A 3-shifted symplectic Lie 2-group}
			Let $G$ be a Lie group with Lie algebra $\g$ and denote by $(\g[1], Q=d_{CE})$ the corresponding degree $1$ $Q$-manifold. We saw in $\S \ref{sec:higcou}$ that the degree $3$ symplectic $Q$-manifold $(T^*[3]\g[1],\omega_\can,  \Lie_Q)$ encodes the higher Courant algebroid $(\g\oplus\wedge^2\g^*,\langle\cdot,\cdot\rangle,[\![\cdot,\cdot]\!])$ with pairing and bracket given as in  \eqref{eq:brahigc}. Here we exhibit the integration found in \cite[\S 5.1]{cue:ve}.
			
			The coadjoint representation of G on $\g^*$ defines a full double Lie groupoid 
			\begin{eqnarray*}
				\xymatrix{G\ltimes\g^*\ar@<-.5ex>[r]\ar@<.5ex>[r]\ar@<-.5ex>[d] \ar@<.5ex>[d]& G\ar@<-.5ex>[d] \ar@<.5ex>[d]\\ pt\ar@<-.5ex>[r]\ar@<.5ex>[r]&pt}
			\end{eqnarray*}
			where the horizontal groupoid is the one associated to the vector bundle $\pr_G:T^*G\to G$ and the vertical groups are: the semi-direct product with the coadjoint representation and the group itself. So Theorem \ref{thm:raj} implies that $\overline{W}_\bullet(G\ltimes\g^*)$ is a Lie $2$-group. 
			
			Denoting the right Maurer-Cartan $1$-form on $G$ by $\theta^r\in\Omega^1(G;\g)$, the tautological 1-form and the canonical symplectic form on $T^*G=\g^*\times G$ can be written as
			$$\theta=\langle \pr_1, \pr_2^*\theta^r\rangle\in\Omega^1(\g^*\times G)\quad\text{and}\quad \omega=-d\theta\in\Omega^2(\g^*\times G).$$
			It was shown in \cite[Proposition 5.1]{cue:ve} that the Lie $2$-group $\overline{W}_\bullet(G\ltimes\g^*)$ carries a $3$-shifted symplectic structure given by  $$p_2^*\omega\in\Omega^2(\overline{W}_3(G\ltimes\g^*)),$$ where 
			$p_{2}:\overline{W}_3(G\ltimes\g^*)\to \g^*\times G,$ $p_{2}(\xi_1, \xi_2,g_2; \xi,g_1; g_0)=(\xi_2,g_1).$ Moreover, the image of $p_2^*\omega$ under the van Est map defined in that article coincides with the canonical symplectic form on the degree $3$ $Q$-manifold $(T^*[3]\g[1], \Lie_Q)$. So we get
			$$\text{Lie}\big(\overline{W}_\bullet(G\ltimes\g^*)\big)=(T^*[3]\g[1], \Lie_Q)\quad \text{and}\quad \text{VE}(p_2^*\omega)=\omega_\can\in\Omega^2_3(T^*[3]\g[1]).$$
			
			Therefore, the $3$-shifted symplectic Lie $2$-group $\big(\overline{W}_\bullet(G\ltimes\g^*), p^*_2\omega\big)$ provides an integration of the higher Courant algebroid $(\g\oplus\wedge^2\g^*,\langle\cdot,\cdot\rangle,[\![\cdot,\cdot]\!])$. Moreover, denote by $\cB G$ the differentiable stack presented by the Lie group $G\rightrightarrows pt$. Then $\overline{W}_\bullet(G\ltimes\g^*)$ provides a presentation for the differentiable stack $T^*[3]\cB G,$ and the $3$-shifted symplectic form $p_2^*\omega$ gives a concrete model for the $3$-shifted symplectic structure on $T^*[3]\cB G$ introduced in \cite{cal:cot}.

			\subsection{On Prequantization}
			
			In Proposition \ref{prop:prequantum} we showed that symplectic $Q$-manifolds admit a prequantum $Q$-bundle. For shifted symplectic Lie $n$-algebroids, as  introduced in $\S~\ref{sec:ssa}$, prequantization  becomes  a significantly more involved problem. Using the language of higher gerbes,  a prequantization scheme  for shifted symplectic structures on stacks was proposed in \cite{saf:pre}.  However, a general prequantization theory for shifted symplectic structures on higher Lie groupoids is still lacking.
			
			Nevertheless, there have been several attempts in the literature to construct meaningful prequantum bundles for specific examples. Some notable contributions include:
			\begin{itemize}
				\item The pioneering works of Alan Weinstein on using symplectic groupoids to quantize Poisson manifolds, see e.g. \cite{bw:geoqua, w2:q, w1:q}.
				\item The work \cite{haw:quan}, which develops a full quantization scheme for symplectic groupoids; see also \cite{cc:jacobi} for the integrability problem and \cite{bon:quan, bon:q2} for further examples.
				\item  The article \cite{xu:quan}, which introduces prequantization of $1$-shifted symplectic Lie $1$-groupoids; see \cite{kre:qua} for a recent treatment and \cite{mtv:shif} for its connection with shifted contact geometry.
				\item As mentioned in in $\S~\ref{sec:bg}$,  there has been extensive study over the past twenty years on the prequantization of the $2$-shifted symplectic structure on the classifying space, e.g. \cite{baez:from, cjmsw, wald:string, mei:bas}.
				\item  For the interplay between field theories, higher gerbes, and prequantization, see  \cite{frs:hig, frs:pre}.
			\end{itemize}

			\section{Shifted Lagrangian structures}\label{sec:7}
			For a symplectic manifold, the natural subobjects of interest are \emph{Lagrangian submanifolds}.   In this section, we adapt the notion of   shifted Lagrangian structures from \cite{ptvv} to the realm of $m$-shifted symplectic Lie $n$-groupoids. As an application, we show how various moment map theories can be reformulated in terms of shifted Lagrangian structures. Building on ideas of Alan Weinstein \cite{wei:symcat}, it was demonstrated in 
			\cite[Thm. 4.4.]{cal:lag} that such structures can be composed. We adapt their work to the smooth setting and describe its connection  with the reduction procedures arising from moment map theories.  Most of the material in this section is based on \cite{bur:lag}.
			
			\subsection{Definition}
			
			Here we introduce shifted Lagrangian structures on Lie $n$-groupoids. The point of view presented here is a straightforward translation of the definitions introduced in \cite{ptvv} and \cite{cal:lag} in the realm of derived algebraic geometry. A precise treatment of this subject will appear in \cite{bur:lag}.
			
			Recall that a submanifold $j:L\to M$ of a symplectic manifold $(M,\omega)$ is called \emph{Lagrangian} if 
			\begin{equation}\label{eq:clalag}
				j^*\omega=0\quad\text{and the map} \quad \nu_L\xrightarrow{T^*j\circ \omega^\flat}T^*L    
			\end{equation}
			is an isomorphism, where $\nu_L=j^*TM/TL$ is the normal bundle.
			
			Let $(K_\bullet, \omega_\bullet)$ be an $m$-shifted symplectic Lie $n$-groupoid. In order to adapt the above definition to the realm of Lie $n$-groupoids we need to introduce the following concepts. The analog of the first condition in \eqref{eq:clalag} is as follows: An \emph{isotropic structure} on a Lie $n$-groupoid morphism $\Phi_\bullet:J_\bullet \to K_\bullet$ is an $(m-1)$-shifted $2$-form $\eta_\bullet$ on $J_\bullet$ such that $$\Phi_\bullet^*\omega_\bullet=D\eta_\bullet.$$
			
			To reproduce the second condition in \eqref{eq:clalag} we need a replacement for the normal bundle of a submanifold. Notice that we can interpret $\nu_L$  as the cohomology of the mapping cone for the morphism $Tj:TL\to j^*TM$. Therefore, we make an analogous definition using the tangent complex: The \emph{normal complex} $\nu_\Phi$ of a Lie  $n$-groupoid morphism $\Phi_\bullet:J_\bullet\to K_\bullet$ is the mapping cone of the induced map $\cT \Phi:\cT^J\to \Phi_0^*\cT^K$ of tangent complexes, i.e. 
			\begin{equation}\label{eq:normalcx}
				\nu_\Phi:=\left(\cT^J\oplus\Phi_0^*\cT^K[-1], \begin{pmatrix}\partial^J& 0\\ \cT\Phi&-\partial^K\end{pmatrix}\right).
			\end{equation}
			An immediate consequence, see \cite{cal:lag}, of the normal being a mapping cone is the next result.
			\begin{proposition}\label{prop:iso} 
				Let $(K_\bullet,\omega_\bullet)$ be an $m$-shifted symplectic Lie $n$-groupoid, and  let $(J_\bullet, \Phi_\bullet, \eta_\bullet)$ be an isotropic morphism into $(K_\bullet,\omega_\bullet)$.
				Then the following is a cochain map 
				\begin{equation*}
					\begin{array}{rccl}
						\lambda^{\Phi,\eta}:&\nu_\Phi&\to& \cT^{*J}[m-1]  \\
						&v\oplus w&\mapsto& \lambda^{\eta_{m-1}}(v,\cdot)+ (\cT\Phi)^*\circ \lambda^{\omega_m} (w,\cdot)
					\end{array}
				\end{equation*}
				where $\lambda^{\eta_{m-1}}$ and $ \lambda^{\omega_m}$ are the IM-pairings as in \eqref{eq:im-form}.
			\end{proposition}
			
			We are ready to translate Definition 2.8 in  \cite{ptvv}, see also \cite[\S~2.1.1]{cal:lag} and \cite[\S~2.1]{saf:qua}, to higher Lie groupoids. An isotropic morphism $(J_\bullet, \Phi_\bullet, \eta_\bullet)$ into an $m$-shifted symplectic Lie $n$-groupoid $(K_\bullet, \omega_\bullet)$ is called a \emph{shifted Lagrangian} if the cochain map 
			\begin{equation}\label{eq:lagcon}
				\lambda^{\Phi,\eta}:\nu_\Phi\to \cT^{*J}[m-1]
			\end{equation}
			is a quasi-isomorphism. In this case, we say that $\eta_\bullet$ is a {\em shifted Lagrangian structure} on the morphism $\Phi_\bullet$. 
			
			\begin{remark}[On differentiation III]
				In Remark \ref{rmk:diff1} we pointed out that a van Est map for higher Lie groupoids will give a way to differentiate $m$-shifted symplectic structures. The same may hold for shifted Lagrangian structures and it is expected to give the shifted Lagrangians introduced in \cite{pym:shift}, see also \S~\ref{sec:ssa}.
			\end{remark}
			
			\begin{remark}[On Morita invariance]
				As currently introduced, shifted Lagrangians are not Morita invariant. A way to overcome this is by replacing the strict morphism $\Phi_\bullet$ by a zig-zag of morphisms with the left leg given by a hypercover. We do not explore such features here.
			\end{remark}
			
			\subsection{First Examples}
			To guide the discussion, we provide some basic examples of shifted Lagrangian structures. 
			
			\subsubsection{Shifted Lagrangians on a point}\label{sec:lagpt}
			
			Recall that a point $pt$, regarded as a constant simplicial manifold, is $m$-shifted symplectic for any $m\in\bN$ with the zero $2$-form, see \S~\ref{subsub:sss0}. Then, a shifted Lagrangian morphism into a point, viewed as $m$-shifted symplectic, is the same as a Lie $n$-groupoid $J_\bullet$ equipped with an $(m-1)$-shifted $2$-form $\eta_\bullet$ such that $D\eta_\bullet=0 $ and $\lambda^{\Phi,\eta}=\lambda^{\eta_{m-1}}:\cT^J\to \cT^{*J}[m-1]$ is a quasi-isomorphism (note that in this case $ \nu_\Phi=\cT^J$), i.e., $(J_\bullet,\eta_\bullet)$ is an $(m-1)$-shifted symplectic Lie $n$-groupoid. 
			
			\subsubsection{The diagonal and symplectic Morita equivalences}\label{sec:diag}
			Let $(K_\bullet, \omega_\bullet)$ be an $m$-shifted symplectic Lie $n$-groupoid and denote by $\overline{K}$ the $m$-shifted symplectic groupoid defined by $-\omega_\bullet$. Then the zero $2$-form is a shifted Lagrangian structure on the diagonal morphism $\text{Diag}:K_\bullet\to \overline{K}_\bullet\times K_\bullet$. Indeed the isotropic condition follows from $\text{Diag}^*(\text{pr}_1^*\omega_\bullet-\text{pr}_2^*\omega_\bullet)=0$. To check that $(\text{Diag},0)$ is shifted Lagrangian, consider the normal complex
			$\nu_{\text{Diag}}=\cT^{K}\oplus \cT^K [-1]\oplus \cT^K[-1]$ as in \eqref{eq:normalcx}, and  note that the map
			$\nu_\text{Diag}\to \cT^K[-1]$, $(u,v,w)\mapsto v-w$ is a quasi-isomorphism. Then, the map $\lambda^{\text{Diag},0}: \nu_\text{Diag} \to \cT^{*K}[m-1]$, which is given by
			$(u,v,w) \mapsto \lambda^{\omega_m}(v)-\lambda^{\omega_m}(w)$, is a quasi-isomorphism since it is the composition of two quasi-isomorphisms:
			\begin{equation*}
				\xymatrix{\nu_{\text{Diag}} \ar[d]_{}\ar[rr]^{\lambda^{\text{Diag},0}}&& \cT^{*K}[m-1]\\ \cT^K[-1]\ar[urr]_{\lambda^{\omega_m}}&}
			\end{equation*}
			
			An application of this example is as follows. Recall from \S~\ref{sec:defss} that two $m$-shifted symplectic Lie $n$-groupoids $(K_\bullet, \omega^K_\bullet)$ and $(J_\bullet, \omega^J_\bullet)$ are symplectic Morita equivalent if there exist another Lie $n$-groupoid $Z_\bullet$ with an $(m-1)$-shifted $2$-form $\eta_\bullet$ and hypercovers $\Phi_\bullet:Z_\bullet \to K_\bullet$ and $\Psi_\bullet:Z_\bullet\to J_\bullet$ satisfying
			$$\Psi_\bullet^*\omega^J_\bullet-\Phi_\bullet^*\omega^K_\bullet=D\eta_\bullet.$$
			
			Using Proposition \ref{prop:hyp} and the diagonal it is easy to show that $\eta_\bullet$ is a shifted Lagrangian structure on the morphism
			$$Z_\bullet\xrightarrow{\text{Diag}}Z_\bullet\times Z_\bullet\xrightarrow{\Phi_\bullet\times \Psi_\bullet} \overline{K}_\bullet\times J_\bullet$$
			where $\overline{K}_\bullet\times J_\bullet$ is endowed with the $m$-shifted symplectic structure $\pr_J^*\omega^J_\bullet-\pr^*_K\omega_\bullet^K.$
			
			\subsubsection{The path object}\label{sec:path}
			
			On a category of fibrant objects, the path object allows us to compute fiber products, see e.g. \cite{gh:simp}. Even though Lie $n$-groupoids do not form a category of fibrant objects, see \cite{rcc:hom}, they nevertheless possess a path object  that allows us to compute fiber products in many situations. Its construction is as follows.
			
			The category of simplicial sets has a \emph{cartesian product} given by  $$(X\times Y)_l=X_l\times Y_l,$$ where $X_\bullet$ and $Y_\bullet$ are simplicial sets. In particular, the non-degenerate $(n+1)$-simplices of $\Delta[n]\times \Delta[1]$ are in one-to-one correspondence with 
			\begin{equation}\label{eq:EZ-map}
				\varrho_a: \Delta[n+1]\to \Delta[n]\times \Delta[1],
			\end{equation}
			morphism of simplicial sets, such that, as permutations of $(n+1)$-elements, they are $(n,1)$-shuffles. Therefore on $\Delta[n]\times\Delta[1]$ we have exactly $n+1$ non-degenerate simplices.
			
			On simplicial sets one introduces a mapping space as the adjoint to the cartesian product, see e.g. \cite[\S~I.5]{gh:simp} and compare with the mapping space for $\bZ$-graded manifolds defined in \S~\ref{sec:akszT}. More concretely, given simplicial sets $X_\bullet, Y_\bullet$ the \emph{mapping space} is the simplicial set given by $$\text{Maps}_l(X, Y)=\Hom(\Delta[l]\times X, Y).$$ 
			\begin{proposition}\label{prop:trans}
				Let $K_\bullet$ be a Lie $n$-groupoid. The following assertions hold:
				\begin{enumerate}
					\item[1.] There is an evaluation map $ev:\Delta[1]\times\text{Maps}(\Delta[1],K_\bullet)\to K_\bullet,\  ev_l(u,f)=f(0\cdots l, u)$.
					\item[2.] The simplicial set $K^{\Delta[1]}_\bullet:=\text{Maps}(\Delta[1],K_\bullet)$ is a Lie $n$-groupoid. Moreover, it is Morita equivalent to $K_\bullet=\text{Maps}(\Delta[0],K)$ via the  weak equivalence $$(\sigma^0)^*:\text{Maps}(\Delta[0],K)\to \text{Maps}(\Delta[1],K).$$
					\item[3.] The map $(ev_0,ev_1)=(\delta^0,\delta^1)^*:\text{Maps}(\Delta[1],K)\to \text{Maps}(\partial\Delta[1],K)=K_\bullet\times K_\bullet$ is a Kan fibration. Moreover it satisfies $$\text{Diag}=(ev_0,ev_1)\circ (\sigma^0)^*:K_\bullet\to K_\bullet\times K_\bullet.$$
					\item[4.] The maps in \eqref{eq:EZ-map} induce a transgression morphism $\bT_{\Delta[1]}:\Omega^j(K_{i+1})\to \Omega^j(K^{\Delta[1]}_i),$ for $\omega\in \Omega^j(K_{i+1})$ then $$\bT_{\Delta[1]}(\omega)=\sum_{a=1}^{i+1}sign(\varrho_a)\ \varrho_a^*\omega\quad \text{satisfying}\quad \bT_{\Delta[1]} D\omega+D\bT_{\Delta[1]}\omega=ev_1^*\omega-ev_0^*\omega.$$
					Indeed  we have $\bT_{\Delta[1]}(\omega)=\int_{\Delta[1]}ev^*\omega$ where $\int_{\Delta[1]}$ is the Eilenberg-Zilber map. 
				\end{enumerate}
			\end{proposition}
			\begin{proof}
				The first item is standard from being a mapping space, see e.g. \cite[Sec. I.9]{gh:simp}. The second and the third are proven in  \cite[Prop 7.1 and 7.2]{rcc:hom}. The last item is the usual simplicial decomposition of the cylinders as in \cite[Prop 2.10]{hat:alg}.
			\end{proof}

			\begin{corollary}
				Let $(K_\bullet,\omega_\bullet)$ be an $m$-shifted symplectic Lie $n$-groupoid. Then $\bT_{\Delta[1]}(\omega_\bullet)$ is a shifted Lagrangian structure on the morphism 
				$$(ev_0, ev_1):K^{\Delta[1]}_\bullet\to \overline{K}_\bullet\times K_\bullet.$$
			\end{corollary}
			\begin{proof}
				The isotropic condition holds because $D(\omega_\bullet)=0$ and hence 
				$$D\bT_{\Delta[1]}(\omega_\bullet)=ev_1^*\omega_\bullet-ev_0^*\omega_\bullet$$
				by the fourth item in Proposition \ref{prop:trans}. While the Lagrangian condition is a consequence of the following facts: the morphism $\text{Diag}:K_\bullet\to \overline{K}_\bullet\times K_\bullet$ is shifted Lagrangian with the zero $2$-form; the diagonal can be written as $\text{Diag}=(ev_0,ev_1)\circ (\sigma^0)^*$ and the map $ (\sigma^0)^*$ is a weak equivalence. 
			\end{proof}
			
			In other words,  Proposition \ref{prop:trans} shows that the \emph{path object} $K^{\Delta[1]}_\bullet=\text{Maps}(\Delta[1],K)$ is a fiber replacement of the diagonal and the above corollary indicates that the transgression defines an explicit shifted Lagrangian structure on $(ev_0,ev_1):K^{\Delta[1]}_\bullet\to \overline{K}_\bullet\times K_\bullet$. Thus, from the perspective of Morita theory $$(K_\bullet,\text{Diag}, 0)\quad \text{and}\quad \big(K_\bullet^{\Delta[1]},(ev_0,ev_1), \bT_{\Delta[1]}(\omega_\bullet)\big)$$ give the same shifted Lagrangian on $\overline{K}_\bullet\times K_\bullet$.
			
			\subsection{The symplectic ``category"}\label{sec:symcatlag}
			Inspired by ideas from geometric quantization, it was proposed in \cite{wei:symcat}  that symplectic geometry should be encoded in a ``category" $\cS ymp$ whose objects are symplectic manifolds and morphisms are \emph{Lagrangian relations}, i.e.
			$$\Hom_{\cS ymp}\big((M,\omega),(M',\omega')\big)=\{ L\subseteq\overline{M}\times M' \ | \ L \text{ is a Lagrangian submanifold} \}.$$ The problem with this idea  is that the composition of morphisms is not well defined. One of the main achievements   in \cite{cal:lag, ptvv} was to define an honest symplectic category whose objects are shifted symplectic derived Artin stacks and whose morphisms are shifted Lagrangian relations; see \cite[Thm. 4.4]{cal:lag} and \cite[Thm. 2.9]{ptvv}  (see also \cite{cal:aksz} for a higher categorical statement).
			
			Recall that Lie $n$-groupoids just model higher stacks and therefore do not contain derived directions. In consequence, the symplectic category constructed in \cite{cal:lag} is still not available. Nevertheless, it remains a shadow of the above construction as the next result shows.
			
			\begin{claim}\label{thm:lagint}
				Let $(K^i_\bullet,\omega_\bullet^i)$ be $m$-shifted symplectic Lie $n$-groupoids, $i=1,2,3$, and consider shifted Lagrangian morphisms $(J^1_\bullet,\Phi_\bullet,\eta^1_\bullet)$ into $\overline{K}^1_\bullet\times K_\bullet^2$ and  $(J_\bullet^2,\Psi_\bullet,\eta^2_\bullet)$ into $\overline{K}_\bullet^2\times K_\bullet^3$. 
				If $\Phi^2_\bullet: J_\bullet^1\to K_\bullet^2$ and $\Psi^2_\bullet: J_\bullet^2\to K_\bullet^2$ are transverse maps, then
				$$
				(J_\bullet^1 \times_{K_\bullet^2} J_\bullet^2, (\Phi^1_\bullet,\Psi^3_\bullet), \pr_1^*\eta^1_\bullet+\pr_2^*\eta^2_\bullet)
				$$
				is a shifted Lagrangian morphism into $\overline{K}^1_\bullet\times K_\bullet^3$. 
			\end{claim}
			
			The following diagram gives a visual interpretation for the previous claim
			\begin{equation*}
				\xymatrix{&&J_\bullet^1 \times_{K_\bullet^2} J_\bullet^2\ar[dl]\ar[dr]&&\\
					& J^1_\bullet\ar[dl]_{\Phi^1_\bullet}\ar[dr]^{\Phi^2_\bullet}&& J_\bullet^2\ar[dl]_{\Psi^2_\bullet}\ar[dr]^{\Psi^3_\bullet}&\\
					K_\bullet^1&&K_\bullet^2&&K_\bullet^3.}
			\end{equation*}
			
			For $m=n=1$, this result can be found e.g. in \cite[Thm 4.1]{max:coi}, for $m=2$ and $n=1$ see \cite{bur:lag}, for higher cases,  we do not know any references in the literature.  The transversality condition in the above result  ensures that the fiber product exists as a Lie $n$-groupoid and that the resulting $(m-1)$-shifted $2$-form satisfies the shifted Lagrangian condition. 
			
			An important consequence of Claim \ref{thm:lagint}, when $K_\bullet^1=K_\bullet^3=pt$, together with the fact that shifted Lagrangians on a point are shifted symplectic, see \S~\ref{sec:lagpt}, is the next result that provides a far-reaching generalization of the symplectic reduction procedure on several moment map theories.  
			
			\begin{corollary}\label{cor:red}
				For $i=1,2$,  let  $(J^i_\bullet,\Phi^i_\bullet,\eta^i_\bullet)$ be shifted Lagrangian structures  on the $m$-shifted symplectic Lie $n$-groupoid $(K_\bullet, \omega_\bullet)$. If $\Phi^1_\bullet$ and $\Phi^2_\bullet$ are transverse then 
				$$
				(J_\bullet^1 \times_{K_\bullet} J_\bullet^2, \pr_1^*\eta^1_\bullet+\pr_2^*\eta^2_\bullet)
				$$
				is an $(m-1)$-shifted symplectic Lie $n$-groupoid.
			\end{corollary}

			\subsection{Moment maps, reduction and shifted Lagrangians}\label{sec:mome}
			Here we illustrate how the theory of shifted Lagrangian structures both unifies and generalizes many of the moment map theories and their associated reduction procedures.
			
			\subsubsection{Hamiltonian actions and reduction}
			Let $G$ be a Lie group with Lie algebra $\g$. Then $G$ acts on $\g^*$ via the coadjoint representation, and its action groupoid $G\ltimes\g^*\rightrightarrows\g^*$ is a symplectic groupoid with respect to  $\omega_\can\in\Omega^2(T^*G=G\times \g^*),$ see \cite[\S II.4]{cdw:sym}.

			On the one hand, to a given point $\xi\in\g^*$ we can associate the following shifted Lagrangian morphisms on $(G\ltimes\g^*\rightrightarrows\g^*, \omega_{\can})$:
			\begin{enumerate}
				\item[(L1)] The subgroupoid $\cL\rightrightarrows\mathcal{O}_\xi$, where $\mathcal{O}_\xi$ is the coadjoint orbit through $\xi$ and $$\cL=\{(g,x)\in G\times \g^*=T^*G\ | \ g\in G_x, x\in\mathcal{O}_\xi\}$$ with shifted Lagrangian structure $\eta=0$.
				\item [(L2)] The isotropy subgroup $G_\xi\rightrightarrows\xi$ with $\eta=0$.
				\item [(L3)] The action groupoid $G\ltimes\mathcal{O}_\xi\rightrightarrows\mathcal{O}_\xi$ obtained by restricting $T^*G$ to the coadjoint orbit, and $\eta=\omega_{\text{kks}}\in\Omega^2(\mathcal{O}_\xi)$ the \emph{Kirillov-Kostant-Souriau} symplectic form.
			\end{enumerate}
			The fact that (L1) and (L2) give a shifted Lagrangian morphism is immediate, while (L3) will follow from Proposition \ref{prop:hamlag} below.
			
			On the other hand, recall that given a symplectic manifold  $(M,\omega)$, we say that an action $\mathtt{a}:G\times M\to M$ is \emph{Hamiltonian} if there exists a map
			$$\mu:M\to \g^*$$
			satisfying:
			\begin{enumerate}
				\item[1.] The symplectic form is invariant
				$$\mathtt{a}_g^*\omega=\omega\quad \forall g\in G.$$ 
				\item[2.] The moment map is Hamiltonian, i.e.   $$\langle d\mu, X\rangle=-\iota_{\rho(X)}\omega,$$
				where $\rho(X)\in\fX^1(M)$ is the generating vector field  corresponding to $X\in\g$.
				\item[3.] The moment map is equivariant with respect to the coadjoint action, i.e. 
				$$\mu\circ\mathtt{a}_g=Ad_g^*\circ \mu\quad \forall g\in G.$$
			\end{enumerate}
			
			One can summarize all the above information by using the shifted symplectic geometry language as follows.
			
			\begin{proposition}\label{prop:hamlag}
				Let $G$ be a Lie group and $(M,\omega,\mathtt{a}, \mu)$ be the data defining a Hamiltonian action. Then  $(G\ltimes_{\mathtt{a}}M\rightrightarrows M, \id\times\mu,\omega) $ is a shifted Lagrangian structure for the $1$-shifted symplectic groupoid $(T^*G\rightrightarrows \g^*, \omega_\can).$
			\end{proposition}
			\begin{proof}
				Let $(M,\omega,\mathtt{a}, \mu)$ be the above data defining a Hamiltonian action and consider the action groupoid $G\ltimes_{\mathtt{a}}M\rightrightarrows M$. The third item implies that the moment map defines a Lie groupoid morphism
				\begin{equation*}
					\xymatrix{G\ltimes_{\mathtt{a}}M\ar[r]^{\id\times \mu}\ar@<-.5ex>[d] \ar@<.5ex>[d]& G\ltimes \g^*\ar@<-.5ex>[d] \ar@<.5ex>[d]\\ M\ar[r]^{\mu}&\g^*}
				\end{equation*}
				while the first and second items, together with the fact that $d\omega=0$, imply that 
				$$(\id\times\mu)^*\omega_{\can}=\pr_M^*\omega-\mathtt{a}^*\omega=\delta(\omega)+d\omega=D(\omega).$$
				
				Therefore, $(G\ltimes_{\mathtt{a}}M\rightrightarrows M, \id\times\mu,\omega) $ is an isotropic structure for the $1$-shifted symplectic groupoid $(T^*G\rightrightarrows \g^*, \omega_\can).$
				
				To check the quasi-isomorphism condition, consider $m\in M$. Then the map in \eqref{eq:lagcon} at the point $m$ becomes 
				\begin{equation*}
					\xymatrix{\g\ar[d]\ar[rr]^{(\rho_m, \id)\qquad}&& T_mM\oplus \g\ar[d]^{\omega^\flat_m+T_m^*\mu\circ \langle\cdot,\mu(m)\rangle} \ar[rr]^{\qquad T_m\mu-ad^*_\cdot\mu(m)}&& \g^*\ar[d]^{\id}\\ 
						0\ar[rr]&& T_m^*M\ar[rr]_{-\rho^*_m}&& \g^*
					}
				\end{equation*}
				and it defines a quasi-ismorphism because $\omega^\flat$ is an isomorphism and it satisfies the second item of a Hamiltonian action. 
			\end{proof}
			
			A direct consequence of Corollary \ref{cor:red} is the symplectic reduction obtained in \cite[Thm. 1]{mw:red}.
			
			\begin{corollary}
				Let $G$ be a Lie group with a Hamiltonian action $(M,\omega, \mathtt{a}, \mu)$. Assume that $\xi\in\g^*$ is a regular value of $\mu$ and suppose that
				$G_\xi$ acts freely and properly on the manifold $\mu^{-1}(\xi)$. Then
				there is a unique symplectic structure $\omega_{red}$ on the reduced phase space $M_\text{red}=\mu^{-1}(\xi)/G_{\xi}
				$, such that $p^*\omega_{red}=i^*\omega$ where $p:\mu^{-1}(\xi)\to M_{\text{red}}$ is the projection and $i:\mu^{-1}(\xi)\to M$ the inclusion.
			\end{corollary}
			\begin{proof}
				Since  $\xi\in\g^*$ is a regular value of $\mu$ and
				$G_\xi$ acts freely and properly on the manifold $\mu^{-1}(\xi)$ then the shifted Lagrangians 
				$$(G\ltimes_{\mathtt{a}}M\rightrightarrows M, \omega)\xrightarrow{\id\times \mu}(T^*G\rightrightarrows\g^*, \omega_\can)\xleftarrow{j}(G_\xi\rightrightarrows\xi, 0)$$
				are tranversal. Hence,  Corollary \ref{cor:red} implies that $(G_\xi\ltimes \mu^{-1}(\xi)\rightrightarrows\mu^{-1}(\xi), i^*\omega)$ is a $0$-shifted symplectic groupoid. Moreover, since $G_\xi$ acts freely and properly it is Morita equivalent to the manifold  $M_\text{red}$ and by \S~\ref{subsub:sss0}, it follows that  it is a symplectic manifold. 
			\end{proof}

			\begin{remark}
				The shift trick in Marsden-Weinstein reduction is explained by the fact that the shifted Lagrangian structures on $(L2)$ and $(L3)$ are Morita equivalent. 
			\end{remark}

			\subsubsection{Quasi-Hamiltonian reduction}
			Let $(G\rightrightarrows M, \omega+H)$ be a $1$-shifted symplectic Lie $1$-groupoid as in \S~\ref{subsec:sss1}. Recall that these are called twisted presymplectic groupoids in \cite{bcwz:dir} and quasi-symplectic groupoids in \cite{xu:mome}.

			Following \cite{xu:mome}, we say that $(\phi, X,\eta)$ is a \emph{Hamiltonian $G$-space} if $G\rightrightarrows M$ acts on $\phi:X\to M$, the graph of the action is isotropic inside $G\times X\times X$ (for the 2-form $\omega\times \eta\times -\eta$), $$\phi^*H=d\eta\quad \text{and}\quad\ker \eta_x=\{ \rho(\xi)_x \ | \ \xi\in (\al_G)_{|\phi(x)} \ \text{and}\ \xi^l_{\phi(x)}\in \ker\omega\},$$ where $\rho$ is the infinitesimal action.
			In this situation, the action groupoid $G\times_MX\rightrightarrows X$ together with the morphism $\Phi=\pr_{G}:G\times_MX\to G$ and $\eta$ constitutes a shifted Lagrangian structure on $(\cG,\omega+H)$, see \cite{cal:lag,saf:qua}. These action groupoids are the global counterparts of \emph{presymplectic realizations} as explained in \cite[\S 7.1]{bcwz:dir}.

			The two most important families of examples are given by: 
			\begin{itemize}
				\item Hamiltonian spaces of Lu-Weinstein symplectic groupoids, i.e. one of the symplectic groupoids in Diagram \eqref{eq:dlg_gamma}. These Hamiltonian spaces are symplectic manifolds equipped with a complete Poisson action \cite{luwei, semdres}.
				\item Quasi-Hamiltonian $G$-spaces \cite{amm:lie}. These are the Hamiltonian spaces of the AMM groupoid explained in \S~\ref{sec:amm}. One of the most important features of this last example is that it allows us to construct the symplectic structure on the moduli space of flat connections as explained in \cite{amm:lie}. Indeed that reduction is another instance of Corollary \ref{cor:red}. 
			\end{itemize} 
			
			\begin{remark}
				It was explained in \cite[Thm. 5.3]{saf:qua} that one can understand the \emph{internal fusion} of quasi-Hamiltonian $G$-spaces as a particular case of Claim \ref{thm:lagint}.
			\end{remark}
			
			\subsection{Classical TFT as a functor to the symplectic category}\label{sec:CTFT}
			
			Recall that in \S~\ref{sec:aksz} we indicated how to construct a topological field theory out of a symplectic $Q$-manifold. That approach is based on the Feynman path integral \eqref{eq:correlation}. However, it was proposed in \cite{ati:tqft} that a topological quantum field theory 
			is a functor of symmetric monoidal categories from the category of cobordism to the category of vector spaces.
			
			Using that idea, \cite{cal:aksz} (see also \cite{cal:lec3, cal:lag}) shows that for each $m$-shifted symplectic derived Artin stack, one can construct a (higher) functor between an extended cobordism (higher) category and a (higher) symplectic category.
			
			Here we give an informal presentation of the main ideas in two cases: Classical $3d$-Chern-Simons theory as explained in \cite{saf:qua}, and a $4d$ field theory with defects that contains the Fock-Rosly Poisson structure \cite{fr:poi}, which will be explained in detail in \cite{bur:lag}.
			
			\subsubsection{Classical Chern-Simons as a functor}\label{sec:csf}
			In this section we mainly follow \cite{saf:qua}. The idea is to construct a functor between the higher category of oriented cobordisms \cite{bd:cob} and the higher symplectic category, whose morphisms are Lagrangian relations and relations among them; see \cite{cal:aksz} for a precise definition. 
			
			Recall that the input for Chern-Simons theory is a Lie group $G$ whose Lie algebra is quadratic. We said in Proposition \ref{prop:BG} that this data makes $( G, \omega_2+\omega_1)$ into a $2$-shifted symplectic Lie $1$-group.
			
			If we think of classical Chern-Simons as a fully extended field theory, see \cite{bd:cob, lur:cla}, the functor $$CS:\cC ob\to \cS ymp$$ must be determined by its image on the point. Therefore, following \cite{saf:qua} we define
			$$CS(pt)=( G, \omega_2+\omega_1).$$
			Now, if we think of an interval $\bullet\to\bullet$ as a trivial cobordism between two points, then 
			$$CS(\bullet\to\bullet)=\id: G\to  G\equiv \text{Diag}: G\to\overline{ G}\times  G.$$
			Observe that the change of sign on the symplectic form corresponds to the fact that the cobordism is oriented. 
			
			The next step is to understand the value of the functor on $S^1$. Observe that $S^1$ is a cobordism from the empty set to itself. Therefore, it must produce a $1$-shifted symplectic groupoid.  On the other hand, by marking two points on $S^1$, we obtain that:  
			\begin{equation*}
				\begin{array}{c}
					\begin{tikzpicture}
						\filldraw[fill=red!0](-5,0) circle (0.25);
						\filldraw  (-4.75,0) circle (0.5pt);
						\filldraw  (-5.25,0) circle (0.5pt);
						\draw (-4.35,0.05) -- (-4.1,0.05);
						\draw (-4.35,-0.05) -- (-4.1,-0.05);
						\draw[->] (-3.25,0.1) arc (0:180:0.25);
						\draw[->] (-3.75,-0.1) arc (180:360:0.25);
						\filldraw  (-3.75,-0.1) circle (0.5pt);
						\filldraw  (-3.75,0.1) circle (0.5pt);
						\filldraw  (-3.25,0.1) circle (0.5pt);
						\filldraw  (-3.25,-0.1) circle (0.5pt);
					\end{tikzpicture}
				\end{array}.
			\end{equation*}
			Hence, the value $CS(S^1)$ must coincide with the homotopic fiber product of two intervals, in other words
			$$CS(S^1)= G\underset{\overline{ G}\times G}{\times^h}  G.$$
			Now, by replacing the diagonal with the path object explained in\S~\ref{sec:path},  we get that 
			$$CS(S^1)= G\underset{\overline{ G}\times  G}{\times^h}  G= G\underset{\overline{G}\times  G}{\times}  G^{\Delta[1]}.$$
			Indeed, it is not difficult to show that the last fiber product has transversal morphisms and, therefore, we can use Corollary \ref{cor:red}. An easy computation using the formulas of the transgression shows that 
			$$CS(S^1)=  G\underset{\overline{G}\times  G}{\times}  G^{\Delta[1]}= G^{AMM}$$
			with $1$-shifted symplectic structure given by \eqref{eq:AMM} and \eqref{eq:c3f}. 
			
			Therefore, we proved Proposition \ref{prop:amm1} using Corollary \ref{cor:red}. Notice also that this gives another explanation for the symplectic Morita equivalence found in \S~\ref{sec:amm}.
			
			Finally, recall that any closed surface $\Sigma$ can be glued from the following two kinds of pieces: disks $D$ and pair of pants $P$. 
			\begin{itemize}
				\item The disk is a cobordism from the empty set into $S^1$. Therefore, $CS(D)$ must be a shifted Lagrangian on the AMM groupoid. It was computed in \cite[\S 5.3]{saf:qua} that this shifted Lagrangian is given by the conjugacy class of the identity, i.e.
				$$CS(D)=e: (G\rightrightarrows pt, 0)\to (G^{AMM},\omega+H), \quad e(pt)=e\quad \text{and}\quad e(g)=(g,e).$$
				\item The pair of pants is a cobordism from $S^1$ into two copies of $S^1$. It was shown in \cite[Thm. 5.3]{saf:qua} that this shifted Lagrangian correspondence is given by the \emph{internal fusion of quasi-Hamiltonian $G$-spaces} \cite{amm:lie}. More precisely, this is a shifted Lagrangian structure on  $K:=\overline{G}\times G\times G$ given as follows: Denote by $J$ the action groupoid of $G$ acting diagonally on $G^2$ by conjugation. The morphism $\Phi:J\to K$ together with   $\eta\in\Omega^2(G^2)$ given by
				\[ 
				(g,(a,b))\xrightarrow{\Phi} ((g,ab),(g,a),(g,b)) \quad \text{and}\quad \eta=\frac{1}{2}\langle \pr_1^*\theta^l, \pr_2^*\theta^r\rangle
				\]
				for  $(g,(a,b))\in G\ltimes G^2$ is a shifted Lagrangian structure.
			\end{itemize}
			Therefore, for any closed surface $\Sigma$, we get that $CS(\Sigma)$ gives the symplectic structure on the moduli space of flat $G$-connections, see \cite{ab:ym}, as a fusion of quasi-Hamiltonian $G$-spaces \cite{amm:lie}.
			
			\subsubsection{A 4d TFT with defects} 
			Here, instead of considering the plain (higher) category of cobordisms  we attach to it one new morphism given by an interval decorated by a marked point on its interior. For more precise mathematical formulations of decorations (also known as defects), see e.g. \cite{cat:def, sev:quil}.
			
			As in the previous section \S~\ref{sec:csf}, the input data is a Lie group $G$ with a quadratic Lie algebra $(\g,[\cdot,\cdot],\langle\cdot,\cdot\rangle),$ so we get a $2$-shifted symplectic structure $(G, \omega_2+\omega_1).$ 
			
			Here we construct a ``functor" $\cF:\cC ob[a]\to \cS ymp$, where the $[a]$  means that we have an extra morphism, namely an interval decorated with a marked point on the interior. For precise mathematical details of this construction we refer to \cite{bur:lag}.  
			Let $Pair(G)$ denote the pair groupoid of $G,$ which, as we know, is Morita equivalent to a point, and we think of it as a $3$-shifted symplectic Lie group. Then 
			$$\cF(pt)=Pair(G)\equiv pt.$$
			Reasoning in the same way as before, one gets that $\cF(\bullet\to\bullet)=pt,$ which is a $3$-shifted Lagrangian on two points and also $2$-shifted symplectic. 
			
			Denote by $a$ a segment with a marked point in its interior. Therefore, we have freedom to send this new morphism to a $3$-shifted Lagrangian on a point, i.e. a $2$-shifted symplectic Lie $n$-groupoid. Let us define
			$$\cF(a)=\overline{G}\times G.$$
			Some easy computations show that 
			\begin{equation*}
				\begin{array}{c}
					\begin{tikzpicture}
						\filldraw[fill=red!0](-5,0) circle (0.25);
						\filldraw  (-5,0.25) circle (0.5pt);
						\filldraw (-5,-0.25) circle (0.5pt);
						\filldraw [fill=red!100] (-4.75,0) circle (1.2pt);
						\draw (-4.35,-0.05) -- (-4.1,-0.05);
						\draw (-4.35,0.05) -- (-4.1,0.05);
						\draw (-3.75,0.25) -- (-3.75,-0.25);
						\filldraw  (-3.75,0.25) circle (0.5pt);
						\filldraw  (-3.75,-0.25) circle (0.5pt);
						\filldraw [fill=red!100] (-3.75,0) circle (1.2pt);
					\end{tikzpicture}
				\end{array}.
			\end{equation*}
			The next task is to compute $\cF(D_2)$, where $D_2$ is a disk with two marked points on the boundary, as shown in the next diagram
			\begin{equation*}
				\begin{array}{c}
					\begin{tikzpicture}
						\filldraw[fill=blue!25](-5,0) circle (0.25);
						\filldraw  (-5,0.25) circle (0.5pt);
						\filldraw (-5,-0.25) circle (0.5pt);
						\filldraw [fill=red!100] (-4.75,0) circle (1.2pt);
						\filldraw [fill=red!100] (-5.25,0) circle (1.2pt);
					\end{tikzpicture}
				\end{array}.
			\end{equation*}
			Since there are two marked points, this means that to $D_2$ one should assign a shifted Lagrangian structure on $(\overline{G}\times G)^2$. This shifted Lagrangian structure is given by $(Pair(G)\times G^2,\Phi,\beta)\to (\overline{G}\times G)^2$, where
			\[ \Phi(u_1,u_2,g_1,g_2)=((u_1g_2u_2^{-1},g_1),(u_1^{-1}g_1u_2,g_2))\in (\overline{G}\times G)^2\quad\text{and}\]
			$$\beta=(u_1,g_2u_2^{-1})^*\omega_2-(g_2,u_2^{-1})^*\omega_2-(u_1^{-1},g_1u_2)^*\omega_2+(g_1,u_2)^*\omega_2\in\Omega^2\big(G^4=(Pair(G)\times G^2)_1\big).$$
			
			Reasoning similarly, we get that a pair of pants with three marked points $P_3$, as drawn in \eqref{fig:p3} also corresponds to a disk with three marked points $P_3=D_3,$ 
			\begin{equation}\label{fig:p3}
				\begin{array}{c}
					\begin{tikzpicture}
						\filldraw[fill=blue!25] (-5,1) -- (-5,0.5) -- (-4.5,0.25) -- (-4.5,-0.25)-- (-5,-0.5)-- (-5,-1)-- (-2,-0.25)-- (-2,0.25)-- cycle;
						\filldraw [fill=red!100] (-5,0.75) circle (1.2pt);
						\filldraw [fill=red!100] (-5,-0.75) circle (1.2pt);
						\filldraw [fill=red!100] (-2,0) circle (1.2pt);
						\draw (-1.75,-0.05) -- (-1.5,-0.05);
						\draw (-1.75,0.05) -- (-1.5,0.05);
						\filldraw[fill=blue!25](-0.75,0) circle (0.5);
						\filldraw [fill=red!100] (-0.75,0.5) circle (1.2pt);
						\filldraw [fill=red!100] (-1.25,0) circle (1.2pt);
						\filldraw [fill=red!100] (-0.25,0) circle (1.2pt);
					\end{tikzpicture}
				\end{array}.
			\end{equation}
			The image of the functor $\cF(\cdot)$ on $D_3$ is the shifted Lagrangian on $(\overline{G}\times G)^3$ given by the multiplication of the pair groupoid. More concretely, the shifted Lagrangian is $\cF(D_3)=(G^3, \Phi, 0)$
			where $$\Phi:G^3\to (\overline{G}\times G)^3\quad \text{is given by}\quad \Phi(z,y, x)=(z,y,y,x ,x,z )\in(\overline{G}\times G)^3.$$
			
			One can show that figure \eqref{fig:p3} gives the \emph{internal fusion for quasi-Poisson $\g_\Delta$-manifolds}, see \cite{bur:lag, sev:quil} for more details. In particular, a given surface with marked points on the boundary $(\Sigma, V)$ is obtained by composing the shifted Lagrangians of $D_2$ and $D_3$.
			
			One of the main relevance of the above field theory is that, given a surface $\Sigma$ with  $V$ marked points on its boundary, then $\Hom(\pi_1(\Sigma, V), G)$ carries a Poisson structure known as the \emph{Fock-Rosly Poisson structure} \cite{fr:poi}.  In \cite{bur:lag}, an integration of such a Poisson structure was provided using the shifted Lagrangian structures associated with  $D_2$ and $D_3$,  and another one coming from the quasi-triangular $r$-matrix on $\g$.

			\bibliographystyle{siam}
			\bibliography{biblioN}

		\end{document}